\documentclass[ijoc,nonblindrev]{informs3} 

\OneAndAHalfSpacedXII 

\usepackage{mathrsfs}
\usepackage{amsmath}
\usepackage{dsfont}
\usepackage{enumerate}
\usepackage{longtable}
\usepackage{graphicx}
\usepackage{epstopdf}
\usepackage{amsfonts}
\usepackage{nicefrac}
\usepackage{color}
\usepackage{soul}

\usepackage{algorithmic}
\usepackage{algorithm}

\usepackage{subfigure}

\newcommand{\tabincell}[2]{\begin{tabular}{@{}#1@{}}#2\end{tabular}}

\usepackage{natbib}
 \bibpunct[, ]{(}{)}{,}{a}{}{,}%
 \def\bibfont{\small}%
\TheoremsNumberedThrough     

\EquationsNumberedThrough    

\begin{document}


\RUNAUTHOR{Zhang et al.} 

\RUNTITLE{Asymptotically Optimal Sampling Policy for Selecting Top-$m$ Alternatives}

\TITLE{Asymptotically Optimal Sampling Policy for Selecting Top-$m$ Alternatives}

\ARTICLEAUTHORS{
\AUTHOR{Gongbo Zhang}
\AFF{Department of Management Science and Information Systems, Guanghua School of Management, Peking University, Beijing 100871, China, \EMAIL{gongbozhang@pku.edu.cn}} 
\AUTHOR{Yijie Peng}
\AFF{Department of Management Science and Information Systems, Guanghua School of Management, Peking University, Beijing 100871, China, \EMAIL{pengyijie@pku.edu.cn}}
\AUTHOR{Jianghua Zhang}
\AFF{School of Management, Shandong University, Jinan 250100, China, \EMAIL{zhangjianghua@sdu.edu.cn}}
\AUTHOR{Enlu Zhou}
\AFF{School of Industrial and Systems Engineering, Georgia Institute of Technology, GA 30332, USA, \EMAIL{enlu.zhou@isye.gatech.edu}}
} 

\ABSTRACT{%
 We consider selecting the top-$m$ alternatives from a finite number of alternatives via Monte Carlo simulation. Under a Bayesian framework, we formulate the  sampling decision as a stochastic dynamic programming problem, and develop a sequential sampling policy that maximizes a value function approximation one-step look ahead. To show the asymptotic optimality of the proposed procedure, the asymptotically optimal sampling ratios which optimize large deviations rate of the probability of false selection for selecting top-$m$ alternatives  has been rigorously defined. The proposed sampling policy is not only proved to be consistent but also achieve the asymptotically optimal sampling ratios. Numerical experiments demonstrate superiority of the proposed allocation procedure over existing ones.
}%



\maketitle

%


\section{Introduction}

Monte Carlo simulation has been widely used for evaluating and analyzing complex stochastic systems such as manufacturing, communication and service systems, whose performances do not have analytical forms~\citep{cassandras2009introduction}. Since simulation is typically expensive, efficiency becomes a major concern, especially when the number of competing alternatives is large. Ordinal optimization (OO) was introduced in simulation of discrete event dynamic systems (DEDS) to soften hard optimization problems~\citep{ho1992ordinal}. The tenet of OO is that correctly ordering the performance of competing alternatives would be easier than accurately estimating their performance, and selecting a subset of good alternatives would be easier than selecting the best. In this work, we consider selecting alternatives with the top-$m$ largest (or smallest) means from a finite set of $k$ alternatives ($1 \le m < k$, $m,k \in \mathbb{Z}^{+}$). The true mean performance of each alternative is unknown and can only be estimated by Monte Carlo simulation within a finite simulation budget. Our research objective is to efficiently allocate the simulation budget to alternatives so that the probability of correct selection (PCS) for top-$m$ alternatives can be maximized.

Our problem generalizes a classic ranking and selection (R\&S) problem, which is to find the alternative with the largest (or the smallest) mean, i.e., the best alternative. One of the most commonly used metrics for R\&S procedures is PCS, and the other popular metric is expected opportunity cost (EOC). There are two different ways to categorize R\&S procedures. In terms of methodologies, there are frequentist and Bayesian branches, whereas in terms of goals, there are fixed-precision and fixed-budget procedures. See~\citet{kim2006selecting, hunter2017parallel} and~\citet{hong2021review} for overviews. Fixed-precision procedures~\citep[e.g.,][]{rinott1978two,kim2001fully,frazier2014fully} allocate simulation replications to guarantee a pre-specified PCS level, whereas fixed-budget procedures~\citep[e.g.,][]{chen2000simulation,frazier2009knowledge,chick2010sequential} maximize a posterior performance under a fixed simulation budget constraint. Frequentist procedures tend to allocate more simulation replications than necessary to guarantee a PCS level, and Bayesian procedures usually achieve higher PCS under the same simulation budget or achieve the same level PCS using fewer simulation replications~\citep{chen2000simulation}. R\&S shares many similarities with the multi-armed bandit (MAB) problem, however, there are also some differences, e.g., the reward is collected when the total simulation replications are exhausted and depends on the states of all alternatives in the R\&S, whereas the rewards are collected at all steps and each reward depends on the state of each arm in the standard MAB problem.

Many R\&S procedures focus on enhancing the efficiency for finding the best alternative, including optimal computing budget allocation (OCBA)~\citep{chen2006efficient,chen2011stochastic}, expected value of information (EVI)~\citep{chick2001new, chick2010sequential}, knowledge gradient (KG)~\citep{gupta1996bayesian, frazier2008knowledge}, expected improvement (EI)~\citep{jones1998efficient, ryzhov2016convergence} and asymptotically optimal allocation procedure (AOAP)~\citep{peng2018ranking}. OCBA formulates the sampling allocation decision as a static optimization problem, and approximates the objective function with various bounds and relaxations for simplicity of computation, and then solves the approximate problem using gradient ascent or greedy heuristics or by asymptotically optimal in the long run. OCBA is originally a two-stage method, and is then extended to a fully sequential manner by combining with certain heuristic rules, e.g., ``most starving" sequential rule \citep{chen2011stochastic}.~\citet{glynn2004large} provide a rigorous framework for determining the asymptotically optimal sampling ratios for selecting the best by maximizing the large deviations rate of the probability of false selection (PFS). The sampling ratios of OCBA are proved to be an approximation of the asymptotically optimal sampling ratios under normal sampling distributions. EVI is derived based on a Bayesian decision-theoretic approach and it is originally proposed by using approximation techniques such as Bonferroni and Slepian's inequalities to relax the objective function and then derived as a myopic allocation policy. KG and EI are myopic allocation procedures, which look one-step ahead to maximize the posterior information gains. AOAP is a myopic allocation procedure that maximizes a value function approximation one-step look ahead, and has an analytical form reflecting a mean-variance trade-off, which is the hallmark of OCBA. AOAP can achieve the asymptotically optimal sampling ratios defined in \citet{glynn2004large} under normal sampling distributions, whereas neither  KG nor EI achieves the asymptotically optimal sampling ratios \citep{peng2018review}.

Some sequential sampling procedures are developed based on the rationale that certain desirable asymptotic property can be achieved. For example, ~\citet{xiao2013optimal} and~\citet{xiao2017simulation} develop sequential sampling procedures guided by estimated ratios derived from certain asymptotic optimality conditions, and~\citet{gao2016new,gao2017new,chen2019balancing} and~\citet{xiao2020optimal} develop sampling procedures by sequentially balancing the certain asymptotic optimality conditions. However, recent studies show that the asymptotic property is inadequate to capture the finite-sample behavior of sequential sampling policies. Specifically, the performance of OCBA under perfect information (assuming true parameters are known) could be far inferior to that of a sequential OCBA algorithm with imperfect information \citep{chen2006efficient}. This means that the asymptotic optimality conditions per se could lead to poor finite-sample performance \citep{peng2018review}. A sequential sampling policy guided by the asymptotic optimality conditions could happen to achieve a desirable finite-sample performance in certain scenarios, but this is not always the case. For selecting the alternative with the optimal mean under normal sampling distributions, \cite{peng2017gradient} consider a low-confident scenario, where the means of alternatives are close, while their variances are relatively large and simulation budget is relatively small. In low-confident scenarios, many sequential sampling policies guided by the asymptotic optimality conditions decrease PCS, whereas equal allocation can increase the PCS \citep{peng2018review}. For exponential sampling distributions and selecting the optimal quantile, sequential sampling policies guided by the asymptotic optimality conditions are reported to perform poorly \citep{zhang2020wsc,shin2016tractable}. Therefore, developing a sequential sampling procedure solely by asymptotic analysis may lead to misleading results for more general and complicated problems, and a stochastic dynamic programming formulation rigorously capturing dynamic sampling and selection decisions would be advantageous \citep{peng2014dynamic,peng2018ranking}.

Instead of focusing on selection of the best, selecting top-$m$ alternatives provides more choices for decision makers, and more choices sometimes can lead to a more resilient design. For example, in an emergency evacuation problem, the evacuees are evacuated to safe areas in response to disasters such as floods, earthquakes, hurricanes and stampedes. A critical issue is route planning which uses network flow and routing algorithms to schedule evacuees on each route. In the process of evacuation, the travel time and capacity of a route are usually uncertain since the impact of a disaster on the evacuation network is unknown in advance. A metric for the emergency response to a disaster is clearance time, which is the duration between the start of evacuation and the time when the last evacuee reaches safe areas. Evacuation plans with different expected clearance time can be regarded as alternatives with different performances. Some routes may become unavailable when a disaster happens, making some prepared evacuation plans infeasible, and therefore, several well-chosen evacuation plans would improve resilience in response to a disaster compared with finding a single best plan~\citep{cova2003network,pillac2016conflict,wang2016evacuation,zhang2018emergency}.

We briefly review the literature for selecting top-$m$ alternatives.~\citet{koenig1985procedure} adopt an indifference zone paradigm in the frequentist branch, and propose a two-stage procedure which provides a PCS guarantee. \citet{chen2008efficient} generalize OCBA to selection of top-$m$ by using a parameter to separate the top-$m$ and other alternatives, however, the proposed allocation procedure is highly sensitive to a separating parameter, and cannot be reduced to OCBA when $m=1$. \citet{zhang2012improved, zhang2015simulation} improve the approximation in \citet{chen2008efficient} by considering two sub-optimization problems to avoid using the separating parameter, however, when $m=1$, the proposed sampling allocation procedure can only be reduced to OCBA in a special case where the variances of all alternatives are equal. \citet{gao2015efficient} and \citet{gao2015note} propose sequential sampling procedures to asymptotically optimize the EOC and PCS metrics, respectively.~\citet{gao2016new} further extend the results in~\citet{gao2015note} to general sampling distributions by using a large deviations principle and Bonferroni inequality to approximate PCS. Both sampling procedures in \citet{gao2015note} and \citet{gao2016new} are developed by sequentially balancing certain asymptotic conditions derived under a static optimization problem, without considering finite-sample behavior in their formulation.

In this work, we  extend the results in both~\citet{peng2018ranking} and~\citet{glynn2004large} to the problem of selecting top-$m$ alternatives. We rigorously formulate the dynamic sampling decision for selecting top-$m$ alternatives as a stochastic dynamic programming problem under a Bayesian framework, which captures the finite-sample behavior. The optimal sampling policy can be governed by the Bellman equations of a Markov decision process (MDP). However, solving a MDP typically suffers from curse-of-dimensionality. In principle, some state-of-art simulation-based reinforcement learning techniques could be used to address the difficulty. Yet as in the R\&S-related literature, we try to derive an efficient sampling procedure by approximations so that our proposed procedure is comparable to the existing sampling procedures for selecting top-$m$ alternatives. We focus on a commonly assumed normal sampling distributions, and propose a sequential allocation policy that maximizes a value function approximation one-step look ahead, which has an analytical form. The proposed sampling procedure is derived from approximating the value function of the MDP rather than any asymptotic property of certain static optimization problem, and it is reduced to AOAP when $m=1$.  Yet our procedure will be proved to achieve asymptotic optimality for selecting top-$m$ alternatives. To show this, we rigorously define the asymptotically optimal sampling ratios that optimize large deviations rate of the PFS. The asymptotic optimality condition for the top-$m$ selection is a generalization and much more complicated to derive, relative to that for selection of the best. Our asymptotic analysis buttresses the asymptotic solutions of approximate static optimization problems in \citet{gao2015note} and \citet{gao2016new}, but also reveals that their results are incomplete for asymptotic optimality conditions in the top-$m$ selection problem. Neither the asymptotic sampling ratios of the allocation procedures in \citet{chen2008efficient} nor those in \citet{zhang2012improved, zhang2015simulation} guarantee asymptotic optimality. The efficient sampling procedure proposed in our work is not only proved to be consistent, i.e., the top-$m$ alternatives can be definitely selected as the simulation budget goes to infinity, but also leads to the asymptotically optimal sampling ratios. Moreover, numerical results show that our proposed procedure leads to superior performance. The contributions of our work are summarized as follows:
\begin{itemize}
    \item We provide a stochastic dynamic programming framework to capture the finite-sample behavior of sequential sampling decision for selecting top-$m$ alternatives.
    \item We rigorously define the asymptotically optimal sampling ratios which optimize large deviations rate of PFS for selecting top-$m$ alternatives.
    \item We propose an efficient sequential allocation procedure for selecting top-$m$ alternatives, which is proved to be asymptotically optimal and leads to superior performance.
\end{itemize}

The rest of the paper is organized as follows. Section~\ref{SectionPF} formulates the sequential sampling problem for selecting top-$m$ alternatives. Section~\ref{SectionAOA} proposes a sequential sampling procedure. Section~\ref{Sectionopt} derives the asymptotically optimal sampling ratios for selecting top-$m$ alternatives and proves the proposed sequential allocation procedure to be asymptotically optimal. Section~\ref{SectionNuEx} presents numerical examples, and Section~\ref{SectionCon} concludes the paper. The proofs of the theorems and lemmas in the paper are relegated to the online appendices, which also include additional numerical experiments.

\section{Problem Formulation}\label{SectionPF}

In this section, we formulate the problem under a Bayesian framework, and define the optimal sampling policy in a stochastic dynamic programming formulation. Then we introduce the conjugate prior of normal sampling distributions with unknown mean and known variance for computational convenience.

\subsection{Bayesian Framework}

Suppose there are $k$ alternatives and the performance of each alternative is measured by unknown mean $\mu_i \in \mathbb{R}$, $i = 1,2, \cdots, k$, where $\mu_i$ can be estimated by Monte Carlo simulation. We assume that the performance of each alternative can be distinguished, i.e., $\mu_i \ne \mu_j$, $i,j=1,2,\cdots,k$, $i \ne j$. The objective is to find the alternatives with the top-$m$ largest means. Denote the true set of top-$m$ alternatives as ${\mathcal{F}^m}\mathop  = \limits^\Delta \left\{ {\left\langle 1 \right\rangle ,\left\langle 2 \right\rangle, \cdots ,\left\langle m \right\rangle } \right\}$, where $\left\langle i \right\rangle$, $i=1,2,\cdots,k$ are indices ranked by system performances such that ${\mu _{\left\langle 1 \right\rangle }} > {\mu _{\left\langle 2 \right\rangle }} >  \cdots  > {\mu _{\left\langle m \right\rangle }} > {\mu _{\left\langle {m + 1} \right\rangle }} >  \cdots  > {\mu _{\left\langle k \right\rangle }}$.

Let $X_{i,t}$, $t \le T$, $t \in \mathbb{Z}^+$ be the $t$-th independent and identically distributed (i.i.d.) replication for alternative $i$, where $T < \infty$ is the number of total simulation replications. We assume that ${X_t} \mathop  = \limits^\Delta \left( {{X_{1,t}},{X_{2,t}}, \cdots ,{X_{k,t}}} \right)$ follows a joint distribution $Q\left( { \cdot ;\theta } \right)$, and replications across different alternatives are independent, where $\theta \mathop  = \limits^\Delta \left( {{\theta _1},{\theta _2}, \cdots ,{\theta _k}} \right) \in \Theta $ comprises all unknown parameters of alternatives, including the unknown means, i.e., $\mu \subset \theta$, $\mu \mathop  = \limits^\Delta \left\{\mu_{\left\langle 1 \right\rangle},\mu_{\left\langle 2 \right\rangle},\cdots,\mu_{\left\langle k \right\rangle}\right\}$. Let ${Q_i}\left( { \cdot ;{\theta _i}} \right)$ be the marginal distribution of alternative $i$, where ${\theta _i}$, $i=1,2,\cdots,k$ comprises all unknown parameters in the marginal distribution, and $Q\left( { \cdot ;\theta } \right) = \prod\nolimits_{i = 1}^k {{Q_i}\left( { \cdot ;{\theta _i}} \right)}$. The framework includes non-Gaussian and correlated cases since no assumptions are made for $Q\left( { \cdot ;\theta } \right)$ until later sections where a specific analytical form of the sampling allocation decision is derived.

Denote the set of selected top-$m$ alternatives after allocating $T$ simulation replications as ${\widehat {\mathcal{F}}_T^m} \mathop  = \limits^\Delta  \left\{ {{\left\langle 1 \right\rangle}_T ,{\left\langle 2 \right\rangle}_T , \cdots ,{\left\langle m \right\rangle}_T } \right\}$. A correct selection occurs when the selected set ${\widehat {\mathcal{F}}_T^m}$ is the true set ${\mathcal{F}^m}$, i.e., ${\widehat {\mathcal{F}}_T^m} = {\mathcal{F}^m}$, that is, the worst performance of alternatives in ${\widehat {\mathcal{F}}_T^m}$ should be better than the best performance of other $\left(k-m\right)$ alternatives.

A possible sampling allocation for selecting the top-$m$ alternatives can be given by the solution of a static optimization problem from a frequentist perspective:
\begin{alignat}{2}\label{staticopt}
\mathop {\max }\limits_{{r_1},{r_2}, \cdots ,{r_k}} \quad & \Pr \left\{ \left. \bigcap\limits_{i = 1}^m {\bigcap\limits_{j = m + 1}^k {\left( {\bar{X}_{\langle i\rangle}\left(r_{\langle i\rangle} T\right)>\bar{X}_{\langle j\rangle}\left(r_{\langle j \rangle} T\right)} \right)} }  \right|\theta \right\} \\
\mbox{s.t.}\quad
& \sum\limits_{i = 1}^k {{r_{\left\langle i \right\rangle}} = 1}, \quad r_{\left\langle i \right\rangle} \ge 0 \nonumber~,
\end{alignat}
where $r_i$ is the sampling ratio. $t_i = r_i T$ is the number of simulation replications allocated to alternative $i$, where a technicality that the number of simulation replications needs to be integer is ignored, and ${\bar{X}_i}\left(r_i T\right) = { \left. \sum\nolimits_{\ell = 1}^{t_i} {{X_{ i,\ell}}} \right/ t_i}$, $i=1,2,\cdots,k$ is sample mean. In general, the optimal solution of problem (\ref{staticopt}) does not have an analytical form and is expensive to calculate, and an alternative static optimization problem (\ref{opt0}) based on large deviations rate of PFS in Section 4 is much simpler. Even though the optimal solution of (\ref{staticopt}) or (\ref{opt0}) can be computed exactly, it depends on the unknown parameters  which need to be estimated sequentially from samples in implementation. Thus, the sampling allocation decision essentially becomes a policy adaptive to the information filtration generated by allocated samples. For analyzing the efficiency of a sampling policy, we introduce a Bayesian framework to quantify the uncertainty of unknown parameter $\theta$.

Suppose the unknown parameter $\theta$ follows a prior distribution $F\left( { \cdot ;{\zeta _0}} \right)$, where $\zeta_0$ contains all hyper-parameters for the parametric family of the prior distribution. Let $X_i^{\left( t \right)} \mathop  = \limits^\Delta \left( {{X_{i,1}},{X_{i,2}}, \cdots , {X_{i,{t_i}}}} \right)$ be the sample observations of alternative $i$ throughout $t$ steps, $i=1,2,\cdots,k$, which are ordered in a chronological arrangement, and then we have sample information set ${\mathcal{E} _t} \mathop  = \limits^\Delta \{ {X_1^{\left( t \right)},X_2^{\left( t \right)}, \cdots , X_k^{\left( t \right)}},{\zeta _0} \}$ for $0<t \le T$ and ${\mathcal{E} _0} \mathop  = \limits^\Delta {\zeta _0}$. Following the Bayes' rule, the posterior distribution of $\theta$ conditional on the sample information set $\mathcal{E}_t$ can be calculated as
$$F\left( {\left. {d\theta } \right|{\mathcal{E} _t}} \right) = \frac{{L\left( {{\mathcal{E} _t};\theta } \right)F\left( {d\theta ;{\zeta _0}} \right)}}{{\int_{\theta  \in \Theta } {L\left( {{\mathcal{E} _t};\theta } \right)F\left( {d\theta ;{\zeta _0}} \right)} }}~,$$
where $L\left(  \cdot; \theta \right)$ is a likelihood function that evaluates joint density of simulation replications at their observations and $d \cdot$ stands for Lebesgue and counting measure for continuous and discrete distributions, respectively. The predictive distribution for ${X_{i,{t_i} + 1}}$ conditioned on the sample information set $\mathcal{E}_t$ can be calculated as
$${Q_i}\left( {\left. {d{x_{i,{t_i} + 1}}} \right|{\mathcal{E} _{t}}} \right) = \frac{{\int_{\theta  \in \Theta } {{Q_i}\left( {d{x_{i,{t_i} + 1}};{\theta _i}} \right)L\left( {{\mathcal{E} _{t}};\theta } \right)F\left( {d\theta ;{\zeta _0}} \right)} }}{{\int_{\theta  \in \Theta } {L\left( {{\mathcal{E} _{t}};\theta } \right)F\left( {d\theta ;{\zeta _0}} \right)} }}~.$$

In the case where the prior distribution $F\left( { \cdot ;{\zeta _0}} \right)$ is a conjugate prior of the sampling distribution $Q\left( { \cdot ;\theta } \right)$, the posterior distribution lies in the same parametric family of $F\left( { \cdot ;{\zeta _0}} \right)$, i.e., $F\left( { \cdot ;{\zeta _t}} \right)$, where ${\zeta _t}$ is the posterior hyperparameter, and the sample information set $\mathcal{E}_t$ can be completely determined by ${\zeta _t}$, i.e., $\mathcal{E}_t = \zeta _t$. Under the Bayesian framework, the indices in $\widehat {\mathcal{F}}_T^m$ are ranked by posterior means, and the posterior PCS for selecting top-$m$ alternatives can be expressed as
$$\begin{aligned}
{\rm {PCS}_T}  & = \Pr \left\{ \left. \widehat {\mathcal{F}}_T^m = {\mathcal{F}^m} \right|\mathcal{E}_T \right\} = \Pr \left\{ \left. {{\bigcap\limits_{i = 1}^m {\bigcap\limits_{j = m + 1}^k {\left( {{\mu _{{{\left\langle i \right\rangle }_T}}} > {\mu _{{{\left\langle j \right\rangle }_T}}}} \right)} } } } \right|\mathcal{E}_T \right\} \\
& = \Pr \left\{ {\left. {\mathop {\min }\limits_{i = 1,2, \cdots ,m} {\mu _{{{\left\langle i \right\rangle }_T}}} > \mathop {\max }\limits_{j = m + 1,m + 2, \cdots ,k} {\mu _{{{\left\langle j \right\rangle }_T}}}} \right|{\mathcal{E}_T}} \right\}~.
\end{aligned}$$

\subsection{Optimal Sampling Policy}

In order to develop an optimal sequential sampling policy such that the posterior PCS for selecting the top-$m$ alternatives can be maximized, we formulate the dynamic sampling decision as a stochastic dynamic programming problem under the Bayesian framework. The allocation policy ${\mathcal{A}_T}\left(  \cdot  \right) \mathop  = \limits^\Delta \left( {{A_1}\left(  \cdot  \right),{A_2}\left(  \cdot  \right), \cdots ,{A_T}\left(  \cdot  \right)} \right)$ is a sequence of mappings, where ${A_t}\left( \cdot \right) \in \left\{ {1,2, \cdots ,k} \right\}$ allocates the $t$-th simulation replication to an alternative based on $\mathcal{E} _{t - 1}$ collected throughout the first $\left( {t - 1} \right)$ steps. Notice that the allocation decision and the sample information set are nested in each other as $t$ evolves. Define ${A_{i,t}}\left( {{\mathcal{E} _{t - 1}}} \right) \mathop  = \limits^\Delta \mathds{1}\left( {{A_t}\left( {{\mathcal{E} _{t - 1}}} \right) = i} \right)$, where $\mathds{1}\left(\cdot\right)$ is an indicator function that equals 1 when the event in the bracket is true and equals 0 otherwise, and then the number of simulation replications allocated to alternative $i$ throughout $t$ steps can be expressed as ${t_i} \mathop = \sum\nolimits_{\ell = 1}^t {{A_{i,\ell}{\left( \mathcal{E}_{\ell-1} \right)}}}$, $i=1,2,\cdots,k$, where we suppress the $\mathcal{E}_t$ of $t_i$ for notation simplicity. We can also consider the selection policy at the end after allocating all simulation budget \citep{peng2014dynamic}. For simplicity, we fix the selection policy as to select alternatives with top-$m$ posterior means, i.e., $\widehat{\mathcal{F}}_T^m$.

Let ${R_T} ( \theta; {\widehat{\mathcal{F}}_T^m})\mathop = \limits^\Delta \mathds{1} ( {\widehat{\mathcal{F}}_T^m} = {\mathcal{F}^m} )$ be the reward for correctly determining the top-$m$ alternatives. Under the Bayesian setting, we can recursively define the expected payoff for a sequential sampling policy $\mathcal{A}_T$ in the stochastic dynamic programming problem by
\begin{align}\label{valuefunction}
{V_T}\left( {{\mathcal{E} _T};{\mathcal{A}_T}} \right) &\mathop  = \limits^\Delta \mathbb{E} [ {{R_T} ( \theta; {\widehat{\mathcal{F}}_T^m} )|{\mathcal{E} _T}} ] =  \int_{\theta \in \Theta}  {{R_T}( \theta; {\widehat{\mathcal{F}}_T^m} ) F\left( {d\theta |{\mathcal{E} _T}} \right)}  \nonumber \\
& = \left. \Pr \left\{ {{\bigcap\limits_{i = 1}^m {\bigcap\limits_{j = m + 1}^k {\left( {{\mu _{{{\left\langle i \right\rangle }_T}}} > {\mu _{{{\left\langle j \right\rangle }_T}}}} \right)} } } } \right|\mathcal{E}_T \right\}  \nonumber ~,
\end{align}
which is posterior integrated PCS (IPCS), since it is the weighted average of the PCS under different parameters values, and for $0 \le t < T$,
$$\begin{aligned}
{V_t}\left( {{{\mathcal{E}} _t};\mathcal{A}_T} \right) &\mathop  = \limits^\Delta \mathbb{E}{\left. {\left[ {\left. {{V_{t + 1}}\left( {{\mathcal{E} _t} \cup \left\{ {{X_{i,{t_{i + 1}}}}} \right\};{\mathcal{A}_T\left(\cdot\right)}} \right)} \right|{\mathcal{E} _t}} \right]} \right|_{i = {A_{t + 1}}\left( {{\mathcal{E} _t}} \right)}} \\
&= \int_{{\mathcal{X}_i}} {{{\left. {{V_{t + 1}}\left( {{\mathcal{E} _t} \cup \left\{ {{x_{i,{t_i} + 1}}} \right\};\mathcal{A}_T\left(\cdot\right)} \right){Q_i}\left( {\left. {d{x_{i,{t_i} + 1}}} \right|{\mathcal{E} _t}} \right)} \right|}_{i = {A_{t + 1}}\left( {{\mathcal{E} _t}} \right)}}}~,
\end{aligned}$$
where ${X_{i,{t_i} + 1}}$ is the $\left({t_i} + 1\right)$-th simulation replication for the allocated alternative $i$, $i = A_{t+1}\left(\mathcal{E}_t\right)$, and $\mathcal{X}_i$ is the support of ${X_{i,{t_i} + 1}}$. Then an optimal sequential sampling policy can be defined as the solution of the stochastic dynamic programming problem:
\begin{equation}\label{SDPsolution}
\mathcal{A}_T^*\mathop  = \limits^\Delta \arg \mathop {\max }\limits_{\mathcal{A}_T, \mathcal{S}} {V_0}\left( {{\zeta _0};\mathcal{A}_T} \right)~.
\end{equation}

Such a stochastic dynamic programming problem can be viewed as a MDP with $\left(T+1\right)$ stages, where $\zeta_0$ is the state at stage 0, and the sample information set $\mathcal{E}_{t}$ is the state at stage $t$, $0<t \le T$. The action is $A_{t+1}$ with the transition ${\mathcal{E} _t} \to {\mathcal{E} _{t + 1}} \mathop  = \limits^\Delta {\left. {\left\{ \mathcal{E}_{t} \cup {\left\{X_{i,{t_i} + 1}\right\}} \right\}} \right|_{i = {A_{t + 1}}}}$, $0 \le t < T$, where ${X_{i,{t_i} + 1}} \sim {Q_i}\left( {\left.  \cdot  \right|{\mathcal{E}_t}} \right)$, $i = A_{t+1}$, and final reward ${V_T}\left( {{{\mathcal{E}} _T};\mathcal{A}_T} \right)$. Note that the reward for the MDP only occurs at final step $T$ and the dimension of state space $\mathcal{E}_t$ grows as the step grows. Then we can recursively compute the optimal sampling policy in (\ref{SDPsolution}) by solving the following Bellman equation:
\begin{equation*}
V_T\left(\mathcal{E}_T\right) \mathop  = \limits^\Delta \mathbb{E}[ {{R_T} ( \theta; {\widehat{\mathcal{F}}_T^m} )|{\mathcal{E} _T}} ]~,
\end{equation*}
and
for $0 \le t < T$,
\begin{equation}\label{BM2}
{V_t}\left( {{\mathcal{E} _t}} \right) \mathop  = \limits^\Delta {\left. {\mathbb{E}\left[ {\left. {{V_{t + 1}}\left( {{\mathcal{E} _t} \cup {\left\{X_{i,{t_i} + 1}\right\} }} \right)} \right|{\mathcal{E} _t}} \right]} \right|_{i = A_{t + 1}^*\left( {{\mathcal{E}_t}} \right)}}~,
\end{equation}
where $A_{t + 1}^*\left( {{\mathcal{E} _t}} \right) = \arg \mathop {\max }\limits_{i = 1,2, \cdots ,k} \mathbb{E}\left[ {\left. {{V_{t + 1}}\left( {{\mathcal{E}_t}\cup {\left\{ {{X_{i,{t_i} + 1}}} \right\}} } \right)} \right|{\mathcal{E}_t}} \right]$.

The optimality condition of a sequential sampling policy for selecting top-$m$ alternatives is governed by the Bellman equation of the MDP. None of the previous sequential sampling procedures under the Bayesian framework for selecting the top-$m$ alternatives are derived from this MDP perspective.

\subsection{Normal Sampling Distributions}

Suppose $Q_i\left(\cdot;\theta_i\right)$ is normal with unknown mean $\mu_i$ and known variance $\sigma_i^2$, that is, ${X_{i,t}}\mathop \sim \limits^{i.i.d.} N\left( {{\mu _i},\sigma _i^2} \right)$, $i=1,2,\cdots,k$. The assumption is commonly used in the R\&S problem, since simulation outputs are obtained from average performances or batched means such that a central limit theorem holds~\citep{chen2011stochastic}. Under the Bayesian setting, we can obtain the posterior estimates of the unknown performances $\mu_i$, $i=1,2,\cdots,k$. We suppose a conjugate prior for each $\mu_i$, which is also a normal distribution $N ( {\mu _i^{\left( 0 \right)},{( {\sigma _i^{\left( 0 \right)}} )^2}} )$. By conjugacy~\citep{degroot2005optimal}, the posterior distribution of $\mu_i$ is $N( {\mu _i^{\left( t \right)},{( {\sigma _i^{\left( t \right)}} )^2}} )$, $i=1,2,\cdots,k$ with posterior variance
\begin{equation}\label{pv}
{( {\sigma _i^{\left( t \right)}} )^2}{\rm{ = }}{\left( {\frac{1}{{{( {\sigma _i^{\left( 0 \right)}} )^2}}} + \frac{{{t_i}}}{{\sigma _i^2}}} \right)^{ - 1}}~,
\end{equation}
and posterior mean
\begin{equation}\label{pm}
\mu _i^{\left( t \right)} = {( {\sigma _i^{\left( t \right)}} )^2}\left( {\frac{{\mu _i^{\left( 0 \right)}}}{{{( {\sigma _i^{\left( 0 \right)}} )^2}}} + \frac{{{t_i}{{\bar{X}_i}\left(t_i\right)}}}{{\sigma _i^2}}} \right)~.
\end{equation}

Note that if $\sigma _i^{\left( 0 \right)} \to \infty$, then $\mu _i^{\left( t \right)} ={\bar{X}_i}\left(t_i\right)$ and ${({\sigma _i^{\left( t \right)}})^2} = {{\sigma _i^2} \mathord{\left/
 {\vphantom {{\sigma _i^2} {{t_i}}}} \right.\kern-\nulldelimiterspace} {{t_i}}}$, $i=1,2,\cdots,k$, which only includes sample information without any hyper-parameters. Such case is called uninformative. Under the conjugate prior, the sample information set $\mathcal{E}_t$ can be completely determined by the posterior parameters $\mu _i^{\left( t \right)}$ and ${( {\sigma _i^{\left( t \right)}} )^2}$, i.e., $\mathcal{E}_t = \{\mu _1^{\left( t \right)},\mu _2^{\left( t \right)},\cdots,\mu _k^{\left( t \right)},{( {\sigma _1^{\left( t \right)}} )^2},{( {\sigma _2^{\left( t \right)}} )^2,\cdots,{( {\sigma _k^{\left( t \right)}} )^2},\mathcal{E}_0 } \}$, and then the dimension of $\mathcal{E}_t$ is fixed at each step $t$. In addition, the predictive distribution of $X_{i,{t_i+1}}$ is $N ( {\mu _i^{\left( t \right)},\sigma _i^2 + {( {\sigma _i^{\left( t \right)}})^2}} )$. In practice, $\sigma_i^2$, $i=1,2,\cdots,k$ are usually unknown and can be estimated by their sample variances, i.e., $\widehat{\sigma}_i^2 = \sum\nolimits_{\ell = 1}^{{t_i}} {{{( {{X_{i,\ell}} - {\bar{X}_i}\left(t_i\right)} )^2} / {\left({t_{i}}-1\right)}}}$.  For a normal distribution with unknown variance, there is a normal-gamma conjugate prior~\citep{degroot2005optimal}.

\section{Efficient Dynamic Allocation Policy}\label{SectionAOA}

Since the dimension of the state space is high, using backward induction to solve the Bellman equation  (\ref{BM2}) is computationally intractable. Under normal sampling distributions, we adopt an approximate dynamic programming (ADP) paradigm to address this difficulty. The ADP approach makes dynamic decisions and keeps learning the value function approximation with decisions moving forward.

We suppose any step $t$ could be the last step. We use a single feature of the value function one-step look ahead for allocating the next simulation replication and certainty equivalence to approximate the value function. Specifically, conditioned on $\mathcal{E}_t$, $\mu_i$ follows a normal distribution with mean $\mu _i^{\left( t \right)}$ and variance $( {\sigma _i^{\left( t \right)}})^2$, $i=1,2, \cdots, k$. The joint distribution of vector
$$( {{\mu _{{{\left\langle 1 \right\rangle }_t}}} - {\mu _{{{\left\langle {m + 1} \right\rangle }_t}}}, \cdots ,{\mu _{{{\left\langle 1 \right\rangle }_t}}} - {\mu _{{{\left\langle k \right\rangle }_t}}}, \cdots ,{\mu _{{{\left\langle m \right\rangle }_t}}} - {\mu _{{{\left\langle {m + 1} \right\rangle }_t}}}, \cdots ,{\mu _{{{\left\langle m \right\rangle }_t}}} - {\mu _{{{\left\langle k \right\rangle }_t}}}} )~,$$
follows a joint normal distribution with mean vector
\begin{equation}\label{jointvector}
\left( {\mu _{_{{{\left\langle 1 \right\rangle }_t}}}^{\left( t \right)} - \mu _{{{\left\langle {m + 1} \right\rangle }_t}}^{\left( t \right)}, \cdots , \mu _{_{{{\left\langle 1 \right\rangle }_t}}}^{\left( t \right)} - \mu _{{{\left\langle k \right\rangle }_t}}^{\left( t \right)}, \cdots ,\mu _{_{{{\left\langle m \right\rangle }_t}}}^{\left( t \right)} - \mu _{{{\left\langle {m + 1} \right\rangle }_t}}^{\left( t \right)}, \cdots ,\mu _{{{\left\langle m \right\rangle }_t}}^{\left( t \right)} - \mu _{{{\left\langle k \right\rangle }_t}}^{\left( t \right)}} \right)~,
\end{equation}
and covariance matrix $\Sigma =\Gamma^{\prime} \Lambda \Gamma$, where $^\prime$ denotes the transpose operation of the matrix, $\Lambda \in \mathbb{R}^{k \times k }$, $\Lambda \mathop  = \limits^\Delta diag ({{( {\sigma _{{{\left\langle 1 \right\rangle }_t}}^{\left( t \right)}} )^2}, \cdots ,{( {\sigma _{{{\left\langle m \right\rangle }_t}}^{\left( t \right)}})^2},{( {\sigma _{{{\left\langle {m + 1} \right\rangle }_t}}^{\left( t \right)}})^2},\cdots,{( {\sigma _{{{\left\langle k \right\rangle }_t}}^{\left( t \right)}})^2}} )$ is a diagonal matrix, and $\Gamma \in \mathbb{R}^{k \times \left[ m \times \left( {k - m} \right)\right]}$ is a block matrix,
$$\Gamma \mathop  = \limits^\Delta {\left( {\begin{array}{*{20}{c}}
{{R_1}}&{{R_2}}& \cdots &{{R_m}}\\
{{S_1}}&{{S_2}}& \cdots &{{S_m}}
\end{array}} \right)_{k \times \left[ m \times \left( {k - m} \right)\right] }}~,$$
where $R_i \in \mathbb{R}^{m \times \left( {k - m} \right)}$, $i=1,2, \cdots, m$ is a matrix in which $i$-th row is 1 and the remaining rows are 0, and $S_i \in \mathbb{R}^{\left(k-m\right) \times \left(k-m\right)}$, $i=1,2,\cdots,m$ is a diagonal matrix, i.e.,
$${R_i} \mathop  = \limits^\Delta \left( {\begin{array}{*{20}{c}}
0&0& \cdots &0\\
0&0& \cdots &0\\
 \vdots & \vdots & \ddots & \vdots \\
1&1& \cdots &1\\
 \vdots & \vdots & \ddots & \vdots \\
0&0& \cdots &0
\end{array}} \right)\begin{array}{*{20}{c}}
{}\\
{ \to i\rm{th}\;row}~,
\end{array}\qquad {S_i} \mathop  = \limits^\Delta {\left( {\begin{array}{*{20}{c}}
{ - 1}&0& \cdots &0\\
0&{ - 1}& \cdots &0\\
 \vdots & \vdots & \ddots & \vdots \\
0&0& \cdots &{ - 1}
\end{array}} \right)_{\left( {k - m} \right) \times \left( {k - m} \right)}}~.$$

\begin{lemma}\label{lemVFAintegral}
The value function $\Pr \{ {{\bigcap\nolimits_{i = 1}^m {\bigcap\nolimits_{j = m + 1}^k {\left( {{\mu _{{{\left\langle i \right\rangle }_t}}} > {\mu _{{{\left\langle j \right\rangle }_t}}}} \right)} } } } |\mathcal{E}_t \}$ for selecting the alternatives in $\widehat{\mathcal{F}}_t^m$ is an integral of the density of a $\left[ m \times \left( {k - m} \right) \right]$-dimensional standard normal distribution over an area covered by hyperplanes ${\sum\nolimits_{\ell = 1}^q {u_{\ell, q} z_{\ell}} = \mu_{\langle j\rangle_t}^{(t)}-\mu_{\langle i\rangle_t}^{(t)}}$, $q=1,2,\cdots,m \times (k-m)$, where $u_{ij}$ are elements of an upper triangular matrix $U \in \mathbb{R}^{\left[ m \times \left( {k - m} \right) \right] \times \left[ m \times \left( {k - m} \right) \right] }$, i.e., ${u_{ij}} = 0$ if $i>j$, $i,j = 1,2,\cdots,m\times\left( {k - m} \right)$, and $z_{\ell} \stackrel{i . i . d .}{\sim} N(0,1)$, $\ell=1,2, \cdots, {m \times \left( {k - m} \right)}$.
\end{lemma}

\begin{figure}[htbp]
 \centering
 \includegraphics[width=0.6\textwidth]{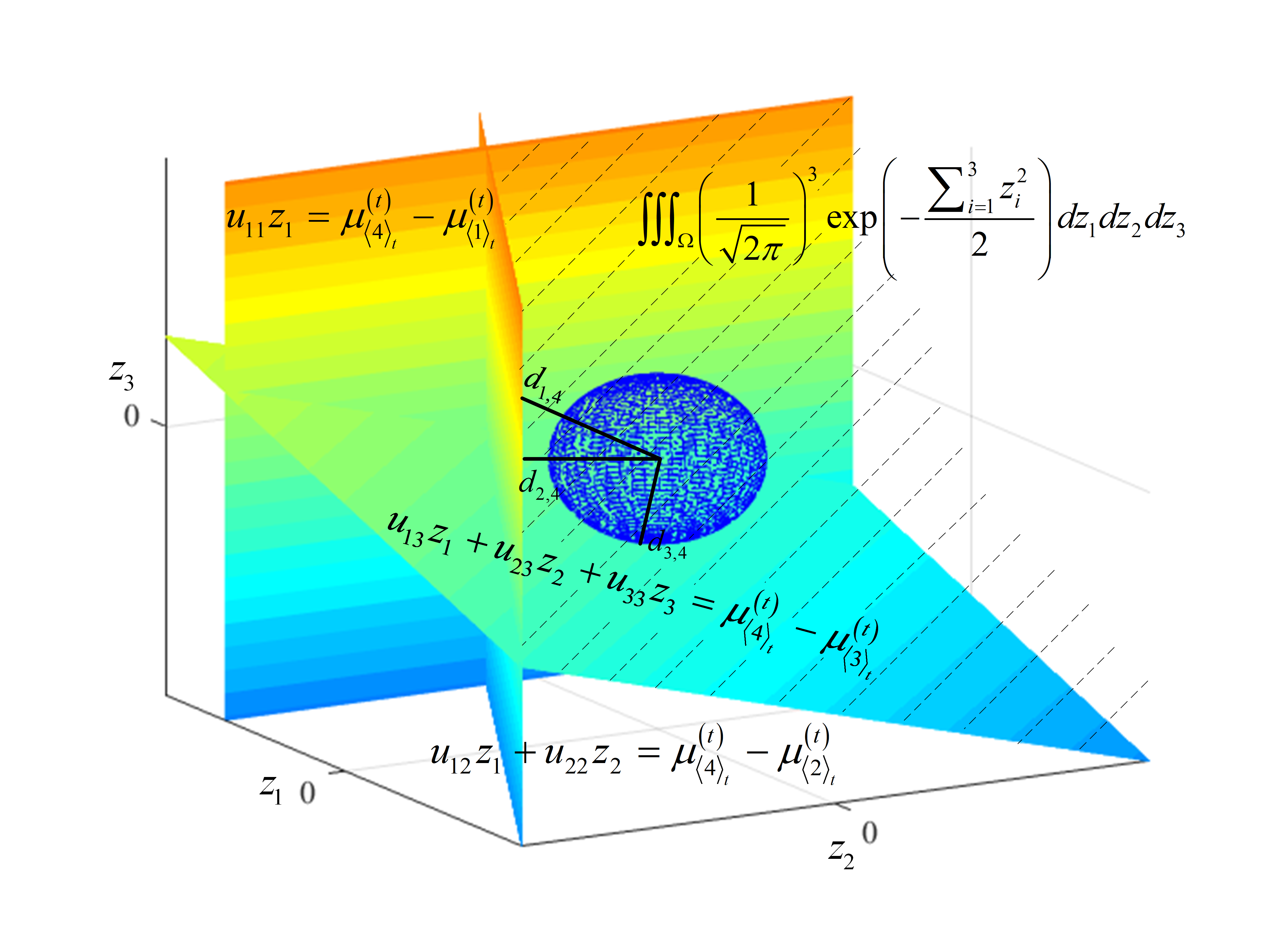}
 \caption{
 Area of integration for approximation is the ball, where dominant values of integrand $\exp \left(  - \left({{z_1^2 + z_2^2 + z_3^2}}\right)/{2}\right)$ are captured.}
 \label{VFAfig}
\end{figure}

We adopt a technique developed in \citet{peng2018ranking} to simplify the value function by an integral over a largest internally tangent ball in the integral region encompassed by the hyperplanes, since the standard normal distribution decays at an exponential rate with the distance from the origin. Then maximizing the integral over the largest inscribed ball is equivalent to maximizing the volume of the ball due to symmetry of normal distribution. The approximation error decreases to 0 at an exponential rate as the radius of the largest inscribed ball goes to infinity \citep{zhang2021efficient}. Let $d\left( {{\mathcal{E} _t}} \right) \mathop  = \limits^\Delta \mathop {\min }\limits_{i = 1,2, \cdots ,m \atop j = m + 1,m + 2, \cdots ,k} {d_{ij}}\left( {{\mathcal{E} _t}} \right)$ be the radius of the largest internally tangent ball of the integration region, where
$${d_{ij}}\left( {{\mathcal{E} _{t}}} \right) = \frac{{\mu _{{{\left\langle i \right\rangle }_{t}}}^{\left( {t} \right)} - \mu _{{{\left\langle j \right\rangle }_{t}}}^{\left( {t} \right)}}}{{\sqrt {{\left( {\sigma _{{{\left\langle i \right\rangle }_{t}}}^{\left( {t} \right)}} \right)^2} + {\left( {\sigma _{{{\left\langle j \right\rangle }_{t}}}^{\left( {t} \right)}} \right)^2}} }},\quad i\in\left\{1,2,\cdots,m\right\},\;j \in \left\{m+1,m+2,\cdots,k\right\}~,$$
is the distance from the origin to each hyperplane. Then an analytical VFA becomes ${\widetilde V_{t}}\left( {{\mathcal{E} _{t}}} \right) \mathop  = \limits^\Delta {d^2}\left( {{\mathcal{E} _{t}}} \right)$. Note that as $t \to \infty$, if every alternative is sampled infinitely often, then by the law of large numbers (LLN), we have $\mathop {\lim }\limits_{t \to \infty } \mu _{{{\left\langle \ell \right\rangle }_t }}^{\left(t\right)}  = {\mu _{\left\langle \ell \right\rangle }},\;a.s.$, $\mathop {\lim }\limits_{t \to \infty } \sigma _{{{\left\langle \ell \right\rangle }_{t} }}^{\left(t\right)}  = 0,\;a.s.$, $\ell=1,2,\cdots,k$, and $\mathop {\lim }\limits_{t \to \infty } {d_{ij}}\left( {{\mathcal{E} _t}} \right) = \infty,\;a.s.$, $i=1,2,\cdots,m$, $j=m+1,m+2,\cdots,k$. We visualize the approximation for  $k=4$ and $m=3$. In Figure~\ref{VFAfig}, the value function is the integration of the three-dimensional standard normal density over the shadowed area, and an approximation is the integration over the largest internally tangent ball with radius $d_{3,4}$.

If the $(t+1)$-th simulation replication is the last one, an approximation looking one step ahead at step $t$ by allocating alternative $i$, $i=1,2,\cdots,k$ can be given as ${\widetilde V_t}\left( {{\mathcal{E} _t};i} \right) \mathop  = \limits^\Delta \mathbb{E} [ { {{{\widetilde V}_{t + 1}}\left( {{\mathcal{E} _t} \cup \{X_{i,{t_i} + 1}}\} \right)} {|\mathcal{E} _t}} ]$. Since the expectation is computationally expensive to calculate, we use the certainty equivalence in~\citet{bertsekas2005dynamic} to further approximate the $\mathbb{E}\left(\cdot\right)$, i.e., replacing stochastic quantities by their expected values, which can be expressed as $\mathbb{E} [ { {{{\widetilde V}_{t + 1}}\left( {{\mathcal{E} _t} \cup \{X_{i,{t_i} + 1}}\} \right)} {|\mathcal{E} _t}} ]  \approx \widetilde{V}_{t+1}\left(\mathcal{E}_{t} \cup \mathbb{E}\left[X_{i, t_{i}+1} | \mathcal{E}_{t}\right]\right)$. Let $\widehat{V}_{t}\left(\mathcal{E}_{t} ; i\right) \mathop  = \limits^\Delta  \widetilde{V}_{t+1}\left(\mathcal{E}_{t} \cup \mathbb{E}\left[X_{i, t_{i}+1} | \mathcal{E}_{t}\right]\right)$, $i=1,2,\cdots,k$. To derive an analytical $\widehat{V}_{t}\left(\cdot\right)$ from $\widetilde{V}_{t}\left(\mathcal{E}_t\right)$, when obtaining one additional simulation replication, we note that $\mu _i^{\left( t \right)}$ remains unchanged since ${X_{i,{t_i} + 1}}$ takes a value of $\mathbb{E}\left[X_{i, t_{i}+1} | \mathcal{E}_{t}\right]$, and only ${( {\sigma _i^{\left( t \right)}} )^2}$ is updated. Then an analytical form of a one-step look ahead dynamic allocation policy can be expressed as: for $i,h = 1,2, \cdots ,m$ and $j,\ell = m+1,m+2, \cdots ,k$,
\begin{equation}\label{equ21}
\widehat{V}_{t}\left(\mathcal{E}_{t} ; i \right) = \min \left\{ {\mathop {\min }\limits_j \frac{{{( {\mu _{{{\left\langle i \right\rangle }_t}}^{\left( t \right)} - \mu _{{{\left\langle j \right\rangle }_t}}^{\left( t \right)}})^2}}}{{{({\sigma _{{{\left\langle i \right\rangle }_t}}^{\left( {t + 1} \right)}})^2} + {( {\sigma _{{{\left\langle j \right\rangle }_t}}^{\left( t \right)}})^2}}},{\mathop {\min }\limits_{\scriptstyle h \ne i\atop \scriptstyle j}} \frac{{{( {\mu _{{{\left\langle h \right\rangle }_t}}^{\left( t \right)} - \mu _{{{\left\langle j \right\rangle }_t}}^{\left( t \right)}})^2}}}{{{( {\sigma _{{{\left\langle h \right\rangle }_t}}^{\left( t \right)}})^2} + {( {\sigma _{{{\left\langle j \right\rangle }_t}}^{\left( t \right)}})^2}}}} \right\}~,
\end{equation}
and
\begin{equation}\label{equ22}
\widehat{V}_{t}\left(\mathcal{E}_{t} ; j \right) = \min \left\{ {\mathop {\min }\limits_i \frac{{{( {\mu _{{{\left\langle i \right\rangle }_t}}^{\left( t \right)} - \mu _{{{\left\langle j \right\rangle }_t}}^{\left( t \right)}} )^2}}}{{{( {\sigma _{{{\left\langle i \right\rangle }_t}}^{\left( {t} \right)}} )^2} + {( {\sigma _{{{\left\langle j \right\rangle }_t}}^{\left( {t + 1} \right)}} )^2}}},{\mathop {\min }\limits_{\scriptstyle \ell \ne j\atop \scriptstyle i}} \frac{{{( {\mu _{{{\left\langle i \right\rangle }_t}}^{\left( t \right)} - \mu _{{{\left\langle \ell \right\rangle }_t}}^{\left( t \right)}} )^2}}}{{{( {\sigma _{{{\left\langle i \right\rangle }_t}}^{\left( t \right)}} )^2} + {( {\sigma _{{{\left\langle \ell \right\rangle }_t}}^{\left( t \right)}} )^2}}}} \right\}~,
\end{equation}
where
$$( \sigma _{{\widetilde \ell}}^{\left( t+1 \right)} )^2 = {\left( {\frac{1}{{{( {\sigma _{{\widetilde \ell}}^{( 0 )}} )^2}}} + \frac{t_{{\widetilde \ell}}+1}{\sigma_{{{{\widetilde \ell}}}}^2}} \right)^{ - 1}}, \quad {\widetilde \ell}=1,2,\cdots,k~.$$

An asymptotically optimal allocation policy for selecting top-$m$ alternatives (AOAm) that optimizes the VFA is given by
\begin{equation}\label{AOAm}
{\widehat{A}_{t + 1}}\left( {{\mathcal{E} _t}} \right) = \arg \mathop {\max }\limits_{i = 1,2, \cdots ,k} \widehat {{V_t}}\left( {{\mathcal{E} _t};i} \right)~.
\end{equation}

The Algorithm 1 in online appendix A.1 presents the proposed sequential AOAm procedure. Based on the posterior information, the proposed sequential AOAm procedure uses the posterior means and posterior variances of the estimated unknown performances, which are calculated by Bayes' rule. Notice that ${( {\mu _{{{\left\langle i \right\rangle }_t}}^{\left( t \right)} - \mu _{{{\left\langle j \right\rangle }_t}}^{\left( t \right)}} )^2}$ and ${{{( {\sigma _{{{\left\langle i \right\rangle }_t}}^{\left( t \right)}} )^2} + {( {\sigma _{{{\left\langle j \right\rangle }_t}}^{\left( t \right)}} )^2}}}$ are the squared mean and variance of the posterior distribution of the differences in performances of alternatives $\left\langle i \right\rangle_t$ and $\left\langle j \right\rangle_t$, respectively, $i=1,2,\cdots,m$, $j=m+1,m+2,\cdots,k$. The VFA can be rewritten as $\mathop {\min }\limits_{\scriptstyle i,j} {1} / {{c_v^2\left( {i,j} \right)}}$, where $c_v\left( {i,j} \right)$ is the  posterior noise-signal ratio (coefficient of variation) of ${\mu _{{{\left\langle i \right\rangle }_t}}} - {\mu _{{{\left\langle j \right\rangle }_t}}}$. The AOAm (\ref{AOAm}) minimizes the maximum noise-signal ratio. The larger is the value of $c_v\left( {i,j} \right)$, the higher the difficulty in comparing ${\mu _{{{\left\langle i \right\rangle }_t}}}$ and ${\mu _{{{\left\langle j \right\rangle }_t}}}$ from posterior information. At each step, the AOAm sequentially allocates a replication to an alternative to reduce the noise-signal ratio of the differences in performances of the pair of alternatives most difficult in comparison among all $\left(m \times \left(k - m \right) \right)$ pairs based on posterior information. In particular, if $m=1$, (\ref{equ21}) and (\ref{equ22}) lead to, for $j,\ell=2,3,\cdots,k$, $j \ne \ell$,
$$\widehat{V}_{t}\left(\mathcal{E}_{t} ; {\left\langle 1 \right\rangle}_t \right) = {\mathop {\min }\limits_{j=2,3,\cdots,k} \frac{{{( {\mu _{{{\left\langle 1 \right\rangle }_t}}^{\left( t \right)} - \mu _{{{\left\langle j \right\rangle }_t}}^{\left( t \right)}})^2}}}{{{({\sigma _{{{\left\langle 1 \right\rangle }_t}}^{\left( {t + 1} \right)}})^2} + {( {\sigma _{{{\left\langle j \right\rangle }_t}}^{\left( t \right)}})^2}}}}~,$$
and
$$\widehat{V}_{t}\left(\mathcal{E}_{t} ; {\left\langle j \right\rangle}_t \right) = \min \left\{ {\mathop \frac{{{( {\mu _{{{\left\langle 1 \right\rangle }_t}}^{\left( t \right)} - \mu _{{{\left\langle j \right\rangle }_t}}^{\left( t \right)}} )^2}}}{{{( {\sigma _{{{\left\langle 1 \right\rangle }_t}}^{\left( {t} \right)}} )^2} + {( {\sigma _{{{\left\langle j \right\rangle }_t}}^{\left( {t + 1} \right)}} )^2}}},{\mathop {\min }\limits_{\scriptstyle \ell \ne j}} \frac{{{( {\mu _{{{\left\langle 1 \right\rangle }_t}}^{\left( t \right)} - \mu _{{{\left\langle \ell \right\rangle }_t}}^{\left( t \right)}} )^2}}}{{{( {\sigma _{{{\left\langle 1 \right\rangle }_t}}^{\left( t \right)}} )^2} + {( {\sigma _{{{\left\langle \ell \right\rangle }_t}}^{\left( t \right)}} )^2}}}} \right\}~,$$
which reduce to the allocation policy for selecting the best alternative derived in~\citet{peng2018ranking}.

The proposed sequential AOAm procedure is proved to be consistent in the following theorem, i.e., as $t \to \infty$, the top-$m$ estimated set $\widehat {\mathcal{F}}_t^{m}$ will be the true top-$m$ set ${\mathcal{F}^m}$.

\begin{theorem}\label{thmconsistent}
The proposed asymptotically optimal allocation policy for selecting top-$m$ alternatives (AOAm) (\ref{AOAm}) is consistent, i.e.,
$$\mathop {\lim }\limits_{t \to \infty } \widehat {\mathcal{F}}_t^{m} = {\mathcal{F}^m},\quad a.s.~.$$
\end{theorem}

In the proof of Theorem \ref{thmconsistent}, we show that every alternative will be sampled infinitely often as the simulation budget goes to infinity. In the next section, we show that the AOAm is not only consistent but also achieves the asymptotically optimal sampling ratios.

\section{Asymptotic Optimality}\label{Sectionopt}

Although \citet{glynn2004large} derive the asymptotically optimal sampling ratios for selecting the best alternative, no previous work rigorously establishes the asymptotically optimal sampling ratios for selecting top-$m$ alternatives under the static optimization problem (\ref{staticopt}), which will be established in the rest of the section.

\subsection{Asymptotically Optimal Sampling Ratios}

Under the static optimization setting, the PFS can be expressed as:
$$\begin{aligned}
{\rm {PFS}} & = 1 - \Pr \left\{ \left. \bigcap\limits_{i = 1}^m {\bigcap\limits_{j = m + 1}^k {\left( {\bar{X}_{\langle i\rangle}\left(r_{\langle i\rangle} T\right)>\bar{X}_{\langle j\rangle}\left(r_{\langle j \rangle} T\right)} \right)} }  \right|\theta \right\} \\
& = \Pr \left\{ \left. {\bigcup\limits_{i = 1}^m {\bigcup\limits_{j = m + 1}^k {\left( {\bar{X}_{\langle i\rangle}\left(r_{\langle i\rangle} T\right) \le \bar{X}_{\langle j\rangle}\left(r_{\langle j \rangle} T\right)} \right)} } } \right|\theta \right\}~,
\end{aligned}$$
which is lower bounded by
$$\mathop {\max }\limits_{\scriptstyle i = 1,2, \cdots ,m\atop
\scriptstyle j = m + 1,m+2, \cdots ,k} \Pr \left\{ \left. {\bar{X}_{\langle i\rangle}\left(r_{\langle i\rangle} T\right) \le \bar{X}_{\langle j\rangle}\left(r_{\langle j \rangle} T\right)}  \right|\theta \right\}~,$$
and is upper bounded by
$$m\times\left( {k - m} \right) \cdot \mathop {\max }\limits_{\scriptstyle i = 1,2, \cdots ,m\atop
\scriptstyle j = m + 1,m+2, \cdots ,k} \Pr \left\{ \left. {\bar{X}_{\langle i\rangle}\left(r_{\langle i\rangle} T\right) \le \bar{X}_{\langle j\rangle}\left(r_{\langle j \rangle} T\right)} \right|\theta \right\}~.$$

We define the cumulant generating function of sample mean $\bar X_i$ to be ${\Lambda_i}\left( \lambda  \right) = {\log M_i}\left( \lambda  \right)$, where $M_i\left( \lambda  \right)= \mathbb{E} ( {e^{\lambda {\bar X_i}}} )$, $i=1,2,\cdots,k$ are exponential moment generating functions. Denote $\Lambda_i^*\left(  x  \right)$ as the Fenchel-Legendre transform of $\Lambda_i\left(  \lambda  \right)$, i.e., $\Lambda_i^*\left(  x  \right) = \mathop {\sup }_{\lambda  \in \mathbb{R}} \left\{ {\lambda x - {\Lambda_i}\left( \lambda  \right)} \right\}$. Some standard assumptions in large deviations contexts \citep{zeitouni2010large} can be found in online appendix A.2.

For $i \in \left\{1,2,\cdots,m\right\}$, $j \in \left\{m+1,m+2,\cdots,k\right\}$ and $r_{\ell} > 0$, $\ell = 1,2,\cdots,k$, if there exists a rate function ${G_{ij}}\left( { \cdot , \cdot } \right)$, such that
$$\mathop {\lim }\limits_{T \to \infty } \frac{1}{T}\log \Pr \left\{ \left. {\bar{X}_{\langle j\rangle}\left(r_{\langle j \rangle} T\right) \ge \bar{X}_{\langle i\rangle}\left(r_{\langle i\rangle} T\right) } \right|\theta \right\} =  -{G_{ij}}\left( r_{\left\langle i \right\rangle }, r_{\left\langle j \right\rangle } \right)~,$$
and then we have
$$\begin{aligned}
& \mathop {\lim }\limits_{T \to \infty } \frac{1}{T}\log \left( {\mathop {\max }\limits_{\scriptstyle i = 1,2, \cdots ,m\atop
\scriptstyle j = m + 1,m+2, \cdots ,k}  \Pr \left\{ \left. {\bar{X}_{\langle j\rangle}\left(r_{\langle j \rangle} T\right) \ge \bar{X}_{\langle i\rangle}\left(r_{\langle i\rangle} T\right) } \right|\theta \right\} } \right) \\
 = & \mathop {\max }\limits_{\scriptstyle i = 1,2, \cdots ,m\atop \scriptstyle j = m + 1,m+2, \cdots ,k} \mathop {\lim }\limits_{T \to \infty } \frac{1}{T} {\log  \Pr \left\{ \left. {\bar{X}_{\langle j\rangle}\left(r_{\langle j \rangle} T\right) \ge \bar{X}_{\langle i\rangle}\left(r_{\langle i\rangle} T\right) } \right|\theta \right\} }  \\
 = & - \mathop {\min }\limits_{\scriptstyle i = 1,2, \cdots ,m\atop
\scriptstyle j = m + 1,m+2, \cdots ,k} {G_{ij}}\left( r_{\left\langle i \right\rangle }, r_{\left\langle j \right\rangle } \right)~,
\end{aligned}$$
where the first equality holds from the monotonicity of function $y = \log \left( x \right)$. Let $r \mathop  = \limits^\Delta \left(r_{\left\langle 1 \right\rangle },r_{\left\langle 2 \right\rangle }, \cdots, r_{\left\langle k \right\rangle }\right)$ and $\widetilde G\left(r\right)\mathop  = \limits^\Delta  \mathop {\min }\limits_{\scriptstyle i = 1,2, \cdots ,m\atop
\scriptstyle j = m + 1,m+2, \cdots ,k} {G_{ij}}\left( r_{\left\langle i \right\rangle} , r_{\left\langle j \right\rangle } \right)$. Then
$$\mathop {\lim }\limits_{T \to \infty } \frac{1}{T}\log {\rm{PFS}} = \mathop {\lim }\limits_{T \to \infty } \frac{1}{T}\log \left( {\mathop {\max }\limits_{\scriptstyle i = 1,2, \cdots ,m\atop
\scriptstyle j = m + 1,m+2, \cdots ,k} \Pr \left\{ \left. {\bar{X}_{\langle j\rangle}\left(r_{\langle j \rangle} T\right) \ge \bar{X}_{\langle i\rangle}\left(r_{\langle i\rangle} T\right) } \right|\theta \right\} } \right) = - \widetilde G \left(r\right)~,$$
which implies that the PFS will decay exponentially with respect to $T$ at a rate function given by $\widetilde G \left(r\right)$. Theorem~\ref{thm1} states that the rate function $G_{ij}\left(\cdot,\cdot\right)$ exist and derives the rate function for PFS.

\begin{theorem}\label{thm1}
Under Assumptions 1-4 in online appendix A.2, the rate function of PFS for selecting top-$m$ alternatives is given by
\begin{align}
 - \mathop {\lim }\limits_{T \to \infty } \frac{1}{T} \rm{log} \rm{PFS} &= \mathop {\min }\limits_{\scriptstyle i = 1,2, \cdots ,m\atop\scriptstyle j = m + 1,m+2, \cdots ,k} {G_{ij}}\left( {{r_{\left\langle i\right\rangle}},{r_{\left\langle j\right\rangle}}} \right)\nonumber\\
& = \mathop {\min }\limits_{\scriptstyle i = 1,2, \cdots ,m\atop\scriptstyle j = m + 1,m+2, \cdots ,k} \mathop {\inf }\limits_x \left( r_{\left\langle i\right\rangle}{\Lambda _{\left\langle i\right\rangle}^*\left(  x  \right)} + r_{\left\langle j\right\rangle}{\Lambda _{\left\langle j\right\rangle}^*\left(  x  \right)} \right)\nonumber~.
\end{align}
\end{theorem}


The rate function of PFS depends on comparisons between pairs of alternatives. To solve the asymptotically optimal sampling ratios, an optimization problem that maximizes the rate at which PFS tends to 0 as a function of the sampling ratios can be expressed as
\begin{alignat}{2}\label{opt0}
\max\quad & \mathop {\min }\limits_{\scriptstyle i = 1,2, \cdots ,m\atop
\scriptstyle j = m + 1,m+2, \cdots ,k} {G_{ij}}\left( {{r_{\left\langle  i \right\rangle}},{r_{\left\langle  j \right\rangle}}} \right) \\
\mbox{s.t.}\quad
& \sum\limits_{\ell = 1}^k {r_{\left\langle  \ell \right\rangle}=1}, \quad r_{{\left\langle  \ell \right\rangle}} \ge 0,\quad {\ell= 1,2,\cdots k} \nonumber~.
\end{alignat}

Some properties of the rate function $G_{ij}\left(\cdot , \cdot\right)$ are shown in Lemma \ref{lem001}-\ref{Gfunction_property}.

\begin{lemma}\label{lem001}
The infimum point of $( r_{\left\langle  i \right\rangle}{\Lambda _{\left\langle i \right\rangle}^*\left(  x  \right)} + r_{\left\langle  j \right\rangle}{\Lambda _{\left\langle j \right\rangle}^*\left(  x  \right)} )$ is unique, $i \in \left\{1,2,\cdots,m\right\}$, $j \in \left\{m+1,m+2,\cdots,k\right\}$.
\end{lemma}


\begin{lemma}\label{lem}
For $\left( {r_{\left\langle  i \right\rangle},r_{\left\langle  j \right\rangle}} \right) > 0$, ${G_{ij}}\left( {{r_{\left\langle  i \right\rangle}},{r_{\left\langle  j \right\rangle}}} \right)$ is strictly increasing in $r_{\left\langle  i \right\rangle}$ and $r_{\left\langle  j \right\rangle}$, respectively. Also, ${G_{ij}}\left( {r_{\left\langle  i \right\rangle}},{r_{\left\langle  j \right\rangle}} \right) = 0$ if $\min \left( {{r_{\left\langle  i \right\rangle}},{r_{\left\langle  j \right\rangle}}} \right) = 0$.
\end{lemma}


\begin{lemma}\label{Gfunction_property}
$G_{ij}\left(r_{\left\langle  i \right\rangle},r_{\left\langle  j \right\rangle}\right)$ is a strictly concave and continuous function of $\left(r_{\left\langle  i \right\rangle},r_{\left\langle  j \right\rangle}\right)$.
\end{lemma}


Denote ${G'_{ij}}\left( r \right)={G_{ij}}\left( {{r_{\left\langle i \right\rangle }},{r_{\left\langle j \right\rangle }}} \right)$, $i=1,2,\cdots,m$, $j=m+1,m+2,\cdots,k$. Following Lemma \ref{Gfunction_property}, $G_{ij}\left(r_{\left\langle  i \right\rangle},r_{\left\langle  j \right\rangle}\right)$ is a concave function of $\left(r_{\left\langle  i \right\rangle},r_{\left\langle  j \right\rangle}\right)$, implying that ${G'_{ij}}\left(r\right)$ is a concave function of $r$. The objective function $\widetilde G\left(r\right)$ in (\ref{opt0}) is concave for $r \ge 0$ since the minimum of concave functions ${G'_{ij}}\left( r \right)$ is also concave. Therefore, (\ref{opt0}) is a concave maximization problem, which implies that it is a convex optimization problem and any locally optimal point is globally optimal (see~\citealp{boyd2009convex}, Page 137).

\begin{remark}
Note that unlike the uniqueness of the asymptotically optimal sampling ratios for selection of the best \citep{zhang2021efficient}, the optimal solution to optimization problem (\ref{opt0}) may not be unique. We caution that establishing this uniqueness result could be challenging in general, in contrast to
some existing work reporting the uniqueness of the asymptotically optimal sampling ratios for their problem. \citet{hunter2013optimal} study the problem of finding the best alternative under stochastic constraints, and report the uniqueness of asymptotic optimal sampling ratios by claiming that the objective function is strictly concave in the proof for Proposition 2 of \citet{hunter2013optimal}. However, $G_{1 \ell}^{\prime}\left(r\right)$, $\ell = 2,3,\cdots,k$ is not strictly concave for $r \ge 0$, even though $G_{1 \ell}\left(r_{\left\langle 1 \right\rangle},r_{\left\langle \ell \right\rangle}\right)$ is strictly concave in $\left(r_{\left\langle 1 \right\rangle},r_{\left\langle \ell \right\rangle}\right)$. \citet{xiao2017simulation} study simultaneously selecting the best and worst alternatives and report the uniqueness of asymptotic optimal sampling ratios by claiming strict concavity without giving a proof.
\end{remark}

Then we solve the optimization problem (\ref{opt0}). Equivalently, (\ref{opt0}) can be rewritten as the following optimization problem:
\begin{alignat}{2}\label{opt1}
\max\quad &z \\
\mbox{s.t.}\quad
&{G_{ij}}\left( {r_{\left\langle  i \right\rangle},r_{\left\langle  j \right\rangle}} \right) - z \ge 0,\quad{i = 1,2,\cdots ,m},\;{j = m + 1,m + 2,\cdots ,k}~,\nonumber\\
&\sum\limits_{\ell = 1}^k {r_{\left\langle \ell \right\rangle}=1},\quad r_{\left\langle \ell \right\rangle} \ge 0,\quad {\ell= 1,2,\cdots ,k} \nonumber~.
\end{alignat}

The optimization problem (\ref{opt1}) is a concave maximization problem with differentiable objective and constraint functions. In addition, there exists a point $r = \left( {{1 \mathord{\left/
 {\vphantom {1 k}} \right.
 \kern-\nulldelimiterspace} k},{1 \mathord{\left/
 {\vphantom {1 k}} \right.
 \kern-\nulldelimiterspace} k}, \cdots ,{1 \mathord{\left/
 {\vphantom {1 k}} \right.
 \kern-\nulldelimiterspace} k}} \right)$ and $z=0$ such that the inequality constraints strictly hold, which implies that Slater's condition holds for (\ref{opt1}), and then the strong duality holds (see~\citealp{boyd2009convex}, Page 226). Therefore, the Karush-Kuhn Tucker (KKT) conditions provide necessary and sufficient conditions for the global optimality of (\ref{opt1}) (see~\citealp{boyd2009convex}, Page 244). We define an optimization problem (\ref{opt2}) by forcing ${r_\ell}$, $\ell=1,2,\cdots,k$ to be strictly positive
\begin{alignat}{2}\label{opt2}
\max\quad &z \\
\mbox{s.t.}\quad
&\mathop {\min }\limits_{j = m + 1,m+2, \cdots ,k} {G_{ij}\left(r_{\left\langle  i \right\rangle},r_{\left\langle  j \right\rangle}\right)} = z,\quad i=1,2,\cdots,m~,\nonumber\\
&\mathop {\min }\limits_{i = 1,2, \cdots,m} {G_{ij}\left(r_{\left\langle  i \right\rangle},r_{\left\langle  j \right\rangle}\right)} = z,\quad j=m+1,m+2,\cdots,k~,\nonumber\\
&\sum\limits_{\ell = 1}^k {r_{\left\langle \ell \right\rangle}=1},\quad r_{\left\langle \ell \right\rangle} > 0,\quad {\ell = 1,2,\cdots k}~. \nonumber
\end{alignat}

\begin{lemma}\label{thm2}
The optimization problem (\ref{opt1}) and (\ref{opt2}) are equivalent, i.e., a solution ${r^*} = ( r_{\left\langle 1 \right\rangle}^*,r_{\left\langle 2 \right\rangle}^*, \cdots ,r_{\left\langle k \right\rangle}^* )$ to (\ref{opt1}) is optimal, if and only if it is also an optimal solution to (\ref{opt2}).
\end{lemma}


\begin{remark}
In particular, when $m=1$, KKT conditions ($\ref{K2}$)-($\ref{K5}$) in online appendix A.2 reduce to $\sum\nolimits_{\ell = 2}^k {{\lambda_{1\ell}}}  = 1$, and
\begin{equation}\label{KK6}
    {\sum\limits_{\ell = 2}^k {{\lambda_{1\ell}}\left. {\frac{{\partial {G_{1\ell}} ( {y,r_{\left\langle  \ell \right\rangle}^*} )}}{{\partial y}}} \right|} _{y = r_{\left\langle 1 \right\rangle}^*}} = \gamma~,
\end{equation}
\begin{equation}\label{KK7}
{\lambda_{1\ell}}{\left. {\frac{{\partial {G_{1\ell}} ( {r_{\left\langle  1 \right\rangle}^*,y} )}}{{\partial y}}} \right|_{y = r_{\left\langle  \ell \right\rangle}^*}} = \gamma,\quad \ell=2,3,\cdots,k~,
\end{equation}
$${\lambda_{1\ell}}\left[ {z - {G_{1\ell}}\left( {r_{\left\langle  1 \right\rangle}^*,r_{\left\langle  \ell \right\rangle}^*} \right)} \right] = 0,\quad \ell=2,3,\cdots,k~,$$
respectively, which correspond to the first order conditions for the optimization problem in~\citet{glynn2004large}. Note that each Lagrangian multipliers $\lambda_{1\ell}$, $\ell=2,3,\cdots,k$ can be directly solved from (\ref{KK7}). However, for selecting top-$m$ alternatives problem (in particular, $m \ne 1$ and $m \ne k-1$), since the summation notations in (\ref{K3}) and (\ref{K4}), each $\lambda_{ij}$, $i\in\left\{1,2,\cdots,m\right\}$, $j \in \left\{m+1,m+2,\cdots,k\right\}$ in KKT conditions in online appendix A.2 cannot be directly solved, which leads to more complicated analyses compared with selecting the best alternative.
\end{remark}

\begin{corollary}\label{thm22}
The asymptotically optimal sampling ratios which attain the optimal large deviations rate of the PFS for selecting top-$m$ alternatives satisfy, for $i,h=1,2,\cdots,m$, $j,\ell=m+1,m+2,\cdots,k$,
\begin{equation}\label{optopt}
\mathop {\min }\limits_{j = m + 1,m+2, \cdots ,k} {G_{hj}\left(r_{\left\langle  h \right\rangle}^*,r_{\left\langle  j \right\rangle}^*\right)} = \mathop {\min }\limits_{i = 1,2, \cdots ,m} {G_{i \ell}\left(r_{\left\langle  i \right\rangle}^*,r_{\left\langle \ell \right\rangle}^*\right)}~,
\end{equation}
where $r_{\widetilde \ell}^*$, ${\widetilde \ell}=1,2,\cdots,k$ are asymptotically optimal sampling ratios.
\end{corollary}

\proof{Proof of Corollary~\ref{thm22}}
The Corollary is a direct consequence of Lemma~\ref{thm2}.
\endproof

\begin{remark}
Suppose there are four competing alternatives and we aim to select top-2 alternatives, i.e., $k=4$ and $m=2$. Then (\ref{optopt}) yields
$$\begin{aligned}
& \min \{ {{G_{13} (r_{\left\langle 1 \right\rangle}^*,r_{\left\langle 3 \right\rangle}^* )},{G_{14} (r_{\left\langle 1 \right\rangle}^*,r_{\left\langle 4 \right\rangle}^* )}} \} = \min \{ {{G_{23} (r_{\left\langle 2 \right\rangle}^*,r_{\left\langle 3 \right\rangle}^* )},{G_{24} (r_{\left\langle 2 \right\rangle}^*,r_{\left\langle 4 \right\rangle}^* )}} \} \\
= & \min \{ {{G_{13} (r_{\left\langle 1 \right\rangle}^*,r_{\left\langle 3 \right\rangle}^*)},{G_{23} (r_{\left\langle 2 \right\rangle}^*,r_{\left\langle 3 \right\rangle}^* )}} \} = \min \{ {{G_{14} (r_{\left\langle 1 \right\rangle}^*,r_{\left\langle 4 \right\rangle}^* )},{G_{24}(r_{\left\langle 2 \right\rangle}^*,r_{\left\langle 4 \right\rangle}^*)}} \}~,
\end{aligned}$$
which contains 7 different cases, corresponding to different values of Lagrangian multipliers $\lambda_{ij}$ in the KKT conditions in online appendix A.2. Table~\ref{table1} shows the 7 cases in detail. Note that once values of Lagrangian multipliers are given, an optimal decreasing rate can be determined correspondingly. In the above example, the number of competing alternatives is small and the different cases of optimal decreasing rate can be easily enumerated. However, the number of cases in (\ref{optopt}) increases as the number of competing alternatives increases.
By considering the possible number of Lagrangian multipliers equals $0$, the number of cases in (\ref{optopt}) can be extremely large as the number of competing alternatives increases, making it difficult to derive allocation procedures that optimize or asymptotically optimize the rate directly.

\begin{table}[htbp]
\caption{An Example of Optimal Decreasing Rate of the Large Deviations of PFS for Selecting Top-$2$ Alternatives From 4 Alternatives}
\label{table1}
\centering
\begin{tabular}{ccc}
\hline
   Cases & Lagrangian Multipliers & Optimal Decreasing Rate \\
\hline
   1 & $\lambda_{14},\lambda_{23},\lambda_{24} >0,\lambda_{13}=0$ & ${G_{14} (r_{\left\langle 1 \right\rangle}^*,r_{\left\langle 4 \right\rangle}^* )} = {G_{23} (r_{\left\langle 2 \right\rangle}^*,r_{\left\langle 3 \right\rangle}^* )} = {G_{24} (r_{\left\langle 2 \right\rangle}^*,r_{\left\langle 4 \right\rangle}^* )}$ \\
   2 & $\lambda_{13},\lambda_{23},\lambda_{24} >0,\lambda_{14}=0$ & ${G_{13} (r_{\left\langle 1 \right\rangle}^*,r_{\left\langle 3 \right\rangle}^* )} = {G_{23} (r_{\left\langle 2 \right\rangle}^*,r_{\left\langle 3 \right\rangle}^* )} = {G_{24} (r_{\left\langle 2 \right\rangle}^*,r_{\left\langle 4 \right\rangle}^*\ )}$ \\
   3 & $\lambda_{13},\lambda_{14},\lambda_{24} >0,\lambda_{23}=0$ & ${G_{13} (r_{\left\langle 1 \right\rangle}^*,r_{\left\langle 3 \right\rangle}^*\ )} = {G_{14} (r_{\left\langle 1 \right\rangle}^*,r_{\left\langle 4 \right\rangle}^* )} = {G_{24} (r_{\left\langle 2 \right\rangle}^*,r_{\left\langle 4 \right\rangle}^* )}$ \\
   4 & $\lambda_{13},\lambda_{14},\lambda_{23} >0,\lambda_{24}=0$ & ${G_{13} (r_{\left\langle 1 \right\rangle}^*,r_{\left\langle 3 \right\rangle}^* )} = {G_{14} (r_{\left\langle 1 \right\rangle}^*,r_{\left\langle 4 \right\rangle}^* )} = {G_{23} (r_{\left\langle 2 \right\rangle}^*,r_{\left\langle 3 \right\rangle}^* )}$ \\
   5 & $\lambda_{14},\lambda_{23} >0,\lambda_{13}=\lambda_{24}=0$ & ${G_{14} (r_{\left\langle 1 \right\rangle}^*,r_{\left\langle 4 \right\rangle}^* )} = {G_{23} (r_{\left\langle 2 \right\rangle}^*,r_{\left\langle 3 \right\rangle}^* )}$ \\
   6 & $\lambda_{13},\lambda_{24} >0,\lambda_{14}=\lambda_{23}=0$ & ${G_{13} (r_{\left\langle 1 \right\rangle}^*,r_{\left\langle 3 \right\rangle}^* )} = {G_{24} (r_{\left\langle 2 \right\rangle}^*,r_{\left\langle 4 \right\rangle}^* )}$ \\
   7 & $\lambda_{13},\lambda_{14},\lambda_{23},\lambda_{24} >0$ & ${G_{13} (r_{\left\langle 1 \right\rangle}^*,r_{\left\langle 3 \right\rangle}^* )} = {G_{14} (r_{\left\langle 1 \right\rangle}^*,r_{\left\langle 4 \right\rangle}^* )} = {G_{23} (r_{\left\langle 2 \right\rangle}^*,r_{\left\langle 3 \right\rangle}^* )} = {G_{24} (r_{\left\langle 2 \right\rangle}^*,r_{\left\langle 4 \right\rangle}^* )}$ \\
\hline
\end{tabular}
\end{table}
\end{remark}

For convenience of analysis, let ${J^{\left(i\right)}} \mathop  = \limits^\Delta \arg \mathop {\min }\limits_{j = m + 1, m+2, \cdots ,k} {G_{ij} (r_{\left\langle  i \right\rangle}^*,r_{\left\langle  j \right\rangle}^*)}$, $i \in \left\{1,2,\cdots,m\right\}$ be a set containing alternatives $\langle j\rangle$ that achieve the minimum value of $G_{ij} ({r_{\left\langle  i \right\rangle}^*,r_{\left\langle  j \right\rangle}^*} )$ for a fixed alternative $\langle i \rangle$, and ${I^{\left(j\right)}} \mathop  = \limits^\Delta \arg \mathop {\min }\limits_{i = 1,2, \cdots ,m} {G_{ij} (r_{\left\langle  i \right\rangle}^*,r_{\left\langle  j \right\rangle}^* )}$, $j \in \left\{ m+1,m+2,\cdots,k\right\}$, be a set containing alternatives $\langle i \rangle$, which achieve the minimum value of $G_{ij} ({r_{\left\langle  i \right\rangle}^*,r_{\left\langle  j \right\rangle}^*} )$ for a fixed alternative $\langle j \rangle$. Corresponding to different cases in (\ref{optopt}), the number of alternatives contained in sets ${J^{\left(i\right)}}$, $i=1,2,\cdots,m$ takes possible values in $\left\{1,2,\cdots,k-m\right\}$, and the number of alternatives contained in sets ${I^{\left(j\right)}}$, $j=m+1,m+2,\cdots,k$ takes possible values in $\left\{1,2,\cdots,m\right\}$. For example, for Case 1 in Table 1, ${\left|J^{\left( 1 \right)}\right|} = 1$, ${\left|J^{\left( 2 \right)}\right|} = 2$, ${\left|I^{\left( 3 \right)}\right|} = 1$, ${\left|I^{\left( 4 \right)}\right|} = 2$, and for Case 5 in Table 1, ${\left|J^{\left( 1 \right)}\right|} = 1$, ${\left|J^{\left( 2 \right)}\right|} = 1$, ${\left|I^{\left( 3 \right)}\right|} = 1$, ${\left|I^{\left( 4 \right)}\right|} = 1$, where $\left| A \right|$ denotes the cardinality of a set $A$.
With the definition of ${J^{\left( i \right)}}$ and ${I^{\left( j \right)}}$, Lagrangian multipliers $\lambda_{ij}$ in the KKT conditions in online appendix A.2 satisfy $\lambda_{ij} > 0$ if $i \in {I^{\left(j\right)}}$ and $j \in {J^{\left(i\right)}}$. In addition, for $h \in {I^{\left(j\right)}}$ and $\ell \in {J^{\left(i\right)}}$, we have $G_{hj} (r_{\left\langle  h \right\rangle}^*,r_{\left\langle  j \right\rangle}^*) = G_{i\ell}(r_{\left\langle  i \right\rangle}^*,r_{\left\langle \ell \right\rangle}^*) = z$.

\begin{remark}
$G_{ij}\left(\cdot,\cdot\right)$ is a rate function of a wrong comparison between alternatives $\langle i\rangle$ and $\langle j\rangle$, $i \in \left\{1,2,\cdots,m\right\}$, $j \in \left\{m+1,m+2,\cdots,k\right\}$, then the larger $G_{ij}\left(\cdot,\cdot\right)$ is, the less likely we make a wrong comparison between alternatives $\langle i\rangle$ and $\langle j\rangle$. Each of set $J^{\left(i\right)}$ contains alternatives that are most likely to lead to wrong comparisons with alternative $\langle i\rangle$. Similarly, each of set $I^{\left(j\right)}$ contains alternatives that are most likely to lead to wrong comparisons with alternative $\langle j\rangle$. Then $\mathop {\min }\limits_{j = m + 1,m+2, \cdots ,k} {G_{ij} (r_{\left\langle  i \right\rangle}^*,r_{\left\langle  j \right\rangle}^* )}$ and $\mathop {\min }\limits_{i = 1,2, \cdots ,m} {G_{ij} (r_{\left\langle  i \right\rangle}^*,r_{\left\langle  j \right\rangle}^* )}$ represent the the difficulty of correctly selecting alternatives $\langle i\rangle$ and $\langle j\rangle$, respectively. (\ref{optopt}) shows that the difficulty of correctly selecting alternatives $\langle i\rangle$ and $\langle j\rangle$, $i=1,2,\cdots,m$, $j=m+1,m+2,\cdots,k$ should be equal.
\end{remark}

\begin{remark}\label{inadeGao}
(\ref{optopt}) coincides with the asymptotic solutions of an approximate static optimization problem in~\citet{gao2016new}, which are derived by asymptotically optimizing the PCS approximated by a large deviations principle and Bonferroni inequality. In the analysis of their asymptotic solutions, by deeming that each of set $J^{\left(i\right)}$, $i = 1,2,\cdots,m$ and $I^{\left(j\right)}$, $j = m+1,m+2,\cdots,k$ contain only one alternative, i.e., $\left| {{J^{\left( i \right)}}} \right| = 1$ and $\left| {{I^{\left( j \right)}}} \right| = 1$, \citet{gao2016new} only consider some special cases in (\ref{optopt}). For the example in Table~\ref{table1}, only Cases 5 and 6 are considered. In this work, we show that (\ref{optopt}) contains an increasing number of cases as the number of competing alternatives grows.
\end{remark}

\subsection{Asymptotic Optimality of AOAm}

We consider that the sampling distributions of all alternatives are normal. If ${X_{\ell,{t_\ell}}}$ follows i.i.d. normal distribution $N\left( {{\mu _\ell},\sigma _\ell^2} \right)$, $\ell=1,2\cdots,k$, we have ${M_{\left\langle \ell \right\rangle}}\left( y \right) = \exp ( {y{\mu _{\left\langle \ell \right\rangle}} + \frac{1}{2}\sigma _{\left\langle \ell \right\rangle}^2{y^2}} )$, and ${\Lambda_{\left\langle \ell \right\rangle}^*}\left( x \right) = {{{\left( x - {\mu _{\left\langle \ell \right\rangle}} \right)^2}} / (2\sigma _{\left\langle \ell \right\rangle}^2)}$. With Lemma~\ref{lem001},
$${r_{\left\langle  i \right\rangle}}{\left. {\frac{{d\Lambda_{\left\langle  i \right\rangle}^*\left( x \right)}}{{dx}}} \right|_{x = x( {{r_{\left\langle  i \right\rangle}},{r_{\left\langle  j \right\rangle}}} )}} + {r_{\left\langle  j \right\rangle}}{\left. {\frac{{d\Lambda_{\left\langle  j \right\rangle}^*\left( x \right)}}{{dx}}} \right|_{x = x( {{r_{\left\langle  i \right\rangle}},{r_{\left\langle  j \right\rangle}}})}} = 0,\quad i=1,\cdots,m,\;j=m+1,\cdots,k~,$$
we have
$${r_{\left\langle  i \right\rangle}}\frac{{x\left( {{r_{\left\langle  i \right\rangle}},{r_{\left\langle  j \right\rangle}}} \right) - {\mu _{\left\langle i \right\rangle}}}}{{\sigma _{\left\langle i \right\rangle}^2}} + {r_{\left\langle  j \right\rangle}}\frac{{x\left( {{r_{\left\langle  i \right\rangle}},{r_{\left\langle  j \right\rangle}}} \right) - {\mu _{\left\langle j \right\rangle}}}}{{\sigma _{\left\langle j \right\rangle}^2}} = 0,\quad i=1,\cdots,m,\;j=m+1,\cdots,k~,$$
and the infimum point $x\left(r_{\left\langle  i \right\rangle},r_{\left\langle  j \right\rangle}\right)$ can be expressed as
$$x\left( {{r_{\left\langle  i \right\rangle}},{r_{\left\langle  j \right\rangle}}} \right) = \left( {\frac{{{{r_{\left\langle  i \right\rangle}} \mathord{\left/
 {\vphantom {{r_{\left\langle  i \right\rangle}} {\sigma _{\left\langle i \right\rangle}^2}}} \right.
 \kern-\nulldelimiterspace} {\sigma _{\left\langle i \right\rangle}^2}}}}{{{{r_{\left\langle  i \right\rangle}} \mathord{\left/
 {\vphantom {{r_{\left\langle  i \right\rangle}} {\sigma _{\left\langle i \right\rangle}^2}}} \right.
 \kern-\nulldelimiterspace} {\sigma _{\left\langle i \right\rangle}^2}} + {{r_{\left\langle  j \right\rangle}} \mathord{\left/
 {\vphantom {{r_{\left\langle  j \right\rangle}} {\sigma _{\left\langle j \right\rangle}^2}}} \right.
 \kern-\nulldelimiterspace} {\sigma _{\left\langle j \right\rangle}^2}}}}} \right){\mu _{\left\langle i \right\rangle}} + \left( {\frac{{{{r_{\left\langle  j \right\rangle}} \mathord{\left/
 {\vphantom {{r_{\left\langle  j \right\rangle}} {\sigma _{\left\langle j \right\rangle}^2}}} \right.
 \kern-\nulldelimiterspace} {\sigma _{\left\langle j \right\rangle}^2}}}}{{{{r_{\left\langle  i \right\rangle}} \mathord{\left/
 {\vphantom {{r_{\left\langle  i \right\rangle}} {\sigma _{\left\langle i \right\rangle}^2}}} \right.
 \kern-\nulldelimiterspace} {\sigma _{\left\langle i \right\rangle}^2}} + {{r_{\left\langle  j \right\rangle}} \mathord{\left/
 {\vphantom {{r_{\left\langle  j \right\rangle}} {\sigma _{\left\langle j \right\rangle}^2}}} \right.
 \kern-\nulldelimiterspace} {\sigma _{\left\langle j \right\rangle}^2}}}}} \right){\mu _{\left\langle j \right\rangle}}~.$$

From ${G_{ij}}\left( {{r_{\left\langle  i \right\rangle}},{r_{\left\langle  j \right\rangle}}} \right) = {r_{\left\langle  i \right\rangle}}\Lambda_{\left\langle  i \right\rangle}^*\left( {x\left( {{r_{\left\langle  i \right\rangle}},{r_{\left\langle  j \right\rangle}}} \right)} \right) + {r_{\left\langle  j \right\rangle}}\Lambda_{\left\langle  j \right\rangle}^*\left( {x\left( {{r_{\left\langle  i \right\rangle}},{r_{\left\langle  j \right\rangle}}} \right)} \right)$, we obtain
 \begin{equation}\label{normalrate}
 {G_{ij}}\left( {r_{\left\langle  i \right\rangle},r_{\left\langle  j \right\rangle}} \right) = \frac{{{{\left( {{\mu _{\left\langle i \right\rangle}} - {\mu _{\left\langle j \right\rangle}}} \right)}^2}}}{{2\left( {{{\sigma _{\left\langle i \right\rangle}^2} \mathord{\left/
 {\vphantom {{\sigma _{\left\langle i \right\rangle}^2} {r_{\left\langle  i \right\rangle}}}} \right.
 \kern-\nulldelimiterspace} {r_{\left\langle  i \right\rangle}}} + {{\sigma _{\left\langle j \right\rangle}^2} \mathord{\left/
 {\vphantom {{\sigma _{\left\langle j \right\rangle}^2} {r_{\left\langle  j \right\rangle}}}} \right.
 \kern-\nulldelimiterspace} {r_{\left\langle  j \right\rangle}}}} \right)}},\quad i=1,2,\cdots,m,\;j=m+1,m+2,\cdots,k~,
 \end{equation}

\begin{figure}[htbp]
\center{\includegraphics[width=0.5\textwidth]{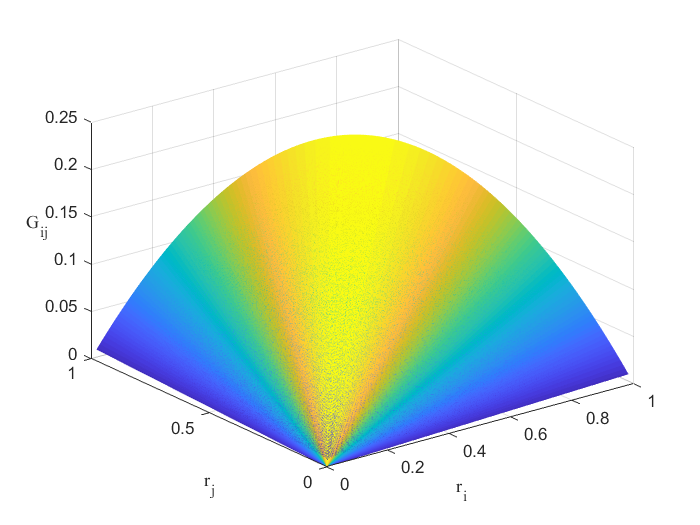}}
\caption{Rate function $G_{ij}\left(r_{\left\langle  i \right\rangle},r_{\left\langle  j \right\rangle}\right)$ under normal assumptions when $\left(\mu_{\left\langle j \right\rangle}-\mu_{\left\langle i \right\rangle}\right)=1$ and $\sigma_{\left\langle i \right\rangle}=\sigma_{\left\langle j \right\rangle}=1$, $i \in \left\{1,2,\cdots,m\right\}$, $j \in \left\{m+1,m+2,\cdots,k\right\}$.}
\label{fig1}
\end{figure}

We visualize a rate function $G_{ij}\left(r_{\left\langle  i \right\rangle},r_{\left\langle  j \right\rangle}\right)$ under normal sampling distributions when $\left(\mu_{\left\langle j \right\rangle}-\mu_{\left\langle i \right\rangle}\right)=1$, and $\sigma_{\left\langle i \right\rangle}=\sigma_{\left\langle j \right\rangle}=1$ in Figure~\ref{fig1}, where we can see that $G_{ij}\left(r_{\left\langle  i \right\rangle},r_{\left\langle  j \right\rangle}\right)$ is a strictly concave function of $\left(r_{\left\langle  i \right\rangle},r_{\left\langle  j \right\rangle}\right)$ and is strictly increasing in $r_{\left\langle  i \right\rangle}$ and $r_{\left\langle  j \right\rangle}$, respectively, $i \in \left\{1,2,\cdots,m\right\}$, $j \in \left\{m+1,m+2,\cdots,k\right\}$.

\begin{lemma}\label{normal}
Consider that $X_{{\widetilde{\ell}},{t_{\widetilde{\ell}}}}$ follows i.i.d. normal distribution $N ( {{\mu_{\widetilde{\ell}}},\sigma _{\widetilde{\ell}}^2} )$, ${\widetilde{\ell}}=1,2,\cdots,k$. If there exist $\left( {{r_{\left\langle  i \right\rangle}},{r_{\left\langle  j \right\rangle}}} \right) > 0$ and $\left( {{r_{\left\langle  h \right\rangle}},{r_{\left\langle  j \right\rangle}}} \right) > 0$, such that ${G_{ij}}\left( {{r_{\left\langle  i \right\rangle}},{r_{\left\langle  j \right\rangle}}} \right) = {G_{hj}}\left( {{r_{\left\langle  h \right\rangle}},{r_{\left\langle  j \right\rangle}}} \right)$, $i,h = 1,2, \cdots ,m$, $j = m + 1,m + 2, \cdots ,k$, then we have
$$\frac{{{{\left. {{{\partial {G_{ij}}\left( {{r_{\left\langle  i \right\rangle}},y} \right)} \mathord{\left/
 {\vphantom {{\partial {G_{ij}}\left( {{r_{\left\langle  i \right\rangle}},y} \right)} {\partial y}}} \right.
 \kern-\nulldelimiterspace} {\partial y}}} \right|}_{y = {r_{\left\langle  j \right\rangle}}}}}}{{{{\left. {{{\partial {G_{hj}}\left( {{r_{\left\langle h \right\rangle}},y} \right)} \mathord{\left/
 {\vphantom {{\partial {G_{hj}}\left( {{r_{\left\langle h \right\rangle}},y} \right)} {\partial y}}} \right.
 \kern-\nulldelimiterspace} {\partial y}}} \right|}_{y = {r_{\left\langle  j \right\rangle}}}}}} = \frac{{{{\sigma _{\left\langle h \right\rangle}^2} \mathord{\left/
 {\vphantom {{\sigma _{\left\langle h \right\rangle}^2} {{r_h}}}} \right.
 \kern-\nulldelimiterspace} {{r_{\left\langle h \right\rangle}}}} + {{\sigma _{\left\langle j \right\rangle}^2} \mathord{\left/
 {\vphantom {{\sigma _{\left\langle j \right\rangle}^2} {{r_{\left\langle  j \right\rangle}}}}} \right.
 \kern-\nulldelimiterspace} {{r_{\left\langle  j \right\rangle}}}}}}{{{{\sigma _{\left\langle i \right\rangle}^2} \mathord{\left/
 {\vphantom {{\sigma _{\left\langle i \right\rangle}^2} {{r_{\left\langle  i \right\rangle}}}}} \right.
 \kern-\nulldelimiterspace} {{r_{\left\langle  i \right\rangle}}}} + {{\sigma _{\left\langle j \right\rangle}^2} \mathord{\left/
 {\vphantom {{\sigma _{\left\langle j \right\rangle}^2} {{r_{\left\langle  j \right\rangle}}}}} \right.
 \kern-\nulldelimiterspace} {{r_{\left\langle  j \right\rangle}}}}}}~,$$
and similarly, if there exist $\left( {{r_{\left\langle  i \right\rangle}},{r_{\left\langle  j \right\rangle}}} \right) > 0$ and $\left( {{r_{\left\langle  i \right\rangle}},{r_{\left\langle \ell \right\rangle}}} \right) > 0$, which satisfy ${G_{ij}}\left( {{r_{\left\langle  i \right\rangle}},{r_{\left\langle  j \right\rangle}}} \right) = {G_{i \ell}}\left( {{r_{\left\langle  i \right\rangle}},{r_{\left\langle \ell \right\rangle}}} \right)$, $i = 1,2, \cdots ,m$, $j,\ell = m + 1,m + 2, \cdots ,k$, then we have
$$\frac{{{{\left. {{{\partial {G_{ij}}\left( {y,{r_{\left\langle  j \right\rangle}}} \right)} \mathord{\left/
 {\vphantom {{\partial {G_{ij}}\left( {y,{r_{\left\langle  j \right\rangle}}} \right)} {\partial y}}} \right.
 \kern-\nulldelimiterspace} {\partial y}}} \right|}_{y = {r_{\left\langle  i \right\rangle}}}}}}{{{{\left. {{{\partial {G_{i\ell}}\left( {y,{r_{\left\langle \ell \right\rangle}}} \right)} \mathord{\left/
 {\vphantom {{\partial {G_{i\ell}}\left( {y,{r_{\left\langle \ell \right\rangle}}} \right)} {\partial y}}} \right.
 \kern-\nulldelimiterspace} {\partial y}}} \right|}_{y = {r_{\left\langle  i \right\rangle}}}}}} = \frac{{{{\sigma _{\left\langle i \right\rangle}^2} \mathord{\left/
 {\vphantom {{\sigma _{\left\langle i \right\rangle}^2} {{r_{\left\langle  i \right\rangle}}}}} \right.
 \kern-\nulldelimiterspace} {{r_{\left\langle  i \right\rangle}}}} + {{\sigma _{\left\langle \ell \right\rangle}^2} \mathord{\left/
 {\vphantom {{\sigma _{\left\langle \ell \right\rangle}^2} {{r_{\left\langle \ell \right\rangle}}}}} \right.
 \kern-\nulldelimiterspace} {{r_{\left\langle \ell \right\rangle}}}}}}{{{{\sigma _{\left\langle i \right\rangle}^2} \mathord{\left/
 {\vphantom {{\sigma _{\left\langle i \right\rangle}^2} {{r_{\left\langle  i \right\rangle}}}}} \right.
 \kern-\nulldelimiterspace} {{r_{\left\langle  i \right\rangle}}}} + {{\sigma _{\left\langle j \right\rangle}^2} \mathord{\left/
 {\vphantom {{\sigma _{\left\langle j \right\rangle}^2} {{r_{\left\langle j \right\rangle}}}}} \right.
 \kern-\nulldelimiterspace} {{r_{\left\langle  j \right\rangle}}}}}}~.$$
\end{lemma}



\begin{theorem}\label{thmLDTnormal}
The asymptotically optimal sampling ratios that optimize large deviations rate of PFS for selecting top-$m$ alternatives with normal sampling distributions satisfy, for $i,h = 1,2,\cdots,m$, $j,\ell=m+1,m+2,\cdots,k$,
\begin{equation}\label{optrate}
 \mathop {\min }\limits_{j = m + 1,m+2, \cdots ,k} \frac{{{{\left( {{\mu _{\left\langle h \right\rangle}} - {\mu _{\left\langle j \right\rangle}}} \right)}^2}}}{ {{{\sigma _{\left\langle h \right\rangle}^2} \mathord{\left/
 {\vphantom {{\sigma _{\left\langle h \right\rangle}^2} {r_{\left\langle h \right\rangle}^*}}} \right.
 \kern-\nulldelimiterspace} {r_{\left\langle h \right\rangle}^*}} + {{\sigma _{\left\langle j \right\rangle}^2} \mathord{\left/
 {\vphantom {{\sigma _{\left\langle j \right\rangle}^2} {r_{\left\langle  j \right\rangle}^*}}} \right.
 \kern-\nulldelimiterspace} {r_{\left\langle  j \right\rangle}^*}}} } = \mathop {\min }\limits_{i = 1,2, \cdots ,m} \frac{{{{\left( {{\mu _{\left\langle i \right\rangle}} - {\mu _{\left\langle \ell \right\rangle}}} \right)}^2}}}{ {{{\sigma _{\left\langle i \right\rangle}^2} \mathord{\left/
 {\vphantom {{\sigma _{\left\langle i \right\rangle}^2} {r_{\left\langle  i \right\rangle}^*}}} \right.
 \kern-\nulldelimiterspace} {r_{\left\langle i \right\rangle}^*}} + {{\sigma _{\left\langle \ell \right\rangle}^2} \mathord{\left/
 {\vphantom {{\sigma _{\left\langle j \right\rangle}^2} {r_{\left\langle \ell \right\rangle}^*}}} \right.
 \kern-\nulldelimiterspace} {r_{\left\langle \ell \right\rangle}^*}}}}~,
 \end{equation}
 \begin{equation}\label{thmbalance}
\sum\limits_{i = 1}^m {\frac{{{{\left( {r_{\left\langle  i \right\rangle}^*} \right)}^2}}}{{\sigma _{\left\langle i \right\rangle}^2}}}  = \sum\limits_{j = m + 1}^k {\frac{{{{\left( {r_{\left\langle  j \right\rangle}^*} \right)}^2}}}{{\sigma _{\left\langle j \right\rangle}^2}}}~,
\end{equation}
where $\sum\nolimits_{{\widetilde {\ell}} = 1}^k {r_{ {\widetilde {\ell}}}^*}  = 1$, ${r_{{\widetilde {\ell}}}^*} > 0$, ${\widetilde {\ell}}=1,2,\cdots,k$.
\end{theorem}

\begin{remark}\label{rmkspecificequal}
In the proof of Theorem \ref{thmLDTnormal} in online appendix A.2, we show that (\ref{thmbalance}) holds for all cases in (\ref{optrate}). In the following, we consider a special case in (\ref{optopt}), where for $i,h=1,2,\cdots,m$ and $j,\ell=m+1,m+2,\cdots,k$, ${G_{hj}} ( {r_{\left\langle h \right\rangle}^*,r_{\left\langle  j \right\rangle}^*} ) = {G_{i{\ell}}} ( {r_{\left\langle  i \right\rangle}^*,r_{\left\langle \ell \right\rangle}^*} )$. With the definition of $J^{\left(i\right)}$ and $I^{\left(j\right)}$, the special case implies $J^{\left(i\right)}=\left\{m+1,m+2,\cdots,k\right\}$, $\left| {{J^{\left( i \right)}}} \right| = k - m$, and $I^{\left(j\right)}=\left\{1,2,\cdots,m\right\}$, $\left| {{I^{\left( j \right)}}} \right| = m$. Correspondingly, Lagrangian multipliers satisfy $\lambda_{ij} > 0$ in the KKT conditions of optimization problem (\ref{opt1}). The special case requires all optimal rates $G_{ij} (r_{\left\langle  i \right\rangle}^*,r_{\left\langle  j \right\rangle}^*)$ to be equal, which is stronger than  (\ref{optopt}), so the special case leads to a sufficient condition rather than a necessary one. Following normal sampling distributions, the special case yields $\frac{{{( {{\mu _{\left\langle h \right\rangle}} - {\mu _{\left\langle j \right\rangle}}} )^2}}}{{{{\sigma _{\left\langle h \right\rangle}^2} / {r_{\left\langle h \right\rangle}^* + {{\sigma _{\left\langle j \right\rangle}^2} / {r_{\left\langle  j \right\rangle}^*}}}}}} =  \frac{{{( {{\mu _{\left\langle i \right\rangle}} - {\mu _{\left\langle \ell \right\rangle}}} )^2}}}{{{{\sigma _{\left\langle i \right\rangle}^2} / {r_{\left\langle  i \right\rangle}^* + {{\sigma _{\left\langle \ell \right\rangle}^2} / {r_{\left\langle \ell \right\rangle}^*}}}}}}$.

For the special case, the second term in the denominator of (\ref{toolong}) in online appendix A.2 vanishes, since $\left| {{J^{\left( i \right)}}} \right| = k - m$, $i=1,2\cdots,m$, and $\widehat \ell \notin {J^{\left( i \right)}}$ is empty, which leads to
\begin{equation}\label{topmratio}
{\sum\limits_{i = 1}^m \frac{1}{{\sum\limits_{j = m + 1}^k {\frac{{{{\left. {{{\partial {G_{ij}}\left( {y,r_{\left\langle  j \right\rangle}^*} \right)} \mathord{\left/{\vphantom {{\partial {G_{ij}}\left( {y,r_{\left\langle  j \right\rangle}^*} \right)} {\partial y}}} \right.\kern-\nulldelimiterspace} {\partial y}}} \right|}_{y = r_{\left\langle  i \right\rangle}^*}}}}{{{{\left. {{{\partial {G_{ij}}\left( {r_{\left\langle  i \right\rangle}^*,y} \right)} \mathord{\left/{\vphantom {{\partial {G_{ij}}\left( {r_{\left\langle  i \right\rangle}^*,y} \right)} {\partial y}}} \right.
 \kern-\nulldelimiterspace} {\partial y}}} \right|}_{y = r_{\left\langle  j \right\rangle}^*}}}}} }}} = 1~.
\end{equation}

In particular, when $m=1$, the optimal decreasing rate of the PFS for selecting the best alternative in \citet{glynn2004large} satisfies, for $j,\ell=2,3,\cdots,k$, ${G_{1j}} ( {r_{\left\langle 1 \right\rangle }^*,r_{\left\langle j \right\rangle }^*} ) = {G_{1 \ell}} ( {r_{\left\langle 1 \right\rangle }^*,r_{\left\langle \ell \right\rangle }^*} )$, which belongs to the special case above. Then (\ref{optrate})-(\ref{topmratio}) lead to
\begin{equation}\label{glynn}
 \frac{{{{\left( {{\mu _{\left\langle 1 \right\rangle}} - {\mu _{\left\langle j \right\rangle}}} \right)}^2}}}{{ {{{\sigma _{\left\langle 1 \right\rangle}^2} \mathord{\left/
 {\vphantom {{\sigma _{\left\langle 1 \right\rangle}^2} {r_{\left\langle 1 \right\rangle}^*}}} \right.
 \kern-\nulldelimiterspace} {r_{\left\langle 1 \right\rangle}^*}} + {{\sigma _{\left\langle j \right\rangle}^2} \mathord{\left/
 {\vphantom {{\sigma _{\left\langle j \right\rangle}^2} {r_{\left\langle  j \right\rangle}^*}}} \right.
 \kern-\nulldelimiterspace} {r_{\left\langle  j \right\rangle}^*}}}}} = \frac{{{{\left( {{\mu _{\left\langle 1 \right\rangle}} - {\mu _{\left\langle \ell \right\rangle}}} \right)}^2}}}{{ {{{\sigma _{\left\langle 1 \right\rangle}^2} \mathord{\left/
 {\vphantom {{\sigma _{\left\langle 1 \right\rangle}^2} {r_{\left\langle 1 \right\rangle}^*}}} \right.
 \kern-\nulldelimiterspace} {r_{\left\langle 1 \right\rangle}^*}} + {{\sigma _{\left\langle \ell \right\rangle}^2} \mathord{\left/
 {\vphantom {{\sigma _{\left\langle \ell \right\rangle}^2} {r_{\left\langle \ell \right\rangle}^*}}} \right.
 \kern-\nulldelimiterspace} {r_{\left\langle \ell \right\rangle}^*}}}}}~,
\end{equation}
\begin{equation}\label{wscequ2}
r_{\left\langle 1 \right\rangle}^* = {\sigma _{\left\langle 1 \right\rangle}}\sqrt {\sum\nolimits_{j = 2}^k {{{{( {r_{\left\langle  j \right\rangle}^*} )^2}} \mathord{\left/
 {\vphantom {{{( {r_{\left\langle  j \right\rangle}^*} )^2}} {{{\left( {{\sigma _{\left\langle j \right\rangle}}} \right)}^2}}}} \right.
 \kern-\nulldelimiterspace} {{\sigma_{\left\langle j \right\rangle}^2}}}} }~,
\end{equation}
and
\begin{equation}\label{wscequ3}
\sum\limits_{j = 2}^k {\frac{{{{\left. {{{\partial {G_{1j}}( {y,r_{\left\langle  j \right\rangle}^*} )} \mathord{\left/
 {\vphantom {{\partial {G_{1j}}\left( {y,r_{\left\langle  j \right\rangle}^*} \right)} {\partial y}}} \right.
 \kern-\nulldelimiterspace} {\partial y}}} \right|}_{y = r_{\left\langle  1 \right\rangle}^*}}}}{{{{\left. {{{\partial {G_{1j}} ( {r_{\left\langle  1 \right\rangle}^*,y} )} \mathord{\left/
 {\vphantom {{\partial {G_{1j}}\left( {r_{\left\langle  1 \right\rangle}^*,y} \right)} {\partial y}}} \right.
 \kern-\nulldelimiterspace} {\partial y}}} \right|}_{y = r_{\left\langle  j \right\rangle}^*}}}}}  = 1~,
\end{equation}
respectively, which correspond to the results for selecting the best alternative with normal sampling distributions derived in~\citet{glynn2004large}. Similarly, $m=k-1$ also belongs to the above special case. (\ref{optopt}) for selecting top-$m$ alternatives contains different cases, making $\widehat \ell \notin {J^{\left( i \right)}}$, $i=1,2,\cdots,m$ in the second term in the denominator of (\ref{toolong}) in the online appendix A.2 not empty for some cases; and therefore, leads to more complex analyses compared with selecting the best.
\end{remark}

\begin{remark}
(\ref{thmbalance}) uses the variance information of the alternatives, and provides a certain balance between $r_{\left\langle  i \right\rangle}^*$, $i=1,2,\cdots,m$ and $r_{\left\langle  j \right\rangle}^*$, $j=m+1,m+2,\cdots,k$, i.e., the ratios of simulation replications allocated to the top-$m$ alternatives and other $\left(k-m\right)$ alternatives. (\ref{optrate}) further adjusts the ratios of the simulation replications allocated to some pair-wise comparison alternatives. (\ref{optrate}) and (\ref{thmbalance}) coincide with the asymptotic solutions of an approximate static optimization problem in \citet{gao2015note}, which are derived by maximizing the PCS approximated by Bonferroni inequality. In addition, (\ref{thmbalance}) is missed in~\citet{gao2016new} for normal sampling distributions. \citet{zhang2012improved,zhang2015simulation} use the rate function of PFS for selecting top-$m$ alternatives to compare the efficiency of several allocation rules in a simple case where the variances of all alternatives are equal, i.e., $\sigma _{\left\langle 1 \right\rangle }^2 = \sigma _{\left\langle 2 \right\rangle }^2 =  \cdots  = \sigma _{\left\langle k \right\rangle }^2$. Based on the rate function of the allocation procedures shown in online appendix A.0, neither the asymptotic sampling ratios of the allocation procedure in \citet{chen2008efficient} nor \citet{zhang2012improved,zhang2015simulation} guarantee asymptotic optimality.
\end{remark}

\begin{remark}
(\ref{optopt}) applies to not only normal sampling distributions. If simulation replications follow i.i.d. exponential sampling distributions with rate parameters $\widetilde {\lambda_{\widetilde \ell}}$, ${\widetilde \ell}=1,2,\cdots,k$, with rate functions derived in~\citet{gao2016optimal}, (\ref{optopt}) leads to, for $i,h = 1,2,\cdots,m$, $j,\ell = m+1,m+2,\cdots,k$,
$$\begin{aligned}
& \mathop {\min }\limits_{j = m + 1,m + 2, \cdots ,k} r_{\left\langle  h \right\rangle}^*\log \frac{{{{\widetilde \lambda }_{\left\langle  h \right\rangle}}\left( {r_{\left\langle  h \right\rangle}^* + r_{\left\langle  j \right\rangle}^*} \right)}}{{r_{\left\langle  h \right\rangle}^*{{\widetilde \lambda }_{\left\langle  h \right\rangle}} + r_{\left\langle  j \right\rangle}^*{{\widetilde \lambda }_{\left\langle  j \right\rangle}}}} + r_{\left\langle  j \right\rangle}^*\log \frac{{{{\widetilde \lambda }_{\left\langle  j \right\rangle}}\left( {r_{\left\langle  h \right\rangle}^* + r_{\left\langle  j \right\rangle}^*} \right)}}{{r_{\left\langle  h \right\rangle}^*{{\widetilde \lambda }_{\left\langle  h \right\rangle}} + r_{\left\langle  j \right\rangle}^*{{\widetilde \lambda }_{\left\langle  j \right\rangle}}}} \\
= & \mathop {\min }\limits_{i = 1,2, \cdots ,m} r_{\left\langle  i \right\rangle}^*\log \frac{{{{\widetilde \lambda }_{\left\langle  i \right\rangle}}\left( {r_{\left\langle  i \right\rangle}^* + r_{\left\langle  \ell \right\rangle}^*} \right)}}{{r_{\left\langle  i \right\rangle}^*{{\widetilde \lambda }_{\left\langle  i \right\rangle}} + r_{\left\langle  \ell \right\rangle}^*{{\widetilde \lambda }_{\left\langle  \ell \right\rangle}}}} + r_{\left\langle  \ell \right\rangle}^*\log \frac{{{{\widetilde \lambda }_{\left\langle  j \right\rangle}}\left( {r_{\left\langle  i \right\rangle}^* + r_{\left\langle  \ell \right\rangle}^*} \right)}}{{r_{\left\langle  i \right\rangle}^*{{\widetilde \lambda }_{\left\langle  i \right\rangle}} + r_{\left\langle  \ell \right\rangle}^*{{\widetilde \lambda }_{\left\langle  \ell \right\rangle}}}}~,
\end{aligned}$$
and if simulation replications follow i.i.d. Bernoulli sampling distributions with success probability $q_{\widetilde \ell}$, ${\widetilde \ell}=1,2,\cdots,k$, with rate functions derived in~\citet{glynn2004large}, (\ref{optopt}) yields, for $i,h=1,2,\cdots,m$, $j,\ell=m+1,m+2,\cdots,k$,
$$\begin{aligned}
& \mathop {\min }\limits_{j = m + 1,m + 2, \cdots ,k}  - \left( {r_{\left\langle  h \right\rangle}^* + r_{\left\langle  j \right\rangle}^*} \right)\log \left[ {{{\left( {1 - {q_{\left\langle  h \right\rangle}}} \right)}^{\frac{{r_{\left\langle  h \right\rangle}^*}}{{r_{\left\langle  h \right\rangle}^* + r_{\left\langle  j \right\rangle}^*}}}}{{\left( {1 - {q_{\left\langle  j \right\rangle}}} \right)}^{\frac{{r_{\left\langle  j \right\rangle}^*}}{{r_{\left\langle  h \right\rangle}^* + r_{\left\langle  j \right\rangle}^*}}}} + q_{\left\langle  h \right\rangle}^{\frac{{r_{\left\langle  h \right\rangle}^*}}{{r_{\left\langle  h \right\rangle}^* + r_{\left\langle  j \right\rangle}^*}}}q_{\left\langle  j \right\rangle}^{\frac{{r_{\left\langle  j \right\rangle}^*}}{{r_{\left\langle  h \right\rangle}^* + r_{\left\langle  j \right\rangle}^*}}}} \right]\\
= & \mathop {\min }\limits_{i = 1,2, \cdots ,m}  - \left( {r_{\left\langle  i \right\rangle}^* + r_{\left\langle  \ell \right\rangle}^*} \right)\log \left[ {{{\left( {1 - {q_{\left\langle  i \right\rangle}}} \right)}^{\frac{{r_{\left\langle  i \right\rangle}^*}}{{r_{\left\langle  i \right\rangle}^* + r_{\left\langle  \ell \right\rangle}^*}}}}{{\left( {1 - {q_{\left\langle  \ell \right\rangle}}} \right)}^{\frac{{r_{\left\langle  \ell \right\rangle}^*}}{{r_{\left\langle  i \right\rangle}^* + r_{\left\langle  \ell \right\rangle}^*}}}} + q_{\left\langle  i \right\rangle}^{\frac{{r_{\left\langle  i \right\rangle}^*}}{{r_{\left\langle  i \right\rangle}^* + r_{\left\langle  \ell \right\rangle}^*}}}q_{\left\langle  j \right\rangle}^{\frac{{r_{\left\langle  \ell \right\rangle}^*}}{{r_{\left\langle  i \right\rangle}^* + r_{\left\langle  \ell \right\rangle}^*}}}} \right]~.
\end{aligned}$$

Though (\ref{optopt}) is derived for general sampling distributions, computational complexity could arise in calculating  $x\left(r_{\left\langle  i \right\rangle},r_{\left\langle  j \right\rangle}\right)$ and $G_{ij}\left(r_{\left\langle  i \right\rangle},r_{\left\langle  j \right\rangle}\right)$, $i=1,2,\cdots,m$, $j=m+1,m+2,\cdots,k$. Moreover, for heavy-tailed distributions, estimating $\Lambda_{\left\langle {\widetilde \ell} \right\rangle}^*\left( x \right)$, ${\widetilde \ell}=1,2,\cdots,k$ in  $G_{ij}\left(r_{\left\langle  i \right\rangle},r_{\left\langle  j \right\rangle}\right)$ would be statistically inefficient~\citep{glynn2011ordinal}.
\end{remark}

We then establish the asymptotic optimality of the sampling ratios of the proposed AOAm (\ref{AOAm}) in the following theorem.

\begin{theorem}\label{thmasymptotic}
Suppose (\ref{optrate}) and (\ref{thmbalance}) determine a unique solution. As $t \to \infty$, with the proposed AOAm procedure, the sampling ratio of each alternative sequentially achieves the asymptotically optimal sampling ratios, i.e.,
$$\mathop {\lim }\limits_{t \to \infty } r_{\left\langle  i \right\rangle}^{\left( t \right)} = r_{\left\langle  i \right\rangle}^*,\quad a.s.,\quad i=1,2,\cdots,k~,$$
where $r_{\left\langle  i \right\rangle}^{\left( t \right)} = {{{t_{\left\langle  i \right\rangle}}} \mathord{\left/{\vphantom {{{t_i}} t}} \right.\kern-\nulldelimiterspace} t}$, $\sum\nolimits_{i = 1}^k {r_{\left\langle  i \right\rangle}^{\left( t \right)}}  = 1$, ${r_{\left\langle  i \right\rangle}^*} > 0$, and $r^*= (r_{\left\langle 1 \right\rangle}^*,r_{\left\langle 2 \right\rangle}^*,\cdots,r_{\left\langle k \right\rangle}^*)$ satisfies (\ref{optrate}) and (\ref{thmbalance}).
\end{theorem}

\begin{remark}
If the solution of  (\ref{optrate}) and (\ref{thmbalance}) is not unique, the proof of Theorem~\ref{thmasymptotic} can be adapted to prove that the clustering points of the sampling ratios following the proposed AOAm (\ref{AOAm}) procedure satisfy (\ref{optrate}) and (\ref{thmbalance}). For the special case of selecting the best alternative, i.e., $m=1$, the solution of  (\ref{glynn}) and (\ref{wscequ2}) can be proved to be unique.
\end{remark}

\section{Numerical Experiments}\label{SectionNuEx}
In this section, we conduct numerical experiments to demonstrate the efficiency of the proposed allocation policy. The proposed sequential AOAm procedure is tested with equal allocation (EA), optimal computing budget allocation for selecting top-$m$ alternatives (OCBAm) in \citet{chen2008efficient}, improved optimal computing budget allocation for selecting top-$m$ alternatives (OCBAm$+$) in \citet{zhang2012improved,zhang2015simulation}, optimal computing budget allocation for subset selection (OCBAss) in \citet{gao2015note}, and improved optimal computing budget allocation for subset selection (OCBASS) in \citet{gao2016optimal}. EA equally allocates simulation budget to estimate the performance of each alternative (roughly $T / k$ simulation replications for each alternative). OCBAm allocates simulation replications to alternatives proportional to ${r_{\left\langle i \right\rangle }}/{r_{\left\langle j \right\rangle }} = \sigma _{\left\langle i \right\rangle }^2{ ( {{\mu _{\left\langle j \right\rangle }} - c} )^2}/ ( {\sigma _{\left\langle j \right\rangle }^2{ ( {{\mu _{\left\langle i \right\rangle }} - c} )^2}} )$, where $c$ separates $\mu_{\left\langle m \right\rangle}$ and $\mu_{\left\langle {m+1}\right\rangle}$. Remark \ref{OCBAmparac} in online appendix A.0 shows that the suggested $c$ in \citet{chen2008efficient} does not perform well. We choose a better parameter $c$ in comparison. More details of the compared algorithms can be found in online appendix A.0. In all numerical experiments, we set $n_0=10$, and the statistical efficiency of the sampling procedure is measured by the IPCS and EOC, respectively, estimated by 100,000 independent macro experiments. The IPCS and EOC are reported as functions of sampling budget $T$ in each numerical experiment, i.e., ${\rm{IPCS}}_T =  \mathbb{E} [ {{\mathds{1} ( {\widehat {\mathcal{F}}_T^m = {\mathcal{F}^m}} )}} ]$ and ${\rm{EOC}}_T = \mathbb{E} [ {\sum\nolimits_{\ell \in {\widehat {\mathcal{F}}_T^m}} {{\mu _\ell}}  - \sum\nolimits_{\ell \in {\mathcal{F}^m}} {{\mu _\ell}}} ]$. The code for the numerical experiments in this paper can be found in https://github.com/gongbozhang-pku/Top-m-Selection.

\subsection{Synthetic Examples}

\textit{Experiment 1: 20 alternatives with equal variances}. We test our proposed AOAm in a synthetic example with $k=20$ and $m=5$. In each macro experiment, the performance of each alternative is generated from $\mu_{i} \sim {N\left(0,{1^2}\right)}$, $i=1,2,\cdots,20$, and the samples are drawn independently from a normal distribution $N\left(\mu_i, \sigma_i^2\right)$, where $\sigma_i^2=1$. The total simulation budget is $T=5000$.

\begin{figure}[htbp]
\centering
\subfigure[IPCS]{\includegraphics[width=0.49\textwidth]{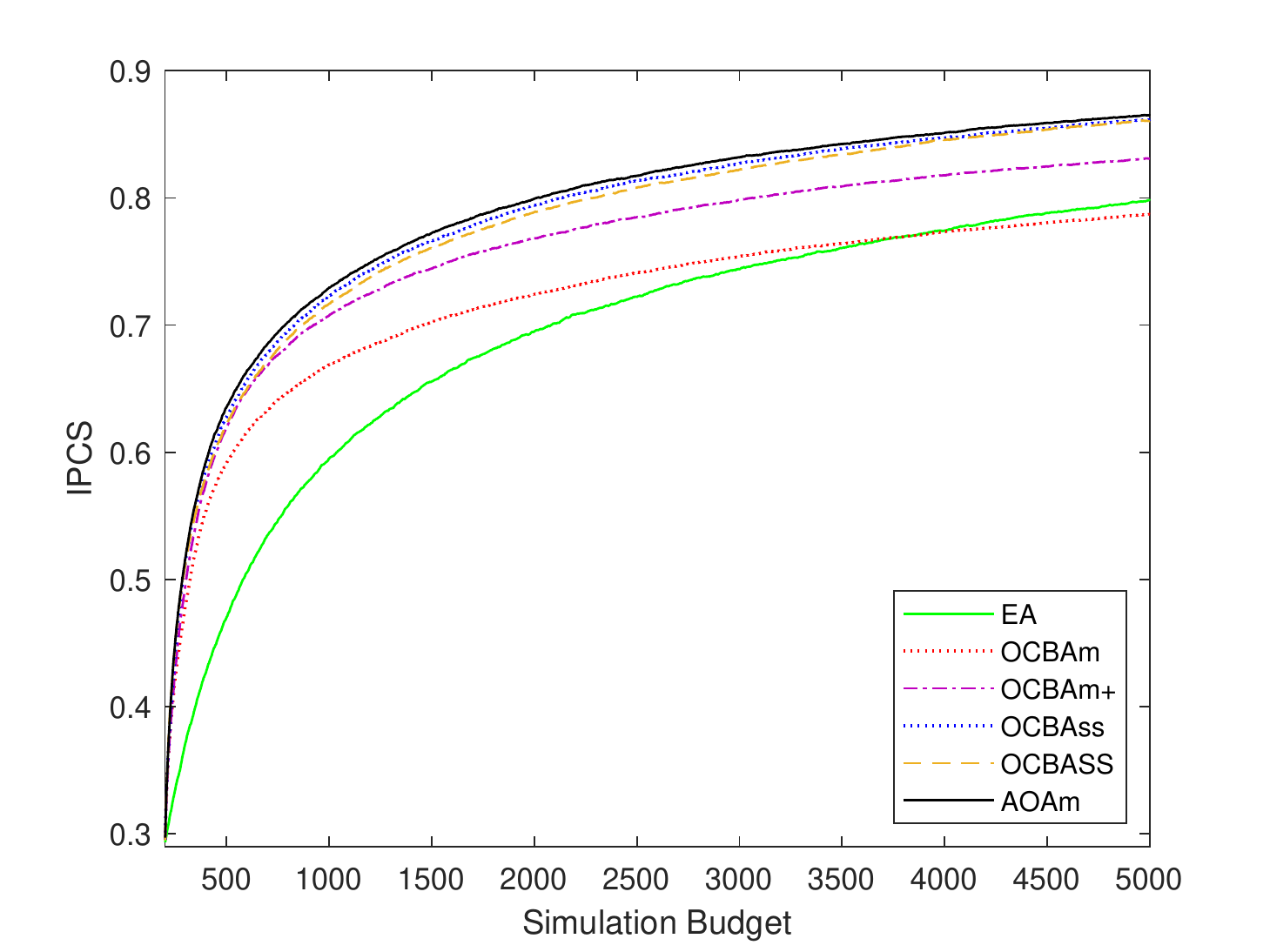}}
\subfigure[EOC]{\includegraphics[width=0.49\textwidth]{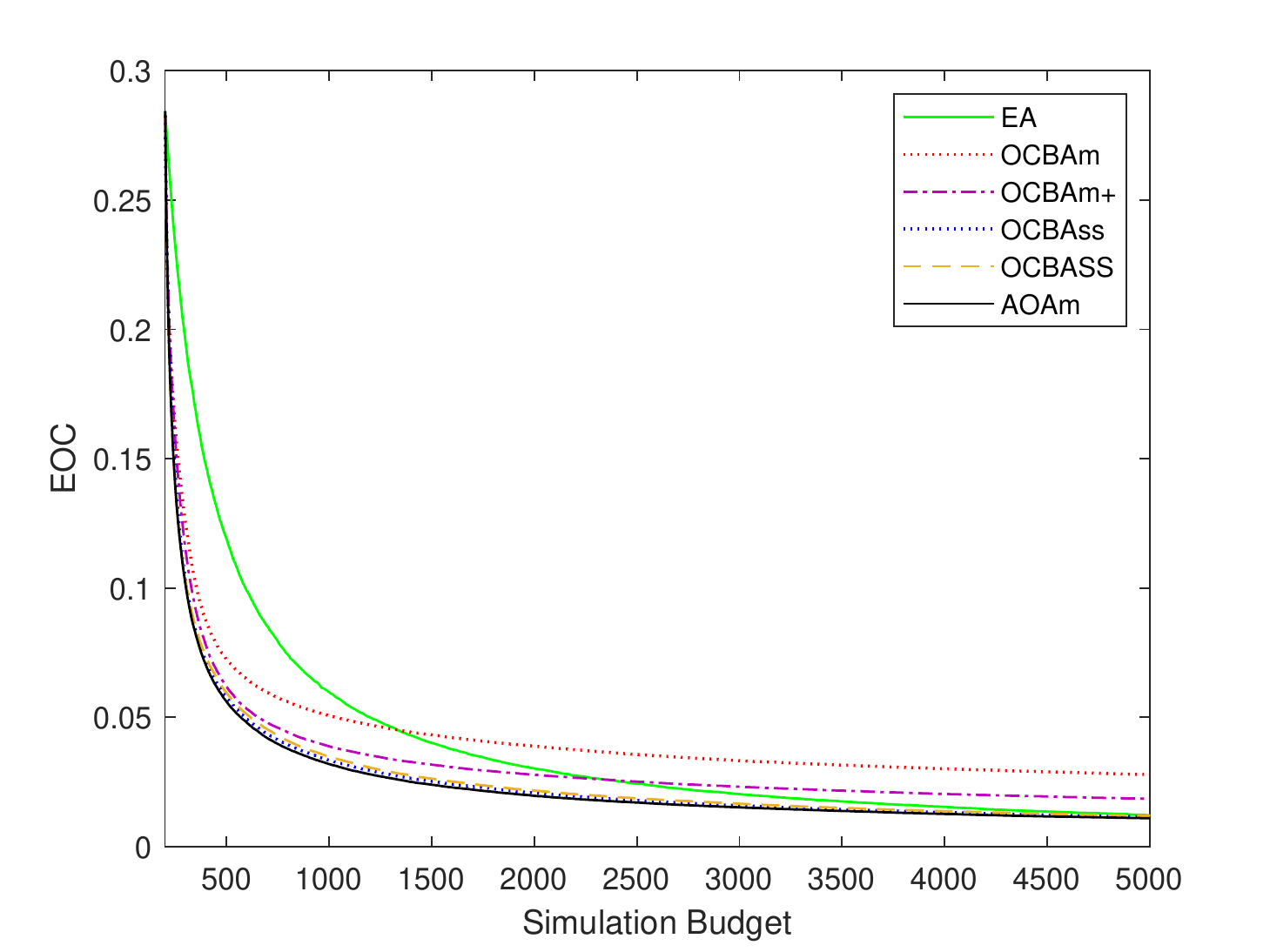}}
\caption{Comparison of IPCS and EOC for 6 sampling allocation procedures in \textit{Experiment 1}.}
\label{fig2}
\end{figure}

From Figure~\ref{fig2} (a), we can see that the IPCS of EA increases at a slow pace at the beginning and it surpasses the IPCS of OCBAm when the simulation budget reaches around 3755.  AOAm, OCBAss and OCBASS have a comparable performance and perform better than OCBAm$+$ as the simulation budget increases.  The performances of different allocation procedures in terms of EOC shown in Figure~\ref{fig2} (b) exhibit similar behaviors as measured by IPCS, and this observation is consistent in all experiments. Online Appendix A.3.1 includes the performances of OCBAss and OCBASS under perfect information (assuming parameters are known), which illustrate a significant difference between the asymptotic property and finite-sample performance for a sampling procedure.

\textit{Experiment 2: 50 alternatives with decreasing variances}. The proposed AOAm is tested in a synthetic example with $k=50$ and $m=15$, where the differences in the true means of alternatives (controlled by standard deviation) would be relatively small compared with the sampling variances. Specifically, the performance of each alternative is generated from $\mu_i \sim N (0, ({\left( {51 - i} \right)} / {\sqrt {10} } )^2 )$, $i=1,2,\cdots,50$, and the samples are drawn independently from a normal distribution $N\left(\mu_i,\sigma_i^2\right)$, where $\sigma_i^2 = \left(51 - i\right)^2$ is one order of magnitude larger than the variances of the true means. The total simulation budget is $T=12000$.

\begin{figure}[htbp]
\centering
\subfigure[IPCS]{\includegraphics[width=0.49\textwidth]{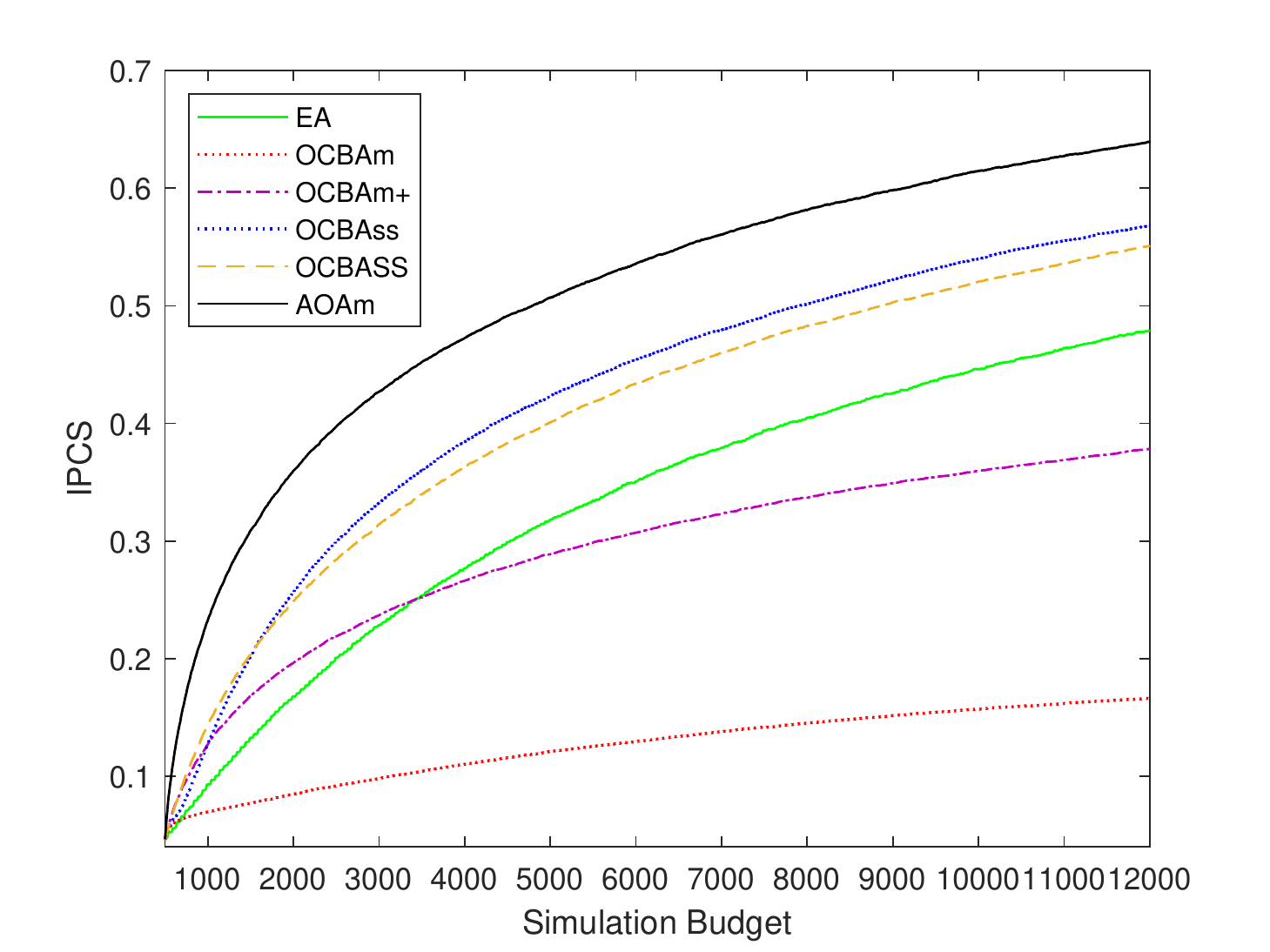}}
\subfigure[EOC]{\includegraphics[width=0.49\textwidth]{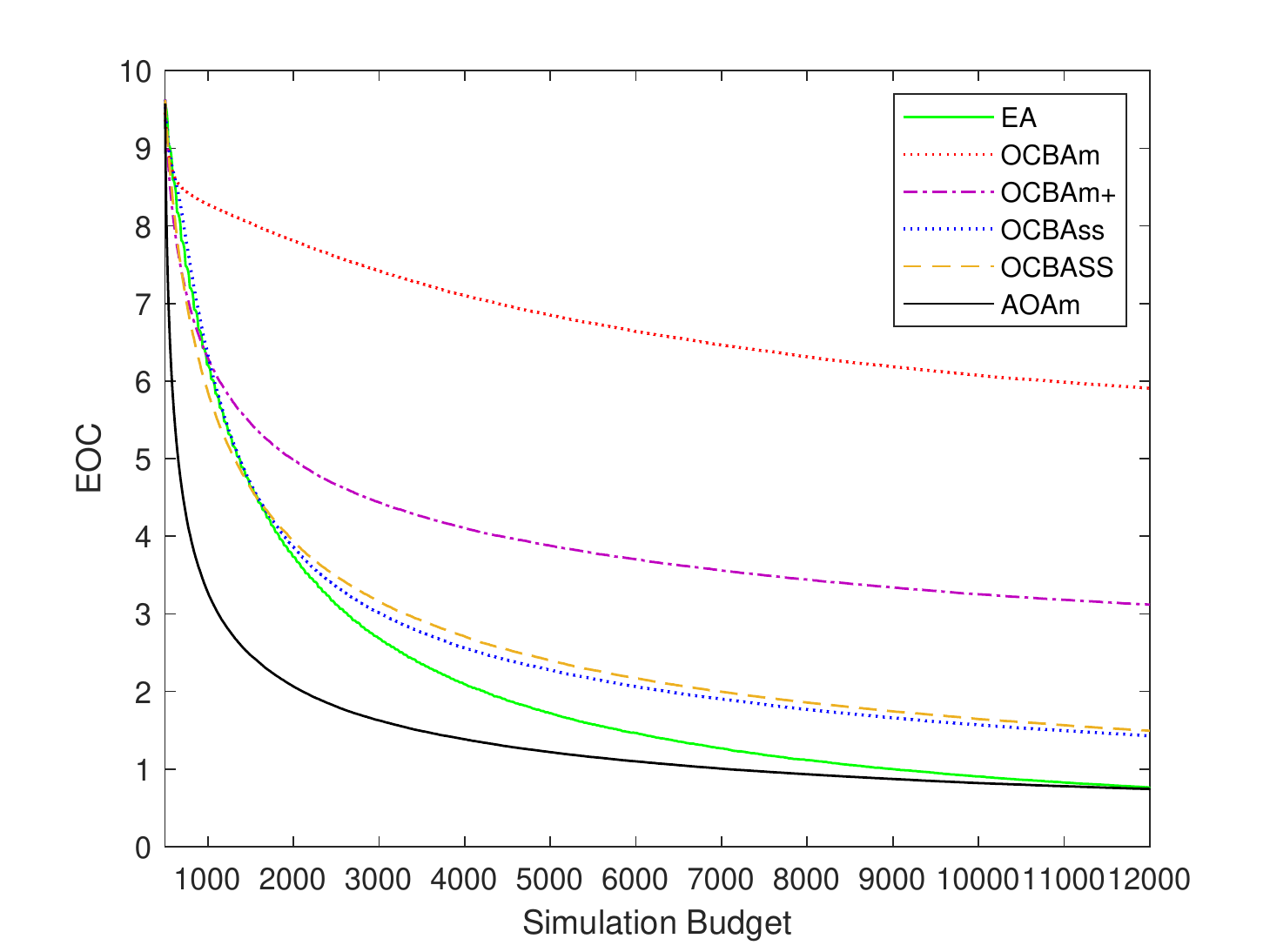}}
\caption{Comparison of IPCS and EOC for 6 sampling allocation procedures in \textit{Experiment 2}.}
\label{fig4}
\end{figure}

From Figure~\ref{fig4} (a), we can see that the IPCS of OCBAm increases at a slow pace and becomes the worst as the simulation budget grows. The IPCS of EA is the lowest at the beginning, and it surpasses OCBAm and OCBAm+ with the number of allocated simulation replications growing up to around 690 and 3300, respectively. The IPCS of OCBAss increases at a slower pace than the IPCS of OCBASS at the beginning, and it surpasses OCBAm+ and OCBASS when the simulation budget reaches around 990 and 1650, respectively. AOAm performs the best among all allocation procedures in terms of both IPCS and EOC performances. In order to attain IPCS=0.55, OCBAss requires more than 11900 simulation replications, whereas AOAm consumes 6530 simulation replications, that is, AOAm reduces the simulation budget by more than $45\%$. The performance enhancement of AOAm could be attributed to the benefit of using a stochastic dynamic programming framework and its asymptotic optimality.

\textit{Experiment 3: 50 alternatives with increasing variances}. We test our proposed AOAm in a synthetic example with $k=50$ and $m=15$, where the differences in the true means of alternatives would be relatively smaller compared with the sampling variances. Specifically, the performance of each alternaitve is generated from $\mu_i \sim N\left(0,(i/10)^2\right)$, $i=1,2,\cdots,50$, and the samples are drawn independently from a normal distribution $N\left(\mu_i,\sigma_i^2\right)$, where $\sigma_i^2 = i^2$ is two orders of magnitude larger than the variances of the true means. The total simulation budget is $T=12000$.

\begin{figure}[htbp]
\centering
\subfigure[IPCS]{\includegraphics[width=0.49\textwidth]{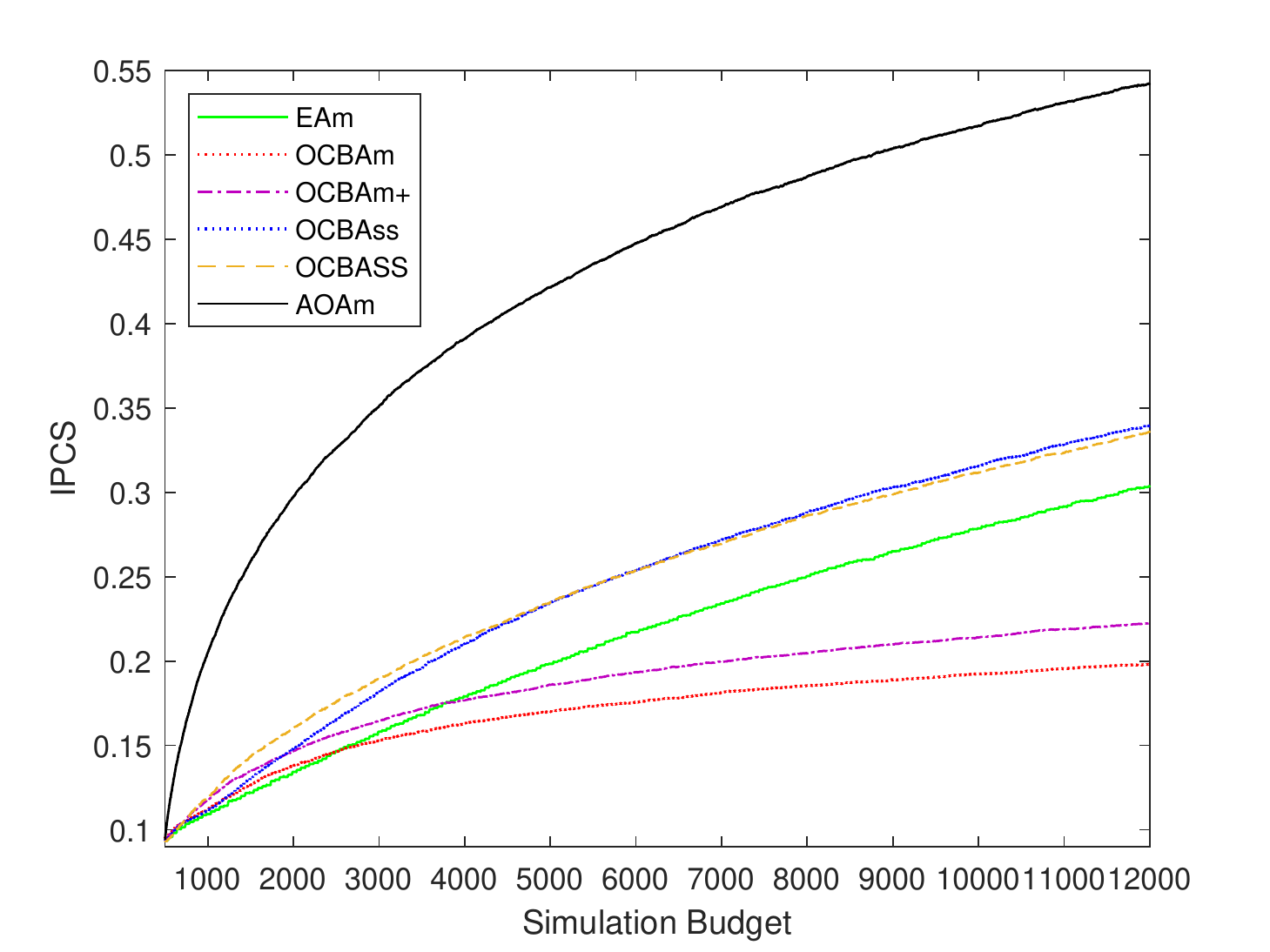}}
\subfigure[EOC]{\includegraphics[width=0.49\textwidth]{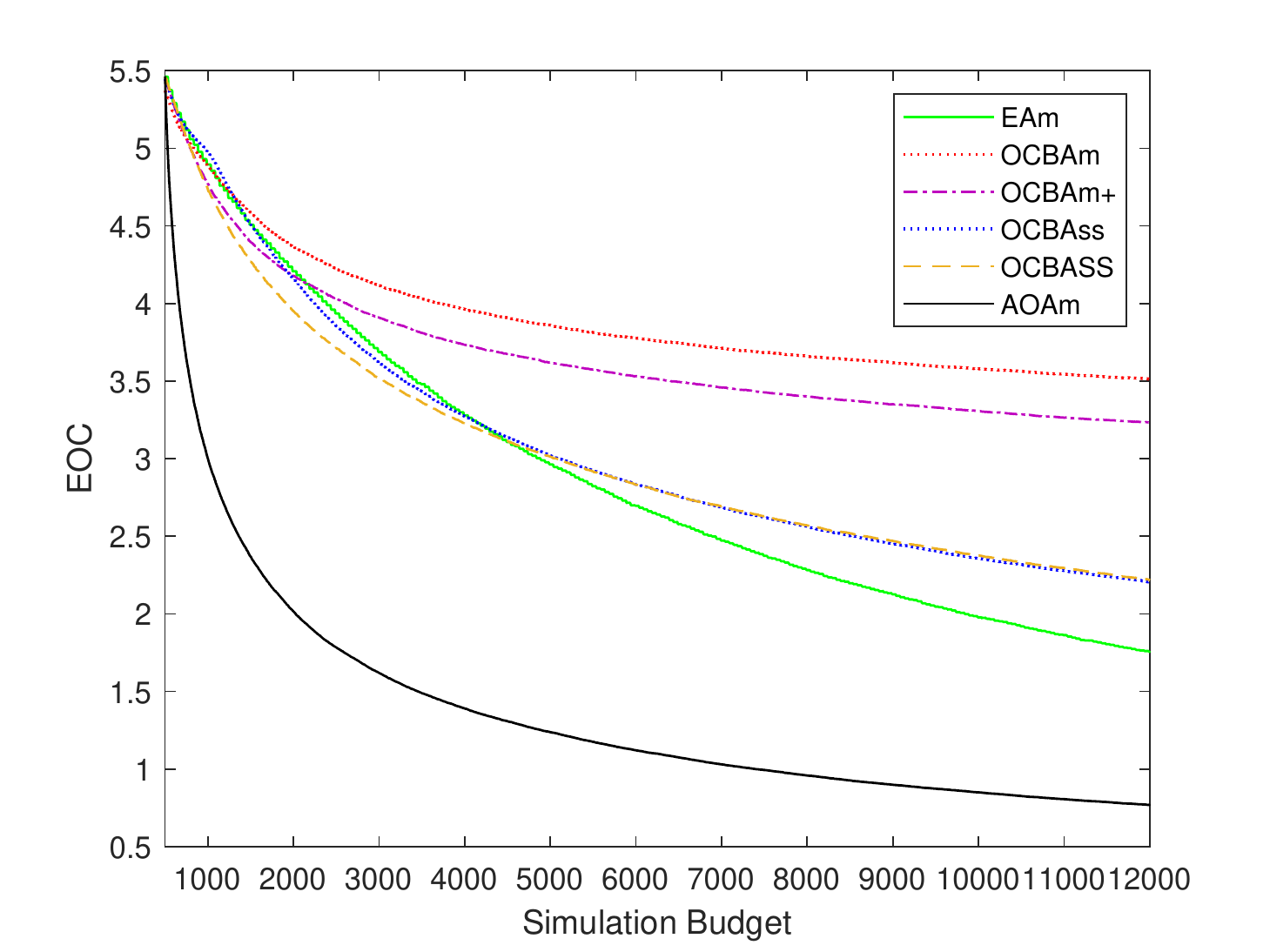}}
\caption{Comparison of IPCS and EOC for 6 sampling allocation procedures in \textit{Experiment 3}.}
\label{fig6}
\end{figure}

From Figure~\ref{fig6} (a), we can see that the IPCS of EA increases at a slow pace at the beginning, and it surpasses OCBAm and OCBAm+ when the simulation budget reaches around 2440 and 3790, respectively. OCBAm+ has an edge over OCBAm, which increases at a slow pace as the simulation budget grows. The performance of OCBAss is comparable to OCBAm at the beginning, and it surpasses OCBASS with the number of allocated simulation replications growing up to around 5275. AOAm performs significantly better than the other 5 allocation procedures. In order to attain IPCS = 0.3, AOAm consumes 2030 simulation replications, whereas OCBAss requires more than 9070 simulation replications, that is, AOAm reduces the simulation budget by more than $77\%$. It can be noticed that the advantage of AOAm is more significant when the difficulty of the problem become higher, i.e., the differences between competing alternatives are small and variances are large.

To illustrate the effect of parameter $m$ on the performance of each allocation procedure, under the same numerical settings as \textit{Experiment 3}, Table \ref{tableA21} reports the IPCS and EOC obtained by each sampling procedure using nine different $m$ values: $5,10,15,20,25,30,35,40,45$. From Tabel~\ref{tableA21}, we can see that the performance of AOAm is superior to other sampling procedures for all $m$ values. The IPCS of OCBAm surpasses EA as the value of parameter $m$ increases. The performance of OCBAm+ is better than the performance of OCBAm. The IPCS of OCBAss has a slight edge over OCBASS that is better than OCBAm+. In addition, numerical results show that, for a fixed simulation budget $T$, the IPCS for each sampling procedure does not monotonically increase with the value of parameter $m$ grows, which could be attributed to the reason that the difficulty of correctly determining the top-$m$ alternatives increases at first and then decreases as the value of parameter $m$ grows.

\begin{table}[htbp]
\caption{IPCS and EOC of different procedures with various $m$  values (simulation budget $T=12000$)}
\label{tableA21}
\centering
\begin{tabular}{cccccccc}
\hline
   $m$ & Metric & EA & OCBAm & OCBAm+ & OCBAss & OCBASS & AOAm  \\
\hline
   5 & \tabincell{c}{IPCS \\ EOC}  & \tabincell{c}{0.2887 \\ 2.3019} & \tabincell{c}{0.1959 \\ 4.5360} & \tabincell{c}{0.2189 \\ 4.1712} & \tabincell{c}{0.3454 \\ 2.7269} & \tabincell{c}{0.3442 \\ 2.7104} & \tabincell{c}{0.5489 \\ 0.9034} \\
   10 & \tabincell{c}{IPCS \\ EOC} & \tabincell{c}{0.2587 \\ 2.2409} & \tabincell{c}{0.1560 \\ 4.6196} & \tabincell{c}{0.1699 \\ 4.3281} & \tabincell{c}{0.2849 \\ 2.9221} & \tabincell{c}{0.2836 \\ 2.9514} & \tabincell{c}{0.4903 \\ 1.0396} \\
   15 & \tabincell{c}{IPCS \\ EOC} & \tabincell{c}{0.3036 \\ 1.7562} & \tabincell{c}{0.1982 \\ 3.5127} & \tabincell{c}{0.2224 \\ 3.2348} & \tabincell{c}{0.3393 \\ 2.2053} & \tabincell{c}{0.3361 \\ 2.2176} & \tabincell{c}{0.5421 \\ 0.7689} \\
   20 & \tabincell{c}{IPCS \\ EOC} & \tabincell{c}{0.3716 \\ 1.2707} & \tabincell{c}{0.2654 \\ 2.4875} & \tabincell{c}{0.3060 \\ 2.1947} & \tabincell{c}{0.4268 \\ 1.4669} & \tabincell{c}{0.4221 \\ 1.4838} & \tabincell{c}{0.6215 \\ 0.4951} \\
   25 & \tabincell{c}{IPCS \\ EOC} & \tabincell{c}{0.4467 \\ 0.8780} & \tabincell{c}{0.3595 \\ 1.6007} & \tabincell{c}{0.4183 \\ 1.3388} & \tabincell{c}{0.5381 \\ 0.8662} & \tabincell{c}{0.5290 \\ 0.8830} & \tabincell{c}{0.6978 \\ 0.2979} \\
   30 & \tabincell{c}{IPCS \\ EOC} & \tabincell{c}{0.5418 \\ 0.5459} & \tabincell{c}{0.4800 \\ 0.9188} & \tabincell{c}{0.5546 \\ 0.7132} & \tabincell{c}{0.6618 \\ 0.4405} & \tabincell{c}{0.6529 \\ 0.4504} & \tabincell{c}{0.7769 \\ 0.1516} \\
   35 & \tabincell{c}{IPCS \\ EOC} & \tabincell{c}{0.6487 \\ 0.2954} & \tabincell{c}{0.6311 \\ 0.4293} & \tabincell{c}{0.7221 \\ 0.2806} & \tabincell{c}{0.7857 \\ 0.1733} & \tabincell{c}{0.7816 \\ 0.1791} & \tabincell{c}{0.8536 \\ 0.0658} \\
   40 & \tabincell{c}{IPCS \\ EOC} & \tabincell{c}{0.7760 \\ 0.1200} & \tabincell{c}{0.8066 \\ 0.1253} & \tabincell{c}{0.8785 \\ 0.0645} & \tabincell{c}{0.8989 \\ 0.0416} & \tabincell{c}{0.8968 \\ 0.0437} & \tabincell{c}{0.9249 \\ 0.0177} \\
   45 & \tabincell{c}{IPCS \\ EOC} & \tabincell{c}{0.9212 \\ 0.0205} & \tabincell{c}{0.9596 \\ 0.0095} & \tabincell{c}{0.9759 \\ 0.0029} & \tabincell{c}{0.9772 \\ 0.0024} & \tabincell{c}{0.9761 \\ 0.0026} & \tabincell{c}{0.9782 \\ 0.0018} \\
\hline
\end{tabular}
\end{table}

\textit{Experiment 4: 100 alternatives with increasing variances}. In this example, our proposed AOAm is tested in a larger synthetic case with $k=100$ and $m=15$. The numerical settings are the same as~\textit{Experiment 3}. The total simulation budget is $T=200,000$.

\begin{figure}[htbp]
\centering
\subfigure[IPCS]{\includegraphics[width=0.49\textwidth]{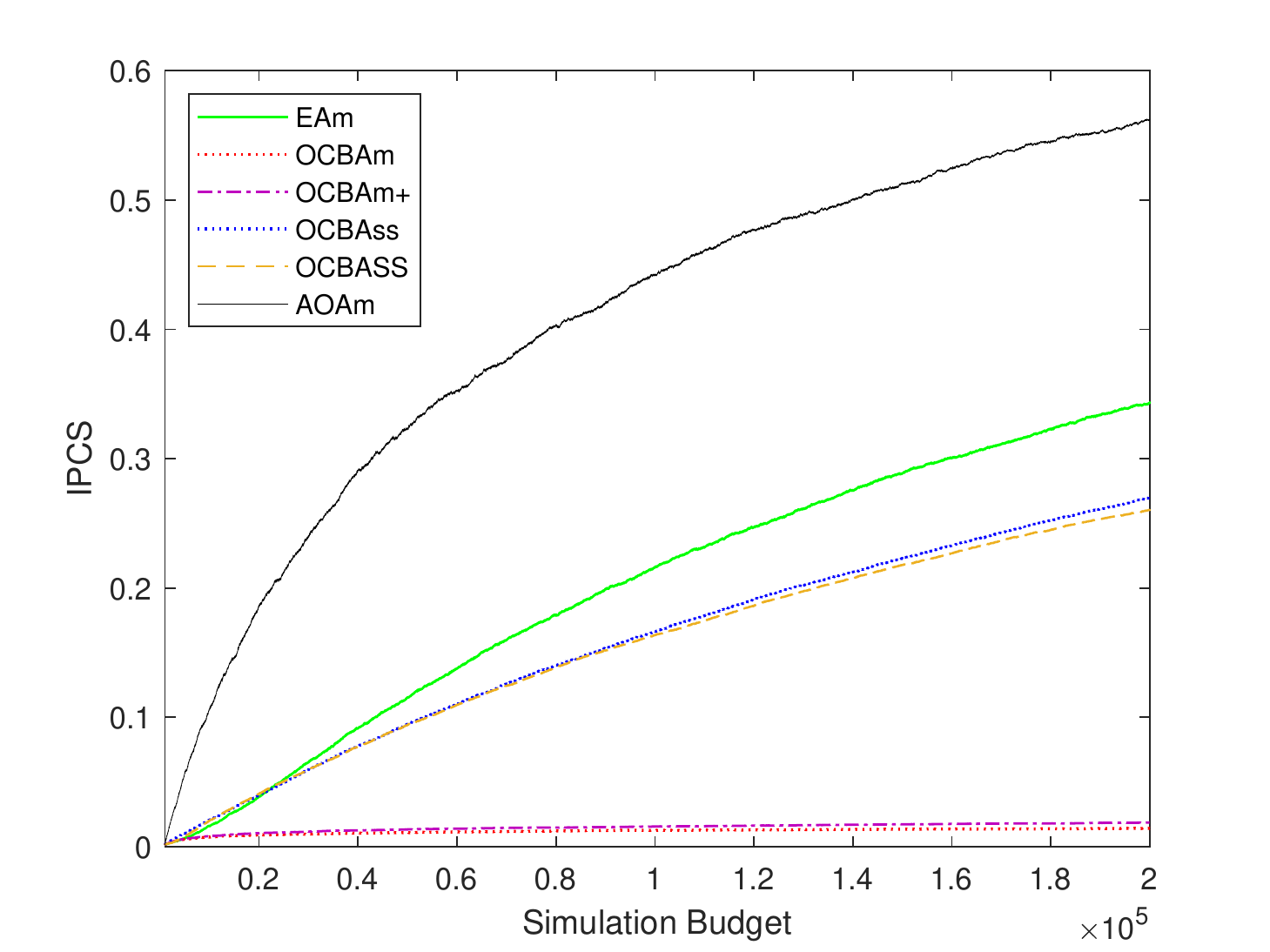}}
\subfigure[EOC]{\includegraphics[width=0.49\textwidth]{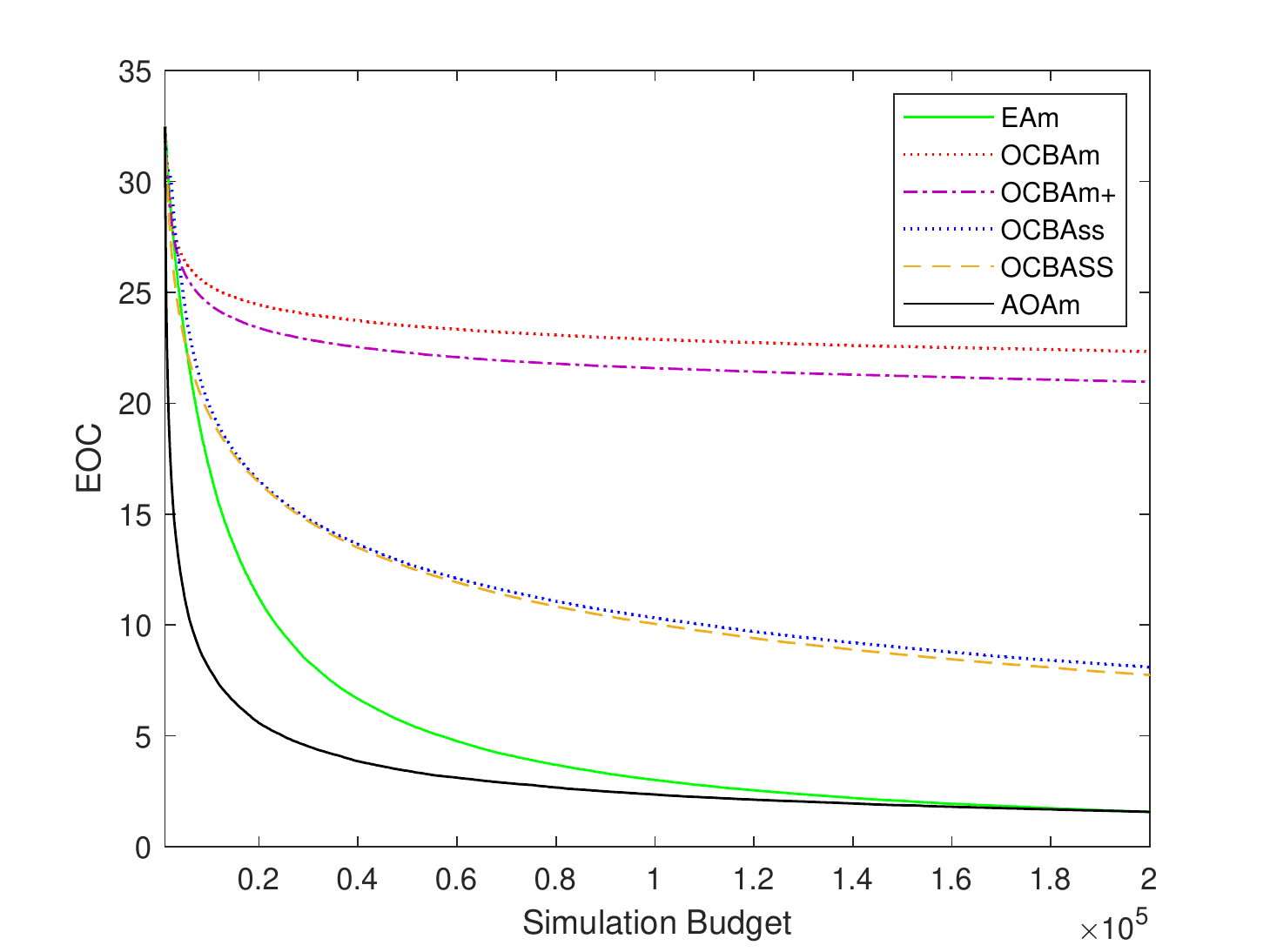}}
\caption{Comparison of IPCS and EOC for 6 sampling allocation procedures in \textit{Experiment 4}.}
\label{fig8}
\end{figure}

From Figure~\ref{fig8} (a), we can see that the IPCS of EA increases at a slow pace at the beginning, and it surpasses OCBAm, OCBAm+, OCBAss and OCBASS when the simulation budget reaches around 2380, 3800, 22180 and 22185, respectively. OCBAm+ has a slight edge over OCBAm, which flattens as the simulation budget grows. OCBAss and OCBASS have a comparable performance at the beginning, and the former surpasses the latter when the simulation budget reaches 46337. AOAm is superior to other allocation procedures in comparison. AOAm requires 31808 simulation replications to attain IPCS = 0.25, whereas OCBAss requires 177315 simulation replications, that is, AOAm reduces the simulation budget by more than $82\%$. AOAm requires around $1.39 \times 10^{5}$ simulation replications to attain IPCS = 0.5, whereas EAm and OCBAss cannot achieve the same IPCS level even when the simulation budget reaches $2 \times 10^5$. Compared with~\textit{Experiment 3}, the advantage of AOAm appears to be more significant as the number of competing alternatives increases.

\subsection{Emergency Evacuation Planning Example}

\begin{figure}[htbp]
\center{\includegraphics[width=0.8\textwidth]{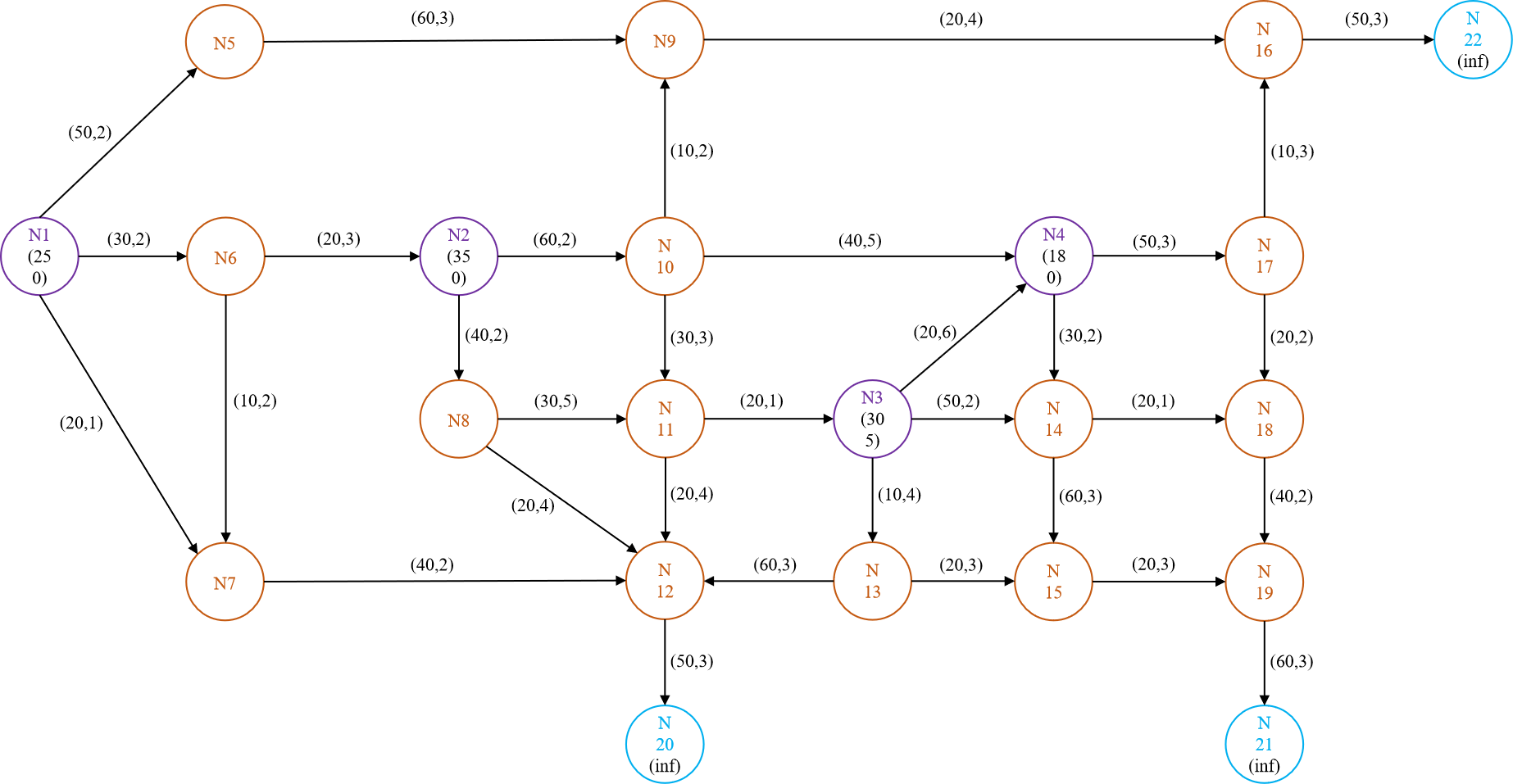}}
\caption{Evacuation network based on a real map of a region in Shanghai, China.  The bivariate vector on the top of each edge comprises of expected capacity $\mathbb{E}\left( {c_{ij}} \right)$ and expected travel time $\mathbb{E}\left( {t_{ij}} \right)$, respectively.    }
\label{figrealnetwork}
\end{figure}

When disasters occur, evacuating evacuees to safe areas is often the top priority. Preparing some evacuation plans in advance can make evacuation be executed orderly. The goal of an evacuation is to transfer all evacuees to safe areas using less time. The evacuation network is modeled as a directed graph $G^{\prime} = \left\langle {V^{\prime},E^{\prime}} \right\rangle$, where nodes $v_{i} \in {V^{\prime}}$, $i=1,2,\cdots,\left| {V'} \right|$ represent evacuation areas and edges $e_{ij} \in E^{\prime}$, $i,j\in \left\{1,2,\cdots,\left| {V'} \right|\right\}$ represent evacuation routes. $V^{\prime}$ contains three types of nodes: source nodes which represent affected areas, destination nodes which represent safe areas and intermediate nodes which represent intermediate evacuation areas. The capacity of each node is unlimited and each edge has two attributes: capacity $c_{ij}$ and travel time $t_{ij}$, which are uncertain. An evacuation network is shown in Figure~\ref{figrealnetwork}, which is based on a transportation network of a region in  Shanghai, China, shown in online Appendix A.3.2. The evacuation network consists of 22 nodes and 33 edges, with 4 source nodes: $N_{1},N_{2},N_{3},N_{4}$ and 3 destination nodes: $N_{20},N_{21},N_{22}$. The number of evacuees at source nodes are $NES = (250, 350, 305, 180)$. The capacity and travel time of each route are assumed to follow i.i.d. log-normal distribution with means $\mathbb{E}\left( {c_{ij}} \right)$ and $\mathbb{E}\left( {t_{ij}} \right)$, respectively, reported on the top of each edge in Figure~\ref{figrealnetwork}, and variances $Var\left( {{c_{ij}}} \right) = 1$ and $Var\left( {{t_{ij}}} \right) = 0.01$, respectively.

A feasible evacuation plan aims to schedule all evacuees to any of destination nodes under the capacity constraints. The evacuees at each source node are sent by groups. Let $P_{s}^r\mathop  = \limits^\Delta \{ {n_0^r = {N_s},n_1^r,n_2^r, \cdots ,n_{h_{sd}^r}^r,n_{h_{sd}^r + 1}^r = {N_d}} \}$, $s \in \left\{1,2,3,4\right\}$, $d \in \left\{20,21,22\right\}$, $r=1,2,\cdots,R_{s}$ denotes the $r$-th evacuation path of a source node $s$, where $h_{sd}^r$  is the number of intermediate nodes in the evacuation path and $R_{s}$ is the number of total evacuation paths from a source node $s$. The total travel time and maximum flow of an evacuation path are $T_s^r = \sum\nolimits_{i = 0}^{h_{sd}^r} {{t_{n_i^r,n_{i + 1}^r}}}$ and $F_s^r = \mathop {\min }\limits_{i = 0,1, \cdots ,h_{sd}^r} {c_{n_i^r,n_{i + 1}^r}}$, respectively. For each source node $s$, the path with shorter travel time $T_s^r$ has higher priority to evacuate evacuees. The route congestion is considered when evacuation paths have common routes, and evacuation order of evacuees when route congestion exists is determined by the evacuation time of each source node,
$$T_{s} = \left\lceil {\frac{{NE{S_s} + \sum\nolimits_{r = 1}^{{R_s}} {\left( {T_s^r - 1} \right) \cdot F_s^r} }}{{\sum\nolimits_{r = 1}^{{R_s}} {F_s^r} }}} \right\rceil,\quad s=1,2,3,4,\;r=1,2,\cdots,R_s~,$$
where $NE{S_s}$ is the $s$-th element of $NES$ and $\left\lceil a \right\rceil$ is a ceiling function mapping $a$ to the least integer greater than or equal to $a$. The evacuees from a source node $s$ with longer $T_{s}$ have higher priority to be evacuated when route congestion exists. The evacuation plan uses a user equilibrium model, which balances the evacuation time of each evacuation path $P_s^r$ under an assumption that evacuees have full information on capacity and travel time of each route. The clearance time of each evacuation plan is obtained by the evacuation times that the evacuees reaching any safe area through different evacuation paths. The example has also been considered in \citet{zhang2018emergency}, and the parameters are set based on \citet{zhang2018emergency}.

The expected clearance time of each evacuation plan is calculated by the average of $10^3$ clearance time. We use 2 evacuation paths simultaneously for each source node, i.e., $R_{s} = 2$, $s = 1,2,3,4$, which could reduce the clearance time compared with evacuating evacuees via a single evacuation path. 81 alternatives are summarized in online appendix A.3.2, and we aim to find the top-5 evacuation plans with the shortest expected clearance time, which are shown in Table~\ref{table2}, determined by $10^6$ macro experiments for estimating the expected clearance time of each alternative. Since the performance of each alternative is deterministic, the effectiveness of each allocation procedure is measured by the classical PCS and EOC estimated by 100,000 independent macro experiments. The total simulation budget is $T = 1500$.

\begin{figure}[htbp]
\centering
\subfigure[PCS]{\includegraphics[width=0.49\textwidth]{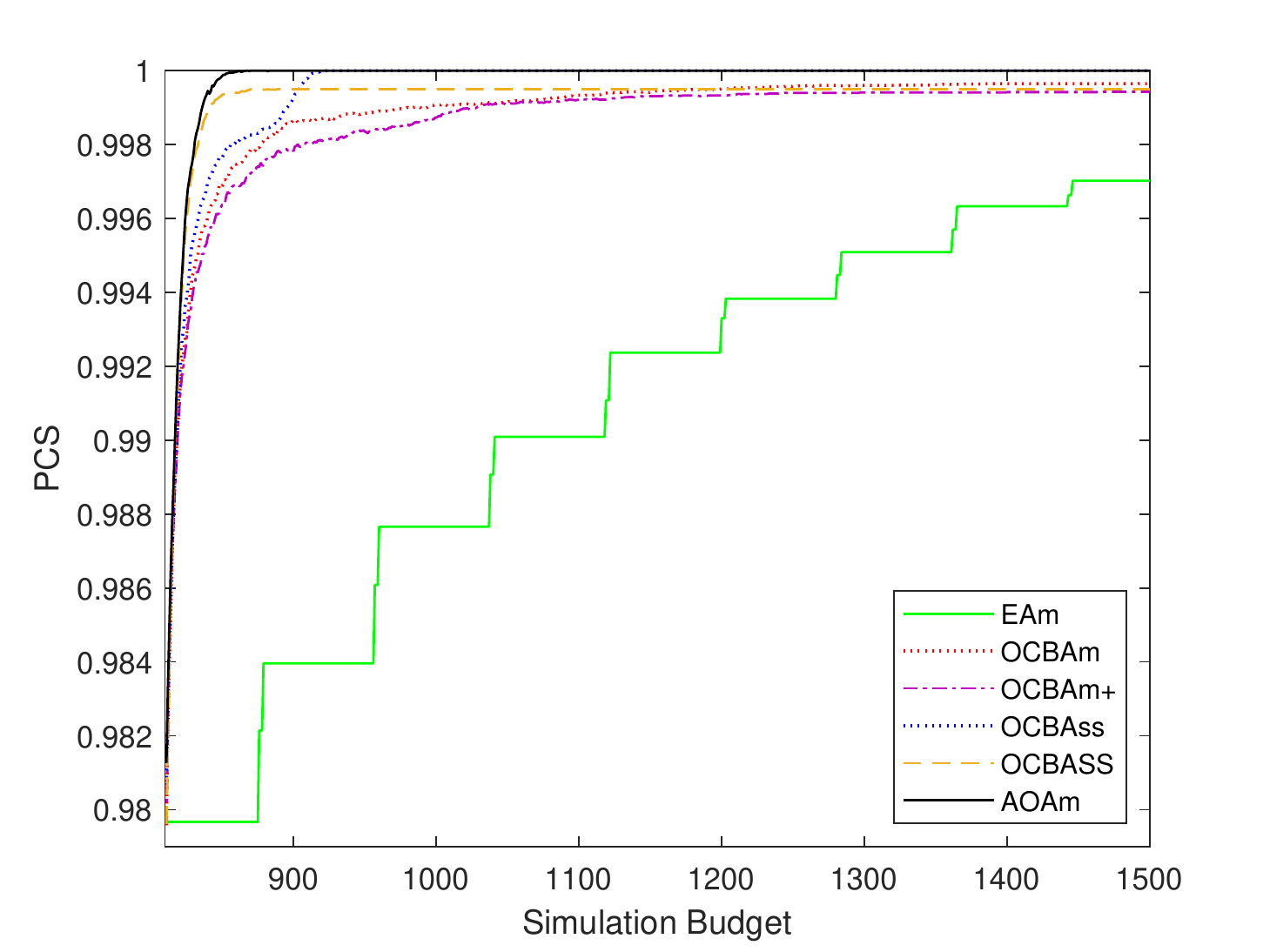}}
\subfigure[EOC]{\includegraphics[width=0.49\textwidth]{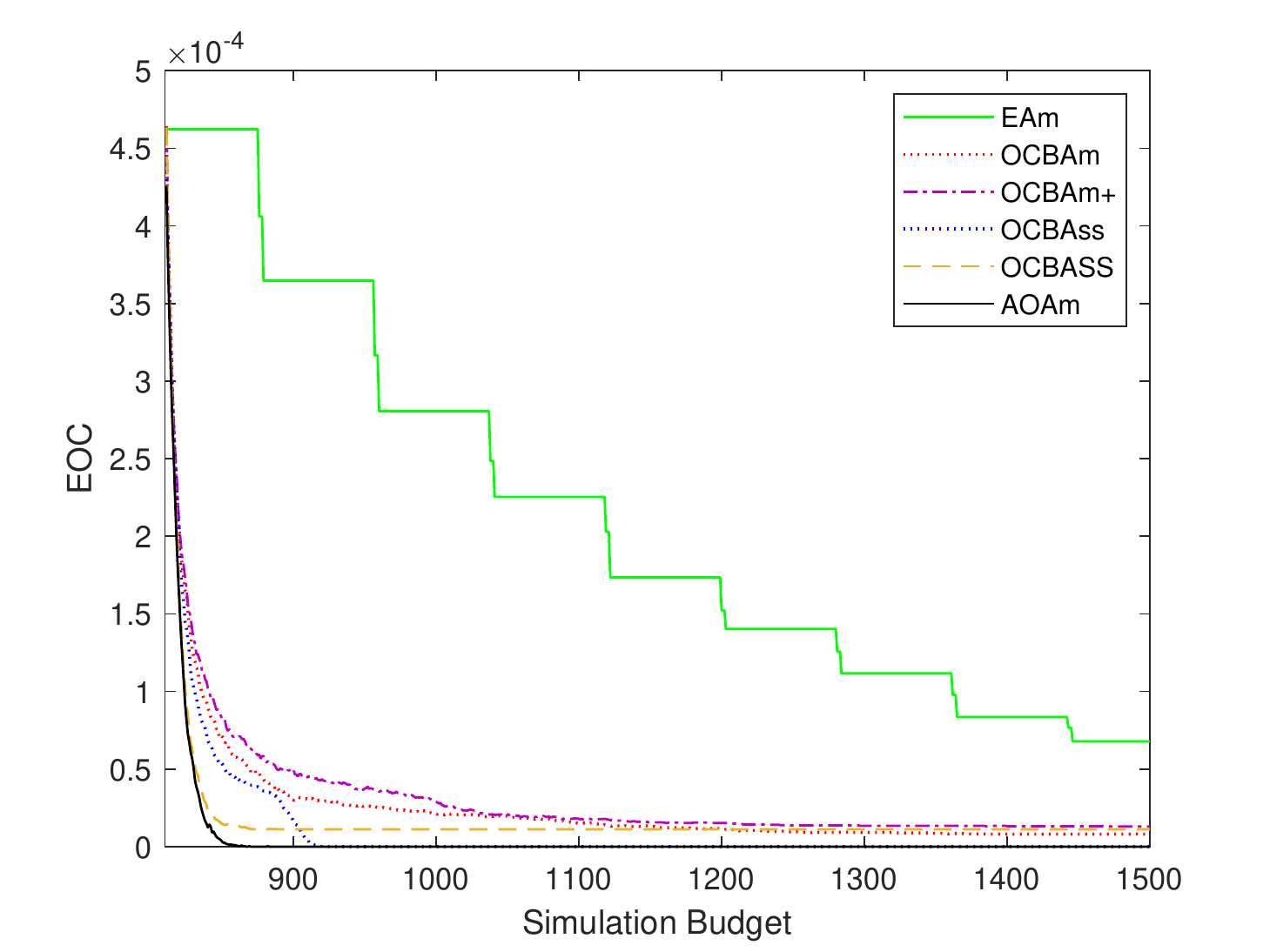}}
\caption{Comparison of PCS and EOC for 6 sampling allocation procedures in the emergency evacuation planning example.}
\label{fig9}
\end{figure}

From Figure~\ref{fig9} (a), we can see that EA performs the worst among all allocation procedures. OCBAm has an edge over OCBAm+ at the beginning, and the latter has a comparable performance with the the former as the simulation budget grows. The PCS of OCBASS increases at a faster pace than the PCS of OCBAss at the beginning, and the latter surpasses the former when the simulation budget reaches around 905. AOAm achieves the best performance among
all allocation procedures, 
and the PCS of AOAm approaches 1 when the simulation budget
reaches around 870.

\section{Conclusions}\label{SectionCon}

We study a dynamic sampling problem for selecting the top-$m$ alternatives from a finite set of alternatives. We rigorously define the asymptotically optimal sampling ratios which optimize large deviations rate of PFS for selecting top-$m$ alternatives. Under a Bayesian framework, a stochastic dynamic programming problem is proposed to capture the dynamic sampling decision, and an efficient sampling procedure named as AOAm is derived by maximizing a VFA one-step look ahead, which is proved to be consistent and asymptotically optimal. Numerical results show that the proposed AOAm is superior to the existing methods developed for the problem setting in this work. The superiority of AOAm is more pronounced when the differences between competing alternatives are small and variances are large. How to develop fast learning schemes when the simulation replications are obtained in parallel deserves further study.

%
%
%

\ACKNOWLEDGMENT{This work was supported in part by the National Science Foundation of China (NSFC) under Grants 72022001 and 71901003, and by the US National Science Foundation (NSF) under Grant DMS-2053489. A preliminary version of this work has been published in Proceedings of 2021 Winter Simulation Conference~\citep{zhang2021wsc}.}

\bibliographystyle{informs2014}
\bibliography{myref}

\clearpage

\section*{Online Supplements}

\subsection*{A.0. OCBA Policies for Selecting top-$m$ Alternatives }

In this section, we summarize some details of the existing OCBA policies for selecting top-$m$ alternatives, including OCBAm, OCBAm$+$, OCBAss and OCBASS.
\begin{enumerate}[(1)]
\item OCBAm: \citet{chen2008efficient} approximate the posterior PCS with a pre-specified parameter $c$, which lies between $\mu _{{\left\langle m \right\rangle }}$ and $\mu _{\left\langle m+1 \right\rangle }$, i.e.,
$$\mathop {\max }\limits_{{r_1},{r_2}, \cdots ,{r_k}} \prod\limits_{i = 1,2, \cdots ,m} {\Pr \left\{ {\left. {{\mu _{{{\left\langle i \right\rangle }_T}}} > c} \right|{\mathcal{E}_T}} \right\}} \prod\limits_{j = m + 1,m+2, \cdots ,k} {\Pr \left\{ {\left. {{\mu _{{{\left\langle j \right\rangle }_T}}} < c} \right|{\mathcal{E}_T}} \right\}}~,$$
and derive an asymptotic solution of the approximated optimization problem optimal in the limit. A sequential allocation procedure is heuristically proposed based on the asymptotic solution, where the unknown parameters are sequentially estimated from samples. Specifically, at each step, OCBAm allocates a simulation replication to an alternative according to the allocation rule:
$$\frac{{\bar r}_{i}^{\left( t \right)}}{{\bar r}_{j}^{\left( t \right)}} = {\left( {\frac{{{{\widehat{\sigma}_{i}} / {({\bar{X}_{i}} - c )}}}}{{{{{\widehat{\sigma} _{j}}} / {({{\bar X}_{j}} - c)}}}}}\right)^2},\quad i,j = 1,2,\cdots,k~,$$
where
\begin{equation}\label{corg}
c = \frac{{{{\bar \sigma }_{\left\langle {m + 1} \right\rangle }}{\bar{X} _{{\left\langle m \right\rangle }}} + {{\bar \sigma }_{\left\langle m \right\rangle }}{\bar{X}_{\left\langle {m + 1} \right\rangle }}}}{{{{\bar \sigma }_{\left\langle m \right\rangle }} + {{\bar \sigma }_{\left\langle {m + 1} \right\rangle }}}}~,
\end{equation}
where ${{\bar \sigma }_{\ell}} = {{{\widehat {\sigma}_{\ell}}} / {\sqrt {{t_{\ell}}} }}$, and the allocated alternative is chosen by the ``most starving" sequential rule, i.e., ${A_{t + 1}}\left( {{\mathcal{E} _t}} \right) = \arg \mathop {\max }\limits_{i = 1,2, \cdots ,k} \{ {\left( {t + 1} \right){\bar r}_i^{\left( t \right)} - {t_i}} \}$. Typically, the performances of alternatives that are close to $c$ will receive more simulation replications. The performance of OCBAm is highly sensitive to the value of $c$, and OCBAm cannot reduce to OCBA in \citet{chen2011stochastic} for selecting the best alternative, when $m=1$. \citet{zhang2015simulation} shows that when the variances of all alternatives are equal, i.e., $\sigma _{\left\langle 1 \right\rangle }^2 = \sigma _{\left\langle 2 \right\rangle }^2 =  \cdots  = \sigma _{\left\langle k \right\rangle }^2 = {\sigma ^2}$, the large deviations rate of PFS of OCBAm is
$$\min \left\{ {\frac{{{{\left( {{\mu _{\left\langle m \right\rangle }} - {\mu _{\left\langle k \right\rangle }}} \right)}^2}}}{{2{\sigma ^2}\left( {{1 / {{{\bar r}_{\left\langle m \right\rangle}}}} + {1 / {{{\bar r}_{\left\langle k \right\rangle}}}}} \right)}},\frac{{{{\left( {{\mu _{\left\langle 1 \right\rangle }} - {\mu _{\left\langle {m + 1} \right\rangle }}} \right)}^2}}}{{2{\sigma ^2}\left( {{1 / {{{\bar r}_{\left\langle 1 \right\rangle}}}} + {1 / {{{\bar r}_{\left\langle m+1 \right\rangle}}}}} \right)}}} \right\}~,$$
where ${\bar r}_{\left\langle i \right\rangle}$, $i=1,2,\cdots,k$ is the allocation ratio of the OCBAm rule.

\begin{remark}\label{OCBAmparac}
Note that ${\bar \sigma}_\ell$ in (\ref{corg}) depends on the number of allocated simulation replications $t_\ell$, and alternatives $\left\langle m \right\rangle$ and $\left\langle{m+1} \right\rangle$ will receive most of simulation replications when sequentially implementing the OCBAm. We suggest ${\bar \sigma}_\ell$ in $c$ are independent of $t_\ell$ at each step, i.e.,
\begin{equation}\label{cadp}
c = \frac{{{{\widehat \sigma }_{\left\langle {m + 1} \right\rangle }}{\bar{X} _{{\left\langle m \right\rangle }}} + {{\widehat \sigma }_{\left\langle m \right\rangle }}{\bar{X}_{\left\langle {m + 1} \right\rangle }}}}{{{{\widehat \sigma }_{\left\langle m \right\rangle }} + {{\widehat \sigma }_{\left\langle {m + 1} \right\rangle }}}}~.
\end{equation}

We then test our suggested $c$ in (\ref{cadp}) with some numerical experiments. The following allocation procedures are tested: EA; sequential OCBAm with parameter $c$ in (\ref{corg}) (OCBAm(sequential)), which allocates simulation replications based on OCBAm rule with parameters estimated sequentially by available sample information; sequential OCBAm with parameter $c$ in (\ref{cadp}) (OCBAm), which allocates simulation replications based on OCBAm rule with parameters estimated sequentially by available sample information; two-stage OCBAm with parameter $c$ in (\ref{corg}) (OCBAm(two-stage)), where the initial simulation replications are equally allocated to estimate unknown parameters and then the rest simulation replications are allocated according to OCBAm rule. We set $n_0=20$. PCS and EOC are estimated based on 100,000 independent macro experiments.

\textit{Experiment A.0.1: 10 alternatives with equal variances.}  We test the allocation procedures with $k=10$ and $m=3$. The performance of each alternative is $\mu_i = i$, $i=1,2,\cdots,10$, and the samples are drawn independently from a normal distribution $N\left(\mu_i,\sigma_i^2\right)$, where $\sigma_i=6$. The total simulation budget is $T=8000$. The numerical settings are the same as \textit{Example 1} in~\citet{chen2008efficient}.

\begin{figure}[htbp]
\centering
\subfigure[PCS]{\includegraphics[width=0.49\textwidth]{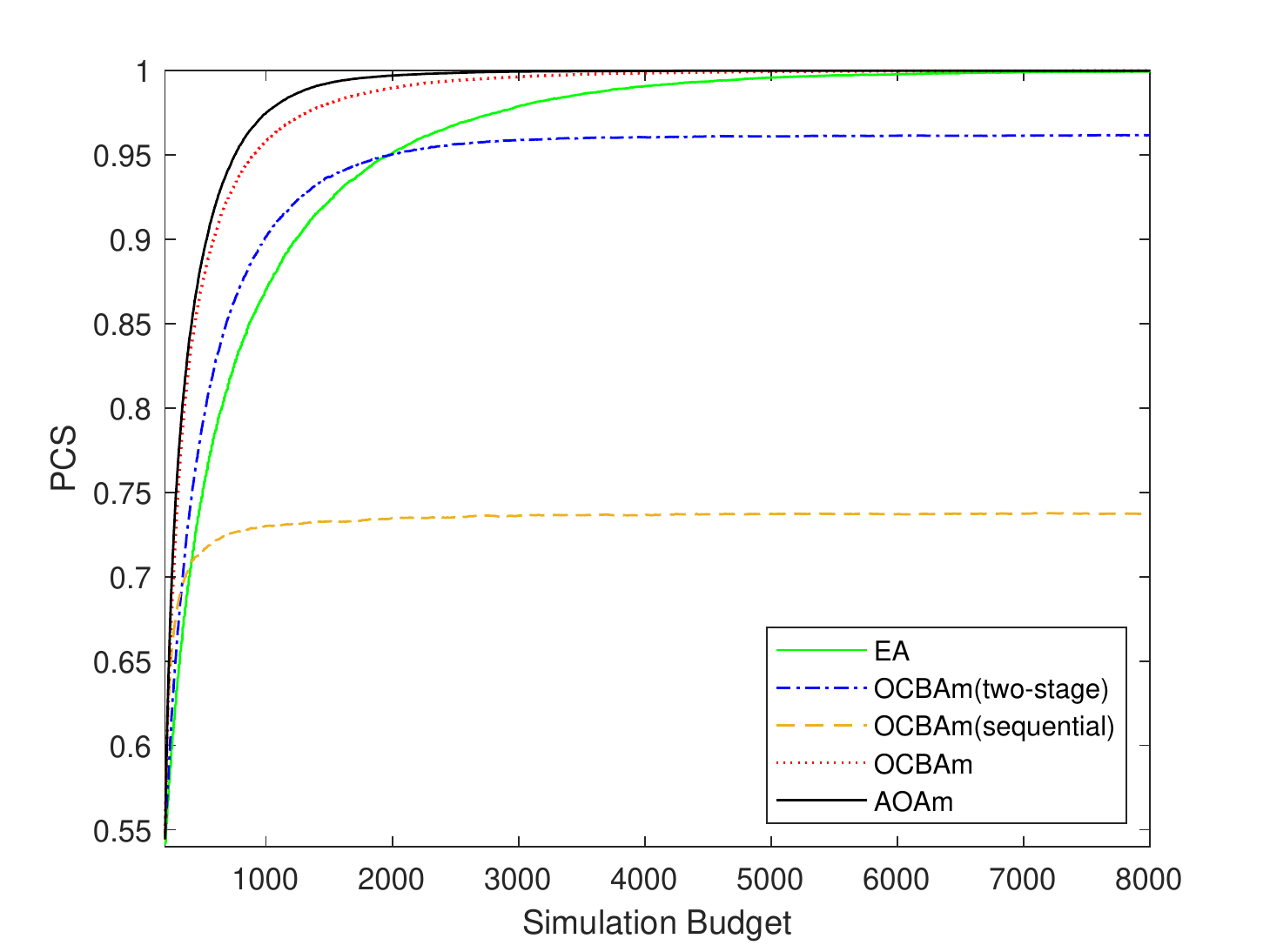}}
\subfigure[EOC]{\includegraphics[width=0.49\textwidth]{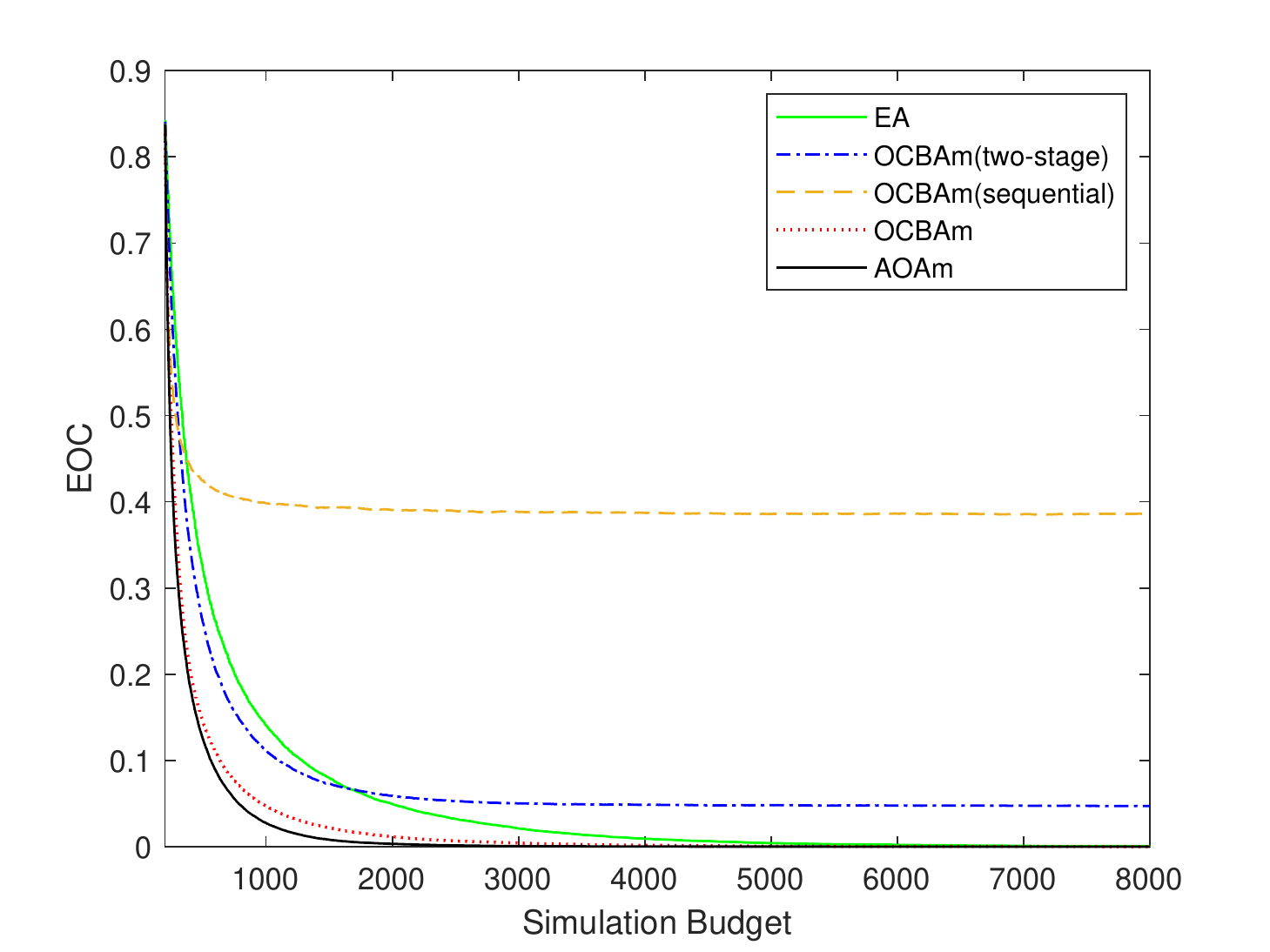}}
\caption{Comparison of PCSs and EOCs for 5 sampling allocation procedures in Experiment A.0.1.}
\label{figA221}
\end{figure}

From Figure~\ref{figA221} (a), we can see that the PCS of OCBAm(sequential) flattens and becomes the smallest as the simulation budget grows. The PCS of EA increases at a slow rate and surpasses the PCS of OCBAm(two-stage) when the simulation budget reaches around 2000. The PCS of OCBAm(two-stage) increases at a fast pace at the beginning and flattens as the simulation budget increases. The performance of OCBAm is significantly better than the performance of OCBAm(sequential). AOAm stands out among all allocation procedures in terms of both PCS and EOC. The performances of different allocation procedures in terms of EOC shown in Figure~\ref{figA221} (b) exhibit similar behaviors as measured by PCS, and this observation is consistent in all experiments.

\textit{Experiment A.0.2: 10 alternatives with increasing variances.} In this example, the sampling procedures are tested with $k=10$ and $m=3$. The performance of each alternative is $\mu_i = i$, $i=1,2,\cdots,10$, and the samples are drawn independently from a normal distribution $N\left(\mu_i,\sigma_i^2\right)$, where $\sigma_i=i$. The total simulation budget is $T=7000$. The numerical settings are the same as \textit{Example 2} in \citet{chen2008efficient}.

\begin{figure}[htbp]
\centering
\subfigure[PCS]{\includegraphics[width=0.49\textwidth]{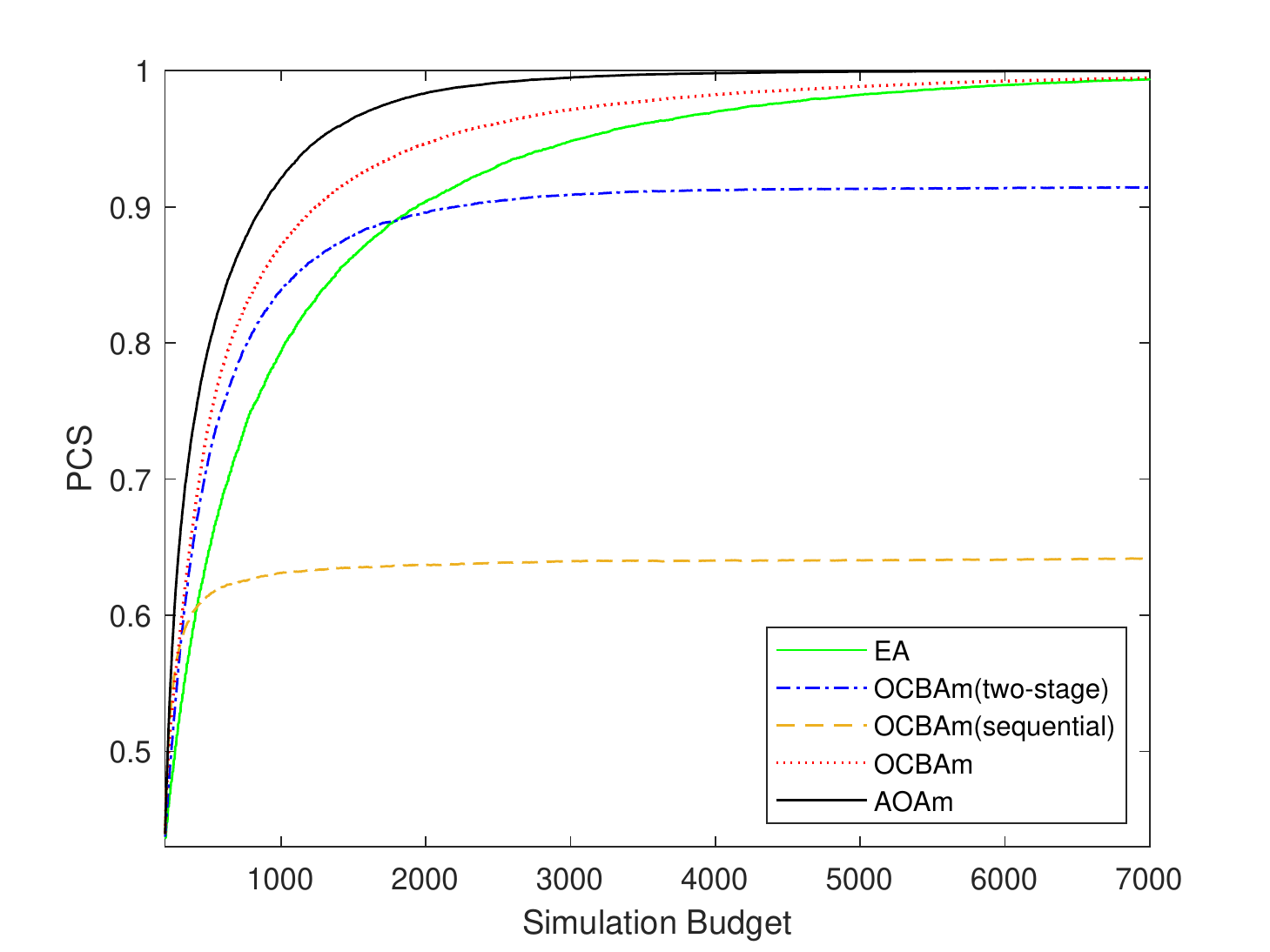}}
\subfigure[EOC]{\includegraphics[width=0.49\textwidth]{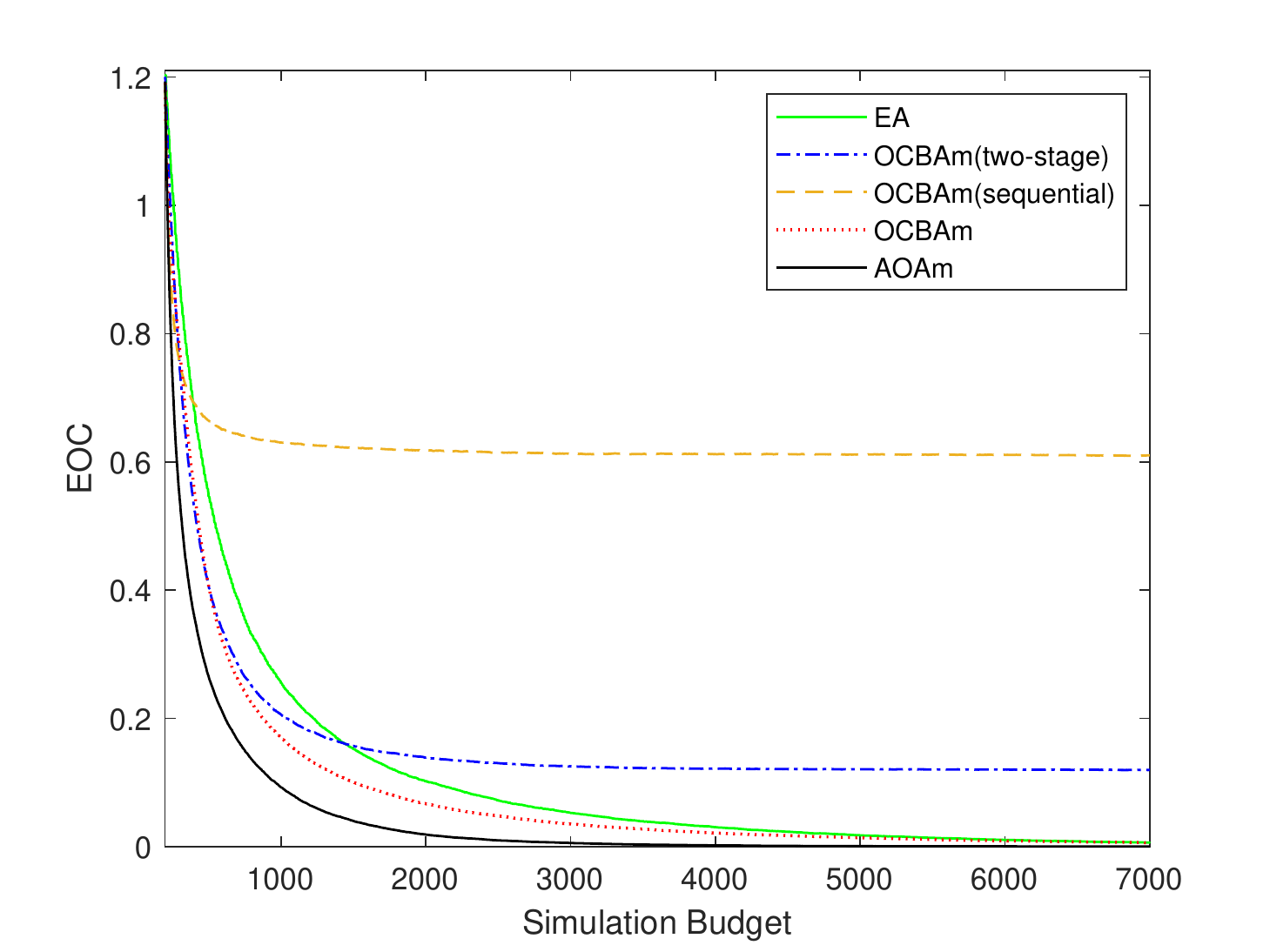}}
\caption{Comparison of PCSs and EOCs for 5 sampling allocation procedures in Experiment A.0.2.}
\label{figA222}
\end{figure}

From Figure~\ref{figA222} (a), we can see that the PCS of OCBAm(sequential) flattens as the simulation budget grows large. The performance of OCBAm(two-stage) is better than EA at the beginning, and the PCS of the latter surpasses the PCS of the former when the simulation budget reaches around 1900. The PCS of OCBAm increases at a fast pace, and the performance of OCBAm is significantly better than the performance of OCBAm(sequential). AOAm performs the best among all allocation procedures.

\textit{Experiment A.0.3: 20 alternatives with equal variances.} The numerical settings are the same as \textit{Experiment 1} in the main body of the paper.

\begin{figure}[htbp]
\centering
\subfigure[IPCS]{\includegraphics[width=0.49\textwidth]{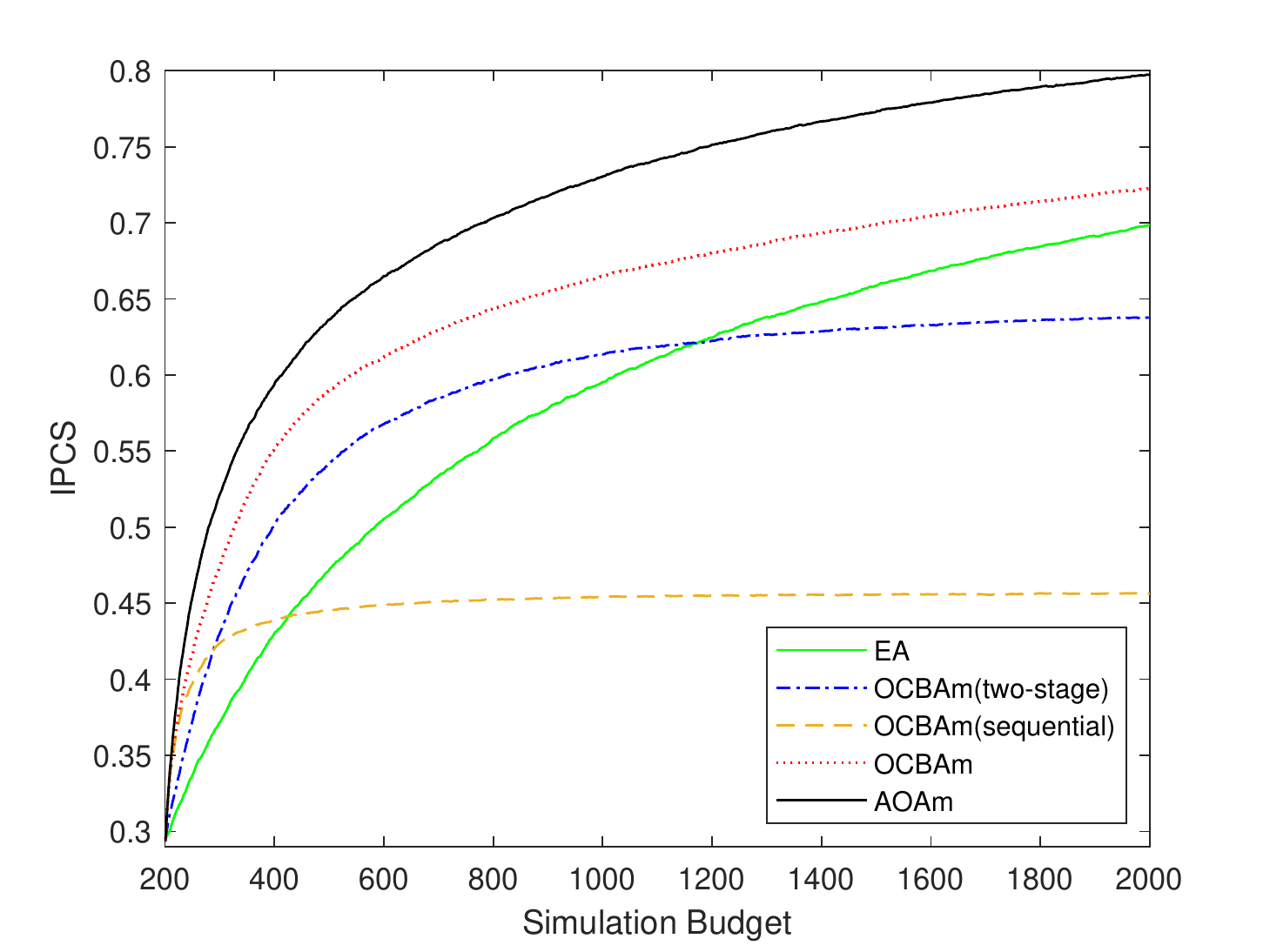}}
\subfigure[EOC]{\includegraphics[width=0.49\textwidth]{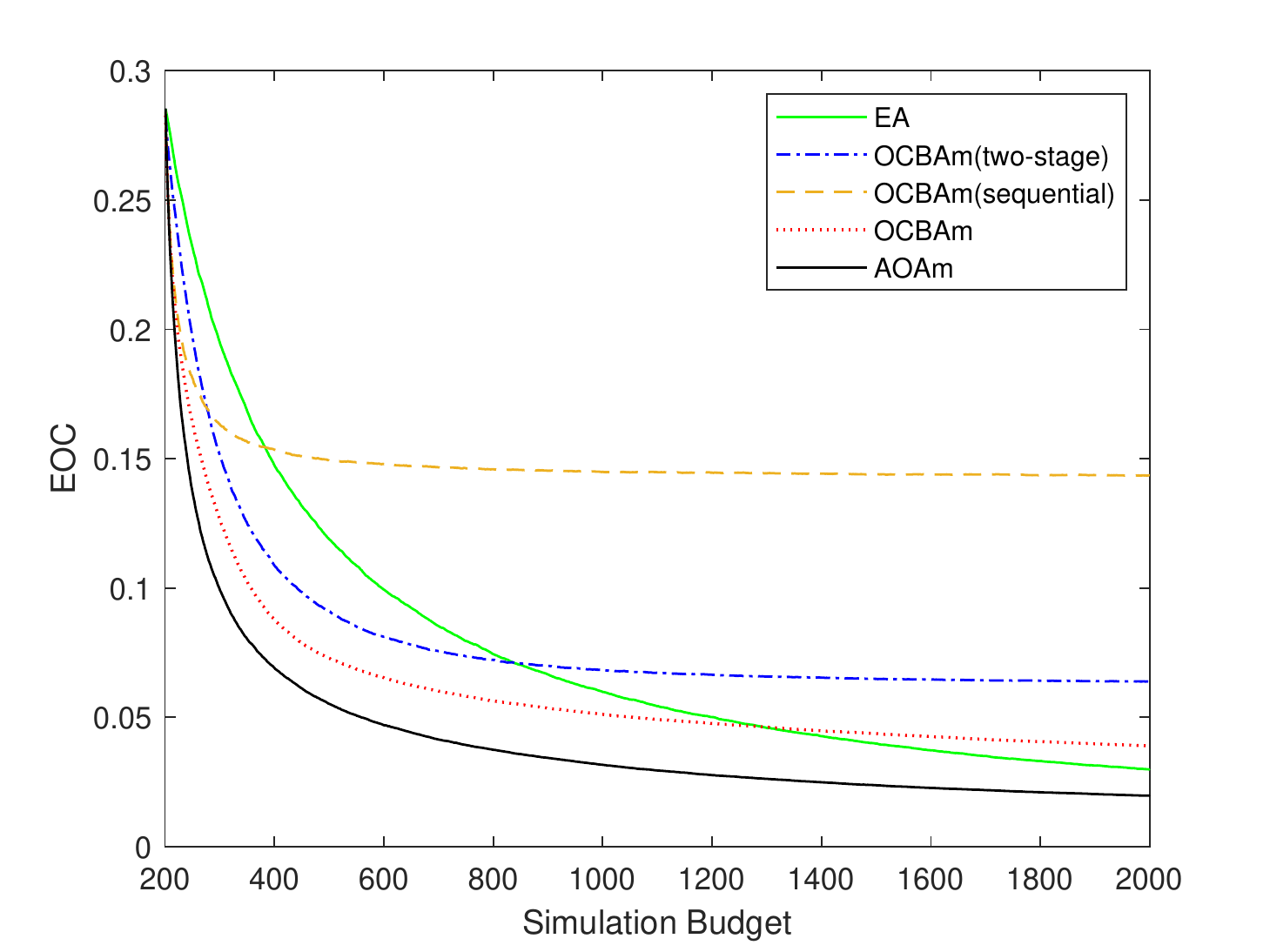}}
\caption{Comparison of IPCSs and EOCs for 5 sampling allocation procedures in Experiment A.0.3.}
\label{figA23}
\end{figure}

From Figure~\ref{figA23}, we can see that the IPCS of OCBAm(sequential) increases at a slow pace and becomes the worst as the simulation budget grows. The IPCS of EA is the smallest at the beginning, and it surpasses OCBAm(sequential) and OCBAm(two-stage) with the number of allocated simulation replications up to around 425 and 1175, respectively. The IPCS of OCBAm(two-stage) increases at a slower pace than the IPCS of OCBAm(sequential) at the beginning, and surpasses the latter when the simulation budget reaches around 290. The performance of OCBAm is significantly better than the performance of OCBAm(sequential). AOAm stands out among all allocation procedures in terms of both IPCS and EOC metrics.

Numerical results show that AOAm performs the best among all allocation procedures. The sequential OCBAm with parameter $c$ in (\ref{corg}) performs the worst as the simulation budget grows. The performance of sequential OCBAm with parameter $c$ in (\ref{cadp}) is significantly better than the performance of sequential OCBAm with parameter $c$ in (\ref{corg}). The suggested choice of parameter $c$ in~\citet{chen2008efficient} could not lead to a good performance when sequentially implementing the OCBAm sampling rule.
\end{remark}

In the numerical experiments in the main text, we use the OCBAm rule with $c$ in (\ref{cadp}), which could result in a better performance of OCBAm for comparison.

\item OCBAm$+$: \citet{zhang2012improved,zhang2015simulation} decompose static optimization problem (\ref{staticopt}) into two static sub-problems by focusing on alternatives $\left\langle m \right\rangle$ and $\left\langle {m+1} \right\rangle$, respectively, i.e.,
$$\mathop {\max }\limits_{{r_1},{r_2}, \cdots ,{r_k}} 1 - \sum\limits_{i = 1}^{m - 1} {\Pr \left\{ {\left. {{{\bar X}_{\left\langle i \right\rangle}}\left( {{r_{\left\langle  i \right\rangle}}T} \right) > {{\bar X}_{\left\langle m \right\rangle}}\left( {{r_{\left\langle m \right\rangle}}T} \right)} \right|\theta } \right\}}  - \sum\limits_{j = m + 1}^k {\Pr \left\{ {\left. {{{\bar X}_{\left\langle m \right\rangle}}\left( {{r_{\left\langle m \right\rangle}}T} \right) > {{\bar X}_{\left\langle j \right\rangle}}\left( {{r_{\left\langle  j \right\rangle}}T} \right)} \right|\theta } \right\}}~,$$
and
$$\begin{aligned}
\mathop {\max }\limits_{{r_1},{r_2}, \cdots ,{r_k}} 1 - & \sum\limits_{i = 1}^m {\Pr \left\{ {\left. {{{\bar X}_{\left\langle i \right\rangle}}\left( {{r_{\left\langle  i \right\rangle}}T} \right) > {{\bar X}_{\left\langle m+1 \right\rangle}}\left( {{r_{\left\langle m+1 \right\rangle}}T} \right)} \right|\theta } \right\}}  \\
& - \sum\limits_{j = m + 2}^k {\Pr \left\{ {\left. {{{\bar X}_{\left\langle m+1 \right\rangle}}\left( {{r_{{\left\langle m+1 \right\rangle}}}T} \right) > {{\bar X}_{\left\langle j \right\rangle}}\left( {{r_{\left\langle  j \right\rangle}}T} \right)} \right|\theta } \right\}}~,
\end{aligned}$$
and then derive asymptotic solution of each approximated optimization problem optimal in the limit. A sequential sampling procedure is proposed based on those asymptotic solutions, where unknown parameters are sequentially estimated from samples. Specifically, at each step, OCBAm$+$ allocates a simulation replication to an alternative according to the rule:
$${{\mathord{\buildrel{\lower3pt\hbox{$\scriptscriptstyle\smile$}} \over r} }^{\left(t\right)}} = \left\{
\begin{array}{rcl}
{\mathord{\buildrel{\lower3pt\hbox{$\scriptscriptstyle\smile$}} \over r} }^{{\left(t\right)}^{'}}, & & {\rm{if}\;\; {\bar{X}_{\left\langle m \right\rangle }} - {\bar{X}_{\left\langle {m - 1} \right\rangle }} \ge {\bar{X}_{\left\langle {m + 2} \right\rangle }} - {\bar{X}_{\left\langle {m + 1} \right\rangle }}}\\
{\mathord{\buildrel{\lower3pt\hbox{$\scriptscriptstyle\smile$}} \over r} }^{{\left(t\right)}^{''}}, & & {\rm{if}\;\; {\bar{X}_{\left\langle m \right\rangle }} - {\bar{X}_{\left\langle {m - 1} \right\rangle }} < {\bar{X}_{\left\langle {m + 2} \right\rangle }} - {\bar{X}_{\left\langle {m + 1} \right\rangle }}}\\
\end{array} \right.~,$$
where for $i,j=1,2,\cdots,k$, $i,j \ne m$,
$$\frac{{\mathord{\buildrel{\lower3pt\hbox{$\scriptscriptstyle\smile$}} \over r} }_{\left\langle i \right\rangle}^{{\left(t\right)}^{'}}}{{\mathord{\buildrel{\lower3pt\hbox{$\scriptscriptstyle\smile$}} \over r} }_{\left\langle j \right\rangle}^{{\left(t\right)}^{'}}} = {\left( {\frac{{{{{\widehat \sigma _{\left\langle i \right\rangle}}} / {({\bar X _{\left\langle i \right\rangle}} - {\bar X _{\left\langle m \right\rangle }})}}}}{{{{{\widehat \sigma _{\left\langle j \right\rangle}}} / {({\bar X _{\left\langle j \right\rangle}} - {\bar X _{\left\langle m \right\rangle }})}}}}} \right)^2},\quad {\mathord{\buildrel{\lower3pt\hbox{$\scriptscriptstyle\smile$}} \over r} }_{\left\langle m \right\rangle}^{{\left(t\right)}^{'}} = {\widehat \sigma _{\left\langle m \right\rangle }}\sqrt {\sum\limits_{i \ne m }^k {\frac{{{{\left( {{\mathord{\buildrel{\lower3pt\hbox{$\scriptscriptstyle\smile$}} \over r} }_{\left\langle i \right\rangle}^{{\left(t\right)}^{'}}}\right)}^2}}}{{\widehat \sigma _{\left\langle i \right\rangle}^2}}} }~,$$
and for $i,j=1,2,\cdots,k$, $i,j \ne m+1$,
$$\frac{{\mathord{\buildrel{\lower3pt\hbox{$\scriptscriptstyle\smile$}} \over r} }_{\left\langle i \right\rangle}^{{\left(t\right)}^{''}}}{{\mathord{\buildrel{\lower3pt\hbox{$\scriptscriptstyle\smile$}} \over r} }_{\left\langle j \right\rangle}^{{\left(t\right)}^{''}}} = {\left( {\frac{{{{{\widehat \sigma _{\left\langle i \right\rangle}}} / {({\bar X _{\left\langle i \right\rangle}} - {\bar X _{\left\langle {m+1}\right\rangle }})}}}}{{{{{\widehat \sigma _{\left\langle j \right\rangle}}} / {({\bar X _{\left\langle j \right\rangle}} - {\bar X _{\left\langle {m+1}\right\rangle }})}}}}} \right)^2},\quad {\mathord{\buildrel{\lower3pt\hbox{$\scriptscriptstyle\smile$}} \over r} }_{\left\langle {m+1}\right\rangle}^{{\left(t\right)}^{''}} = {\widehat \sigma _{\left\langle {m+1}\right\rangle }}\sqrt {\sum\limits_{i \ne m+1 }^k {\frac{{{{\left( {{\mathord{\buildrel{\lower3pt\hbox{$\scriptscriptstyle\smile$}} \over r} }_{\left\langle i \right\rangle}^{{\left(t\right)}^{''}}} \right)}^2}}}{{\widehat \sigma _{\left\langle i \right\rangle}^2}}} }~,$$
and the allocated alternative is chosen by the ``most starving" sequential rule ${A_{t + 1}}\left( {{\mathcal{E} _t}} \right) = \arg \mathop {\max }\limits_{i = 1,2, \cdots ,k} \{ {\left( {t + 1} \right){\mathord{\buildrel{\lower3pt\hbox{$\scriptscriptstyle\smile$}}\over r}}_i^{\left( t \right)} - {t_i}} \}$ in \citet{zhang2012improved} or with probability $\Pr \left\{ {{A_{t + 1}}\left( {{\mathcal{E} _t}}\right) = i} \right\} = {\mathord{\buildrel{\lower3pt\hbox{$\scriptscriptstyle\smile$}}\over r}}_i^{\left( t \right)}$ in \citet{zhang2015simulation}. OCBAm$+$ can be reduced to OCBA in \citet{chen2011stochastic} for selecting the best alternative only when the variances of all alternatives are equal, i.e., $\sigma_{\left\langle 1 \right\rangle}^2 = \sigma_{\left\langle 2 \right\rangle}^2 = \cdots = \sigma_{\left\langle k \right\rangle}^2 = \sigma^2$. \citet{zhang2015simulation} shows that when the variances of all alternatives are equal, the large deviations rate of PFS of OCBAm$+$ is
$$\frac{{{{\left( {{\mu _{\left\langle m \right\rangle }} - {\mu _{\left\langle {m + 1} \right\rangle }}} \right)}^2}}}{{2{\sigma ^2}\left( {{1 \mathord{\left/{\vphantom {1{{{\mathord{\buildrel{\lower3pt\hbox{$\scriptscriptstyle\smile$}}\over r} }_{\left\langle m \right\rangle }}}}} \right.\kern-\nulldelimiterspace} {{{\mathord{\buildrel{\lower3pt\hbox{$\scriptscriptstyle\smile$}}\over r} }_{\left\langle m \right\rangle }}}} + {1 \mathord{\left/{\vphantom{1{{{\mathord{\buildrel{\lower3pt\hbox{$\scriptscriptstyle\smile$}}\over r} }_{\left\langle {m + 1} \right\rangle }}}}}\right.\kern-\nulldelimiterspace}{{{\mathord{\buildrel{\lower3pt\hbox{$\scriptscriptstyle\smile$}}\over r} }_{\left\langle {m + 1} \right\rangle }}}}} \right)}}~,$$
where ${\mathord{\buildrel{\lower3pt\hbox{$\scriptscriptstyle\smile$}}\over r} }_{\left\langle i \right\rangle }$, $i=1,2,\cdots,k$ is the allocation ratio of OCBAm+ rule. In the numerical experiments in the main body of the paper, the allocated alternative for OCBAm$+$ rule is chosen by the ``most starving" seuential rule.

\item OCBAss: \citet{gao2015note} approximate the objective function in the static optimization problem (\ref{staticopt}) with Bonferroni inequality, i.e.,
$$\mathop {\max }\limits_{{r_1},{r_2}, \cdots ,{r_k}} 1 - \sum\limits_{i = 1}^m {\sum\limits_{j = m + 1}^k {\Pr \left\{ {\left. {{{\bar X}_{\left\langle i \right\rangle}}\left( {{r_{\left\langle  i \right\rangle}}T} \right) > {{\bar X}_{\left\langle j \right\rangle}}\left( {{r_{\left\langle  j \right\rangle}}T} \right)} \right|\theta } \right\}} }~,$$
and then derive the asymptotic solutions of the approximated optimization problem optimal in the limit. A sequential sampling procedure is proposed to achieve the asymptotic solutions by balancing the equations in the asymptotic solutions, where unknown parameters are sequentially estimated from samples. Specifically, at each step, OCBAss allocates a simulation replication to an alternative according to the  rule:
$$A_{t+1}\left(\mathcal{E}_t\right) =  \left\{
\begin{array}{rcl}
\arg \mathop {\min }\limits_{i = 1,2, \cdots ,m} {I_{i,{\widetilde{j}_i}}}, & & {{\rm{if}}\;\; \sum\limits_{i = 1}^m {\frac{{t_{\left\langle i \right\rangle}^2}}{{\widehat \sigma _{\left\langle i \right\rangle}^2}} < } \sum\limits_{j = m + 1}^k {\frac{{t_{\left\langle j \right\rangle}^2}}{{\widehat \sigma _{\left\langle j \right\rangle}^2}}}}\\
\arg \mathop {\min }\limits_{j = m+1,m+2, \cdots,k} {I_{{\widetilde{i}_j},{j}}}, & & {{\rm{if}}\;\; \sum\limits_{i = 1}^m {\frac{{t_{\left\langle i \right\rangle}^2}}{{\widehat \sigma _{\left\langle i \right\rangle}^2}} \ge } \sum\limits_{j = m + 1}^k {\frac{{t_{\left\langle j \right\rangle}^2}}{{\widehat \sigma _{\left\langle j \right\rangle}^2}}}}\\
\end{array} \right.~,$$
where ${\widetilde{j}_i}\mathop  = \limits^\Delta  \arg \mathop {\min }\limits_{j = m + 1,m + 2, \cdots ,k} {I_{i,j}}$, ${\widetilde{i}_j}\mathop  = \limits^\Delta  \arg \mathop {\min }\limits_{i = 1,2, \cdots ,k} {I_{i,j}}$, and ${I_{i,j}} = \frac{\left(\bar X_{\left\langle i \right\rangle} - \bar X_{\left\langle j \right\rangle}\right)^2}{{{{\widehat \sigma _{\left\langle i \right\rangle}^2} / {{r_{\left\langle i \right\rangle}^{\left(t\right)}}}} + {{\widehat \sigma _{\left\langle j \right\rangle}^2} / {{r_{\left\langle j \right\rangle}^{\left(t\right)}}}}}}$.

\item OCBASS: \citet{gao2016new} extend the results in~\citet{gao2015note} by using a large deviations principle to further approximate the objective function, i.e.,
$$\mathop {\max }\limits_{{r_1},{r_2}, \cdots ,{r_k}} 1 - \sum\limits_{i = 1}^m {\sum\limits_{j = m + 1}^k {\exp \left( { - T{G_{ij}}\left( {{r_{\left\langle  i \right\rangle}},{r_{\left\langle  j \right\rangle}}} \right)} \right)} }~,$$
where ${\exp \left( { - T{G_{ij}}\left( {{r_{\left\langle  i \right\rangle}},{r_{\left\langle  j \right\rangle}}} \right)} \right)}$ is an upper bound for $\Pr \left\{ {\left. {{{\bar X}_{\left\langle i \right\rangle}}\left( {{r_{\left\langle  i \right\rangle}}T} \right) \le {{\bar X}_{\left\langle j \right\rangle}}\left( {{r_{\left\langle  j \right\rangle}}T} \right)} \right|\theta } \right\}$, and then derive the asymptotic solutions of the approximated optimization problem optimal in the limit. A sequential sampling procedure is proposed to achieve the asymptotic solutions, where unknown parameters are sequentially estimated from samples. Specifically, at each step, OCBASS allocates a simulation replication to an alternative according to the rule: $A_{t+1}\left(\mathcal{E}_t\right) = \arg \mathop {\min }\limits_{\ell = 1,2, \cdots ,k} {I_{\ell,{\ell_0}}}$, where
$${\ell_0} \mathop  = \limits^\Delta  \left\{
\begin{array}{rcl}
\arg \mathop {\min }\limits_{j = m + 1,m + 2, \cdots ,k} {I_{\ell,j}}, & & {\rm{if}\;\; \ell=1,2,\cdots,m}\\
\arg \mathop {\min }\limits_{i = 1,2, \cdots,m} {I_{i,\ell}}, & & {\rm{if}\;\; \ell=m+1,m+2,\cdots,k}\\
\end{array} \right.~,$$
however, note that $A_{t+1}\left(\mathcal{E}_t\right)$ is not unique when determining the allocated alternative with OCBASS rule, since $\arg \mathop {\min }\limits_{i,j} {I_{ij}}$, $i = 1,2, \cdots ,m$, $j = m + 1,m + 2, \cdots ,k$ includes two alternatives. In the numerical experiments in the main body of the paper, we randomly choose an alternative among the two alternatives as the allocated alternative.
\end{enumerate}

We caution that the approximated objective functions could be loose when the simulation budget $T$ is small. For example, Bonferroni inequality could lead to a negative lower bound for the objective function in (\ref{staticopt}) which is a probability.

\subsection*{A.1. Algorithm}

The proposed sequential AOAm procedure is shown in Algorithm 1. The initial $n_0$ simulation replications are allocated to each alternative for estimating their unknown parameters. Additional simulation replications are allocated incrementally with one replication in each step to the alternative with the largest VFA until the total simulation replications are exhausted. The posterior information is updated iteratively in the allocation procedure.

\begin{algorithm}[htbp]
 \caption{Asymptotically Optimal Allocation Procedure for Selecting Top-$m$ Alternatives (AOAm)}
 \begin{algorithmic}[1]
 \renewcommand{\algorithmicrequire}{\textbf{Input:}}
 \REQUIRE $k$, $m$, $n_0$, $T$. \\
 \textbf{INITIALIZE:} $t=0$ and perform $n_0$ simulation replications for each alternative. \\
 \textbf{WHILE} $t < \left(T-k \times n_0\right)$, \textbf{DO:}\\
 \quad \quad \textbf{UPDATE:} Calculate the posterior variances ${( {\sigma _{\rm{i}}^{\left( t \right)}})^2}$ and posterior means $\mu_{i}^{\left(t\right)}$ with sample means ${\bar{X}_i}\left(t_i\right)$ and sample variances $\widehat{\sigma}_{i}^2$, $i=1,2,\cdots,k$ according to (\ref{pv}) and (\ref{pm}).\\
 \quad \quad \textbf{ALLOCATE:}  Calculate ${\widehat {V}_t}\left( {{\mathcal{E}_t};i} \right)$, $i=1,2,\cdots,k$ according to (\ref{equ21}) and (\ref{equ22}). \\
 Determine the allocated alternative based on allocation rule (\ref{AOAm}). \\
 \quad \quad \textbf{SIMULATE:} Run an additional simulation for the allocated alternative. \\
 \quad \quad $t=t+1$.\\
 \textbf{END WHILE.}
\end{algorithmic}
\end{algorithm}

\subsection*{A.2. Theoretical Supplements}

\subsubsection*{Proof of Lemma \ref{lemVFAintegral}}
We apply the Cholesky decomposition to the covariance matrix $\Sigma  = \Gamma '\Lambda \Gamma  = U'U$. Then, the value function for selecting the alternatives in $\widehat{\mathcal{F}}_t^m$ can be expressed as
$$\begin{aligned}
& \Pr \left\{ \left. {{\bigcap\limits_{i = 1}^m {\bigcap\limits_{j = m + 1}^k {\left( {{\mu _{{{\left\langle i \right\rangle }_t}}} > {\mu _{{{\left\langle j \right\rangle }_t}}}} \right)} } } } \right|\mathcal{E}_t \right\} \\
= & \int_{{\chi _{1,\left( {m + 1} \right) > 0}}, \cdots, {\chi _{m,k} > 0}} {f\left( {{\chi _{1,\left( {m + 1} \right)}}, \cdots, {\chi _{1,k}}, \cdots, {\chi _{m,\left(m+1\right)}}, \cdots,{\chi _{m,k}}} \right)}d{\chi _{1,\left( {m + 1} \right)}} \cdots d{\chi _{m,k}} \\
 = & \Pr \left\{{\sum\limits_{\ell = 1}^q {{u_{\ell,q}}{z_\ell}}  > \mu _{{{\left\langle j \right\rangle }_t}}^{\left( t \right)} - \mu _{{{\left\langle i \right\rangle }_t}}^{\left( t \right)},\;i = 1,2, \cdots ,m,\;j = m + 1,m + 2, \cdots ,k} \right\} \\
 = & \frac{1}{\left( {\sqrt {2\pi } } \right)^{m \times\left(k - m\right)}} \times \int_{\sum_{\ell = 1}^q {u_{\ell, q} z_{\ell}} \geq \mu_{\langle j\rangle_t}^{(t)}-\mu_{\langle i\rangle_t}^{(t)}} \prod_{h = 1}^{m\times\left( {k - m} \right)} \exp \left(-\frac{z_{h}^{2}}{2}\right) d z_{1} \cdots d z_{m \times \left( {k - m} \right)}~,
\end{aligned}$$
where $\chi_{i,j}\mathop  = \limits^\Delta {{\mu _{{{\left\langle {i} \right\rangle }_t}}} - {\mu _{{{\left\langle {j} \right\rangle }_t}}}}$, $i = 1,2,\cdots,m$, $j=m+1,m+2,\cdots,k$, and the second equality holds by a change of variables ${\chi _{ij}} = \mu _{{{\left\langle {i} \right\rangle }_t}}^{\left( t \right)} - \mu _{{{\left\langle j \right\rangle }_t}}^{\left( t \right)} + \sum\nolimits_{\ell = 1}^{q} {{u_{\ell,q}}} {z_\ell}$, where $q=1,2,\cdots,m\times\left(k-m\right)$ are ranked indices of $( \mu _{{{\left\langle {i} \right\rangle }_t}}^{\left( t \right)} - \mu _{{{\left\langle j \right\rangle }_t}}^{\left( t \right)} )$ in the joint distribution vector (\ref{jointvector}). Summarizing the above, the Lemma is proved.

\subsubsection*{Proof of Theorem~\ref{thmconsistent}}
We only need to prove that every alternative will be sampled infinitely often almost surely, i.e., $\mathop {\lim }\limits_{t \to \infty } {t_i} = \infty$ for $i=1,2,\cdots,k$. Following AOAm (\ref{AOAm}), the consistency will follow by the law of large numbers (LLN).

Define $$\mathcal{I}\mathop  = \limits^\Delta \{ i\in\{1,2,\cdots,k\}:~ \mbox{alternative $\langle i\rangle$ is sampled infinitely often a.s.} \}~.$$

If $\left\{1,2,\cdots,m\right\}\cap\mathcal{I}=\emptyset$ and $\left\{m+1,m+2,\cdots,k\right\}\cap \mathcal{I}\neq \emptyset$, then $\forall~ i \in \left\{1,2,\cdots,m \right\}$, $\mathop {\lim }\limits_{t \to \infty } {( {\sigma _{\langle i\rangle}^{\left( t \right)}})^2} > 0$ and $\mathop {\lim }\limits_{t \to \infty } [{{( {\sigma _{\langle i\rangle}^{\left( t \right)}})^2} - ( {\sigma _{\langle i\rangle}^{\left( t+1 \right)}})^2}]> 0$. $\exists~ h \in \left\{1,2,\cdots,m \right\}$ such that
\begin{equation*}
\mathop {\lim}\limits_{t \to \infty } \left[{{\widehat{V}_t} \left( {{\mathcal{E} _t};h} \right) - \widetilde{V}_t\left( {{\mathcal{E} _t}} \right)} \right] > 0,\quad a.s.~,
\end{equation*}
and $\exists~ j \in \left\{m+1,m+2,\cdots,k \right\}$ such that $\mathop {\lim }\limits_{t \to \infty } {( {\sigma _{\langle j\rangle}^{\left( t \right)}})^2} =0$, $\mathop {\lim }\limits_{t \to \infty } \left[{{( {\sigma _{\langle j\rangle}^{\left( t \right)}})^2} - ( {\sigma _{\langle j\rangle}^{\left( t+1 \right)}})^2}\right]=0$, and
\begin{equation*}
\mathop {\lim}\limits_{t \to \infty } \left[{{\widehat{V}_t} \left( {{\mathcal{E} _t};j} \right) - \widetilde{V}_t\left( {{\mathcal{E} _t}} \right)} \right] =0,\quad a.s.~,
\end{equation*}
which leads to a contradiction to the sampling rule (\ref{AOAm}) that the alternative with the largest VFA is sampled. Therefore, $\left\{1,2,\cdots,m\right\}\cap\mathcal{I}\neq \emptyset$.

If $\left\{m+1,m+2,\cdots,k\right\}\cap\mathcal{I}=\emptyset$, then $\forall~ j \in \left\{m+1,m+2,\cdots,k \right\}$, $\mathop {\lim }\limits_{t \to \infty } {( {\sigma _{\langle j\rangle}^{\left( t \right)}})^2} > 0$ and $\mathop {\lim }\limits_{t \to \infty } \left[{{( {\sigma _{\langle j\rangle}^{\left( t \right)}})^2} - ( {\sigma _{\langle j\rangle}^{\left( t+1 \right)}})^2}\right]> 0$,   $\exists~ j \in  \left\{m+1,m+2,\cdots,k\right\}$ such that
\begin{equation*}
\mathop {\lim}\limits_{t \to \infty } \left[{{\widehat{V}_t} \left( {{\mathcal{E} _t};j} \right) - \widetilde{V}_t\left( {{\mathcal{E} _t}} \right)} \right] >0,\quad a.s.~,
\end{equation*}
and for $i\in \mathcal{I}\cap\left\{1,2,\cdots,m\right\}$,
\begin{equation*}
\mathop {\lim}\limits_{t \to \infty } \left[{{\widehat{V}_t} \left( {{\mathcal{E} _t};i} \right) - \widetilde{V}_t\left( {{\mathcal{E} _t}} \right)} \right] =0,\quad a.s.~,
\end{equation*}
which leads to a contradiction to the sampling rule (\ref{AOAm}). Therefore, $\left\{m+1,m+2,\cdots,k\right\}\cap\mathcal{I}\neq \emptyset$. In addition, if $\left\{m+1,m+2,\cdots,k\right\}\setminus\mathcal{I}\neq \emptyset$, where $\setminus$ denotes set minus, then for $i\in\left\{1,2,\cdots,m\right\}\cap\mathcal{I}$, $j\in\left\{m+1,m+2,\cdots,k\right\}\cap\mathcal{I}$, and  $\ell \in \left\{m+1,m+2,\cdots,k\right\}\setminus\mathcal{I}$,
$$\lim\limits_{t\to \infty}d_{ij}(\mathcal{E}_t)=\infty,\quad \lim\limits_{t\to \infty}d_{i\ell}(\mathcal{E}_t)<\infty~.$$

Thus, for $j\in\left\{m+1,m+2,\cdots,k\right\}\cap\mathcal{I}$,
\begin{align*}
\mathop {\lim}\limits_{t \to \infty } \left[{{\widehat{V}_t} \left( {{\mathcal{E} _t};j} \right) - \widetilde{V}_t\left( {{\mathcal{E} _t}} \right)} \right] =0,\quad a.s.~,
\end{align*}
and $\exists ~q \in \left\{m+1,m+2,\cdots,k\right\}\setminus\mathcal{I}$ such that
\begin{align*}
\mathop {\lim}\limits_{t \to \infty } \left[{{\widehat{V}_t} \left( {{\mathcal{E} _t};q} \right) - \widetilde{V}_t\left( {{\mathcal{E} _t}} \right)} \right] >0,\quad a.s.~,
\end{align*}
which leads to a contradiction to the sampling rule (\ref{AOAm}). Therefore, $\left\{m+1,m+2,\cdots,k\right\}\subset\mathcal{I}$. Finally, if $\left\{1,2,\cdots,m\right\}\setminus\mathcal{I}\neq \emptyset$, then  for $j\in\left\{m+1,m+2,\cdots,k\right\}$, $i\in\left\{1,2,\cdots,m\right\}\cap\mathcal{I}$,  $h\in \left\{1,2,\cdots,m\right\}\setminus\mathcal{I}$,
$$\lim\limits_{t\to \infty}d_{ij}(\mathcal{E}_t)=\infty,\quad \lim\limits_{t\to \infty}d_{hj}(\mathcal{E}_t)<\infty~.$$

Thus, for $i\in\left\{1,2,\cdots,m\right\}\cap\mathcal{I}$,
\begin{align*}
\mathop {\lim}\limits_{t \to \infty } \left[{{\widehat{V}_t} \left( {{\mathcal{E} _t};i} \right) - \widetilde{V}_t\left( {{\mathcal{E} _t}} \right)} \right] =0,\quad a.s.~,
\end{align*}
and $\exists ~h \in \left\{1,2,\cdots,m\right\}\setminus\mathcal{I}$ such that
\begin{align*}
\mathop {\lim}\limits_{t \to \infty } \left[{{\widehat{V}_t} \left( {{\mathcal{E} _t};h} \right) - \widetilde{V}_t\left( {{\mathcal{E} _t}} \right)} \right] >0,\quad a.s.~,
\end{align*}
which leads to a contradiction to the sampling rule (\ref{AOAm}). Therefore, $\left\{1,2,\cdots,m\right\}\subset\mathcal{I}$. Summarizing the above, the theorem is proved.

\subsubsection*{Assumptions}
Let ${\mathcal{D}_{{\Lambda_i}}} \mathop  = \limits^\Delta \left\{ {\lambda  \in \mathbb{R}:{\Lambda_i}\left( \lambda  \right) < \infty } \right\}$ be the effective domain of the function ${\Lambda_i}\left(  \cdot  \right)$, $i=1,2,\cdots,k$, and let ${\mathcal{F}_i} \mathop  = \limits^\Delta  \left\{ {{{d{\Lambda_i}\left( \lambda  \right)} \mathord{\left/{\vphantom {{d{\Lambda_i}\left( \lambda  \right)} {d\lambda }}} \right.\kern-\nulldelimiterspace} {d\lambda }}:\lambda  \in {\mathcal{D}_{\Lambda_i}^{o}}} \right\}$, $i=1,2,\cdots,k$, where $\mathbb{A}^o$ denotes the interior for the set $\mathbb{A}$. The following assumptions are standard.

\begin{assumption}\label{ass1}
The limit ${\Lambda_i}\left( \lambda  \right) = \mathop {\lim }\limits_{T \to \infty } \frac{1}{T}{\Lambda_i}\left( {T\lambda } \right)$, $i=1,2,\cdots,k$ is well defined as an extended real number for all $\lambda$.
\end{assumption}

\begin{assumption}\label{ass2}
The origin belongs to ${\mathcal{D}_{\Lambda_i}^{o}}$, that is, $0 \in {\mathcal{D}_{\Lambda_i}^{o}}$, $i=1,2,\cdots,k$.
\end{assumption}

\begin{assumption}\label{ass3}
${\Lambda_i}\left( \lambda \right)$ are strictly convex, continuous on ${\mathcal{D}_{\Lambda_i}^o}$ and steep, i.e., $\mathop {\lim }\limits_{T \to \infty } \left| {\Lambda_i^\prime \left( {{\lambda _T}} \right)} \right| = \infty$, $i=1,2,\cdots,k$, where $\{\lambda_T\} \in \mathcal{D}_{\Lambda_i}$ is a sequence converging to a boundary point of $\mathcal{D}_{\Lambda_i}^o$.
\end{assumption}

\begin{assumption}\label{ass4}
The interval $\left[ {{\mu _ {\left\langle k \right\rangle}},{\mu _{\left\langle 1 \right\rangle}}} \right] \subset \bigcap\nolimits_{i = 1}^k {\mathcal{F}_i^o}$.
\end{assumption}

Note that Assumption~\ref{ass2} and Assumption~\ref{ass3} are stronger than what is needed for the G$\rm{\ddot a}$rtner-Ellis theorem to hold (see~\citealp{zeitouni2010large}, Definition 2.3.5). By the G$\rm{\ddot a}$rtner-Ellis theorem (see~\citealp{zeitouni2010large}, Theorem 2.3.6), Assumptions~\ref{ass1},~\ref{ass2} and~\ref{ass3} indicate that the sample mean ${{\bar X}_i}\left( {{r_i}T} \right)$ satisfy the large deviations principle with rate functions $\Lambda_{i}^*\left(\cdot\right)$, $i=1,2,\cdots,k$. Assumption~\ref{ass3} indicates that $\Lambda_i^*\left(x\right)$ are strictly convex and continuous for $x\in{\mathcal{F}_i^o}$, where ${\mathcal{F}_i^o}$ is the interior of the set ${\mathcal{F}_i}$, $i=1,2,\cdots,k$ (see~\citealp{zeitouni2010large}, Page 35). Assumption~\ref{ass4} ensures that the sample mean of each alternative can take any value in $\left[ {\mu _{\left\langle k \right\rangle}},{\mu _{\left\langle 1 \right\rangle}} \right]$, and $\Pr \left\{ \left. {\bar{X}_{\langle j\rangle}\left(r_{\langle j \rangle} T\right) \ge \bar{X}_{\langle i\rangle}\left(r_{\langle i\rangle} T\right) } \right|\theta \right\} >0$ for $i = 1,2, \cdots ,m$ and $j = m+1,m+2, \cdots ,k$. The assumptions hold for common distributions that are light-tailed with overlapping support, e.g., normal distribution, Bernoulli and Poisson distribution.

\subsubsection*{Proof of Theorem~\ref{thm1}}
Let ${Z_T} = \left( {{{\bar X}_{\left\langle i \right\rangle}}\left( {{r_{\left\langle  i \right\rangle}}T} \right),{{\bar X}_{\left\langle j \right\rangle}}\left( {{r_{\left\langle  j \right\rangle}}T} \right)} \right)$, $T=1,2,\cdots$, $i=1,2,\cdots,m$ and $j=m+1,m+2,\cdots,k$. The cumulant moment generating function of $Z_T$ is
\begin{align}
{\Lambda _T}\left( {{\lambda _{\left\langle i \right\rangle}},{\lambda _{\left\langle j \right\rangle}}} \right) &= \log \mathbb{E}\left[ {\exp \left( {{\lambda _{\left\langle i \right\rangle}}{{\bar X}_{\left\langle i \right\rangle}}\left( {{r_{\left\langle  i \right\rangle}}T} \right) + {\lambda _{\left\langle j \right\rangle}}{{\bar X}_{\left\langle j \right\rangle}}\left( {{r_{\left\langle  j \right\rangle}}T} \right)} \right)} \right] \nonumber \\
&= \Lambda _{\left\langle  i \right\rangle}^{\left( {{r_{\left\langle  i \right\rangle}}T} \right)}\left( {\frac{\lambda _{\left\langle i \right\rangle}}{t_{\left\langle  i \right\rangle}}} \right) + \Lambda_{\left\langle  j \right\rangle}^{\left( {r_{\left\langle  j \right\rangle}T} \right)}\left( {\frac{\lambda _{\left\langle j \right\rangle}}{t_{\left\langle  j \right\rangle}}} \right) \nonumber~,
\end{align}
under Assumption~\ref{ass1}, we have
$$\mathop {\lim }\limits_{T \to \infty } \frac{1}{T}{\Lambda_T}\left(\lambda_{\left\langle i \right\rangle}T,\lambda_{\left\langle j \right\rangle}T\right) = r_{\left\langle  i \right\rangle}{\Lambda _{\left\langle  i \right\rangle}}\left( {\frac{\lambda _{\left\langle i \right\rangle}}{r_{\left\langle  i \right\rangle}}} \right) + r_{\left\langle  j \right\rangle}{\Lambda _{\left\langle  j \right\rangle}}\left( {\frac{\lambda _{\left\langle j \right\rangle}}{r_{\left\langle  j \right\rangle}}} \right)~.$$

By G$\rm{\ddot a}$rtner-Ellis Theorem~\citep{zeitouni2010large}, $\left\{ {{Z_T},\;T = 1,2,\cdots } \right\}$ has large deviations with a rate function
\begin{align}
\Lambda^*\left( {x_{\left\langle i \right\rangle},x_{\left\langle j \right\rangle}} \right) &= \mathop {\sup }\limits_{{\lambda _{\left\langle i \right\rangle}},{\lambda_{\left\langle j \right\rangle}}} \left( {{\lambda _{\left\langle i \right\rangle}}{x_{\left\langle i \right\rangle}} + {\lambda _{\left\langle j \right\rangle}}{x_{\left\langle j \right\rangle}}-r_{\left\langle  i \right\rangle}{\Lambda _{\left\langle  i \right\rangle}}\left( {\frac{\lambda _{\left\langle i \right\rangle}}{r_{\left\langle  i \right\rangle}}} \right) - r_{\left\langle  j \right\rangle}{\Lambda _{\left\langle  j \right\rangle}}\left( {\frac{\lambda _{\left\langle j \right\rangle}}{r_{\left\langle  j \right\rangle}}} \right)} \right)\nonumber\\
&= \mathop {\sup }\limits_{\lambda _{\left\langle i \right\rangle}} \left( {\lambda_{\left\langle i \right\rangle}{x_{\left\langle i \right\rangle}} - r_{\left\langle  i \right\rangle}{\Lambda _{\left\langle  i \right\rangle}}\left( {\frac{\lambda _{\left\langle i \right\rangle}}{r_{\left\langle  i \right\rangle}}} \right)} \right) + \mathop {\sup }\limits_{\lambda _{\left\langle j \right\rangle}} \left( {\lambda_{\left\langle j \right\rangle}{x_{\left\langle j \right\rangle}} - r_{\left\langle  j \right\rangle}{\Lambda _{\left\langle  j \right\rangle}}\left( {\frac{\lambda _{\left\langle j \right\rangle}}{r_{\left\langle  j \right\rangle}}} \right)} \right)\nonumber\\
&= r_{\left\langle  i \right\rangle}\mathop {\sup }\limits_{{{\lambda_{\left\langle i \right\rangle}} / {r_{\left\langle  i \right\rangle}}}} \left( {\left( {\frac{{{\lambda _{\left\langle i \right\rangle}}}}{{r_{\left\langle  i \right\rangle}}}} \right){x_{\left\langle i \right\rangle}} - {\Lambda _{\left\langle i \right\rangle}}\left( {\frac{{{\lambda _{\left\langle i \right\rangle}}}}{{r_{\left\langle  i \right\rangle}}}} \right)} \right)+ r_{\left\langle  j \right\rangle}\mathop {\sup }\limits_{{{{\lambda _{\left\langle j \right\rangle}}} / {r_{\left\langle  j \right\rangle}}}} \left( {\left( {\frac{{{\lambda _{\left\langle j \right\rangle}}}}{{r_{\left\langle  j \right\rangle}}}} \right){x_{\left\langle j \right\rangle}} - {\Lambda _{\left\langle j \right\rangle}}\left( {\frac{{{\lambda _{\left\langle j \right\rangle}}}}{{r_{\left\langle  j \right\rangle}}}} \right)} \right)\nonumber\\
&= r_{\left\langle i \right\rangle}{\Lambda_{\left\langle i \right\rangle}^*}\left( {{x_{\left\langle i \right\rangle}}} \right) + r_{\left\langle  j \right\rangle}{\Lambda_{\left\langle j \right\rangle}^*}\left( {{x_{\left\langle j \right\rangle}}} \right) \nonumber~.
\end{align}

Therefore, from large deviations theorem, for $i = 1,2,\cdots, m$ and $j = m+1,m+2,\cdots,k$, we have ${G_{ij}}\left( {r_{\left\langle i \right\rangle},r_{\left\langle j \right\rangle}} \right) = \mathop {\inf }\limits_{{x_{\left\langle j \right\rangle}} \ge {x_{\left\langle i \right\rangle}}} \left( {r_{\left\langle i \right\rangle}{\Lambda_{\left\langle i \right\rangle}^*}\left( {{x_{\left\langle i \right\rangle}}} \right) + r_{\left\langle j \right\rangle}{\Lambda_{\left\langle j \right\rangle}^*}\left( {{x_{\left\langle j \right\rangle}}} \right)} \right)$. With Assumption 3 and (\citealp{zeitouni2010large}, Lemma 2.2.5), $\Lambda^*\left( x \right)$ is a strictly convex rate function and $\Lambda^*\left( x \right) \ge 0$, where the equality holds when $x=\mu$. In addition, both ${\Lambda_{\left\langle i \right\rangle}^*}\left( x \right)$ and ${\Lambda_{\left\langle j \right\rangle}^*}\left( x \right)$ are non-decreasing with $x$ for $x>\mu_{\left\langle i \right\rangle}$ and non-increasing with $x$ for $x<\mu_{\left\langle j \right\rangle}$. Therefore, for $x \in \left[\mu_{\left\langle j \right\rangle},\mu_{\left\langle i \right\rangle}\right]$, ${\Lambda_{\left\langle j \right\rangle}^*}\left( x \right)$ is non-decreasing with $x$ and ${\Lambda_{\left\langle i \right\rangle}^*}\left( x \right)$ is non-increasing with $x$, and we can search for the infimum for ${\mu _{\left\langle j \right\rangle}} \le {x_{\left\langle i \right\rangle}} \le {x_{\left\langle j \right\rangle}} \le {\mu _{\left\langle i \right\rangle}}$. Then for $i=1,2,\cdots,m$ and $j=m+1,m+2,\cdots,k$,
$$\begin{aligned}
{G_{ij}}\left( {r_{\left\langle i \right\rangle},r_{\left\langle j \right\rangle}} \right) & = \mathop {\inf }\limits_{x \in \left[ {{\mu _{\left\langle j \right\rangle}},{\mu _{\left\langle i \right\rangle}}} \right]} \left( {r_{\left\langle i \right\rangle}{\Lambda_{\left\langle i \right\rangle}^*}\left( x \right) + r_{\left\langle j \right\rangle}{\Lambda_{\left\langle j \right\rangle}^*}\left( x \right)} \right) \\
& = \mathop {\inf}\limits_x \left( {r_{\left\langle i \right\rangle}{\Lambda_{\left\langle i \right\rangle}^*}\left( x \right) + r_{\left\langle j \right\rangle}{\Lambda_{\left\langle j \right\rangle}^*}\left( x \right)} \right)~.
\end{aligned}$$

\subsubsection*{Proof of Lemma~\ref{lem001}}
 $\Lambda_{\ell}^{*}\left( x \right)$, $\ell = 1,2,\cdots,k$ are strictly convex rate functions and $\Lambda_{\left\langle \ell \right\rangle}^{*}\left( x \right) \ge 0$, where the equality holds when $x = \mu_{\left\langle \ell \right\rangle}$ (see~\citealp{zeitouni2010large}, Lemma 2.2.5), we have ${\left. {\frac{{d{\Lambda_{\left\langle \ell \right\rangle}^*}\left( x \right)}}{{dx}}} \right|_{x = \mu_{\left\langle \ell \right\rangle} }} = 0$, $\ell=1,2,\cdots,k$. In addition, since both $\Lambda _{\left\langle i \right\rangle}^*\left( x \right)$ and $\Lambda _{\left\langle j \right\rangle}^*\left( x \right)$ are non-decreasing with $x$ for $x >\mu_{\left\langle i \right\rangle}$ and non-increasing with $x$ for $x < \mu_{\left\langle j \right\rangle}$, we have ${\left. {\frac{{d\Lambda _{\left\langle j \right\rangle}^*\left( x \right)}}{{dx}}} \right|_{x = {\mu _{\left\langle i \right\rangle}}}} > 0$ and ${\left. {\frac{{d\Lambda _{\left\langle i \right\rangle}^*\left( x \right)}}{{dx}}} \right|_{x = {\mu _{\left\langle j \right\rangle}}}} < 0$, $i = 1,2,\cdots,m$, $j = m+1,m+2,\cdots,k$.
	
For $x \in \left[\mu_{\left\langle j \right\rangle},\mu_{\left\langle i \right\rangle}\right]$, $\frac{{d\Lambda _{\left\langle i \right\rangle}^*\left( x \right)}}{{dx}}$ and $\frac{{d\Lambda _{\left\langle j \right\rangle}^*\left( x \right)}}{{dx}}$ are increasing functions, $i=1,2,\cdots,m$, $j=m+1,m+2,\cdots,k$, since $\Lambda_{\ell}^{*}\left(x\right)$,  $\ell = 1,2,\cdots,k$ are strictly convex functions. With ${r_{\left\langle  i \right\rangle}}{\left. {\frac{{d\Lambda _{\left\langle i \right\rangle}^*\left( x \right)}}{{dx}}} \right|_{x = {\mu _{\left\langle j \right\rangle}}}} + {r_{\left\langle  j \right\rangle}}{\left. {\frac{{d\Lambda _{\left\langle j \right\rangle}^*\left( x \right)}}{{dx}}} \right|_{x = {\mu _{\left\langle j \right\rangle}}}} < 0$, and ${r_{\left\langle  i \right\rangle}}{\left. {\frac{{d\Lambda _{\left\langle i \right\rangle}^*\left( x \right)}}{{dx}}} \right|_{x = {\mu _{\left\langle i \right\rangle}}}} + {r_{\left\langle  j \right\rangle}}{\left. {\frac{{d\Lambda _{\left\langle j \right\rangle}^*\left( x \right)}}{{dx}}} \right|_{x = {\mu _{\left\langle i \right\rangle}}}} > 0$,  there exists a unique solution to ${r_{\left\langle  i \right\rangle}}\frac{{d\Lambda _{\left\langle i \right\rangle}^*\left( x \right)}}{{dx}} + {r_{\left\langle  j \right\rangle}}\frac{{d\Lambda _{\left\langle j \right\rangle}^*\left( x \right)}}{{dx}} = 0$. Therefore, the infimum point of $\left( r_{\left\langle  i \right\rangle}{\Lambda _{\left\langle i \right\rangle}^*\left(  x  \right)} + r_{\left\langle  j \right\rangle}{\Lambda _{\left\langle j \right\rangle}^*\left(  x  \right)} \right)$ is unique.

\subsubsection*{Proof of Lemma~\ref{lem}}
Following Lemma \ref{lem001}, let $x\left( {r_{\left\langle i \right\rangle},r_{\left\langle j \right\rangle}} \right)$ be the unique solution to minimize $\mathop {\inf }\limits_x \left( {{r_{\left\langle  i \right\rangle}}\Lambda {\left\langle i \right\rangle}^*\left( x \right) + {r_{\left\langle  j \right\rangle}}\Lambda _{\left\langle j \right\rangle}^*\left( x \right)} \right)$, and then ${G_{ij}}\left( {{r_{\left\langle  i \right\rangle}},{r_{\left\langle  j \right\rangle}}} \right) = {r_{\left\langle  i \right\rangle}}\Lambda _{\left\langle i \right\rangle}^*\left( {x\left( {{r_{\left\langle  i \right\rangle}},{r_{\left\langle  j \right\rangle}}} \right)} \right) + {r_{\left\langle  j \right\rangle}}\Lambda _{\left\langle j \right\rangle}^*\left( {x\left( {{r_{\left\langle  i \right\rangle}},{r_{\left\langle  j \right\rangle}}} \right)} \right)$, $i=1,2,\cdots,m$, $j=m+1,m+2,\cdots,k$. The partial derivative shows that \begin{align}
{\left. {\frac{{\partial {G_{ij}}\left( {y,{r_{\left\langle  j \right\rangle}}} \right)}}{{\partial y}}} \right|_{y = {r_{\left\langle  i \right\rangle}}}} &= {\left. {\frac{{\partial \left( {y{\Lambda _{\left\langle  i \right\rangle}^*}\left( {x\left( {y,{r_{\left\langle  j \right\rangle}}} \right)} \right) + {r_{\left\langle  j \right\rangle}}{\Lambda_{\left\langle  j \right\rangle}^*}\left( {x\left( {y,{r_{\left\langle  j \right\rangle}}} \right)} \right)} \right)}}{{\partial y}}} \right|_{y = {r_{\left\langle  i \right\rangle}}}}\nonumber\\
&= {\Lambda _{\left\langle  i \right\rangle}^*}\left( {x\left( {{r_{\left\langle  i \right\rangle}},{r_{\left\langle  j \right\rangle}}} \right)} \right) + {r_{\left\langle  i \right\rangle}}{\left. {\frac{{d\Lambda _{\left\langle  i \right\rangle}^*\left( x \right)}}{{dx}}} \right|_{x = x( {{r_{\left\langle  i \right\rangle}},{r_{\left\langle  j \right\rangle}}})}}{\left. {\frac{{\partial x\left( {y,{r_{\left\langle  j \right\rangle}}} \right)}}{{\partial y}}} \right|_{y = {r_{\left\langle  i \right\rangle}}}}\nonumber\\
&+ {r_{\left\langle  j \right\rangle}}{\left. {\frac{{d\Lambda _{\left\langle  j \right\rangle}^*\left( x \right)}}{{dx}}} \right|_{x = x( {{r_{\left\langle  i \right\rangle}},{r_{\left\langle  j \right\rangle}}})}}{\left. {\frac{{\partial x\left( {y,{r_{\left\langle  j \right\rangle}}} \right)}}{{\partial y}}} \right|_{y = {r_{\left\langle  i \right\rangle}}}}\nonumber= {\Lambda _{\left\langle  i \right\rangle}^*}\left( {x\left( {{r_{\left\langle  i \right\rangle}},{r_{\left\langle  j \right\rangle}}} \right)} \right)\nonumber~,
\end{align}
where the third equality holds since ${r_{\left\langle  i \right\rangle}}{\left. {\frac{{d\Lambda _{\left\langle  i \right\rangle}^*\left( x \right)}}{{dx}}} \right|_{x = x( {{r_{\left\langle  i \right\rangle}},{r_{\left\langle  j \right\rangle}}} )}} + {r_{\left\langle  j \right\rangle}}{\left. {\frac{{d\Lambda _{\left\langle  j \right\rangle}^*\left( x \right)}}{{dx}}} \right|_{x = x( {{r_{\left\langle  i \right\rangle}},{r_{\left\langle  j \right\rangle}}})}} = 0$ by the definition of $x\left(r_{\left\langle  i \right\rangle},r_{\left\langle  j \right\rangle}\right)$. By a similar argument, we have ${\left. {\frac{{\partial {G_{ij}}\left( {{r_{\left\langle  i \right\rangle}},y} \right)}}{{\partial y}}} \right|_{y = {r_{\left\langle  j \right\rangle}}}} = {\Lambda_{\left\langle  j \right\rangle}^{*}}\left( {x\left( {r_{\left\langle  i \right\rangle},r_{\left\langle  j \right\rangle}} \right)} \right)$.

First, we prove that
$${\left. {\frac{{\partial {G_{ij}}\left( {y,{r_{\left\langle  j \right\rangle}}} \right)}}{{\partial y}}} \right|_{y = {r_{\left\langle  i \right\rangle}}}} = {\Lambda_i^{*}}\left( {x\left( {r_{\left\langle  i \right\rangle},r_{\left\langle  j \right\rangle}} \right)} \right) > 0,\quad i=1,2,\cdots,m~,$$
and
$${\left. {\frac{{\partial {G_{ij}}\left( {{r_{\left\langle  i \right\rangle}},y} \right)}}{{\partial y}}} \right|_{y = {r_{\left\langle  j \right\rangle}}}} = {\Lambda_j^{*}}\left( {x\left( {r_{\left\langle  i \right\rangle},r_{\left\langle  j \right\rangle}} \right)} \right) > 0,\quad j=m+1,m+2,\cdots,k~,$$
and then we can conclude that $G_{ij}\left(r_{\left\langle  i \right\rangle},r_{\left\langle  j \right\rangle}\right)$ is strictly increasing with $r_{\left\langle  i \right\rangle}$ and $r_{\left\langle  j \right\rangle}$, respectively. Following $\Lambda_{\ell}^{*}\left( x \right)$ are strictly convex rate functions and $\Lambda_{\left\langle \ell \right\rangle}^{*}\left( x \right) \ge 0$, where the equality holds when $x = \mu_{\left\langle \ell \right\rangle}$, $\ell = 1,2,\cdots,k$, the problem boils down to verifying ${\Lambda_{\left\langle  i \right\rangle}^{*}}\left( {x\left( {r_{\left\langle  i \right\rangle},r_{\left\langle  j \right\rangle}} \right)} \right) \ne 0$, $i=1,2,\cdots,m$ and ${\Lambda_{\left\langle  j \right\rangle}^{*}}\left( {x\left( {r_{\left\langle  i \right\rangle},r_{\left\langle  j \right\rangle}} \right)} \right) \ne 0$, $j=m+1,m+2,\cdots,k$.
	
We assume that ${\Lambda_{\left\langle  i \right\rangle}^{*}}\left( {x\left( {r_{\left\langle  i \right\rangle},r_{\left\langle  j \right\rangle}} \right)} \right) = 0$, $i=1,2,\cdots,m$, which implies ${x\left( {r_{\left\langle  i \right\rangle},r_{\left\langle  j \right\rangle}} \right)} = \mu_{\left\langle i \right\rangle}$ and ${\left. {\frac{{d\Lambda _{\left\langle  i \right\rangle}^*\left( x \right)}}{{dx}}} \right|_{x = x\left( {{r_{\left\langle  i \right\rangle}},{r_{\left\langle  j \right\rangle}}} \right)}} = 0$. Since ${r_{\left\langle  i \right\rangle}}{\left. {\frac{{d\Lambda _{\left\langle  i \right\rangle}^*\left( x \right)}}{{dx}}} \right|_{x = x\left( {{r_{\left\langle  i \right\rangle}},{r_{\left\langle  j \right\rangle}}} \right)}} + {r_{\left\langle  j \right\rangle}}{\left. {\frac{{d\Lambda _{\left\langle  j \right\rangle}^*\left( x \right)}}{{dx}}} \right|_{x = x\left( {{r_{\left\langle  i \right\rangle}},{r_{\left\langle  j \right\rangle}}} \right)}} = 0$ by the definition of $x\left(r_{\left\langle  i \right\rangle},r_{\left\langle  j \right\rangle}\right)$, and ${\left. {\frac{{d\Lambda _{\left\langle  j \right\rangle}^*\left( x \right)}}{{dx}}} \right|_{x = x\left( {{r_{\left\langle  i \right\rangle}},{r_{\left\langle  j \right\rangle}}} \right)}} > 0$ due to ${x\left( {r_{\left\langle  i \right\rangle},r_{\left\langle  j \right\rangle}} \right)} = \mu_{\left\langle i \right\rangle}$, we have ${r_{\left\langle  j \right\rangle}} = 0$, which contradicts the condition $\left( {r_{\left\langle  i \right\rangle},r_{\left\langle  j \right\rangle}} \right) > 0$, $i=1,2,\cdots,m$, $j=m+1,m+2,\cdots,k$. Therefore, ${\left. {\frac{{\partial {G_{ij}}\left( {y,{r_{\left\langle  j \right\rangle}}} \right)}}{{\partial y}}} \right|_{y = {r_{\left\langle  i \right\rangle}}}} = {\Lambda_{\left\langle  i \right\rangle}^{*}}\left( {x\left( {r_{\left\langle  i \right\rangle},r_{\left\langle  j \right\rangle}} \right)} \right) > 0$. By a similar argument, we have ${\left. {\frac{{\partial {G_{ij}}\left( {{r_{\left\langle  i \right\rangle}},y} \right)}}{{\partial y}}} \right|_{y = {r_{\left\langle  j \right\rangle}}}} = {\Lambda_{\left\langle  j \right\rangle}^{*}}\left( {x\left( {r_{\left\langle  i \right\rangle},r_{\left\langle  j \right\rangle}} \right)}\right) > 0$.
	
Then we prove that $G_{ij}\left(r_{\left\langle  i \right\rangle},r_{\left\langle  j \right\rangle}\right) = 0$ if $\min\left(r_{\left\langle  i \right\rangle},r_{\left\langle  j \right\rangle}\right) = 0$. If $\min \left( {{r_{\left\langle  i \right\rangle}},{r_{\left\langle  j \right\rangle}}} \right) = 0$ and we assume that $r_{\left\langle  i \right\rangle} = 0$, then $r_{\left\langle  j \right\rangle} \ge 0$ and ${G_{ij}}\left( {{r_{\left\langle  i \right\rangle}},{r_{\left\langle  j \right\rangle}}} \right) = \mathop {\inf }\limits_x \left( {{r_{\left\langle  j \right\rangle}}{\Lambda_{\left\langle  j \right\rangle}^*}\left( x \right)} \right) = 0$. Similarly, if $r_{\left\langle  j \right\rangle} = 0$, then $r_{\left\langle  i \right\rangle} \ge 0$ and ${G_{ij}}\left( {{r_{\left\langle  i \right\rangle}},{r_{\left\langle  j \right\rangle}}} \right) = \mathop {\inf }\limits_x \left( {{r_{\left\langle  i \right\rangle}}{\Lambda_{\left\langle  i \right\rangle}^*}\left( x \right)} \right) = 0$. Therefore, ${G_{ij}}\left( {{r_{\left\langle  i \right\rangle}},{r_{\left\langle  j \right\rangle}}} \right) = 0$ if $\min \left( {{r_{\left\langle  i \right\rangle}},{r_{\left\langle  j \right\rangle}}} \right) = 0$, $i=1,2,\cdots,m$, $j=m+1,m+2,\cdots,k$.

\subsubsection*{ Proof of Lemma~\ref{Gfunction_property}}
For each $x$, $\left( {r_{\left\langle  i \right\rangle}{\Lambda_{\left\langle  i \right\rangle}^*\left(  x  \right)} + r_{\left\langle  j \right\rangle}{\Lambda_{\left\langle  j \right\rangle}^*\left(  x  \right)}} \right)$ is a linear combination of $\left( {r_{\left\langle  i \right\rangle},r_{\left\langle  j \right\rangle}} \right)$, and therefore, it is a concave and continuous function of $\left( {r_{\left\langle  i \right\rangle},r_{\left\langle  j \right\rangle}} \right)$, $i=1,2,\cdots,m$, $j=m+1,m+2,\cdots,k$. Since the infimum of a concave function is also a concave function,  $G_{ij}\left(r_{\left\langle  i \right\rangle},r_{\left\langle  j \right\rangle}\right)$ is a concave and continuous function of $\left(r_{\left\langle  i \right\rangle}, r_{\left\langle  j \right\rangle}\right)$, $i \in \left\{1,2,\cdots,m\right\}$, $j \in \left\{m+1,m+2,\cdots,k\right\}$.

Consider any two points $\bar{P}_{ij} \mathop = \limits^\Delta  \left( {{\bar {r}_{\left\langle  i \right\rangle}},{\bar {r}_{\left\langle  j \right\rangle}}} \right)$ and $\widetilde P_{ij} \mathop {\rm{ = }}\limits^\Delta \left( {{{\widetilde r}_{\left\langle  i \right\rangle}},{{\widetilde r}_{\left\langle  j \right\rangle}}} \right)$ for the rate function $G_{ij}\left(\cdot,\cdot\right)$, where $\left( {{\bar {r}_{\left\langle  i \right\rangle}},{\bar {r}_{\left\langle  j \right\rangle}}} \right) \ne \left( {{{\widetilde r}_{\left\langle  i \right\rangle}},{{\widetilde r}_{\left\langle  j \right\rangle}}} \right)$,  $i \in \left\{1,2,\cdots,m\right\}$, $j \in \left\{m+1,m+2,\cdots,k\right\}$, and then we have
$$\begin{aligned}
& {\left( {\nabla {{\left. {{G_{ij}}\left( {{r_{\left\langle i \right\rangle}},{r_{\left\langle j \right\rangle}}} \right)} \right|}_{\left( {{{\widetilde r}_{\left\langle  i \right\rangle}},{{\widetilde r}_{\left\langle  j \right\rangle}}} \right)}} - \nabla {{\left. {{G_{ij}}\left( {{r_{\left\langle i \right\rangle}},{r_{\left\langle j \right\rangle}}} \right)} \right|}_{\left( {{{\bar r}_{\left\langle  i \right\rangle}},{{\bar r}_{\left\langle  j \right\rangle}}} \right)}}} \right)^\prime }\left( {{{\left( {{{\widetilde r}_{\left\langle  i \right\rangle}},{{\widetilde r}_{\left\langle  j \right\rangle}}} \right)}^\prime } - {{\left( {{{\bar r}_{\left\langle  i \right\rangle}},{{\bar r}_{\left\langle  j \right\rangle}}} \right)}^\prime }} \right) \\
= & {\left( {\begin{array}{*{20}{c}}
		{{{\left. {\frac{{\partial {G_{ij}}\left( {y,{{\widetilde r}_{\left\langle  j \right\rangle}}} \right)}}{{\partial y}}} \right|}_{y = {{\widetilde r}_{\left\langle  i \right\rangle}}}} - {{\left. {\frac{{\partial {G_{ij}}\left( {y,{{\bar r}_{\left\langle  j \right\rangle}}} \right)}}{{\partial y}}} \right|}_{y = {{\bar r}_{\left\langle  i \right\rangle}}}}}\\
		{{{\left. {\frac{{\partial {G_{ij}}\left( {{{\widetilde r}_{\left\langle  i \right\rangle}},y} \right)}}{{\partial y}}} \right|}_{y = {{\widetilde r}_{\left\langle  j \right\rangle}}}} - {{\left. {\frac{{\partial {G_{ij}}\left( {{{\bar r}_{\left\langle  i \right\rangle}},y} \right)}}{{\partial y}}} \right|}_{y = {{\bar r}_{\left\langle  j \right\rangle}}}}}
		\end{array}} \right)^\prime }\left( {\begin{array}{*{20}{c}}
	{{{\widetilde r}_{\left\langle  i \right\rangle}} - {{\bar r}_{\left\langle  i \right\rangle}}}\\
	{{{\widetilde r}_{\left\langle  j \right\rangle}} - {{\bar r}_{\left\langle  j \right\rangle}}}
	\end{array}} \right) \\
= & \left( {\Lambda_{\left\langle  i \right\rangle}^*\left( {x\left( {{{\widetilde r}_{\left\langle  i \right\rangle}},{{\widetilde r}_{\left\langle  j \right\rangle}}} \right)} \right) - \Lambda_{\left\langle  i \right\rangle}^*\left( {x\left( {{{\bar r}_{\left\langle  i \right\rangle}},{{\bar r}_{\left\langle  j \right\rangle}}} \right)} \right)} \right)\left( {{{\widetilde r}_{\left\langle  i \right\rangle}} - {{\bar r}_{\left\langle  i \right\rangle}}} \right) + \\
& \left( {\Lambda_{\left\langle  j \right\rangle}^*\left( {x\left( {{{\widetilde r}_{\left\langle  i \right\rangle}},{{\widetilde r}_{\left\langle  j \right\rangle}}} \right)} \right) - \Lambda_{\left\langle  j \right\rangle}^*\left( {x\left( {{{\bar r}_{\left\langle  i \right\rangle}},{{\bar r}_{\left\langle  j \right\rangle}}} \right)} \right)} \right)\left( {{{\widetilde r}_{\left\langle  j \right\rangle}} - {{\bar r}_{\left\langle  j \right\rangle}}} \right) \\
= & \left( {{{\bar r}_{\left\langle  i \right\rangle}}\Lambda_{\left\langle  i \right\rangle}^*\left( {x\left( {{{\bar r}_{\left\langle  i \right\rangle}},{{\bar r}_{\left\langle  j \right\rangle}}} \right)} \right) + {{\bar r}_{\left\langle  j \right\rangle}}\Lambda_{\left\langle  j \right\rangle}^*\left( {x\left( {{{\bar r}_{\left\langle  i \right\rangle}},{{\bar r}_{\left\langle  j \right\rangle}}} \right)} \right)} \right) - \left( {{{\widetilde r}_{\left\langle  i \right\rangle}}\Lambda_{\left\langle  i \right\rangle}^*\left( {x\left( {{{\bar r}_{\left\langle  i \right\rangle}},{{\bar r}_{\left\langle  j \right\rangle}}} \right)} \right) + {{\widetilde r}_{\left\langle  j \right\rangle}}\Lambda_{\left\langle  j \right\rangle}^*\left( {x\left( {{{\bar r}_{\left\langle  i \right\rangle}},{{\bar r}_{\left\langle  j \right\rangle}}} \right)} \right)} \right)\\
+ & \left( {{{\widetilde r}_{\left\langle  i \right\rangle}}\Lambda_{\left\langle  i \right\rangle}^*\left( {x\left( {{{\widetilde r}_{\left\langle  i \right\rangle}},{{\tilde r}_{\left\langle  j \right\rangle}}} \right)} \right) + {{\widetilde r}_{\left\langle  j \right\rangle}}\Lambda_{\left\langle  j \right\rangle}^*\left( {x\left( {{{\widetilde r}_{\left\langle  i \right\rangle}},{{\widetilde r}_{\left\langle  j \right\rangle}}} \right)} \right)} \right) - \left( {{{\bar r}_{\left\langle  i \right\rangle}}\Lambda_{\left\langle  i \right\rangle}^*\left( {x\left( {{{\widetilde r}_{\left\langle  i \right\rangle}},{{\widetilde r}_{\left\langle  j \right\rangle}}} \right)} \right) + {{\bar r}_{\left\langle  j \right\rangle}}\Lambda_{\left\langle  j \right\rangle}^*\left( {x\left( {{{\widetilde r}_{\left\langle  i \right\rangle}},{{\widetilde r}_{\left\langle  j \right\rangle}}} \right)} \right)} \right) \\
< &  0~,
\end{aligned}$$
where $^\prime$ denotes the transpose operation of the matrix, and the last inequality strictly holds since as shown in Lemma~\ref{lem001}, $x\left(r_{\left\langle  i \right\rangle},r_{\left\langle  j \right\rangle}\right)$ uniquely minimizes $\mathop {\inf }\limits_x \left( {{r_{\left\langle  i \right\rangle}}\Lambda_{\left\langle  i \right\rangle}^*\left( x \right) + {r_{\left\langle  j \right\rangle}}\Lambda_{\left\langle  j \right\rangle}^*\left( x \right)} \right)$, and $x\left({\widetilde r}_{\left\langle  i \right\rangle},{\widetilde r}_{\left\langle  j \right\rangle}\right)$ uniquely minimizes $\mathop {\inf }\limits_x \left( {{{\widetilde r}_{\left\langle  i \right\rangle}}\Lambda_{\left\langle  i \right\rangle}^*\left( x \right) + {{\widetilde r}_{\left\langle  j \right\rangle}}\Lambda_{\left\langle  j \right\rangle}^*\left( x \right)} \right)$. Thus, the function $G_{ij}\left(r_{\left\langle  i \right\rangle},r_{\left\langle  j \right\rangle}\right)$ restricted to the line segment connecting $\widetilde P_{ij}$ and $\bar P_{ij}$, has strictly decreasing gradient. Since $\widetilde P_{ij}$ and $\bar P_{ij}$ are two  arbitrary points for the rate function $G_{ij}\left(r_{\left\langle  i \right\rangle},r_{\left\langle  j \right\rangle}\right)$, we have $G_{ij}\left(r_{\left\langle  i \right\rangle},r_{\left\langle  j \right\rangle}\right)$ is strictly concave in $\left(r_{\left\langle  i \right\rangle},r_{\left\langle  j \right\rangle}\right)$, $i \in \left\{1,2,\cdots,m\right\}$, $j \in \left\{m+1,m+2,\cdots,k\right\}$ by~\citet{bertsekas2003convex}. Summarizing the above, the Lemma is proved.

\subsubsection*{Proof of Lemma~\ref{thm2}}
First, we show that $r_i^*>0$, $i=1,2,\cdots,k$. If there exists $r_{\left\langle \ell \right\rangle}^* = 0$, $\ell \in \left\{1,2,\cdots,k\right\}$, then from Lemma~\ref{lem}, we have ${G_{\ell j}}( {r_{\left\langle \ell \right\rangle}^*,r_{\left\langle  j \right\rangle}^*}) = 0$ or ${G_{i \ell}}( {r_{\left\langle  i \right\rangle}^*,r_{\left\langle \ell \right\rangle}^*} ) = 0$, and $z = 0$, $i=1,2,
\cdots,m$, $j=m+1,m+2,\cdots,k$. However, for ${r_i^*} = {1 \mathord{\left/{\vphantom {1 k}} \right. \kern-\nulldelimiterspace} k}$, $i=1,2,\cdots,k$, we have $z > 0$. Therefore, the optimal solution to optimization problem (\ref{opt1}) must satisfy $r_i^* > 0$, $i = 1,2,\cdots,k$.

Notice that the feasible region of optimization problem (\ref{opt2}) is a subset of that of optimization problem (\ref{opt1}), and the objective functions of problems (\ref{opt1}) and (\ref{opt2}) are identical. Then we show that the equalities hold for $G_{ij}\left(r_{\left\langle  i \right\rangle},r_{\left\langle  j \right\rangle}\right)-z \ge 0$, $i=1,2,\cdots,m$, $j=m+1,m+2,\cdots,k$ in optimization problem (\ref{opt1}) if and only if the equalities $\mathop {\min }\limits_{j = m + 1,m + 2, \cdots ,k} {G_{ij}}\left( {{r_{\left\langle  i \right\rangle}},{r_{\left\langle  j \right\rangle}}} \right) = z$, $i=1,2,\cdots,m$ and $\mathop {\min }\limits_{i = 1,2, \cdots ,m} {G_{ij}}\left( {{r_{\left\langle  i \right\rangle}},{r_{\left\langle  j \right\rangle}}} \right) = z$, $j=m+1,m+2\cdots,k$ hold. Since KKT conditions are necessary and sufficient conditions for optimization problem (\ref{opt1}), there exist multipliers $\gamma$ and $\lambda_{ij} \ge 0$, $i=1,2,\cdots,m$, $j=m+1,m+2,\cdots,k$ such that
\begin{equation}\label{K1}
\sum\limits_{i = 1}^k {r_{\left\langle  i \right\rangle}^*}  = 1~,
\end{equation}
\begin{equation}\label{K2}
\sum\limits_{i = 1}^m {\sum\limits_{j = m + 1}^k {{\lambda _{ij}} = 1} }~,
\end{equation}
\begin{equation}\label{K3}
\sum\limits_{j = m + 1}^k {{\lambda _{ij}}{\left. {\frac{{\partial {G_{ij}}\left( {y,r_{\left\langle  j \right\rangle}^*} \right)}}{{\partial y}}} \right|_{y = r_{\left\langle  i \right\rangle}^*}}}  = \gamma,\quad i = 1,2,\cdots ,m~,
\end{equation}
\begin{equation}\label{K4}
\sum\limits_{i = 1}^m {{\lambda _{ij}}{\left. {\frac{{\partial {G_{ij}}\left( {r_{\left\langle  i \right\rangle}^*,y} \right)}}{{\partial y}}} \right|_{y = r_{\left\langle  j \right\rangle}^*}}}  = \gamma,\quad j = m+1,m+2, \cdots ,k~,
\end{equation}
\begin{equation}\label{K5}
{\lambda _{ij}}\left[ {z - {G_{ij}}\left( {r_{\left\langle  i \right\rangle}^*,r_{\left\langle  j \right\rangle}^*} \right)} \right] = 0,\quad i = 1,2, \cdots, m,\;j = m + 1, m + 2,\cdots ,k~,
\end{equation}
where equations in (\ref{K5}) are complementary slackness conditions. (\ref{K2}) implies at least one $\lambda_{ij} > 0$, that is,
$$\mathop {\max }\limits_{\scriptstyle i = 1,2, \cdots ,m\atop \scriptstyle j = m + 1,m + 2, \cdots ,k} {\lambda _{ij}} > 0~.$$

Then with Lemma~\ref{lem}, ${\left. {\frac{{\partial {G_{ij}}\left( {y,r_{\left\langle  j \right\rangle}^*} \right)}}{{\partial y}}} \right|_{y = r_{\left\langle  i \right\rangle}^*}} > 0$ and ${\left. {\frac{{\partial {G_{ij}}\left( {r_{\left\langle  i \right\rangle}^*,y} \right)}}{{\partial y}}} \right|_{y = r_{\left\langle  j \right\rangle}^*}} > 0$, it follows that $\gamma  > 0$. With (\ref{K3}) and (\ref{K4}), we can conclude that
$$\mathop {\max }\limits_{j = m + 1,m+2, \cdots ,k}  {{\lambda _{ij}}}  > 0,\quad i=1,2,\cdots,m~,$$
and
$$\mathop {\max }\limits_{i = 1,2, \cdots ,m} {{\lambda _{ij}}}  > 0,\quad j=m+1,m+2,\cdots,k~.$$

Therefore, the complementary slackness implies that each of these constraints is binding in (\ref{K5}), yielding
\begin{equation}
\mathop {\min }\limits_{j = m + 1,m+2, \cdots ,k} {G_{ij}}\left( {r_{\left\langle  i \right\rangle}^*,r_{\left\langle  j \right\rangle}^*} \right) = z,\quad i=1,2,\cdots,m~,
\end{equation}

\begin{equation}
\mathop {\min }\limits_{i = 1,2, \cdots ,m} {G_{ij}}\left( {r_{\left\langle  i \right\rangle}^*,r_{\left\langle  j \right\rangle}^*} \right) = z,\quad j=m+1,m+2,\cdots,k~.
\end{equation}

Therefore, a solution satisfying the KKT conditions of optimization problem (\ref{opt1}) must satisfy the equality constraints in optimization problem (\ref{opt2}). In addition, if ${G_{ij}} ( {r_{\left\langle  i \right\rangle}^*,r_{\left\langle  j \right\rangle}^*} ) \ne z$, the inequality ${G_{ij}} ( {r_{\left\langle  i \right\rangle}^*,r_{\left\langle  j \right\rangle}^*} ) > z$, $i \in \left\{1,2,\cdots,m\right\}$, $j \in \left\{m+1,m+2,\cdots,k\right\}$ strictly holds in optimization problem (\ref{opt1}), which proves the theorem.

\subsubsection*{Proof of Lemma~\ref{normal}}
Since $X_{\widetilde{\ell},t_{\widetilde{\ell}}}$ follows i.i.d. normal distribution $N(\mu_{\widetilde{\ell}},\sigma_{\widetilde{\ell}}^2 )$, ${\widetilde{\ell}}=1,2,\cdots,k$, we have, for $i=1,2,\cdots,m$, $j=m+1,m+2,\cdots,k$,

\begin{equation}\label{lammaequ1}
\begin{aligned}
{\left. {\frac{{\partial {G_{ij}}\left( {y,{r_{\left\langle  j \right\rangle}}} \right)}}{{\partial y}}} \right|_{y = {r_{\left\langle  i \right\rangle}}}} & = \frac{{{{\left( {{\mu _{\left\langle i \right\rangle}} - {\mu _{\left\langle j \right\rangle}}} \right)}^2}\sigma _{\left\langle i \right\rangle}^2}}{{2{{\left( {{{\sigma _{\left\langle i \right\rangle}^2} \mathord{\left/{\vphantom {{\sigma _{\left\langle i \right\rangle}^2} {r_{\left\langle  i \right\rangle}}}} \right.\kern-\nulldelimiterspace} {r_{\left\langle  i \right\rangle}}} + {{\sigma _{\left\langle j \right\rangle}^2} \mathord{\left/{\vphantom {{\sigma _{\left\langle j \right\rangle}^2} {r_{\left\langle  j \right\rangle}}}} \right.\kern-\nulldelimiterspace} {r_{\left\langle  j \right\rangle}}}} \right)}^2}{r_{\left\langle  i \right\rangle}^2}}}~,
\end{aligned}
\end{equation}

\begin{equation}\label{lammaequ2}
\begin{aligned}
{\left. {\frac{{\partial {G_{ij}}\left( {{r_{\left\langle  i \right\rangle}},y} \right)}}{{\partial y}}} \right|_{y = {r_{\left\langle  j \right\rangle}}}} &= \frac{{{{\left( {{\mu _{\left\langle i \right\rangle}} - {\mu _{\left\langle j \right\rangle}}} \right)}^2}\sigma _{\left\langle j \right\rangle}^2}}{{2{{\left( {{{\sigma _{\left\langle i \right\rangle}^2} \mathord{\left/{\vphantom {{\sigma _{\left\langle i \right\rangle}^2} {r_{\left\langle  i \right\rangle}}}} \right.\kern-\nulldelimiterspace} {r_{\left\langle  i \right\rangle}}} + {{\sigma _{\left\langle j \right\rangle}^2} \mathord{\left/{\vphantom {{\sigma _{\left\langle j \right\rangle}^2} {r_{\left\langle  j \right\rangle}}}} \right.\kern-\nulldelimiterspace} {r_{\left\langle  j \right\rangle}}}} \right)}^2}{ r_{\left\langle  j \right\rangle} ^2}}}~.
\end{aligned}
\end{equation}

With ${G_{ij}}\left( {{r_{\left\langle  i \right\rangle}},{r_{\left\langle  j \right\rangle}}} \right) = {G_{hj}}\left( {{r_{\left\langle  h \right\rangle}},{r_{\left\langle  j \right\rangle}}} \right)$, $i,h=1,2,\cdots,m$, $j=m+1,m+2,\cdots,k$, we have
$$\frac{{{{\left( {{\mu _{\left\langle i \right\rangle}} - {\mu _{\left\langle j \right\rangle}}} \right)}^2}}}{{2{\left( {{{\sigma _{\left\langle i \right\rangle}^2} \mathord{\left/{\vphantom {{\sigma _{\left\langle i \right\rangle}^2} {r_{\left\langle  i \right\rangle}}}} \right.\kern-\nulldelimiterspace} {r_{\left\langle  i \right\rangle}}} + {{\sigma _{\left\langle j \right\rangle}^2} \mathord{\left/{\vphantom {{\sigma _{\left\langle j \right\rangle}^2} {r_{\left\langle  j \right\rangle}}}} \right.\kern-\nulldelimiterspace} {r_{\left\langle  j \right\rangle}}}} \right)}}} = \frac{{{{\left( {{\mu _{\left\langle h \right\rangle}} - {\mu _{\left\langle j \right\rangle}}} \right)}^2}}}{{2{\left( {{{\sigma _{\left\langle h \right\rangle}^2} \mathord{\left/{\vphantom {{\sigma _{\left\langle h \right\rangle}^2} {r_{\left\langle h \right\rangle}}}} \right.\kern-\nulldelimiterspace} {r_{\left\langle h \right\rangle}}} + {{\sigma _{\left\langle j \right\rangle}^2} \mathord{\left/{\vphantom {{\sigma _{\left\langle j \right\rangle}^2} {r_{\left\langle  j \right\rangle}}}} \right.\kern-\nulldelimiterspace} {r_{\left\langle  j \right\rangle}}}} \right)}}}~.$$

Therefore,
$$\frac{{{{\left. {{{\partial {G_{ij}}\left( {{r_{\left\langle  i \right\rangle}},y} \right)} \mathord{\left/{\vphantom {{\partial {G_{ij}}\left( {{r_{\left\langle  i \right\rangle}},y} \right)} {\partial y}}} \right.\kern-\nulldelimiterspace} {\partial y}}} \right|}_{y = {r_{\left\langle  j \right\rangle}}}}}}{{{{\left. {{{\partial {G_{hj}}\left( {{r_{\left\langle h \right\rangle}},y} \right)} \mathord{\left/{\vphantom {{\partial {G_{hj}}\left( {{r_{\left\langle h \right\rangle}},y} \right)} {\partial y}}} \right.\kern-\nulldelimiterspace} {\partial y}}} \right|}_{y = {r_{\left\langle  j \right\rangle}}}}}} = \frac{{{{\left( {{\mu _{\left\langle i \right\rangle}} - {\mu _{\left\langle j \right\rangle}}} \right)}^2}}}{{{{\left( {{\mu _{\left\langle h \right\rangle}} - {\mu _{\left\langle j \right\rangle}}} \right)}^2}}} \cdot \frac{{{{\left( {{{\sigma _{\left\langle h \right\rangle}^2} \mathord{\left/{\vphantom {{\sigma _{\left\langle h \right\rangle}^2} {{r_{\left\langle h \right\rangle}}}}} \right.\kern-\nulldelimiterspace} {{r_{\left\langle h \right\rangle}}}} + {{\sigma _{\left\langle j \right\rangle}^2} \mathord{\left/{\vphantom {{\sigma _{\left\langle j \right\rangle}^2} {{r_{\left\langle  j \right\rangle}}}}} \right.\kern-\nulldelimiterspace} {{r_{\left\langle  j \right\rangle}}}}} \right)}^2}}}{{{{\left( {{{\sigma _{\left\langle i \right\rangle}^2} \mathord{\left/{\vphantom {{\sigma _{\left\langle i \right\rangle}^2} {{r_{\left\langle  i \right\rangle}}}}} \right.\kern-\nulldelimiterspace} {{r_{\left\langle  i \right\rangle}}}} + {{\sigma _{\left\langle j \right\rangle}^2} \mathord{\left/{\vphantom {{\sigma _{\left\langle j \right\rangle}^2} {{r_{\left\langle  j \right\rangle}}}}} \right.\kern-\nulldelimiterspace} {{r_{\left\langle  j \right\rangle}}}}} \right)}^2}}} = \frac{{{{\sigma _{\left\langle h \right\rangle}^2} \mathord{\left/{\vphantom {{\sigma _{\left\langle h \right\rangle}^2} {{r_{\left\langle h \right\rangle}}}}} \right.\kern-\nulldelimiterspace} {{r_{\left\langle h \right\rangle}}}} + {{\sigma _{\left\langle j \right\rangle}^2} \mathord{\left/{\vphantom {{\sigma _{\left\langle j \right\rangle}^2} {{r_{\left\langle  j \right\rangle}}}}} \right.\kern-\nulldelimiterspace} {{r_{\left\langle  j \right\rangle}}}}}}{{{{\sigma _{\left\langle i \right\rangle}^2} \mathord{\left/{\vphantom {{\sigma _{\left\langle i \right\rangle}^2} {{r_{\left\langle  i \right\rangle}}}}} \right.\kern-\nulldelimiterspace} {{r_{\left\langle  i \right\rangle}}}} + {{\sigma _{\left\langle j \right\rangle}^2} \mathord{\left/{\vphantom {{\sigma _{\left\langle j \right\rangle}^2} {{r_{\left\langle  j \right\rangle}}}}} \right.\kern-\nulldelimiterspace} {{r_{\left\langle  j \right\rangle}}}}}}~.$$

Similarly, if ${G_{ij}}\left( {{r_{\left\langle  i \right\rangle}},{r_{\left\langle  j \right\rangle}}} \right) = {G_{i\ell}}\left( {{r_{\left\langle  i \right\rangle}},{r_{\left\langle \ell \right\rangle}}} \right)$, $i=1,2,\cdots,m$, $j,\ell=m+1,m+2,\cdots,k$, we have
$$\frac{{{{\left( {{\mu _{\left\langle i \right\rangle}} - {\mu _{\left\langle j \right\rangle}}} \right)}^2}}}{{2{\left( {{{\sigma _{\left\langle i \right\rangle}^2} \mathord{\left/{\vphantom {{\sigma _{\left\langle i \right\rangle}^2} {r_{\left\langle  i \right\rangle}}}} \right.\kern-\nulldelimiterspace} {r_{\left\langle  i \right\rangle}}} + {{\sigma _{\left\langle j \right\rangle}^2} \mathord{\left/{\vphantom {{\sigma _{\left\langle j \right\rangle}^2} {r_{\left\langle  j \right\rangle}}}} \right.\kern-\nulldelimiterspace} {r_{\left\langle  j \right\rangle}}}} \right)}}} = \frac{{{{\left( {{\mu _{\left\langle i \right\rangle}} - {\mu _{\left\langle \ell \right\rangle}}} \right)}^2}}}{{2{\left( {{{\sigma _{\left\langle i \right\rangle}^2} \mathord{\left/{\vphantom {{\sigma _{\left\langle i \right\rangle}^2} {r_{\left\langle  i \right\rangle}}}} \right.\kern-\nulldelimiterspace} {r_{\left\langle  i \right\rangle}}} + {{\sigma _{\left\langle \ell \right\rangle}^2} \mathord{\left/{\vphantom {{\sigma _{\left\langle \ell \right\rangle}^2} {r_{\left\langle \ell \right\rangle}}}} \right.\kern-\nulldelimiterspace} {r_{\left\langle \ell \right\rangle}}}}\right)}}}~.$$

Then
$$\frac{{{{\left. {{{\partial {G_{ij}}\left( {y,{r_{\left\langle  j \right\rangle}}} \right)} \mathord{\left/{\vphantom {{\partial {G_{ij}}\left( {y,{r_{\left\langle  j \right\rangle}}} \right)} {\partial y}}} \right.
\kern-\nulldelimiterspace} {\partial y}}} \right|}_{y = {r_{\left\langle  i \right\rangle}}}}}}{{{{\left. {{{\partial {G_{i\ell}}\left( {y,{r_{\left\langle \ell \right\rangle}}} \right)} \mathord{\left/{\vphantom {{\partial {G_{i\ell}}\left( {y,{{\left\langle \ell \right\rangle}\ell}} \right)} {\partial y}}} \right.\kern-\nulldelimiterspace} {\partial y}}} \right|}_{y = {r_{\left\langle  i \right\rangle}}}}}} = \frac{{{{\left( {{\mu _{\left\langle i \right\rangle}} - {\mu _{\left\langle j \right\rangle}}} \right)}^2}}}{{{{\left( {{\mu _{\left\langle i \right\rangle}} - {\mu _{\left\langle \ell \right\rangle}}} \right)}^2}}} \cdot \frac{{{{\left( {{{\sigma _{\left\langle i \right\rangle}^2} \mathord{\left/{\vphantom {{\sigma _{\left\langle i \right\rangle}^2} {{r_{\left\langle  i \right\rangle}}}}} \right.\kern-\nulldelimiterspace} {{r_{\left\langle  i \right\rangle}}}} + {{\sigma _{\left\langle \ell \right\rangle}^2} \mathord{\left/{\vphantom {{\sigma _{\left\langle \ell \right\rangle}^2} {{r_{\left\langle \ell \right\rangle}}}}} \right.\kern-\nulldelimiterspace} {{r_{\left\langle \ell \right\rangle}}}}} \right)}^2}}}{{{{\left( {{{\sigma _{\left\langle i \right\rangle}^2} \mathord{\left/{\vphantom {{\sigma _{\left\langle i \right\rangle}^2} {{r_{\left\langle  i \right\rangle}}}}} \right.\kern-\nulldelimiterspace} {{r_{\left\langle  i \right\rangle}}}} + {{\sigma _{\left\langle j \right\rangle}^2} \mathord{\left/{\vphantom {{\sigma _{\left\langle j \right\rangle}^2} {{r_{\left\langle  j \right\rangle}}}}} \right.\kern-\nulldelimiterspace} {{r_{\left\langle  j \right\rangle}}}}} \right)}^2}}} = \frac{{{{\sigma _{\left\langle i \right\rangle}^2} \mathord{\left/{\vphantom {{\sigma _{\left\langle i \right\rangle}^2} {{r_{\left\langle  i \right\rangle}}}}} \right.\kern-\nulldelimiterspace} {{r_{\left\langle  i \right\rangle}}}} + {{\sigma _{\left\langle \ell \right\rangle}^2} \mathord{\left/{\vphantom {{\sigma _{\left\langle \ell \right\rangle}^2} {{r_{\left\langle \ell \right\rangle}}}}} \right.\kern-\nulldelimiterspace} {{r_{\left\langle \ell \right\rangle}}}}}}{{{{\sigma _{\left\langle i \right\rangle}^2} \mathord{\left/{\vphantom {{\sigma _{\left\langle i \right\rangle}^2} {{r_{\left\langle  i \right\rangle}}}}} \right.\kern-\nulldelimiterspace} {{r_{\left\langle  i \right\rangle}}}} + {{\sigma _{\left\langle j \right\rangle}^2} \mathord{\left/{\vphantom {{\sigma _{\left\langle j \right\rangle}^2} {{r_{\left\langle  j \right\rangle}}}}} \right.\kern-\nulldelimiterspace} {{r_{\left\langle  j \right\rangle}}}}}}~.$$
			
\subsubsection*{Proof of Theorem~\ref{thmLDTnormal}}
			
Substitute the rate functions $G_{ij}\left(r_{\left\langle  i \right\rangle},r_{\left\langle  j \right\rangle}\right)$, $i=1,2,\cdots,m$, $j=m+1,m+2,\cdots,k$ with normal sampling distributions in (\ref{normalrate}) into (\ref{optopt}), and then we have (\ref{optrate}) holds. Then we prove that (\ref{thmbalance}) holds for all cases in (\ref{optrate}).
			
For $h \in I^{\left(j\right)}$, $\ell \in J^{\left(i\right)}$, $i = 1,2,\cdots,m$, $j = m+1,m+2,\cdots,k$, with Lemma~\ref{normal}, (\ref{K3}) and (\ref{K4}) lead to
\begin{equation}\label{thmnormalbab1}
{\left. {\frac{{\partial {G_{i\ell}}\left( {y,r_{\left\langle \ell \right\rangle}^*} \right)}}{{\partial y}}} \right|_{y = r_{\left\langle  i \right\rangle}^*}} \left( {{\lambda _{i\ell}} + \sum\limits_{\scriptstyle j = m + 1\atop\scriptstyle j \ne \ell}^k {\frac{{{{\sigma _{\left\langle i \right\rangle}^2} \mathord{\left/{\vphantom {{\sigma _{\left\langle i \right\rangle}^2} {{r_{\left\langle  i \right\rangle}^*}}}} \right.\kern-\nulldelimiterspace} {{r_{\left\langle  i \right\rangle}^*}}} + {{\sigma _{\left\langle \ell \right\rangle}^2} \mathord{\left/{\vphantom {{\sigma _{\left\langle \ell \right\rangle}^2} {{r_{\left\langle \ell \right\rangle}^*}}}} \right.\kern-\nulldelimiterspace} {{r_{\left\langle \ell \right\rangle}^*}}}}}{{{{\sigma _{\left\langle i \right\rangle}^2} \mathord{\left/{\vphantom {{\sigma _{\left\langle i \right\rangle}^2} {{r_{\left\langle  i \right\rangle}^*}}}} \right.\kern-\nulldelimiterspace} {{r_{\left\langle  i \right\rangle}^*}}} + {{\sigma _{\left\langle j \right\rangle}^2} \mathord{\left/{\vphantom {{\sigma _{\left\langle j \right\rangle}^2} {{r_{\left\langle  j \right\rangle}^*}}}} \right.\kern-\nulldelimiterspace} {{r_{\left\langle  j \right\rangle}^*}}}}}{\lambda _{ij}}} } \right) = \gamma,\quad i=1,2,\cdots,m~,
\end{equation}
and
\begin{equation}\label{thmnormalbab2}
{\left. {\frac{{\partial {G_{hj}}\left( {r_{\left\langle h \right\rangle}^*,y} \right)}}{{\partial y}}} \right|_{y = r_{\left\langle  j \right\rangle}^*}}\left( {{\lambda _{hj}} + \sum\limits_{\scriptstyle i = 1\atop\scriptstyle i \ne h}^m {\frac{{{{\sigma _{\left\langle h \right\rangle}^2} \mathord{\left/{\vphantom {{\sigma _{\left\langle h \right\rangle}^2} {{r_{\left\langle h \right\rangle}^*}}}} \right.\kern-\nulldelimiterspace} {{r_{\left\langle h \right\rangle}^*}}} + {{\sigma _{\left\langle j \right\rangle}^2} \mathord{\left/{\vphantom {{\sigma _{\left\langle j \right\rangle}^2} {{r_{\left\langle  j \right\rangle}^*}}}} \right.\kern-\nulldelimiterspace} {{r_{\left\langle  j \right\rangle}^*}}}}}{{{{\sigma _{\left\langle i \right\rangle}^2} \mathord{\left/{\vphantom {{\sigma _{\left\langle i \right\rangle}^2} {{r_{\left\langle  i \right\rangle}^*}}}} \right.\kern-\nulldelimiterspace} {{r_{\left\langle  i \right\rangle}^*}}} + {{\sigma _{\left\langle j \right\rangle}^2} \mathord{\left/{\vphantom {{\sigma _{\left\langle j \right\rangle}^2} {{r_{\left\langle  j \right\rangle}^*}}}} \right.\kern-\nulldelimiterspace} {{r_{\left\langle  j \right\rangle}^*}}}}}{\lambda _{ij}}} } \right) = \gamma,\quad j=m+1,m+2,\cdots,k~,
\end{equation}
respectively. For ${\ell \in J^{\left(i\right)}}$, $i=1,2,\cdots,m$, (\ref{thmnormalbab1}) and (\ref{thmnormalbab2}) yield
\begin{equation}\label{thmdif1}
\frac{{{{\left. {{{\partial {G_{i\ell}}\left( {y,r_{\left\langle \ell \right\rangle}^*} \right)} \mathord{\left/{\vphantom {{\partial {G_{i\ell}}\left( {y,r_\ell^*} \right)} {\partial y}}} \right.\kern-\nulldelimiterspace} {\partial y}}} \right|}_{y = r_{\left\langle  i \right\rangle}^*}}}}{{{{\left. {{{\partial {G_{i\ell}}\left( {r_{\left\langle  i \right\rangle}^*,y} \right)} \mathord{\left/{\vphantom {{\partial {G_{i\ell}}\left( {r_{\left\langle  i \right\rangle}^*,y} \right)} {\partial y}}} \right.\kern-\nulldelimiterspace} {\partial y}}} \right|}_{y = r_{\left\langle \ell \right\rangle}^*}}}} = \frac{{{\lambda _{i\ell}} + \sum\limits_{\scriptstyle \widetilde h = 1\atop\scriptstyle \widetilde h \ne i}^m {\frac{{{{\sigma _{\left\langle i \right\rangle}^2} \mathord{\left/{\vphantom {{\sigma _{\left\langle i \right\rangle}^2} {{r_{\left\langle  i \right\rangle}^*}}}} \right.\kern-\nulldelimiterspace} {{r_{\left\langle  i \right\rangle}^*}}} + {{\sigma _{\left\langle \ell \right\rangle}^2} \mathord{\left/{\vphantom {{\sigma _{\left\langle \ell \right\rangle}^2} {{r_{\left\langle \ell \right\rangle}^*}}}} \right.\kern-\nulldelimiterspace} {{r_{\left\langle \ell \right\rangle}^*}}}}}{{{{\sigma _{\left\langle {\widetilde h} \right\rangle}^2} \mathord{\left/{\vphantom {{\sigma _{\left\langle {\widetilde h} \right\rangle}^2} {{r_{\left\langle {\widetilde h} \right\rangle}^*}}}} \right.\kern-\nulldelimiterspace} {{r_{\left\langle {\widetilde h} \right\rangle}^*}}} + {{\sigma _{\left\langle \ell \right\rangle}^2} \mathord{\left/{\vphantom {{\sigma _{\left\langle \ell \right\rangle}^2} {{r_{\left\langle \ell \right\rangle}^*}}}} \right.\kern-\nulldelimiterspace} {{r_{\left\langle \ell \right\rangle}^*}}}}}{\lambda _{\widetilde h \ell}}} }}{{{\lambda _{i\ell}} + \sum\limits_{\scriptstyle j = m + 1\atop\scriptstyle j \ne \ell}^k {\frac{{{{\sigma _{\left\langle i \right\rangle}^2} \mathord{\left/{\vphantom {{\sigma _{\left\langle i \right\rangle}^2} {{r_{\left\langle  i \right\rangle}^*}}}} \right.\kern-\nulldelimiterspace} {{r_{\left\langle  i \right\rangle}^*}}} + {{\sigma _{\left\langle \ell \right\rangle}^2} \mathord{\left/{\vphantom {{\sigma _{\left\langle \ell \right\rangle}^2} {{r_{\left\langle \ell \right\rangle}^*}}}} \right.\kern-\nulldelimiterspace} {{r_{\left\langle \ell \right\rangle}^*}}}}}{{{{\sigma _{\left\langle i \right\rangle}^2} \mathord{\left/{\vphantom {{\sigma _{\left\langle i \right\rangle}^2} {{r_{\left\langle  i \right\rangle}^*}}}} \right.\kern-\nulldelimiterspace} {{r_{\left\langle  i \right\rangle}^*}}} + {{\sigma _{\left\langle j \right\rangle}^2} \mathord{\left/{\vphantom {{\sigma _{\left\langle j \right\rangle}^2} {{r_{\left\langle  j \right\rangle}^*}}}} \right.\kern-\nulldelimiterspace} {{r_{\left\langle  j \right\rangle}^*}}}}}{\lambda _{ij}}} }}{\rm{ = }}\frac{{\widetilde \Phi_{i\ell} }}{{\widetilde \Psi_{i\ell} }},\quad i=1,2,\cdots,m~,
\end{equation}
where
$${{\widetilde \Phi }_{i\ell}}\mathop  = \limits^\Delta  \prod\limits_{\scriptstyle j = m + 1\atop\scriptstyle j \ne \ell}^k {\left( {{{\sigma _{\left\langle i \right\rangle}^2} \mathord{\left/{\vphantom {{\sigma _{\left\langle i \right\rangle}^2} {r_{\left\langle  i \right\rangle}^*}}} \right.\kern-\nulldelimiterspace} {r_{\left\langle  i \right\rangle}^*}} + {{\sigma _{\left\langle j \right\rangle}^2} \mathord{\left/{\vphantom {{\sigma _{\left\langle j \right\rangle}^2} {r_{\left\langle  j \right\rangle}^*}}} \right.\kern-\nulldelimiterspace} {r_{\left\langle  j \right\rangle}^*}}} \right)} \left[ {\sum\limits_{\widetilde h = 1}^m {\prod\limits_{\scriptstyle \widehat h = 1\hfill\atop\scriptstyle \widehat h \ne \widetilde h\hfill}^m {\left( {{{\sigma _{\left\langle {\widehat h} \right\rangle}^2} \mathord{\left/{\vphantom {{\sigma _{\left\langle {\widehat h} \right\rangle}^2} {r_{\left\langle {\widehat h} \right\rangle}^*}}} \right.\kern-\nulldelimiterspace} {r_{\left\langle {\widehat h} \right\rangle}^*}} + {{\sigma _{\left\langle \ell \right\rangle}^2} \mathord{\left/{\vphantom {{\sigma _{\left\langle \ell \right\rangle}^2} {r_{\left\langle \ell \right\rangle}^*}}} \right.\kern-\nulldelimiterspace} {r_{\left\langle \ell \right\rangle}^*}}} \right)} } {\lambda _{\widetilde h\ell}}} \right]~,$$
$${{\widetilde \Psi }_{i\ell}}\mathop  = \limits^\Delta  \prod\limits_{\scriptstyle \widetilde h = 1\atop\scriptstyle \widetilde h \ne i}^m {\left( {{{\sigma _{\left\langle {\widetilde h} \right\rangle}^2} \mathord{\left/{\vphantom {{\sigma _{\left\langle {\widetilde h} \right\rangle}^2} {{r_{\left\langle {\widetilde h} \right\rangle}^*}}}} \right.\kern-\nulldelimiterspace} {{r_{\left\langle {\widetilde h} \right\rangle}^*}}} + {{\sigma _{\left\langle \ell \right\rangle}^2} \mathord{\left/{\vphantom {{\sigma _{\left\langle \ell \right\rangle}^2} {{r_{\left\langle \ell \right\rangle}}^*}}} \right.\kern-\nulldelimiterspace} {{r_{\left\langle \ell \right\rangle}^*}}}} \right)} \left[ {\sum\limits_{j = m + 1}^k {\prod\limits_{\scriptstyle \widetilde \ell = m + 1\atop\scriptstyle \widetilde \ell \ne j}^k {\left( {{{\sigma _{\left\langle i \right\rangle}^2} \mathord{\left/{\vphantom {{\sigma _{\left\langle i \right\rangle}^2} {r_{\left\langle  i \right\rangle}^*}}} \right.\kern-\nulldelimiterspace} {r_{\left\langle  i \right\rangle}^*}} + {{\sigma _{\left\langle {\widetilde \ell} \right\rangle}^2} \mathord{\left/{\vphantom {{\sigma _{\left\langle {\widetilde \ell} \right\rangle}^2} {r_{{\left\langle {\widetilde \ell} \right\rangle}}^*}}} \right.\kern-\nulldelimiterspace}{r_{{\left\langle {\widetilde \ell} \right\rangle}}^*}}} \right)} {\lambda _{ij}}} } \right]~,$$
and the second equality holds by multiplying the numerator and the denominator of the second term in (\ref{thmdif1}) by
$$\prod\limits_{\scriptstyle j = m + 1\atop\scriptstyle j \ne \ell}^k {\left( {{{\sigma _{\left\langle i \right\rangle}^2} \mathord{\left/{\vphantom {{\sigma _{\left\langle i \right\rangle}^2} {r_{\left\langle  i \right\rangle}^*}}} \right.\kern-\nulldelimiterspace} {r_{\left\langle  i \right\rangle}^*}} + {{\sigma _{\left\langle j \right\rangle}^2} \mathord{\left/{\vphantom {{\sigma _{\left\langle j \right\rangle}^2} {r_{\left\langle  j \right\rangle}^*}}} \right.\kern-\nulldelimiterspace} {r_{\left\langle  j \right\rangle}^*}}} \right)} \prod\limits_{\scriptstyle \widetilde h = 1\atop\scriptstyle \widetilde h \ne i}^m {\left( {{{\sigma _{\left\langle {\widetilde h} \right\rangle}^2} \mathord{\left/{\vphantom {{\sigma _{\left\langle {\widetilde h} \right\rangle}^2} {r_{\left\langle {\widetilde h} \right\rangle}^*}}} \right.\kern-\nulldelimiterspace} {r_{\left\langle {\widetilde h} \right\rangle}^*}} + {{\sigma _{\left\langle \ell \right\rangle}^2} \mathord{\left/{\vphantom {{\sigma _{\left\langle \ell \right\rangle}^2} {r_{\left\langle \ell \right\rangle}^*}}} \right.\kern-\nulldelimiterspace} {r_{\left\langle \ell \right\rangle}^*}}} \right)}~,$$
simultaneously. In addition, (\ref{lammaequ1}) and (\ref{lammaequ2}) lead to
\begin{equation}\label{term001}
\begin{aligned}
\frac{{{{\left. {{{\partial {G_{i\ell}}\left( {y,r_{\left\langle \ell \right\rangle}^*} \right)} \mathord{\left/{\vphantom {{\partial {G_{i\ell}}\left( {y,r_{\left\langle \ell \right\rangle}^*} \right)} {\partial y}}} \right.\kern-\nulldelimiterspace} {\partial y}}} \right|}_{y = r_{\left\langle  i \right\rangle}^*}}}}{{{{\left. {{{\partial {G_{i\ell}}\left( {r_{\left\langle  i \right\rangle}^*,y} \right)} \mathord{\left/{\vphantom {{\partial {G_{i\ell}}\left( {r_{\left\langle  i \right\rangle}^*,y} \right)} {\partial y}}} \right.\kern-\nulldelimiterspace} {\partial y}}} \right|}_{y = r_{\left\langle \ell \right\rangle}^*}}}} & = {{\left[ {\frac{{{{\left( {{\mu _{\left\langle i \right\rangle}} - {\mu _{\left\langle \ell \right\rangle}}} \right)}^2}\sigma _{\left\langle i \right\rangle}^2}}{{2\left( {{{\sigma _{\left\langle i \right\rangle}^2} \mathord{\left/{\vphantom {{\sigma _{\left\langle i \right\rangle}^2} {r_{\left\langle  i \right\rangle}^*}}} \right.\kern-\nulldelimiterspace} {r_{\left\langle  i \right\rangle}^*}} + {{\sigma _{\left\langle \ell \right\rangle}^2} \mathord{\left/{\vphantom {{\sigma _{\left\langle \ell \right\rangle}^2} {r_{\left\langle \ell \right\rangle}^*}}} \right.\kern-\nulldelimiterspace} {r_{\left\langle \ell \right\rangle}^*}}} \right)^2{{\left( {r_{\left\langle  i \right\rangle}^*} \right)}^2}}}} \right]} \mathord{\left/{\vphantom {{\left[ {\frac{{{{\left( {{\mu _{\left\langle i \right\rangle}} - {\mu _{\left\langle \ell \right\rangle}}} \right)}^2}\sigma _{\left\langle i \right\rangle}^2}}{{2\left( {{{\sigma _{\left\langle i \right\rangle}^2} \mathord{\left/{\vphantom {{\sigma _{\left\langle i \right\rangle}^2} {r_{\left\langle  i \right\rangle}^*}}} \right.\kern-\nulldelimiterspace} {r_{\left\langle  i \right\rangle}^*}} + {{\sigma_{\left\langle \ell \right\rangle}^2} \mathord{\left/{\vphantom {{\sigma _{\left\langle \ell \right\rangle}^2} {r_{\left\langle \ell \right\rangle}^*}}} \right.\kern-\nulldelimiterspace} {r_{\left\langle \ell \right\rangle}^*}}} \right)^2{{\left( {r_{\left\langle  i \right\rangle}^*} \right)}^2}}}} \right]} {\left[ {\frac{{{{\left( {{\mu _{\left\langle i \right\rangle}} - {\mu _{\left\langle \ell \right\rangle}}} \right)}^2}\sigma _{\left\langle \ell \right\rangle}^2}}{{2\left( {{{\sigma _{\left\langle i \right\rangle}^2} \mathord{\left/{\vphantom {{\sigma _{\left\langle i \right\rangle}^2} {r_{\left\langle  i \right\rangle}^*}}} \right.\kern-\nulldelimiterspace} {r_{\left\langle  i \right\rangle}^*}} + {{\sigma _{\left\langle \ell \right\rangle}^2} \mathord{\left/{\vphantom {{\sigma _{\left\langle \ell \right\rangle}^2} {r_{\left\langle \ell \right\rangle}^*}}} \right.\kern-\nulldelimiterspace} {r_{\left\langle \ell \right\rangle}^*}}} \right){{\left( {r_{\left\langle \ell \right\rangle}^*} \right)}^2}}}} \right]}}} \right.\kern-\nulldelimiterspace} {\left[ {\frac{{{{\left( {{\mu _{\left\langle i \right\rangle}} - {\mu _{\left\langle \ell \right\rangle}}} \right)}^2}\sigma _{\left\langle \ell \right\rangle}^2}}{{2\left( {{{\sigma _{\left\langle i \right\rangle}^2} \mathord{\left/{\vphantom {{\sigma _{\left\langle i \right\rangle}^2} {r_{\left\langle  i \right\rangle}^*}}} \right.\kern-\nulldelimiterspace} {r_{\left\langle  i \right\rangle}^*}} + {{\sigma _{\left\langle \ell \right\rangle}^2} \mathord{\left/{\vphantom {{\sigma _{\left\langle \ell \right\rangle}^2} {r_{\left\langle \ell \right\rangle}^*}}} \right.\kern-\nulldelimiterspace} {r_{\left\langle \ell \right\rangle}^*}}} \right)^2{{\left( {r_{\left\langle \ell \right\rangle}^*} \right)}^2}}}} \right]}} \\
& = \frac{{{{\sigma _{\left\langle i \right\rangle}^2} \mathord{\left/{\vphantom {{\sigma _{\left\langle i \right\rangle}^2} {{{\left( {r_{\left\langle  i \right\rangle}^*} \right)}^2}}}} \right.\kern-\nulldelimiterspace} {{{\left( {r_{\left\langle  i \right\rangle}^*} \right)}^2}}}}}{{{{\sigma _{\left\langle \ell \right\rangle}^2} \mathord{\left/{\vphantom {{\sigma _{\left\langle \ell \right\rangle}^2} {{{\left( {r_{\left\langle \ell \right\rangle}^*} \right)}^2}}}} \right.\kern-\nulldelimiterspace} {{{\left( {r_{\left\langle \ell \right\rangle}^*} \right)}^2}}}}}~.
\end{aligned}
\end{equation}
			
For ${\widehat \ell} \in \left\{m+1,m+2,\cdots,k\right\}$, ${\widehat \ell} \notin {J^{\left( i \right)}}$, $i = 1,2,\cdots,m$, i.e., $\lambda_{i{\widehat \ell}} = 0$ and $G_{i{\widehat \ell}}\left(r_{\left\langle  i \right\rangle}^*,r_{\left\langle {\widehat \ell} \right\rangle}^*\right) > z$ by the definition of $J^{\left(i\right)}$, the term with $\lambda_{i{\widehat \ell}} = 0$ in KKT conditions (\ref{K3}) and (\ref{K4}) vanishes correspondingly. Then for $h \in {I^{(\widehat \ell)}}$, $\ell \in {J^{\left(i\right)}}$, (\ref{thmnormalbab1}) and (\ref{thmnormalbab2}) lead to
\begin{equation}\label{thmdif2}
\begin{aligned}
\frac{{{{\left. {{{\partial {G_{i\ell}}\left( {y,r_{\left\langle \ell \right\rangle}^*} \right)} \mathord{\left/{\vphantom {{\partial {G_{i\ell}}\left( {y,r_{\left\langle \ell \right\rangle}^*} \right)} {\partial y}}} \right.\kern-\nulldelimiterspace} {\partial y}}} \right|}_{y = r_{\left\langle  i \right\rangle}^*}}}}{{{{\left. {{{\partial {G_{h\widehat \ell}}\left( {r_{\left\langle h \right\rangle}^*,y} \right)} \mathord{\left/{\vphantom {{\partial {G_{h\widehat \ell}}\left( {r_{\left\langle h \right\rangle}^*,y} \right)} {\partial y}}} \right.\kern-\nulldelimiterspace} {\partial y}}} \right|}_{y = r_{\left\langle {\widehat \ell} \right\rangle}^*}}}} = \frac{{{\lambda _{h\widehat \ell}} + \sum\limits_{\scriptstyle \widetilde h = 1\atop\scriptstyle \widetilde h \ne h}^m {\frac{{{{\sigma_{\left\langle h \right\rangle}^2} \mathord{\left/{\vphantom {{\sigma _{\left\langle h \right\rangle}^2} {r_{\left\langle h \right\rangle}^*}}} \right.\kern-\nulldelimiterspace} {r_{\left\langle h \right\rangle}^*}} + {{\sigma_{\left\langle {\widehat \ell} \right\rangle}^2} \mathord{\left/{\vphantom {{\sigma _{\left\langle {\widehat \ell} \right\rangle}^2} {r_{\left\langle {\widehat \ell} \right\rangle}^*}}} \right.\kern-\nulldelimiterspace} {r_{\left\langle {\widehat \ell} \right\rangle}^*}}}}{{{{\sigma _{\left\langle {\widetilde h} \right\rangle}^2} \mathord{\left/{\vphantom {{\sigma _{\left\langle {\widetilde h} \right\rangle}^2} {r_{\left\langle {\widetilde h} \right\rangle}^*}}} \right.\kern-\nulldelimiterspace} {r_{\left\langle {\widetilde h} \right\rangle}^*}} + {{\sigma _{\left\langle {\widehat \ell} \right\rangle}^2} \mathord{\left/{\vphantom {{\sigma _{\left\langle {\widehat \ell} \right\rangle}^2} {r_{\left\langle {\widehat \ell} \right\rangle}^*}}} \right.\kern-\nulldelimiterspace} {r_{\left\langle {\widehat \ell} \right\rangle}^*}}}}{\lambda _{\widetilde h\widehat \ell}}} }}{{{\lambda _{i\ell}} + \sum\limits_{\scriptstyle j = m + 1\atop\scriptstyle j \ne \ell}^k {\frac{{{{\sigma _{\left\langle i \right\rangle}^2} \mathord{\left/{\vphantom {{\sigma _{\left\langle i \right\rangle}^2} {r_{\left\langle  i \right\rangle}^*}}} \right.\kern-\nulldelimiterspace} {r_{\left\langle  i \right\rangle}^*}} + {{\sigma _{\left\langle \ell \right\rangle}^2} \mathord{\left/{\vphantom {{\sigma _{\left\langle \ell \right\rangle}^2} {r_{\left\langle \ell \right\rangle}^*}}} \right.\kern-\nulldelimiterspace} {r_{\left\langle \ell \right\rangle}^*}}}}{{{{\sigma _{\left\langle i \right\rangle}^2} \mathord{\left/{\vphantom {{\sigma _{\left\langle i \right\rangle}^2} {r_{\left\langle  i \right\rangle}^*}}} \right.\kern-\nulldelimiterspace} {r_{\left\langle  i \right\rangle}^*}} + {{\sigma _{\left\langle j \right\rangle}^2} \mathord{\left/{\vphantom {{\sigma _{\left\langle j \right\rangle}^2} {r_{\left\langle  j \right\rangle}^*}}} \right.\kern-\nulldelimiterspace} {r_{\left\langle  j \right\rangle}^*}}}}{\lambda _{ij}}} }} = \frac{{{{\widehat \Phi }_{i\widehat \ell}}}}{{{{\widehat \Psi }_{i\widehat \ell}}}},\quad i=1,2,\cdots,m~,
\end{aligned}
\end{equation}
where
$${{\widehat \Phi}_{i\widehat \ell}}\mathop  = \limits^\Delta  \prod\limits_{\scriptstyle j = m + 1\atop\scriptstyle j \ne \ell,{\widehat \ell}}^k {\left( {{{\sigma _{\left\langle i \right\rangle}^2} \mathord{\left/{\vphantom {{\sigma _{\left\langle i \right\rangle}^2} {r_{\left\langle  i \right\rangle}^*}}} \right.\kern-\nulldelimiterspace} {r_{\left\langle  i \right\rangle}^*}} + {{\sigma _{\left\langle j \right\rangle}^2} \mathord{\left/{\vphantom {{\sigma _{\left\langle j \right\rangle}^2} {r_{\left\langle  j \right\rangle}^*}}} \right.\kern-\nulldelimiterspace} {r_{\left\langle  j \right\rangle}^*}}} \right)} \left[ {\sum\limits_{\widetilde h = 1}^m {\prod\limits_{\scriptstyle \widehat h = 1\hfill\atop\scriptstyle \widehat h \ne \widetilde h\hfill}^m {{\left({{\sigma _{\left\langle {\widehat h} \right\rangle}^2} \mathord{\left/{\vphantom {{\sigma _{\left\langle {\widehat h} \right\rangle}^2} {r_{\left\langle {\widehat h} \right\rangle}^*}}} \right.\kern-\nulldelimiterspace} {r_{\left\langle {\widehat h} \right\rangle}^*}} + {{\sigma _{\left\langle {\widehat \ell} \right\rangle}^2} \mathord{\left/{\vphantom {{\sigma _{\left\langle {\widehat \ell} \right\rangle}^2} {r_{\left\langle {\widehat \ell} \right\rangle}^*}}} \right.\kern-\nulldelimiterspace} {r_{\left\langle {\widehat \ell} \right\rangle}^*}}\right)}} } {\lambda _{\widetilde h\widehat \ell}}} \right]~,$$
			
$${{\widehat \Psi }_{i\widehat \ell}}\mathop  = \limits^\Delta  \prod\limits_{\scriptstyle \widetilde h = 1\atop\scriptstyle \widetilde h \ne h,i}^m {\left( {{{\sigma _{\left\langle {\widetilde h} \right\rangle}^2} \mathord{\left/{\vphantom {{\sigma _{\left\langle {\widetilde h} \right\rangle}^2} {r_{\left\langle {\widetilde h} \right\rangle}^*}}} \right.\kern-\nulldelimiterspace} {r_{\left\langle {\widetilde h} \right\rangle}^*}} + {{\sigma _{\left\langle {\widehat \ell} \right\rangle}^2} \mathord{\left/{\vphantom {{\sigma _{\left\langle {\widehat \ell} \right\rangle}^2} {r_{\left\langle {\widehat \ell} \right\rangle}^*}}} \right.\kern-\nulldelimiterspace} {r_{\left\langle {\widehat \ell} \right\rangle}^*}}} \right)} \left[ {\sum\limits_{j = m + 1}^k {\prod\limits_{\scriptstyle \widetilde \ell = m + 1\atop\scriptstyle \widetilde \ell \ne j}^k {\left( {{{\sigma _{\left\langle i \right\rangle}^2} \mathord{\left/{\vphantom {{\sigma _{\left\langle i \right\rangle}^2} {r_{\left\langle  i \right\rangle}^*}}} \right.\kern-\nulldelimiterspace} {r_{\left\langle  i \right\rangle}^*}} + {{\sigma _{\left\langle {\widetilde \ell} \right\rangle}^2} \mathord{\left/{\vphantom {{\sigma _{\left\langle {\widetilde \ell} \right\rangle}^2} {r_{\left\langle {\widetilde \ell} \right\rangle}^*}}} \right.\kern-\nulldelimiterspace} {r_{\left\langle {\widetilde \ell} \right\rangle}^*}}} \right)} {\lambda _{ij}}} } \right]~,$$
and the second equality holds by multiplying the numerator and the denominator of the second term in (\ref{thmdif2}) by
$$\frac{{\prod\limits_{\scriptstyle j = m + 1\atop\scriptstyle j \ne \ell}^k {\left( {{{\sigma _{\left\langle i \right\rangle}^2} \mathord{\left/{\vphantom {{\sigma _{\left\langle i \right\rangle}^2} {r_{\left\langle  i \right\rangle}^*}}} \right.\kern-\nulldelimiterspace} {r_{\left\langle  i \right\rangle}^*}} + {{\sigma _{\left\langle j \right\rangle}^2} \mathord{\left/{\vphantom {{\sigma _{\left\langle j \right\rangle}^2} {r_{\left\langle  j \right\rangle}^*}}} \right.\kern-\nulldelimiterspace} {r_{\left\langle  j \right\rangle}^*}}} \right)} \prod\limits_{\scriptstyle \widetilde h = 1\atop\scriptstyle \widetilde h \ne h}^m {\left( {{{\sigma _{\left\langle {\widetilde h} \right\rangle}^2} \mathord{\left/{\vphantom {{\sigma _{\left\langle {\widetilde h} \right\rangle}^2} {r_{\left\langle {\widetilde h} \right\rangle}^*}}} \right.\kern-\nulldelimiterspace} {r_{\left\langle {\widetilde h} \right\rangle}^*}} + {{\sigma _{\left\langle {\widehat \ell} \right\rangle}^2} \mathord{\left/{\vphantom {{\sigma _{\left\langle {\widehat \ell} \right\rangle}^2} {r_{\left\langle {\widehat \ell} \right\rangle}^*}}} \right.\kern-\nulldelimiterspace} {r_{\left\langle {\widehat \ell} \right\rangle}^*}}} \right)} }}{{{{\sigma _{\left\langle i \right\rangle}^2} \mathord{\left/{\vphantom {{\sigma _{\left\langle i \right\rangle}^2} {r_{\left\langle  i \right\rangle}^*}}} \right.\kern-\nulldelimiterspace} {r_{\left\langle  i \right\rangle}^*}} + {{\sigma _{\left\langle {\widehat \ell} \right\rangle}^2} \mathord{\left/{\vphantom {{\sigma _{\left\langle {\widehat \ell} \right\rangle}^2} {r_{\left\langle {\widehat \ell} \right\rangle}^*}}} \right.\kern-\nulldelimiterspace} {r_{\left\langle {\widehat \ell} \right\rangle}^*}}}}~,$$
simultaneously. In addition, from $G_{i\ell}\left(r_{\left\langle  i \right\rangle}^*,r_{\left\langle \ell \right\rangle}^*\right) = G_{h \widehat \ell}\left(r_{\left\langle h \right\rangle}^*,r_{\left\langle {\widehat \ell} \right\rangle}^*\right)$ by (\ref{optrate}) and the definition of $I^{(\widehat \ell)}$ and $J^{\left(i\right)}$, we have
$$\frac{{{{\left( {{\mu _{\left\langle i \right\rangle}} - {\mu _{\left\langle \ell \right\rangle}}} \right)}^2}}}{{2\left( {{{\sigma _{\left\langle i \right\rangle}^2} \mathord{\left/{\vphantom {{\sigma _{\left\langle i \right\rangle}^2} {r_{\left\langle  i \right\rangle}^*}}} \right.\kern-\nulldelimiterspace} {r_{\left\langle  i \right\rangle}^*}} + {{\sigma _{\left\langle \ell \right\rangle}^2} \mathord{\left/{\vphantom {{\sigma _{\left\langle \ell \right\rangle}^2} {r_{\left\langle \ell \right\rangle}^*}}} \right.\kern-\nulldelimiterspace} {r_{\left\langle \ell \right\rangle}^*}}} \right)}} = \frac{{{{\left( {{\mu _{\left\langle h \right\rangle}} - {\mu _{\left\langle {\widehat \ell} \right\rangle}}} \right)}^2}}}{{2\left( {{{\sigma _{\left\langle h \right\rangle}^2} \mathord{\left/{\vphantom {{\sigma _{\left\langle h \right\rangle}^2} {r_{\left\langle h \right\rangle}^*}}} \right.\kern-\nulldelimiterspace} {r_{\left\langle h \right\rangle}^*}} + {{\sigma _{\left\langle {\widehat \ell} \right\rangle}^2} \mathord{\left/{\vphantom {{\sigma _{\left\langle {\widehat \ell} \right\rangle}^2} {r_{\left\langle {\widehat \ell} \right\rangle}^*}}} \right.\kern-\nulldelimiterspace} {r_{\left\langle {\widehat \ell} \right\rangle}^*}}} \right)}}~,$$
and then (\ref{lammaequ1}) and (\ref{lammaequ2}) lead to
\begin{equation}\label{term002}
\begin{aligned}
\frac{{{{\left. {{{\partial {G_{i\ell}}\left( {y,r_{\left\langle \ell \right\rangle}^*} \right)} \mathord{\left/{\vphantom {{\partial {G_{i\ell}}\left( {y,r_{\left\langle \ell \right\rangle}^*} \right)} {\partial y}}} \right.\kern-\nulldelimiterspace} {\partial y}}} \right|}_{y = r_{\left\langle  i \right\rangle}^*}}}}{{{{\left. {{{\partial {G_{h\widehat \ell}}\left( {r_{\left\langle h \right\rangle}^*,y} \right)} \mathord{\left/{\vphantom {{\partial {G_{h\widehat \ell}}\left( {r_{\left\langle h \right\rangle}^*,y} \right)} {\partial y}}} \right.\kern-\nulldelimiterspace} {\partial y}}} \right|}_{y = r_{\left\langle {\widehat \ell} \right\rangle}^*}}}} & = \frac{{{{\left( {{\mu _{\left\langle i \right\rangle}} - {\mu _{\left\langle \ell \right\rangle}}} \right)}^2}\left( {{{\sigma _{\left\langle i \right\rangle}^2} \mathord{\left/{\vphantom {{\sigma _{\left\langle i \right\rangle}^2} {{{\left( {r_{\left\langle  i \right\rangle}^*} \right)}^2}}}} \right.\kern-\nulldelimiterspace} {{{\left( {r_{\left\langle  i \right\rangle}^*} \right)}^2}}}} \right)}}{{{{\left( {{\mu _{\left\langle h \right\rangle}} - {\mu _{\left\langle {\widehat \ell} \right\rangle}}} \right)}^2}\left( {{{\sigma _{\left\langle {\widehat \ell} \right\rangle}^2} \mathord{\left/{\vphantom {{\sigma _{\left\langle {\widehat \ell} \right\rangle}^2} {{{\left( {r_{\left\langle {\widehat \ell} \right\rangle}^*} \right)}^2}}}}\right.\kern-\nulldelimiterspace} {{\left( {r_{\left\langle {\widehat \ell} \right\rangle}^*} \right)^2}}}} \right)}} \cdot \frac{{{\left( {{{\sigma _{\left\langle h \right\rangle}^2} \mathord{\left/{\vphantom {{\sigma _{\left\langle h \right\rangle}^2} {r_{\left\langle h \right\rangle}^*}}} \right.\kern-\nulldelimiterspace} {r_{\left\langle h \right\rangle}^*}} + {{\sigma _{\left\langle {\widehat \ell} \right\rangle}^2} \mathord{\left/{\vphantom {{\sigma _{\left\langle {\widehat \ell} \right\rangle}^2} {r_{\left\langle {\widehat \ell} \right\rangle}^*}}} \right.\kern-\nulldelimiterspace} {r_{\left\langle {\widehat \ell} \right\rangle}^*}}} \right)^2}}}{{{{\left( {{{\sigma _{\left\langle i \right\rangle}^2} \mathord{\left/{\vphantom {{\sigma _{\left\langle i \right\rangle}^2} {r_{\left\langle  i \right\rangle}^*}}} \right.\kern-\nulldelimiterspace} {r_{\left\langle  i \right\rangle}^*}} + {{\sigma _{\left\langle \ell \right\rangle}^2} \mathord{\left/{\vphantom {{\sigma _{\left\langle \ell \right\rangle}^2} {r_{\left\langle \ell \right\rangle}^*}}} \right.\kern-\nulldelimiterspace} {r_{\left\langle \ell \right\rangle}^*}}} \right)}^2}}} \\
& = \frac{{{{\sigma _{\left\langle h \right\rangle}^2} \mathord{\left/{\vphantom {{\sigma _{\left\langle h \right\rangle}^2} {r_{\left\langle h \right\rangle}^*}}} \right.\kern-\nulldelimiterspace} {r_{\left\langle h \right\rangle}^*}} + {{\sigma _{\left\langle {\widehat \ell} \right\rangle}^2} \mathord{\left/{\vphantom {{\sigma _{\left\langle {\widehat \ell} \right\rangle}^2} {r_{\left\langle {\widehat \ell} \right\rangle}^*}}} \right.\kern-\nulldelimiterspace} {r_{\left\langle {\widehat \ell} \right\rangle}^*}}}}{{{{\sigma _{\left\langle i \right\rangle}^2} \mathord{\left/{\vphantom {{\sigma _{\left\langle i \right\rangle}^2} {r_{\left\langle  i \right\rangle}^*}}} \right.\kern-\nulldelimiterspace} {r_{\left\langle  i \right\rangle}^*}} + {{\sigma _{\left\langle \ell \right\rangle}^2} \mathord{\left/{\vphantom {{\sigma _{\left\langle \ell \right\rangle}^2} {r_{\left\langle \ell \right\rangle}^*}}} \right.\kern-\nulldelimiterspace} {r_{\left\langle \ell \right\rangle}^*}}}} \cdot \frac{{{{\sigma _{\left\langle i \right\rangle}^2} \mathord{\left/{\vphantom {{\sigma _{\left\langle i \right\rangle}^2} {{{\left( {r_{\left\langle  i \right\rangle}^*} \right)}^2}}}} \right.\kern-\nulldelimiterspace} {{{\left( {r_{\left\langle  i \right\rangle}^*} \right)}^2}}}}}{{{{\sigma _{\left\langle {\widehat \ell} \right\rangle}^2} \mathord{\left/{\vphantom {{\sigma _{\left\langle {\widehat \ell} \right\rangle}^2} {{{\left( {r_{\left\langle {\widehat \ell} \right\rangle}^*} \right)}^2}}}} \right.\kern-\nulldelimiterspace} {{\left( {r_{\left\langle {\widehat \ell} \right\rangle}^*} \right)^2}}}}}~.
\end{aligned}
\end{equation}
				
Summarizing the above, for a fixed $i \in \left\{1,2,\cdots,m\right\}$ and $h \in J^{(\widehat \ell)}$, we have
\begin{equation}
\begin{aligned}\label{longforfixedi}
& \sum\limits_{\ell \in {J^{\left( i \right)}}} {\frac{{{{\left. {{{\partial {G_{i\ell}}\left( {y,r_{\left\langle \ell \right\rangle}^*} \right)} \mathord{\left/{\vphantom {{\partial {G_{i\ell}}\left( {y,r_{\left\langle \ell \right\rangle}^*} \right)} {\partial y}}} \right.\kern-\nulldelimiterspace} {\partial y}}} \right|}_{y = r_{\left\langle  i \right\rangle}^*}}}}{{{{\left. {{{\partial {G_{i\ell}}\left( {r_{\left\langle  i \right\rangle}^*,y} \right)} \mathord{\left/{\vphantom {{\partial {G_{i\ell}}\left( {r_{\left\langle  i \right\rangle}^*,y} \right)} {\partial y}}} \right.\kern-\nulldelimiterspace} {\partial y}}} \right|}_{y = r_{\left\langle \ell \right\rangle}^*}}}}}  + \sum\limits_{\scriptstyle \widehat \ell \in \left\{ {m + 1, \cdots ,k} \right\}\atop\scriptstyle \widehat \ell \notin {J^{\left( i \right)}}} {\frac{{{{\sigma _{\left\langle i \right\rangle}^2} \mathord{\left/{\vphantom {{\sigma _{\left\langle i \right\rangle}^2} {r_{\left\langle  i \right\rangle}^*}}} \right.\kern-\nulldelimiterspace} {r_{\left\langle  i \right\rangle}^*}} + {{\sigma _{\left\langle \ell \right\rangle}^2} \mathord{\left/{\vphantom {{\sigma _{\left\langle \ell \right\rangle}^2} {r_{\left\langle \ell \right\rangle}^*}}} \right.\kern-\nulldelimiterspace} {r_{\left\langle \ell \right\rangle}^*}}}}{{{{\sigma _{\left\langle h \right\rangle}^2} \mathord{\left/{\vphantom {{\sigma _{\left\langle h \right\rangle}^2} {r_{\left\langle h \right\rangle}^*}}} \right.\kern-\nulldelimiterspace} {r_{\left\langle h \right\rangle}^*}} + {{\sigma _{\left\langle {\widehat \ell} \right\rangle}^2} \mathord{\left/{\vphantom {{\sigma _{\left\langle {\widehat \ell} \right\rangle}^2} {r_{\left\langle {\widehat \ell} \right\rangle}^*}}} \right.\kern-\nulldelimiterspace} {r_{\left\langle {\widehat \ell} \right\rangle}^*}}}} \cdot \frac{{{{\left. {{{\partial {G_{i\ell}}\left( {y,r_{\left\langle \ell \right\rangle}^*} \right)} \mathord{\left/{\vphantom {{\partial {G_{i\ell}}\left( {y,r_{\left\langle \ell \right\rangle}^*} \right)} {\partial y}}} \right.\kern-\nulldelimiterspace} {\partial y}}} \right|}_{y = r_{\left\langle  i \right\rangle}^*}}}}{{{{\left. {{{\partial {G_{h\widehat \ell}}\left( {r_{\left\langle h \right\rangle}^*,y} \right)} \mathord{\left/{\vphantom {{\partial {G_{h\widehat \ell}}\left( {r_{\left\langle h \right\rangle}^*,y} \right)} {\partial y}}} \right.\kern-\nulldelimiterspace} {\partial y}}} \right|}_{y = r_{\left\langle {\widehat \ell} \right\rangle}^*}}}}} \\
= & \sum\limits_{\ell \in {J^{\left( i \right)}}} {\frac{{{{\widetilde \Phi }_{i\ell}}}}{{{{\widetilde \Psi }_{i\ell}}}}} + \sum\limits_{\scriptstyle \widehat \ell \in \left\{ {m + 1, \cdots ,k} \right\}\atop\scriptstyle \widehat \ell \notin {J^{\left( i \right)}}} {\frac{{{{\left( {{{\sigma _{\left\langle i \right\rangle}^2} \mathord{\left/{\vphantom {{\sigma _{\left\langle i \right\rangle}^2} {r_{\left\langle  i \right\rangle}^*}}} \right.\kern-\nulldelimiterspace} {r_{\left\langle  i \right\rangle}^*}} + {{\sigma _{\left\langle \ell \right\rangle}^2} \mathord{\left/{\vphantom {{\sigma _{\left\langle \ell \right\rangle}^2} {r_{\left\langle \ell \right\rangle}^*}}} \right.\kern-\nulldelimiterspace} {r_{\left\langle \ell \right\rangle}^*}}} \right)}} \cdot {{\widehat \Phi }_{i\widehat \ell}}}}{{{{( {{{\sigma _{\left\langle h \right\rangle}^2} \mathord{\left/{\vphantom {{\sigma _{\left\langle h \right\rangle}^2} {r_{\left\langle h \right\rangle}^*}}} \right.\kern-\nulldelimiterspace} {r_{\left\langle h \right\rangle}^*}} + {{\sigma _{\left\langle {\widehat \ell} \right\rangle}^2} \mathord{/{\vphantom {{\sigma _{\left\langle {\widehat \ell} \right\rangle}^2} {r_{\left\langle {\widehat \ell} \right\rangle}^*}}}\kern-\nulldelimiterspace} {r_{\left\langle {\widehat \ell} \right\rangle}^*}}} )}}\cdot{{\widehat \Psi }_{i\widehat \ell}}}}} \\
= & \sum\limits_{{\widetilde j} = m + 1}^k {\frac{{{{\widetilde \Phi }_{i{\widetilde j}}}}}{{{{\widetilde \Psi }_{i{\widetilde j}}}}}}=  \frac{{\sum\limits_{i = 1}^m {\sum\limits_{j = m + 1}^k {\left[ {\prod\limits_{\widetilde h = 1}^m {\prod\limits_{\widetilde \ell = m + 1}^k {\frac{{\left( {{{\sigma _{\left\langle {\widetilde h} \right\rangle}^2} \mathord{\left/{\vphantom {{\sigma _{\left\langle {\widetilde h} \right\rangle}^2} {r_{\left\langle {\widetilde h} \right\rangle}^*}}} \right.\kern-\nulldelimiterspace} {r_{\left\langle {\widetilde h} \right\rangle}^*}} + {{\sigma _{\left\langle {\widetilde \ell} \right\rangle}^2} \mathord{\left/{\vphantom {{\sigma _{\left\langle {\widetilde \ell} \right\rangle}^2} {r_{\left\langle {\widetilde \ell} \right\rangle}^*}}} \right.\kern-\nulldelimiterspace} {r_{\left\langle {\widetilde \ell} \right\rangle}^*}}} \right)}}{{\left( {{{\sigma _{\left\langle i \right\rangle}^2} \mathord{\left/{\vphantom {{\sigma _{\left\langle i \right\rangle}^2} {r_{\left\langle  i \right\rangle}^*}}} \right.\kern-\nulldelimiterspace} {r_{\left\langle  i \right\rangle}^*}} + {{\sigma _{\left\langle j \right\rangle}^2} \mathord{\left/{\vphantom {{\sigma _{\left\langle j \right\rangle}^2} {r_{\left\langle  j \right\rangle}^*}}} \right.\kern-\nulldelimiterspace} {r_{\left\langle  j \right\rangle}^*}}} \right)}}} } {\lambda _{ij}}} \right]} } }}{{\sum\limits_{j = m + 1}^k {\left[ {\prod\limits_{\widetilde h = 1}^m {\prod\limits_{\widetilde \ell = m + 1}^k {\frac{{\left( {{{\sigma _{\left\langle {\widetilde h} \right\rangle}^2} \mathord{\left/{\vphantom {{\sigma _{\left\langle {\widetilde h} \right\rangle}^2} {r_{\left\langle {\widetilde h} \right\rangle}^*}}} \right.\kern-\nulldelimiterspace} {r_{\left\langle {\widetilde h} \right\rangle}^*}} + {{\sigma _{\left\langle \widetilde \ell \right\rangle}^2} \mathord{\left/{\vphantom {{\sigma _{\left\langle {\widetilde \ell} \right\rangle}^2} {r_{\left\langle {\widetilde \ell} \right\rangle}^*}}} \right.\kern-\nulldelimiterspace} {r_{\left\langle {\widetilde \ell} \right\rangle}^*}}} \right)}}{{\left( {{{\sigma _{\left\langle i \right\rangle}^2} \mathord{\left/{\vphantom {{\sigma _{\left\langle i \right\rangle}^2} {r_{\left\langle  i \right\rangle}^*}}} \right.\kern-\nulldelimiterspace} {r_{\left\langle  i \right\rangle}^*}} + {{\sigma _{\left\langle j \right\rangle}^2} \mathord{\left/{\vphantom {{\sigma_{\left\langle j \right\rangle}^2} {r_{\left\langle  j \right\rangle}^*}}} \right.\kern-\nulldelimiterspace} {r_{\left\langle  j \right\rangle}^*}}} \right)}}} } {\lambda _{ij}}} \right]} }}~,
\end{aligned}
\end{equation}
where the second equality holds since
$${{\widetilde \Phi }_{i\widehat \ell}} = \left( {{{\sigma _{\left\langle i \right\rangle}^2} \mathord{\left/{\vphantom {{\sigma _{\left\langle i \right\rangle}^2} {r_{\left\langle  i \right\rangle}^*}}} \right.\kern-\nulldelimiterspace} {r_{\left\langle  i \right\rangle}^*}} + {{\sigma _{\left\langle \ell \right\rangle}^2} \mathord{\left/{\vphantom {{\sigma _{\left\langle \ell \right\rangle}^2} {r_{\left\langle \ell \right\rangle}^*}}} \right.\kern-\nulldelimiterspace} {r_{\left\langle \ell \right\rangle}^*}}} \right){{\widehat \Phi }_{i\widehat \ell}}~,$$
$${{\widetilde \Psi }_{i\widehat \ell}} = \left( {{{\sigma _{\left\langle  h\right\rangle}^2} \mathord{\left/{\vphantom {{\sigma _{\left\langle h \right\rangle}^2} {r_{\left\langle h \right\rangle}^*}}} \right.\kern-\nulldelimiterspace} {r_{\left\langle h \right\rangle}^*}} + {{\sigma _{\left\langle {\widehat \ell} \right\rangle}^2} \mathord{\left/{\vphantom {{\sigma _{\left\langle {\widehat \ell} \right\rangle}^2} {r_{\left\langle {\widehat \ell} \right\rangle}^*}}} \right.\kern-\nulldelimiterspace} {r_{\left\langle {\widehat \ell} \right\rangle}^*}}} \right){{\widehat \Psi }_{i\widehat \ell}}~,$$
and the third equality holds by multiplying the numerator and the denominator of the second equality by $$\prod\limits_{\scriptstyle \widetilde h = 1\hfill\atop\scriptstyle \widetilde h \ne i\hfill}^m {\prod\limits_{\scriptstyle \widetilde \ell = m + 1\atop\scriptstyle \widetilde \ell \ne \widetilde j}^k {\left( {{{\sigma _{\left\langle {\widetilde h} \right\rangle}^2} \mathord{\left/{\vphantom {{\sigma _{\left\langle {\widetilde h} \right\rangle}^2} {r_{\left\langle {\widetilde h} \right\rangle}^*}}} \right.\kern-\nulldelimiterspace} {r_{\left\langle {\widetilde h} \right\rangle}^*}} + {{\sigma _{\left\langle {\widetilde \ell} \right\rangle}^2} \mathord{\left/{\vphantom {{\sigma _{\left\langle {\widetilde \ell} \right\rangle}^2} {r_{\left\langle {\widetilde \ell} \right\rangle}^*}}} \right.\kern-\nulldelimiterspace} {r_{\left\langle {\widetilde \ell} \right\rangle}^*}}} \right)} }~,$$
simultaneously. Then, (\ref{longforfixedi}) leads to
\begin{equation}\label{toolong}
\sum\limits_{i = 1}^m {\frac{1}{{\sum\limits_{\ell \in {J^{\left( i \right)}}} {\frac{{{{\left. {{{\partial {G_{i \ell}}\left( {y,r_{\left\langle \ell \right\rangle}^*} \right)} \mathord{\left/{\vphantom {{\partial {G_{i \ell}}\left( {y,r_{\left\langle \ell \right\rangle}^*} \right)} {\partial y}}} \right.\kern-\nulldelimiterspace} {\partial y}}} \right|}_{y = r_{\left\langle  i \right\rangle}^*}}}}{{{{\left. {{{\partial {G_{i \ell}}\left( {r_{\left\langle  i \right\rangle}^*,y} \right)} \mathord{\left/{\vphantom {{\partial {G_{i \ell}}\left( {r_{\left\langle  i \right\rangle}^*,y} \right)} {\partial y}}} \right.\kern-\nulldelimiterspace} {\partial y}}} \right|}_{y = r_{\left\langle \ell \right\rangle}^*}}}} + \sum\limits_{\scriptstyle \widehat \ell \in \left\{ {m + 1, \cdots ,k} \right\}\atop\scriptstyle \widehat \ell \notin {J^{\left( i \right)}}} {\frac{{{{\sigma _{\left\langle i \right\rangle}^2} \mathord{\left/{\vphantom {{\sigma _{\left\langle i \right\rangle}^2} {r_{\left\langle  i \right\rangle}^*}}} \right.\kern-\nulldelimiterspace} {r_{\left\langle  i \right\rangle}^*}} + {{\sigma _{\left\langle \ell \right\rangle}^2} \mathord{\left/{\vphantom {{\sigma _{\left\langle \ell \right\rangle}^2} {r_{\left\langle \ell \right\rangle}^*}}} \right.\kern-\nulldelimiterspace} {r_{\left\langle \ell \right\rangle}^*}}}}{{{{\sigma _{\left\langle h \right\rangle}^2} \mathord{\left/{\vphantom {{\sigma _{\left\langle h \right\rangle}^2} {r_{\left\langle h \right\rangle}^*}}} \right.\kern-\nulldelimiterspace} {r_{\left\langle h \right\rangle}^*}} + {{\sigma _{\left\langle {\widehat \ell} \right\rangle}^2} \mathord{\left/{\vphantom {{\sigma _{\left\langle {\widehat \ell} \right\rangle}^2} {r_{\left\langle {\widehat \ell} \right\rangle}^*}}} \right.\kern-\nulldelimiterspace} {r_{\left\langle {\widehat \ell} \right\rangle}^*}}}} \cdot \frac{{{{\left. {{{\partial {G_{i \ell}}\left( {y,r_{\left\langle \ell \right\rangle}^*} \right)} \mathord{\left/{\vphantom {{\partial {G_{i \ell}}\left( {y,r_{\left\langle \ell \right\rangle}^*} \right)} {\partial y}}} \right.\kern-\nulldelimiterspace} {\partial y}}} \right|}_{y = r_{\left\langle  i \right\rangle}^*}}}}{{{{\left. {{{\partial {G_{h \widehat \ell}}\left( {r_{\left\langle h \right\rangle}^*,y} \right)} \mathord{\left/{\vphantom {{\partial {G_{h \widehat \ell}}\left( {r_{\left\langle h \right\rangle}^*,y} \right)} {\partial y}}} \right.\kern-\nulldelimiterspace} {\partial y}}} \right|}_{y = r_{\left\langle {\widehat \ell} \right\rangle}^*}}}}} } }}}  = 1~.
\end{equation}
				
With (\ref{term001}) and (\ref{term002}), (\ref{toolong}) leads to
$$\sum\limits_{i = 1}^m {\frac{1}{{\sum\limits_{\ell \in {J^{\left( i \right)}}} {\frac{{{{\sigma _{\left\langle i \right\rangle}^2} \mathord{\left/{\vphantom {{\sigma _{\left\langle i \right\rangle}^2} {{{\left( {r_{\left\langle  i \right\rangle}^{\rm{*}}} \right)}^2}}}} \right.\kern-\nulldelimiterspace} {{{\left( {r_{\left\langle  i \right\rangle}^{\rm{*}}} \right)}^2}}}}}{{{{\sigma _{\left\langle \ell \right\rangle}^2} \mathord{\left/{\vphantom {{\sigma _{\left\langle \ell \right\rangle}^2} {{{\left( {r_{\left\langle \ell \right\rangle}^{\rm{*}}} \right)}^2}}}} \right.\kern-\nulldelimiterspace} {{{\left( {r_{\left\langle \ell \right\rangle}^{\rm{*}}} \right)}^2}}}}} + \sum\limits_{\scriptstyle \widehat \ell \in \left\{ {m + 1,m+2, \cdots ,k} \right\}\atop\scriptstyle \widehat \ell \notin {J^{\left( i \right)}}} {\frac{{{{\sigma _{\left\langle i \right\rangle}^2} \mathord{\left/{\vphantom {{\sigma _{\left\langle i \right\rangle}^2} {{{\left( {r_{\left\langle  i \right\rangle}^{\rm{*}}} \right)}^2}}}} \right.\kern-\nulldelimiterspace} {{{\left( {r_{\left\langle  i \right\rangle}^{\rm{*}}} \right)}^2}}}}}{{{{\sigma _{\left\langle {\widehat \ell} \right\rangle}^2} \mathord{\left/{\vphantom {{\sigma _{\left\langle {\widehat \ell} \right\rangle}^2} {{{\left( {r_{\left\langle {\widehat \ell} \right\rangle}^{\rm{*}}} \right)}^2}}}} \right.\kern-\nulldelimiterspace} {{{\left( {r_{\left\langle {\widehat \ell} \right\rangle}^{\rm{*}}} \right)}^2}}}}}} } }}}  = \sum\limits_{i = 1}^m {\frac{1}{{\sum\limits_{j = m + 1}^k {\frac{{{{\sigma _{\left\langle i \right\rangle}^2} \mathord{\left/{\vphantom {{\sigma _{\left\langle i \right\rangle}^2} {{{\left( {r_{\left\langle  i \right\rangle}^{\rm{*}}} \right)}^2}}}} \right.\kern-\nulldelimiterspace} {{{\left( {r_{\left\langle  i \right\rangle}^{\rm{*}}} \right)}^2}}}}}{{{{\sigma _{\left\langle j \right\rangle}^2} \mathord{\left/{\vphantom {{\sigma _{\left\langle j \right\rangle}^2} {{{\left( {r_{\left\langle  j \right\rangle}^{\rm{*}}} \right)}^2}}}} \right.\kern-\nulldelimiterspace} {{{\left( {r_{\left\langle  j \right\rangle}^{\rm{*}}} \right)}^2}}}}}} }}}  = 1~,$$
which yields
\begin{equation}\label{thmcalbal}
\sum\limits_{i = 1}^m {\frac{1}{{\sum\limits_{j = m + 1}^k {\frac{{{{\sigma _{\left\langle i \right\rangle}^2} \mathord{\left/{\vphantom {{\sigma _{\left\langle i \right\rangle}^2} {{{\left( {r_{\left\langle  i \right\rangle}^{\rm{*}}} \right)}^2}}}} \right.\kern-\nulldelimiterspace} {{{\left( {r_{\left\langle  i \right\rangle}^{\rm{*}}} \right)}^2}}}}}{{{{\sigma _{\left\langle j \right\rangle}^2} \mathord{\left/{\vphantom {{\sigma _{\left\langle j \right\rangle}^2} {{{\left( {r_{\left\langle  j \right\rangle}^{\rm{*}}} \right)}^2}}}} \right.\kern-\nulldelimiterspace} {{{\left( {r_{\left\langle  j \right\rangle}^{\rm{*}}} \right)}^2}}}}}} }}}  = \sum\limits_{i = 1}^m {\frac{{{{\left( {r_{\left\langle  i \right\rangle}^*} \right)}^2} \cdot \prod\limits_{\ell = m + 1}^k {\sigma _{\left\langle \ell \right\rangle}^2} }}{{\sigma _{\left\langle i \right\rangle}^2\sum\limits_{j = m + 1}^k {\left[ {{\left( {r_{\left\langle  j \right\rangle}^{\rm{*}}} \right)^2} \cdot \prod\limits_{\scriptstyle \ell = m + 1\atop\scriptstyle \ell \ne j}^k {\sigma _{\left\langle \ell \right\rangle}^2} } \right]} }}}  = \frac{{\sum\limits_{i = 1}^m {{{\left( {r_{\left\langle  i \right\rangle}^*} \right)}^2}}  \cdot \prod\limits_{\scriptstyle \ell = 1\atop\scriptstyle \ell \ne i}^k {\sigma _{\left\langle \ell \right\rangle}^2} }}{{\sum\limits_{j = m + 1}^k {{{\left( {r_{\left\langle  j \right\rangle}^*} \right)}^2}}  \cdot \prod\limits_{\scriptstyle \ell = 1\atop\scriptstyle \ell \ne j}^k {\sigma _{\left\langle \ell \right\rangle}^2} }} = 1~.
\end{equation}
			
By simultaneously dividing the numerator and denominator of the last equality in (\ref{thmcalbal}) with $\prod\nolimits_{\ell = 1}^k {\sigma _{\left\langle \ell \right\rangle}^2}$, we have (\ref{thmbalance}) must hold for all cases in (\ref{optrate}). Summarizing the above, the theorem is proved.

\subsubsection*{Proof of Theorem~\ref{thmasymptotic}}

By the law of large numbers (LLN), $\mathop {\lim }\limits_{t \to \infty } \mu _{\left\langle i \right\rangle}^{\left( t \right)} \to {\mu _{\left\langle i \right\rangle}}$, $i = 1,2, \cdots ,k$. For simplicity of analysis, we replace the posterior mean $\mu _{\left\langle i \right\rangle}^{\left( t \right)}$ and posterior variance $( {\sigma _{\left\langle i \right\rangle}^{\left( t \right)}} )^2$ with $\mu _{\left\langle i \right\rangle}$ and ${{\sigma _{\left\langle i \right\rangle}^2} / {t_i}}$ in ${\widehat V}_t \left(\mathcal{E}_t;i\right)$, $i=1,2,\cdots,k$, since the asymptotic sampling ratios will be determined by the increasing order of ${\widehat V}_t \left(\mathcal{E}_t;i\right)$ with respect to $t$.

Note that $0 \le r_{\left\langle  i \right\rangle}^{\left(t\right)} \le 1$ is bounded, $i=1,2,\cdots,k$, and $( {r_1^{\left( t \right)},r_2^{\left( t \right)}, \cdots,r_k^{\left( t \right)}} )$ is a bounded sequence. Following the Bolzano-Weierstrass theorem~\citep{rudin1964principles}, there exists a subsequence of $( {r_1^{\left( t \right)},r_2^{\left( t \right)}, \cdots r_k^{\left( t \right)}} )$ converging to $\left( {{{\widetilde r}_1},{{\widetilde r}_2}, \cdots {{\widetilde r}_k}} \right)$, where $\sum\nolimits_{i = 1}^k {{{\widetilde r}_{\left\langle  i \right\rangle}} = 1}$, ${{\widetilde r}_{\left\langle  i \right\rangle}} \ge 0$. Without loss of generality, we assume that $( {r_1^{\left( t \right)},r_2^{\left( t \right)}, \cdots r_k^{\left( t \right)}})$ converges to $\left( {{{\widetilde r}_1},{{\widetilde r}_2}, \cdots {{\widetilde r}_k}} \right)$; otherwise, the following argument is made over a subsequence. Notice that
$$\begin{aligned}
& \mathop {\lim }\limits_{t \to \infty } \left[ {\frac{{{{\left( {{\mu _{\left\langle i \right\rangle }} - {\mu _{\left\langle j \right\rangle }}} \right)}^2}}}{{{{\sigma _{\left\langle i \right\rangle }^2} \mathord{\left/{\vphantom {{\sigma _{\left\langle i \right\rangle }^2} {t_{\left\langle i \right\rangle}}}} \right.\kern-\nulldelimiterspace} {t_{\left\langle i \right\rangle}}} + {{\sigma _{\left\langle j \right\rangle }^2} \mathord{\left/{\vphantom {{\sigma _{\left\langle j \right\rangle }^2} {\left( {{t_{\left\langle j \right\rangle}} + 1} \right)}}} \right.\kern-\nulldelimiterspace} {\left( {{t_{\left\langle j \right\rangle}} + 1} \right)}}}} - \frac{{{{\left( {{\mu _{\left\langle i \right\rangle }} - {\mu _{\left\langle j \right\rangle }}} \right)}^2}}}{{{{\sigma _{\left\langle i \right\rangle }^2} \mathord{\left/{\vphantom {{\sigma _{\left\langle i \right\rangle }^2} {t_{\left\langle i \right\rangle}}}} \right.\kern-\nulldelimiterspace} {t_{\left\langle i \right\rangle}}} + {{\sigma _{\left\langle j \right\rangle }^2} \mathord{\left/{\vphantom {{\sigma _{\left\langle j \right\rangle }^2} {{t_{\left\langle j \right\rangle}}}}} \right.\kern-\nulldelimiterspace} {{t_{\left\langle j \right\rangle}}}}}}} \right]\\
= & \mathop {\lim }\limits_{t \to \infty } t\left[ {\frac{{{{\left( {{\mu _{\left\langle i \right\rangle }} - {\mu _{\left\langle j \right\rangle }}} \right)}^2}}}{{{{\sigma _{\left\langle i \right\rangle }^2} \mathord{\left/{\vphantom {{\sigma _{\left\langle i \right\rangle }^2} {r_{\left\langle  i \right\rangle}^{\left( t \right)}}}} \right.\kern-\nulldelimiterspace} {r_{\left\langle  i \right\rangle}^{\left( t \right)}}} + {{\sigma _{\left\langle j \right\rangle }^2} \mathord{\left/{\vphantom {{\sigma _{\left\langle j \right\rangle }^2} {\left( {r_{\left\langle  j \right\rangle}^{\left( t \right)} + {1 \mathord{\left/{\vphantom {1 t}} \right.\kern-\nulldelimiterspace} t}} \right)}}} \right.\kern-\nulldelimiterspace} {\left( {r_{\left\langle  j \right\rangle}^{\left( t \right)} + {1 \mathord{\left/{\vphantom {1 t}} \right.\kern-\nulldelimiterspace} t}} \right)}}}} - \frac{{{{\left( {{\mu _{\left\langle i \right\rangle }} - {\mu _{\left\langle j \right\rangle }}} \right)}^2}}}{{{{\sigma _{\left\langle i \right\rangle }^2} \mathord{\left/{\vphantom {{\sigma _{\left\langle i \right\rangle }^2} {r_{\left\langle  i \right\rangle}^{\left( t \right)}}}} \right.\kern-\nulldelimiterspace} {r_{\left\langle  i \right\rangle}^{\left( t \right)}}} + {{\sigma _{\left\langle j \right\rangle }^2} \mathord{\left/{\vphantom {{\sigma _{\left\langle j \right\rangle }^2} {r_{\left\langle  j \right\rangle}^{\left( t \right)}}}} \right.\kern-\nulldelimiterspace} {r_{\left\langle  j \right\rangle}^{\left( t \right)}}}}}} \right]\\
= & \mathop {\lim }\limits_{t \to \infty } {\left. {\frac{{\partial {Q_{ij}}\left( {r_{\left\langle  i \right\rangle}^{\left( t \right)},x} \right)}}{{\partial x}}} \right|_{x = r_{\left\langle  j \right\rangle}^{\left( t \right)}}}\\
= & {\left( {\frac{{{\sigma _{\left\langle j \right\rangle }}}}{{\widetilde{r}_{\left\langle j \right\rangle}}}} \right)^2}\frac{{{{\left( {{\mu _{\left\langle i \right\rangle }} - {\mu _{\left\langle j \right\rangle }}} \right)}^2}}}{{{{\left( {{{\sigma _{\left\langle i \right\rangle }^2} \mathord{\left/{\vphantom {{\sigma _{\left\langle i \right\rangle }^2} {r_{\left\langle i \right\rangle }^{\left( t \right)}}}} \right.\kern-\nulldelimiterspace} {\widetilde{r}_{\left\langle i \right\rangle}}} + {{\sigma _{\left\langle j \right\rangle }^2} \mathord{\left/{\vphantom {{\sigma _{\left\langle j \right\rangle }^2} {\widetilde{r}_{\left\langle j \right\rangle}}}} \right.\kern-\nulldelimiterspace} {\widetilde{r}_{\left\langle j \right\rangle}}}} \right)}^2}}}~,
\end{aligned}$$
where
$${Q_{ij}}\left( {{r_{\left\langle  i \right\rangle}},{r_{\left\langle  j \right\rangle}}} \right) \mathop  = \limits^\Delta 2 \times {G_{ij}}\left( {{r_{\left\langle  i \right\rangle}},{r_{\left\langle  j \right\rangle}}} \right) = \frac{{{{\left( {{\mu _{\left\langle i \right\rangle }} - {\mu _{\left\langle j \right\rangle }}} \right)}^2}}}{{{{\sigma _{\left\langle i \right\rangle }^2} \mathord{\left/{\vphantom {{\sigma _{\left\langle i \right\rangle }^2} {{r_{\left\langle  i \right\rangle}} + {{\sigma _{\left\langle j \right\rangle }^2} \mathord{\left/{\vphantom {{\sigma _{\left\langle j \right\rangle }^2} {{r_{\left\langle  j \right\rangle}}}}} \right.\kern-\nulldelimiterspace} {{r_{\left\langle  j \right\rangle}}}}}}} \right.\kern-\nulldelimiterspace} {{r_{\left\langle  i \right\rangle}} + {{\sigma _{\left\langle j \right\rangle }^2} \mathord{\left/{\vphantom {{\sigma _{\left\langle j \right\rangle }^2} {{r_{\left\langle  j \right\rangle}}}}} \right.\kern-\nulldelimiterspace} {r_{\left\langle  j \right\rangle}}}}}}},\quad i=1,\cdots,m,\;j=m+1,\cdots,k~,$$
and
$$\begin{aligned}
& \mathop {\lim }\limits_{t \to \infty } \left[ {\frac{{{{\left( {{\mu _{\left\langle i \right\rangle }} - {\mu _{\left\langle j \right\rangle }}} \right)}^2}}}{{{{\sigma _{\left\langle i \right\rangle }^2} \mathord{\left/{\vphantom {{\sigma _{\left\langle i \right\rangle }^2} {\left( {{t_{\left\langle i \right\rangle}} + 1} \right) + {{\sigma _{\left\langle j \right\rangle }^2} \mathord{\left/{\vphantom {{\sigma _{\left\langle j \right\rangle }^2} {{t_{\left\langle j \right\rangle}}}}} \right.\kern-\nulldelimiterspace} {{t_{\left\langle j \right\rangle}}}}}}} \right.\kern-\nulldelimiterspace} {\left( {{t_{\left\langle i \right\rangle}} + 1} \right) + {{\sigma _{\left\langle j \right\rangle }^2} \mathord{\left/{\vphantom {{\sigma _{\left\langle j \right\rangle }^2} {{t_{\left\langle j \right\rangle}}}}} \right.\kern-\nulldelimiterspace} {{t_{\left\langle j \right\rangle}}}}}}}} - \frac{{{{\left( {{\mu _{\left\langle i \right\rangle }} - {\mu _{\left\langle j \right\rangle }}} \right)}^2}}}{{{{\sigma _{\left\langle i \right\rangle }^2} \mathord{\left/{\vphantom {{\sigma _{\left\langle i \right\rangle }^2} {{t_{\left\langle i \right\rangle}} + {{\sigma _{\left\langle j \right\rangle }^2} \mathord{\left/{\vphantom {{\sigma _{\left\langle j \right\rangle }^2} {{t_{\left\langle j \right\rangle}}}}} \right.\kern-\nulldelimiterspace} {{t_{\left\langle j \right\rangle}}}}}}} \right.\kern-\nulldelimiterspace} {{t_{\left\langle i \right\rangle}} + {{\sigma _{\left\langle j \right\rangle }^2} \mathord{\left/{\vphantom {{\sigma _{\left\langle j \right\rangle }^2} {{t_{\left\langle j \right\rangle}}}}} \right.\kern-\nulldelimiterspace} {{t_{\left\langle j \right\rangle}}}}}}}}} \right] \\
= & \mathop {\lim }\limits_{t \to \infty } t\left[ {\frac{{{{\left( {{\mu _{\left\langle i \right\rangle }} - {\mu _{\left\langle j \right\rangle }}} \right)}^2}}}{{{{\sigma _{\left\langle i \right\rangle }^2} \mathord{\left/{\vphantom {{\sigma _{\left\langle i \right\rangle }^2} {\left( {r_{\left\langle  i \right\rangle}^{\left( t \right)} + {1 \mathord{\left/{\vphantom {1 t}} \right.\kern-\nulldelimiterspace} t}} \right) + {{\sigma _{\left\langle j \right\rangle }^2} \mathord{\left/{\vphantom {{\sigma _{\left\langle j \right\rangle }^2} {r_{\left\langle  j \right\rangle}^{\left( t \right)}}}} \right.\kern-\nulldelimiterspace} {r_{\left\langle  j \right\rangle}^{\left( t \right)}}}}}} \right.\kern-\nulldelimiterspace} {\left( {r_{\left\langle  i \right\rangle}^{\left( t \right)} + {1 \mathord{\left/{\vphantom {1 t}} \right.\kern-\nulldelimiterspace} t}} \right) + {{\sigma _{\left\langle j \right\rangle }^2} \mathord{\left/{\vphantom {{\sigma _{\left\langle j \right\rangle }^2} {r_{\left\langle  j \right\rangle}^{\left( t \right)}}}} \right.\kern-\nulldelimiterspace} {r_{\left\langle  j \right\rangle}^{\left( t \right)}}}}}}} - \frac{{{{\left( {{\mu _{\left\langle i \right\rangle }} - {\mu _{\left\langle j \right\rangle }}} \right)}^2}}}{{{{\sigma _{\left\langle i \right\rangle }^2} \mathord{\left/{\vphantom {{\sigma _{\left\langle i \right\rangle }^2} {r_{\left\langle  i \right\rangle}^{\left( t \right)} + {{\sigma _{\left\langle j \right\rangle }^2} \mathord{\left/{\vphantom {{\sigma _{\left\langle j \right\rangle }^2} {r_{\left\langle  j \right\rangle}^{\left( t \right)}}}} \right.\kern-\nulldelimiterspace} {r_{\left\langle  j \right\rangle}^{\left( t \right)}}}}}} \right.\kern-\nulldelimiterspace} {r_{\left\langle  i \right\rangle}^{\left( t \right)} + {{\sigma _{\left\langle j \right\rangle }^2} \mathord{\left/{\vphantom {{\sigma _{\left\langle j \right\rangle }^2} {r_{\left\langle  j \right\rangle}^{\left( t \right)}}}} \right.\kern-\nulldelimiterspace} {r_{\left\langle  j \right\rangle}^{\left( t \right)}}}}}}}} \right]\\
= & \mathop {\lim }\limits_{t \to \infty } {\left. {\frac{{\partial {Q_{ij}}\left( {x,r_{\left\langle  j \right\rangle}^{\left( t \right)}} \right)}}{{\partial x}}} \right|_{x = r_{\left\langle  i \right\rangle}^{\left( t \right)}}}\\
= & {\left( {\frac{{{\sigma _{\left\langle i \right\rangle }}}}{{\widetilde{r}_{\left\langle i \right\rangle}}}} \right)^2}\frac{{{{\left( {{\mu _{\left\langle i \right\rangle }} - {\mu _{\left\langle j \right\rangle }}} \right)}^2}}}{{{{\left( {{{\sigma _{\left\langle i \right\rangle }^2} \mathord{\left/{\vphantom {{\sigma _{\left\langle i \right\rangle }^2} {\widetilde{r}_{\left\langle i \right\rangle} + {{\sigma _{\left\langle j \right\rangle }^2} \mathord{\left/{\vphantom {{\sigma _{\left\langle j \right\rangle }^2} {\widetilde{r}_{\left\langle j \right\rangle}}}} \right.\kern-\nulldelimiterspace} {\widetilde{r}_{\left\langle j \right\rangle}}}}}} \right.\kern-\nulldelimiterspace} {\widetilde{r}_{\left\langle i \right\rangle} + {{\sigma _{\left\langle j \right\rangle }^2} \mathord{\left/{\vphantom {{\sigma _{\left\langle j \right\rangle }^2} {\widetilde{r}_{\left\langle j \right\rangle}}}} \right.\kern-\nulldelimiterspace} {\widetilde{r}_{\left\langle j \right\rangle}}}}}} \right)}^2}}}~.
\end{aligned}$$
	
Since $\sum\nolimits_{i = 1}^k {{{\widetilde r}_{\left\langle  i \right\rangle}} = 1}$, ${\widetilde r}_{\left\langle  i \right\rangle} \ge 0$, there exists ${\widetilde r}_{\left\langle h \right\rangle} > 0$, $h \in \left\{1,2,\cdots,k\right\}$. We first claim that $\exists~h \in \left\{1,2,\cdots,m\right\}$ such that $\widetilde{r}_{\left\langle h \right\rangle}>0$; otherwise,  $\exists~j \in \left\{m+1,m+2,\cdots,k\right\}$ such that $\widetilde{r}_{\left\langle j \right\rangle} > 0$, and then for $i \in \left\{1,2,\cdots,m\right\}$,
$$\mathop {\lim }\limits_{t \to \infty } {\left( {\frac{{{\sigma _{\left\langle i \right\rangle }}}}{{r_{\left\langle  i \right\rangle}^{\left( t \right)}}}} \right)^2}\frac{{{{\left( {{\mu _{\left\langle i \right\rangle }} - {\mu _{\left\langle j \right\rangle }}} \right)}^2}}}{{\left( {{{\sigma _{\left\langle i \right\rangle }^2} \mathord{\left/{\vphantom {{\sigma _{\left\langle j \right\rangle }^2} {r_{\left\langle  i \right\rangle}^{\left( t \right)} + }}} \right.\kern-\nulldelimiterspace} {r_{\left\langle  i \right\rangle}^{\left( t \right)} + }}{{\sigma _{\left\langle j \right\rangle }^2} \mathord{\left/{\vphantom {{\sigma _{\left\langle j \right\rangle }^2} {r_{\left\langle  j \right\rangle}^{\left( t \right)}}}} \right.\kern-\nulldelimiterspace} {r_{\left\langle  j \right\rangle}^{\left( t \right)}}}} \right)^2}} >0~,$$
$$\mathop {\lim }\limits_{t \to \infty } {\left( {\frac{{{\sigma _{\left\langle j \right\rangle }}}}{{r_{\left\langle  j \right\rangle}^{\left( t \right)}}}} \right)^2}\frac{{{{\left( {{\mu _{\left\langle i \right\rangle }} - {\mu _{\left\langle j \right\rangle }}} \right)}^2}}}{{\left( {{{\sigma _{\left\langle i \right\rangle }^2} \mathord{\left/{\vphantom {{\sigma _{\left\langle j \right\rangle }^2} {r_{\left\langle  i \right\rangle}^{\left( t \right)} + }}} \right.\kern-\nulldelimiterspace} {r_{\left\langle  i \right\rangle}^{\left( t \right)} + }}{{\sigma _{\left\langle j \right\rangle }^2} \mathord{\left/{\vphantom {{\sigma _{\left\langle j \right\rangle }^2} {r_{\left\langle  j \right\rangle}^{\left( t \right)}}}} \right.\kern-\nulldelimiterspace} {r_{\left\langle  j \right\rangle}^{\left( t \right)}}}} \right)^2}} =0~,$$
which contradicts the sampling rule AOAm (\ref{AOAm}). We then claim that $\exists~\ell\in\left\{m+1,m+2,\cdots,k\right\}$ such that $\widetilde{r}_{\left\langle \ell \right\rangle} >0$; otherwise, $\exists~i \in \left\{1,2,\cdots,m\right\}$ such that $\widetilde{r}_{\left\langle i \right\rangle} > 0$, and then for $j \in \left\{m+1,m+2,\cdots,k\right\}$,
$$\mathop {\lim }\limits_{t \to \infty } {\left( {\frac{{{\sigma _{\left\langle j \right\rangle }}}}{{r_{\left\langle  j \right\rangle}^{\left( t \right)}}}} \right)^2}\frac{{{{\left( {{\mu _{\left\langle i \right\rangle }} - {\mu _{\left\langle j \right\rangle }}} \right)}^2}}}{{\left( {{{\sigma _{\left\langle i \right\rangle }^2} \mathord{\left/{\vphantom {{\sigma _{\left\langle j \right\rangle }^2} {r_{\left\langle  i \right\rangle}^{\left( t \right)} + }}} \right.\kern-\nulldelimiterspace} {r_{\left\langle  i \right\rangle}^{\left( t \right)} + }}{{\sigma _{\left\langle j \right\rangle }^2} \mathord{\left/{\vphantom {{\sigma _{\left\langle j \right\rangle }^2} {r_{\left\langle  j \right\rangle}^{\left( t \right)}}}} \right.\kern-\nulldelimiterspace} {r_{\left\langle  j \right\rangle}^{\left( t \right)}}}} \right)^2}} >0~,$$
$$\mathop {\lim }\limits_{t \to \infty } {\left( {\frac{{{\sigma _{\left\langle i \right\rangle }}}}{{r_{\left\langle  i \right\rangle}^{\left( t \right)}}}} \right)^2}\frac{{{{\left( {{\mu _{\left\langle i \right\rangle }} - {\mu _{\left\langle j \right\rangle }}} \right)}^2}}}{{\left( {{{\sigma _{\left\langle i \right\rangle }^2} \mathord{\left/{\vphantom {{\sigma _{\left\langle j \right\rangle }^2} {r_{\left\langle  i \right\rangle}^{\left( t \right)} + }}} \right.\kern-\nulldelimiterspace} {r_{\left\langle  i \right\rangle}^{\left( t \right)} + }}{{\sigma _{\left\langle j \right\rangle }^2} \mathord{\left/{\vphantom {{\sigma _{\left\langle j \right\rangle }^2} {r_{\left\langle  j \right\rangle}^{\left( t \right)}}}} \right.\kern-\nulldelimiterspace} {r_{\left\langle  j \right\rangle}^{\left( t \right)}}}} \right)^2}} =0~,$$
which contradicts the sampling rule AOAm (\ref{AOAm}). In addition, we claim that for $j = m+1,m+2,\cdots,k$ such that $\widetilde{r}_{\left\langle j \right\rangle}>0$; otherwise, $\exists~\ell, q\in\{m+1,m+2,\cdots,k\}$, $\ell \ne q$ such that $\widetilde{r}_{\left\langle \ell \right\rangle}>0$,  $\widetilde{r}_{\left\langle q \right\rangle}=0$, and for $i\in\{ 1,2,\cdots,m\}$ such that ${\widetilde r_{\left\langle  i \right\rangle}} > 0$,
\begin{align*}
\lim_{t\to\infty}{Q_{i\ell}}\left( {{r_{\left\langle  i \right\rangle}^{(t)}},{r_{\left\langle \ell \right\rangle}^{(t)}}} \right)>0,\quad \lim_{t\to\infty}{Q_{iq}}\left( {{r_{\left\langle  i \right\rangle}^{(t)}},{r_{\left\langle q \right\rangle}^{(t)}}} \right)=0~,
\end{align*}
which contradicts the sampling rule AOAm (\ref{AOAm}). Finally, we claim that for $i = 1,2,\cdots,m$ such that ${\widetilde{r}_{\left\langle i \right\rangle}} > 0$; otherwise, $\exists~h,p\in\left\{1,2,\cdots,m\right\}$, $h \ne p$ such that $\widetilde{r}_{\left\langle h \right\rangle} > 0$, $\widetilde{r}_{\left\langle p \right\rangle} = 0$, and for $j = m+1,m+2,\cdots,k$,
\begin{align*}
\lim_{t\to\infty}{Q_{hj}}\left( {{r_{\left\langle h \right\rangle}^{(t)}},{r_{\left\langle  j \right\rangle}^{(t)}}} \right)>0,\quad \lim_{t\to\infty}{Q_{pj}}\left( {{r_{\left\langle p \right\rangle}^{(t)}},{r_{\left\langle  j \right\rangle}^{(t)}}} \right)=0~,
\end{align*}
which contradicts the sampling rule AOAm (\ref{AOAm}). Summarizing the above, we have $\widetilde{r}_{\ell} >0$, $\ell =1,2,\cdots,k$.

Suppose $\left( {{{\widetilde r}_1},{{\widetilde r}_2}, \cdots ,{{\widetilde r}_k}} \right)$ does not satisfy (\ref{optrate}). If $\exists~i,h \in \left\{1,2, \cdots ,m\right\}$, $i \ne h$ such that the inequality
$$\mathop {\min }\limits_{j = m + 1,m + 2, \cdots ,k} {Q_{hj}} ( {{{\widetilde r}_{\left\langle h \right\rangle}},{{\widetilde r}_{\left\langle  j \right\rangle}}} ) > \mathop {\min }\limits_{j = m + 1,m + 2, \cdots ,k} {Q_{ij}} ( {{{\widetilde r}_{\left\langle  i \right\rangle}},{{\widetilde r}_{\left\langle  j \right\rangle}}} )~,$$
strictly holds, then there exists $T_0>0$ such that $\forall t > {T_0}$,
\begin{equation}\label{equequproof1}
\mathop {\min }\limits_{j = m + 1,m + 2, \cdots ,k}{Q_{hj}} ( {r_{\left\langle h \right\rangle}^{\left( t \right)},r_{\left\langle  j \right\rangle}^{\left( t \right)}} ) > \mathop {\min }\limits_{j = m + 1,m + 2, \cdots ,k}{Q_{ij}} ( {r_{\left\langle  i \right\rangle}^{\left( t \right)},r_{\left\langle  j \right\rangle}^{\left( t \right)}})~,
\end{equation}
due to the continuity of $Q_{hj}$ and $Q_{ij}$ on $\left(0,1\right) \times \left(0,1\right)$, as $G_{hj}$ and $G_{ij}$ are continuous functions shown in Lemma~\ref{Gfunction_property}. By the sampling rule AOAm (\ref{AOAm}), the alternative $\left\langle i \right\rangle$ will be sampled and the alternative $\left\langle h \right\rangle$ will stop receiving simulation replications before the inequality sign in (\ref{equequproof1}) reverses, which leads to a contradiction to that $(r_1^{\left(t\right)},r_2^{\left(t\right)},\cdots,r_k^{\left(t\right)} )$ converges to $\left(\widetilde{r}_1,\widetilde{r}_2,\cdots,\widetilde{r}_k\right)$, and then
$$\mathop {\min }\limits_{j = m + 1,m + 2, \cdots ,k} {Q_{hj}} ( {{{\widetilde r}_{\left\langle h \right\rangle}},{{\widetilde r}_{\left\langle  j \right\rangle}}} ) = \mathop {\min }\limits_{j = m + 1,m + 2, \cdots ,k} {Q_{ij}} ( {{{\widetilde r}_{\left\langle  i \right\rangle}},{{\widetilde r}_{\left\langle  j \right\rangle}}} ),\quad {i,h = 1,2,\cdots,m}~,$$
must hold. If $\exists~j,\ell \in \left\{m+1,m+2, \cdots, k\right\}$, $j \ne \ell$ such that the inequality
$$\mathop {\min }\limits_{i = 1,2, \cdots ,m} {Q_{ij}}( {{{\widetilde r}_{\left\langle i \right\rangle}},{{\widetilde r}_{\left\langle j \right\rangle}}} ) > \mathop {\min }\limits_{i = 1,2, \cdots ,m}{Q_{i \ell}} ( {{{\widetilde r}_{\left\langle i \right\rangle}},{{\widetilde r}_{\left\langle \ell \right\rangle}}} )~,$$
strictly holds, by a similar argument, the alternative $\left\langle {j} \right\rangle$ will stop receiving simulation replications before the inequality sign above reverses, and
$$\mathop {\min }\limits_{i = 1,2, \cdots ,m} {Q_{ij}}( {{{\widetilde r}_{\left\langle i \right\rangle}},{{\widetilde r}_{\left\langle j \right\rangle}}} ) = \mathop {\min }\limits_{i = 1,2, \cdots ,m}{Q_{i \ell}} ( {{{\widetilde r}_{\left\langle i \right\rangle}},{{\widetilde r}_{\left\langle \ell \right\rangle}}} ),\quad {j,\ell=m+1,m+2,\cdots,k}~.$$

Since for $i=1,2,\cdots,m$,
$$\mathop {\min }\limits_{j = m + 1,m + 2, \cdots ,k} {Q_{ij}}\left( {{{\widetilde r}_{\left\langle  i \right\rangle}},{{\widetilde r}_{\left\langle  j \right\rangle}}} \right) = \mathop {\min }\limits_{i = 1,2, \cdots ,m} \mathop {\min }\limits_{j = m + 1,m + 2, \cdots ,k} {Q_{ij}}\left( {{{\widetilde r}_{\left\langle  i \right\rangle}},{{\widetilde r}_{\left\langle  j \right\rangle}}} \right)~,$$
and for $j=m+1,m+2,\cdots,k$,
$$\mathop {\min }\limits_{i = 1,2, \cdots ,m} {Q_{ij}}\left( {{{\widetilde r}_{\left\langle  i \right\rangle}},{{\widetilde r}_{\left\langle  j \right\rangle}}} \right) = \mathop {\min }\limits_{i = 1,2, \cdots ,m} \mathop {\min }\limits_{j = m + 1,m + 2, \cdots ,k} {Q_{ij}}\left( {{{\widetilde r}_{\left\langle  i \right\rangle}},{{\widetilde r}_{\left\langle  j \right\rangle}}} \right)~,$$
summarizing the above, we have (\ref{optrate}) must hold.

If $m \le k-m$, with (\ref{optrate}) and $\sum\nolimits_{i = 1}^k {\widetilde{r}_{\left\langle i \right\rangle}}  = 1$, we have a system of $\left(k-m\right)$ equations, and by the implicit function theorem~\citep{rudin1964principles}, there exist implicit functions ${\left. {{{\widetilde r}_{\left\langle  j \right\rangle}}\left( {{x_1},{x_2}, \cdots ,{x_m}} \right)} \right|_{{x_i} = {{\widetilde r}_{\left\langle  i \right\rangle}}}}$, $i=1,2,\cdots,m$, $j = m+1,m+2,\cdots,k$, since
$${\rm{det}}\left(\widetilde \Phi\right) = \prod\limits_{j = m + 1}^k {{\phi _{{i^{\left( j \right)}},j}}} \left\{ {\sum\limits_{j = m + 1}^k {\phi _{{i^{\left( j \right)}},j}^{ - 1}} } \right\} > 0~,$$
where ${i^{\left( j \right)}} \in {\widetilde I^{\left( j \right)}}$, ${\widetilde I^{\left( j \right)}} \mathop {\rm{ = }}\limits^\Delta  \arg \mathop {\min }\limits_{i = 1,2, \cdots ,m} {Q_{ij}}\left( {{{\widetilde r}_{\left\langle i \right\rangle }},{{\widetilde r}_{\left\langle j \right\rangle }}} \right)$, $j \in \left\{m+1,m+2,\cdots,k\right\}$ be a set containing alternatives ${\left\langle i \right\rangle }$, which achieve the minimum value of ${Q_{ij}}\left( {{{\widetilde r}_{\left\langle i \right\rangle }},{{\widetilde r}_{\left\langle j \right\rangle }}} \right)$ for a fixed alternative ${\left\langle j \right\rangle }$,
$${\phi_{i,j}}\mathop  = \limits^\Delta  {\left. {\frac{{\partial {Q_{ij}}( {{\widetilde r}_{\left\langle i \right\rangle}},y)}}{{\partial y}}} \right|_{y = {\widetilde r}_{\left\langle j \right\rangle}}},\quad i \in \left\{1,2,\cdots,m\right\},\;j \in \left\{m+1,m+2,\cdots,k\right\}~,$$
and $\widetilde{\Phi} \in \mathbb{R}^{\left(k-m\right) \times \left(k-m\right)}$,
$$\widetilde \Phi \mathop  = \limits^\Delta  \left( {\begin{array}{*{20}{c}}
{{\phi _{{i^{\left( {m + 1} \right)}},m + 1}}}&{ - {\phi _{{i^{\left( {m + 2} \right)}},m + 2}}}& \cdots &0&0\\
0&{{\phi _{{i^{\left( {m + 2} \right)}},m + 2}}}& \cdots &0&0\\
 \vdots & \vdots & \ddots & \vdots & \vdots \\
0&0& \cdots &{{\phi _{{i^{\left( {k - 1} \right)}},k - 1}}}&{ - {\phi _{{i^{\left( k \right)}},k}}}\\
1&1& \cdots &1&1
\end{array}} \right)~.$$

By taking derivatives of (\ref{optrate}) and $\sum\nolimits_{\ell = 1}^k {{{\widetilde r}_{\left\langle \ell \right\rangle}}} = 1$ with respect to ${\widetilde r}_{\left\langle i \right\rangle}$, $i = 1,2,\cdots,m$, respectively, we have $\widetilde{\Phi} \cdot \widetilde{R}_i = \Upsilon_i$, $i = 1,2,\cdots,m$, where $\widetilde{R}_i \in \mathbb{R}^{\left(k-m\right) \times 1}$ and $\Upsilon_i \in \mathbb{R}^{\left(k-m\right) \times 1}$,
$${\psi _{i,j}}\mathop  = \limits^\Delta  {\left. {\frac{{\partial {Q_{ij}}\left( {y,{{\widetilde r}_{\left\langle j \right\rangle}}} \right)}}{{\partial y}}} \right|_{y = {{\widetilde r}_{\left\langle i \right\rangle}}}},\quad i \in \left\{1,2,\cdots,m\right\},\;j \in \left\{m+1,m+2,\cdots,k\right\}~,$$
$$\psi _{\ell}^{\left( i \right)}\mathop  = \limits^\Delta\left\{
\begin{array}{ccl}
{\psi _{i,m + \ell + 1}} - {\psi _{i,m + \ell}} & , & {{i^{\left( {m + \ell} \right)}} = {i^{\left( {m + \ell + 1} \right)}} = i}\\
 - {\psi _{i,m + \ell}}& , & {i^{\left( {m + \ell} \right)}} = i,\;{i^{\left( {m + \ell + 1} \right)}} \ne i\\
{\psi _{i,m + \ell + 1}} & , & {i^{\left( {m + \ell + 1} \right)}} = i,\;{i^{\left( {m + \ell} \right)}} \ne i\\
0 & ,& {i^{\left( {m + \ell} \right)}} \ne i,\;{i^{\left( {m + \ell + 1} \right)}} \ne i\
\end{array},\quad {i=1,2,\cdots,m},\;\ell=1,2,\cdots,k-m-1~,\right.$$
$${\widetilde{R}_{i}}\mathop  = \limits^\Delta  {\left. {{{\left( {\frac{{\partial {{\widetilde r}_{\left\langle {m+1} \right\rangle}}\left( y \right)}}{{\partial y}},\frac{{\partial {{\widetilde r}_{\left\langle {m+2} \right\rangle}}\left( y \right)}}{{\partial y}}, \cdots, \frac{{\partial {{\widetilde r}_{\left\langle {k} \right\rangle}}\left( y \right)}}{{\partial y}}} \right)}^\prime }} \right|_{y = {{\widetilde r}_{\left\langle i \right\rangle}}}},\quad i=1,2,\cdots,m~,$$
$${\Upsilon _i}\mathop  = \limits^\Delta  {\left( {\psi _1^{\left( i \right)},\psi _2^{\left( i \right)}, \cdots ,\psi _{k - m - 1}^{\left( i \right)}, - 1} \right)^\prime },\quad i=1,2,\cdots,m~,$$
where $^\prime$ denotes the transpose operation of the matrix. Note that
$${\left. {\frac{{\partial {Q_{ij}}\left( {y,{{\widetilde r}_{\left\langle  j \right\rangle}}} \right)}}{{\partial y}}} \right|_{y = {{\widetilde r}_{\left\langle i^{\prime} \right\rangle}}}} = 0,\quad j=m+1,m+2,\cdots,k,\; i \ne i^{\prime}~,$$
due to ${\widetilde r}_{\left\langle i \right\rangle}$ and ${\widetilde r}_{\left\langle i^{\prime} \right\rangle}$, $i \ne i^{\prime}$, $i,i^{\prime} \in \left\{1,2,\cdots,m\right\}$ are independent. In addition, with the sampling AOAm rule (\ref{AOAm}), we have for $i \in \left\{1,2,\cdots,m\right\}$
$${\left. {\frac{{\partial {Q_{{i^{\left( j \right)}}j}}\left( {y,{{\widetilde r}_{\left\langle  j \right\rangle}}} \right)}}{{\partial y}}} \right|_{y = {{\widetilde r}_{\left\langle  i \right\rangle}}}} + {\left. {\frac{{\partial {Q_{{i^{\left( j \right)}}j}}\left( {{{\widetilde r}_{{i^{\left( j \right)}}}},y} \right)}}{{\partial y}}} \right|_{y = {{\widetilde r}_{\left\langle  j \right\rangle}}}}{\left. {\frac{{\partial {{\widetilde r}_{\left\langle  j \right\rangle}}\left( y \right)}}{{\partial y}}} \right|_{y = {{\widetilde r}_{\left\langle  i \right\rangle}}}} = 0,\quad j=m+1,m+2,\cdots,k~,$$
where ${\left. {\frac{{\partial {Q_{{i^{\left( j \right)}}j}}\left( {y,{{\widetilde r}_{\left\langle  j \right\rangle}}} \right)}}{{\partial y}}} \right|_{y = {{\widetilde r}_{\left\langle  i \right\rangle}}}} = 0$ if ${i^{\left( j \right)}} \ne {\left\langle i \right\rangle}$; otherwise, there exists $\ell \in \left\{m+1,m+2,\cdots,k\right\}$ such that the equality above does not hold, say
$${\left. {\frac{{\partial {Q_{{i^{\left( \ell \right)}}\ell}}\left( {y,{{\widetilde r}_{\left\langle \ell \right\rangle}}} \right)}}{{\partial y}}} \right|_{y = {{\widetilde r}_{\left\langle  i \right\rangle}}}} + {\left. {\frac{{\partial {Q_{{i^{\left( \ell \right)}}\ell}}\left( {{{\widetilde r}_{{i^{\left( \ell \right)}}}},y} \right)}}{{\partial y}}} \right|_{y = {{\widetilde r}_{\left\langle \ell \right\rangle}}}}{\left. {\frac{{\partial {{\widetilde r}_{\left\langle \ell \right\rangle}}\left( y \right)}}{{\partial y}}} \right|_{y = {{\widetilde r}_{\left\langle  i \right\rangle}}}} > 0~,$$
and following the sampling rule AOAm (\ref{AOAm}), alternative ${\left\langle i \right\rangle }$ will be sampled and alternative ${\left\langle \ell \right\rangle }$ will stop receiving simulation replications before the inequality sign above is no longer hold, which contradicts that $( {r_1^{\left( t \right)},r_2^{\left( t \right)}, \cdots r_k^{\left( t \right)}} )$ converges to $\left( {{{\widetilde r}_1},{{\widetilde r}_2}, \cdots {{\widetilde r}_k}} \right)$. Then we have $\widetilde{\Psi} \cdot \widetilde{R}_{i} = G_i$, $i= 1,2,\cdots,m$, where $\widetilde{\Psi} \in {\mathbb{R}^{\left(k-m\right) \times \left(k-m\right)}}$ and $G_i \in {\mathbb{R}^{\left(k-m\right) \times 1}}$,
$${\widetilde{\Psi}}\mathop  = \limits^\Delta  diag\left( {{\phi _{{i^{\left(m+1\right)}},m + 1}},{\phi _{{i^{\left(m+2\right)}},m + 2}}, \cdots ,{\phi _{{i^{\left(k\right)}},k}}} \right)~,$$
is a diagonal matrix, and
$${g_\ell^{\left( i \right)}}\mathop  = \limits^\Delta\left\{
\begin{array}{ccl}
 - {\psi _{i,m + \ell}} & , & {{i^{\left( {m + \ell} \right)}} = {\left\langle i \right\rangle}}\\
0 & ,& {i^{\left( {m + \ell} \right)}} \ne {\left\langle i \right\rangle}
\end{array}~,\quad {i=1,2,\cdots,m},\;\ell=1,2,\cdots,k-m~,\right.$$
$${G_i}\mathop  = \limits^\Delta  {\left( {g_1^{\left( i \right)},g_2^{\left( i \right)}, \cdots ,g_{k - m}^{\left( i \right)}} \right)^\prime }~,$$
where $^\prime$ denotes the transpose operation of the matrix. According to Lemma~\ref{lem}, we have ${{\phi _{i,j}}} > 0$, $i = 1,2,\cdots,m$, $j=m+1,m+2,\cdots,k$, which yields $\det \left( {\widetilde{\Psi}} \right) > 0$ and $\widetilde{\Psi}$ is invertible,
$${{\widetilde{\Psi}}^{-1}} = diag\left( {{\phi _{{i^{\left(m+1\right)}},m + 1}^{-1}},{\phi _{{i^{\left(m+2\right)}},m + 2}^{-1}}, \cdots ,{\phi_{{i^{\left(k\right)}},k}^{-1}}} \right)~.$$

Therefore, we have $\Upsilon_i  = {\widetilde{\Phi}} \cdot {{\widetilde{\Psi}}^{-1}} \cdot G_i$, $i = 1,2,\cdots,m$, which leads to
$$\sum\limits_{\left\{ {\left. j \right|{i^{\left( j \right)}} = {\left\langle i \right\rangle}} \right\}} {\frac{{{{\left. {{{\partial {Q_{{i^{\left( j \right)}} j}}\left( {y,{{\widetilde r}_{\left\langle  j \right\rangle}}} \right)} \mathord{\left/
 {\vphantom {{\partial {Q_{{i^{\left( j \right)}} j }}\left( {y,{{\widetilde r}_{\left\langle  j \right\rangle}}} \right)} {\partial y}}} \right.
 \kern-\nulldelimiterspace} {\partial y}}} \right|}_{y = {{\widetilde r}_{{i^{\left( j \right)}}}}}}}}{{{{\left. {{{\partial {Q_{{i^{\left( j \right)}} j }}\left( {{{\widetilde r}_{{i^{\left( j \right)}}}},y} \right)} \mathord{\left/
 {\vphantom {{\partial {Q_{{i^{\left( j \right)}}\left( j \right)}}\left( {{{\widetilde r}_{{i^{\left( j \right)}}}},y} \right)} {\partial y}}} \right.
 \kern-\nulldelimiterspace} {\partial y}}} \right|}_{y = {{\widetilde r}_{\left\langle  j \right\rangle}}}}}}} = 1,\quad i=1,2,\cdots,m~,$$
and with (\ref{lammaequ1}) and (\ref{lammaequ2}), the equality above yields
$$\frac{{{{\left( {{{\widetilde r}_{\left\langle  i \right\rangle}}} \right)}^2}}}{{\sigma _{\left\langle i \right\rangle }^2}} = \sum\limits_{\left\{ {\left. j \right|{i^{\left( j \right)}} = {\left\langle i \right\rangle}} \right\}} {\frac{{{{\left( {{{\widetilde r}_{\left\langle  j \right\rangle}}} \right)}^2}}}{{\sigma _{\left\langle j \right\rangle }^2}}},\quad i=1,2,\cdots,m~.$$

Since $\bigcup\limits_{i = 1,2, \cdots ,m} {\left\{ {\left. j \right|{i^{\left( j \right)}} = {\left\langle i \right\rangle}} \right\}}  = \left\{ {m + 1,m + 2, \cdots ,k} \right\}$ and $\bigcap\limits_{i = 1,2, \cdots ,m} {\left\{ {\left. j \right|{i^{\left( j \right)}} = {\left\langle i \right\rangle}} \right\} = \emptyset }$, we have (\ref{thmbalance}) must hold.

By a similar argument, if $m > k-m$, with (\ref{optrate}) and $\sum\nolimits_{i = 1}^k {\widetilde{r}_{\left\langle i \right\rangle}}  = 1$, we have a system of $m$ equations, and by the implicit function theorem~\citep{rudin1964principles}, there exist implicit functions ${\left. {{{\widetilde r}_{\left\langle  i \right\rangle}}\left( {{x_{m+1}},{x_{m+2}}, \cdots ,{x_k}} \right)} \right|_{{x_j} = {{\widetilde r}_{\left\langle  j \right\rangle}}}}$, $i = 1,2,\cdots,m$, $j=m+1,m+2,\cdots,k$, since
$${\rm{det}}\left(\widehat \Phi\right) = \prod\limits_{i = 1}^m {{\psi _{i,j^{\left( i \right)}}}} \left\{ {\sum\limits_{i =  1}^m {\psi _{i,j^{\left( i \right)}}^{ - 1}} } \right\} > 0~,$$
where ${j^{\left( i \right)}} \in {\widetilde {J}^{\left( i \right)}}$, ${\widetilde J^{\left( i \right)}} \mathop {\rm{ = }}\limits^\Delta  \arg \mathop {\min }\limits_{j = m+1,m+2, \cdots ,k} {Q_{ij}}\left( {{{\widetilde r}_{\left\langle i \right\rangle }},{{\widetilde r}_{\left\langle j \right\rangle }}} \right)$, $i \in \left\{1,2,\cdots,m\right\}$ be a set containing alternatives ${\left\langle j \right\rangle}$, which achieve the minimum value of ${Q_{ij}}\left( {{{\widetilde r}_{\left\langle i \right\rangle }},{{\widetilde r}_{\left\langle j \right\rangle }}} \right)$ for a fixed alternative ${\left\langle i \right\rangle }$, and $\widehat{\Phi} \in \mathbb{R}^{m \times m}$,
$$\widehat \Phi \mathop  = \limits^\Delta  \left( {\begin{array}{*{20}{c}}
{{\psi _{1,j^{\left( 1\right)}}}}&{ - {\psi _{2,j^{\left( 2 \right)}}}}& \cdots &0&0\\
0&{\psi _{2,j^{\left( 2 \right)}}}& \cdots &0&0\\
 \vdots & \vdots & \ddots & \vdots & \vdots \\
0&0& \cdots &{{\psi _{m-1,j^{\left( {m - 1} \right)}}}}&{-{\psi _{m,j^{\left( {m} \right)}}}}\\
1&1& \cdots &1&1
\end{array}} \right)~.$$

By taking derivatives of (\ref{optrate}) and $\sum\nolimits_{\ell = 1}^k {{{\widetilde r}_{\left\langle \ell \right\rangle}}} = 1$ with respect to ${\widetilde r}_{\left\langle j \right\rangle}$, $j = m+1,m+2,\cdots,k$, respectively, we have $\widehat{\Phi} \cdot \widehat{R}_j = \widehat{\Upsilon}_j$, $j = m+1,m+2,\cdots,k$, where $\widehat{R}_j \in \mathbb{R}^{m \times 1}$ and $\widehat{\Upsilon}_j \in \mathbb{R}^{m \times 1}$,
$$\widehat{\psi} _{h}^{\left( j \right)}\mathop  = \limits^\Delta\left\{
\begin{array}{ccl}
{\phi _{h+1,j}} - {\phi _{h,j}} & , & {{j^{\left( h \right)}} = {j^{\left( {h + 1} \right)}} = j}\\
 - {\phi _{h,j}}& , & {j^{\left( {h} \right)}} = j,\;{j^{\left( {h + 1} \right)}} \ne j\\
{\phi _{h+1,j}} & , & {j^{\left( {h + 1} \right)}} = j,\;{j^{\left( {h} \right)}} \ne j\\
0 & ,& {j^{\left( {h} \right)}} \ne j,\;{j^{\left( {h + 1} \right)}} \ne j\
\end{array},\quad {j=m+1,m+2,\cdots,k},\;h=1,2,\cdots,m-1~,\right.$$
$${\widehat{R}_j}\mathop  = \limits^\Delta  {\left. {{{\left( {\frac{{\partial {{\widetilde r}_{\left\langle 1 \right\rangle}}\left( y \right)}}{{\partial y}},\frac{{\partial {{\widetilde r}_{\left\langle 2 \right\rangle}}\left( y \right)}}{{\partial y}}, \cdots, \frac{{\partial {{\widetilde r}_{\left\langle m \right\rangle}}\left( y \right)}}{{\partial y}}} \right)}^\prime }} \right|_{y = {{\widetilde r}_{\left\langle j \right\rangle}}}},\quad j=m+1,m+2,\cdots,k~,$$
$${\widehat{\Upsilon} _j}\mathop  = \limits^\Delta  {\left( {\widehat{\psi} _1^{\left( j \right)},\widehat{\psi} _2^{\left( j \right)}, \cdots ,\widehat{\psi} _{m - 1}^{\left( j \right)}, - 1} \right)^\prime },\quad j=m+1,m+2,\cdots,k~,$$
where $^\prime$ denotes the transpose operation of the matrix. Note that
$${\left. {\frac{{\partial {Q_{ij}}\left( {{{\widetilde r}_{\left\langle  i \right\rangle}},y} \right)}}{{\partial y}}} \right|_{y = {{\widetilde r}_{\left\langle j^{\prime} \right\rangle}}}} = 0,\quad i=1,2,\cdots,m,\; j \ne j^{\prime}~,$$
due to ${\widetilde r}_{\left\langle j \right\rangle}$ and ${\widetilde r}_{\left\langle j^{\prime} \right\rangle}$, $j \ne j^{\prime}$, $j,j^{\prime} \in \left\{m+1,m+2,\cdots,k\right\}$ are independent. In addition, with the sampling AOAm rule (\ref{AOAm}), we have for $j \in \left\{m+1,m+2,\cdots,k\right\}$
$${\left. {\frac{{\partial {Q_{ij^{\left( i \right)}}}\left( {y,{{\widetilde r}_{j^{\left( i \right)}}}} \right)}}{{\partial y}}} \right|_{y = {{\widetilde r}_{\left\langle  i \right\rangle}}}}{\left. {\frac{{\partial {{\widetilde r}_{\left\langle  i \right\rangle}}\left( y \right)}}{{\partial y}}} \right|_{y = {{\widetilde r}_{\left\langle  j \right\rangle}}}} + {\left. {\frac{{\partial {Q_{i{j^{\left( i \right)}}}}\left( {{{\widetilde r}_{\left\langle i \right\rangle}},y} \right)}}{{\partial y}}} \right|_{y = {{\widetilde r}_{\left\langle  j \right\rangle}}}} = 0,\quad i=1,2,\cdots,m~,$$
where ${\left. {\frac{{\partial {Q_{ij^{\left( i \right)}}}\left( {{\widetilde r}_{\left\langle  i \right\rangle}},y\right)}}{{\partial y}}} \right|_{y = {{\widetilde r}_{\left\langle  j \right\rangle}}}} = 0$ if ${j^{\left( i \right)}} \ne {\left\langle j \right\rangle}$; otherwise, there exists $h \in \left\{1,2,\cdots,m\right\}$ such that the equality above does not hold, say
$${\left. {\frac{{\partial {Q_{hj^{\left( h \right)}}}\left( {y,{{\widetilde r}_{j^{\left( h \right)}}}} \right)}}{{\partial y}}} \right|_{y = {{\widetilde r}_{\left\langle h \right\rangle}}}}{\left. {\frac{{\partial {{\widetilde r}_{\left\langle h \right\rangle}}\left( y \right)}}{{\partial y}}} \right|_{y = {{\widetilde r}_{\left\langle  j \right\rangle}}}} + {\left. {\frac{{\partial {Q_{i{j^{\left( h \right)}}}}\left( {{{\widetilde r}_{\left\langle h \right\rangle}},y} \right)}}{{\partial y}}} \right|_{y = {{\widetilde r}_{\left\langle  j \right\rangle}}}} > 0~,$$
and following the sampling rule AOAm (\ref{AOAm}), alternative ${\left\langle j \right\rangle }$ will be sampled and alternative ${\left\langle h \right\rangle }$ will stop receiving simulation replications before the inequality sign above is no longer hold, which contradicts that $( {r_1^{\left( t \right)},r_2^{\left( t \right)}, \cdots r_k^{\left( t \right)}} )$ converges to $\left( {{{\widetilde r}_1},{{\widetilde r}_2}, \cdots {{\widetilde r}_k}} \right)$. Then we have $\widehat{\Psi} \cdot \widehat{R}_j = \widehat{G}_j$, $j= m+1,m+2,\cdots,k$, where $\widehat{\Psi} \in {\mathbb{R}^{m \times m}}$ and $\widehat{G}_j \in {\mathbb{R}^{m \times 1}}$,
$${\widehat{\Psi}}\mathop  = \limits^\Delta  diag\left( {{\psi _{1,j^{\left(1\right)}}},{\psi _{2,j^{\left(2\right)}}}, \cdots ,{\psi _{m,j^{\left(m\right)}}}} \right)~,$$
is a diagonal matrix, and
$${\widehat{g}_h^{\left( j \right)}}\mathop  = \limits^\Delta\left\{
\begin{array}{ccl}
 - {\phi _{h,j}} & , & {{j^{\left( {h} \right)}} = {\left\langle j \right\rangle}}\\
0 & ,& {j^{\left( h \right)}} \ne {\left\langle j \right\rangle}
\end{array},\quad {j=m+1,m+2,\cdots,k},\;h=1,2,\cdots,m~,\right.$$
$${\widehat{G}_j}\mathop  = \limits^\Delta  {\left( {\widehat{g}_1^{\left( j \right)},\widehat{g}_2^{\left( j \right)}, \cdots ,\widehat{g}_{m}^{\left( j \right)}} \right)^\prime }~,$$
where $^\prime$ denotes the transpose operation of the matrix. According to Lemma~\ref{lem}, we have ${{\psi _{i,j}}} > 0$, $i = 1,2,\cdots,m$, $j=m+1,m+2,\cdots,k$, which yields $\det \left( {\widehat{\Psi}} \right) > 0$ and $\widehat{\Psi}$ is invertible,
$${{\widehat{\Psi}}^{-1}} = diag\left( {{\psi _{1,j^{\left(1\right)}}^{-1}},{\psi _{2,j^{\left(2\right)}}^{-1}}, \cdots ,{\psi _{m,j^{\left(m\right)}}^{-1}}} \right)~.$$

Therefore, we have $\widehat{\Upsilon}_j  = {\widehat{\Phi}} \cdot {{\widehat{\Psi}}^{-1}} \cdot \widehat{G}_j$, $j = m+1,m+2,\cdots,k$, which leads to
$$\sum\limits_{\left\{ {\left. i \right|{j^{\left( i \right)}} = {\left\langle j \right\rangle}} \right\}} {\frac{{{{\left. {{{\partial {Q_{i j^{\left( i \right)}}}\left( {{{\widetilde r}_{\left\langle  i \right\rangle}},y} \right)} \mathord{\left/
 {\vphantom {{\partial {Q_{i j^{\left( i \right)}}}\left( {{{\widetilde r}_{\left\langle  i \right\rangle}},y} \right)} {\partial y}}} \right.
 \kern-\nulldelimiterspace} {\partial y}}} \right|}_{y = {{\widetilde r}_{{j^{\left( i \right)}}}}}}}}{{{{\left. {{{\partial {Q_{i j^{\left( i \right)}}}\left( {y,{{\widetilde r}_{{j^{\left( i \right)}}}}} \right)} \mathord{\left/
 {\vphantom {{\partial {Q_{i j^{\left( i \right)}}}\left( {y,{{\widetilde r}_{{j^{\left( i \right)}}}}} \right)} {\partial y}}} \right.
 \kern-\nulldelimiterspace} {\partial y}}} \right|}_{y = {{\widetilde r}_{\left\langle  i \right\rangle}}}}}}} = 1,\quad j=m+1,m+2,\cdots,k~,$$
and with (\ref{lammaequ1}) and (\ref{lammaequ2}), the equality above yields
$$\frac{{{{\left( {{{\widetilde r}_{\left\langle  j \right\rangle}}} \right)}^2}}}{{\sigma _{\left\langle j \right\rangle }^2}} = \sum\limits_{\left\{ {\left. i \right|{j^{\left( i \right)}} = {\left\langle j \right\rangle}} \right\}} {\frac{{{{\left( {{{\widetilde r}_{\left\langle  i \right\rangle}}} \right)}^2}}}{{\sigma _{\left\langle i \right\rangle }^2}}},\quad j=m+1,m+2,\cdots,k~.$$

Since $\bigcup\limits_{j = m+1,m+2, \cdots ,k} {\left\{ {\left. i \right|{j^{\left( i \right)}} = {\left\langle j \right\rangle}} \right\}}  = \left\{ {1,2, \cdots ,m} \right\}$ and $\bigcap\limits_{j = m+1,m+2, \cdots ,k} {\left\{ {\left. i \right|{j^{\left( i \right)}} = {\left\langle j \right\rangle}} \right\} = \emptyset }$, we have (\ref{thmbalance}) must hold.
Summarizing the above, we have $\left( {{{\widetilde r}_1},{{\widetilde r}_2}, \cdots {{\widetilde r}_k}} \right)$ satisfies both (\ref{optrate}) and (\ref{thmbalance}). Since there is only one solution to (\ref{optrate}) and (\ref{thmbalance}), we have $\mathop {\lim }\limits_{t \to \infty } r_{\left\langle  i \right\rangle}^{\left( t \right)} = r_{\left\langle  i \right\rangle}^*$, $i=1,2,\cdots,k$.


\subsection*{A.3. Supplementary Experiments}

\subsubsection*{A.3.1 Comparisons of OCBAss and OCBASS rules given the true parameters with other allocation procedures in \textit{Experiment 1}} To illustrate the difference between the asymptotic property and finite-sample performance of sequential allocation procedures, in this numerical experiment, the proposed AOAm  is tested with: EA, OCBAm, OCBAm$+$, OCBASSt and OCBASS. Specifically, OCBAm allocates simulation replications based on OCBAm rule with parameters estimated sequentially by available sample information (same as the OCBAm in the main text); OCBAm+ allocates simulation replications based on OCBAm+ rule with parameters estimated sequentially by available sample information (same as the OCBAm$+$ in the main text); OCBAsst allocates simulation replications according to the OCBAss rule under perfect information (assuming parameters are known); OCBASSt allocates simulation replications based on OCBASS rule under perfect information (assuming parameters are known). The numerical settings are the same as \textit{Experiment 1} in the main text.

\begin{figure}[htbp]
\centering
\subfigure[IPCS]{\includegraphics[width=0.49\textwidth]{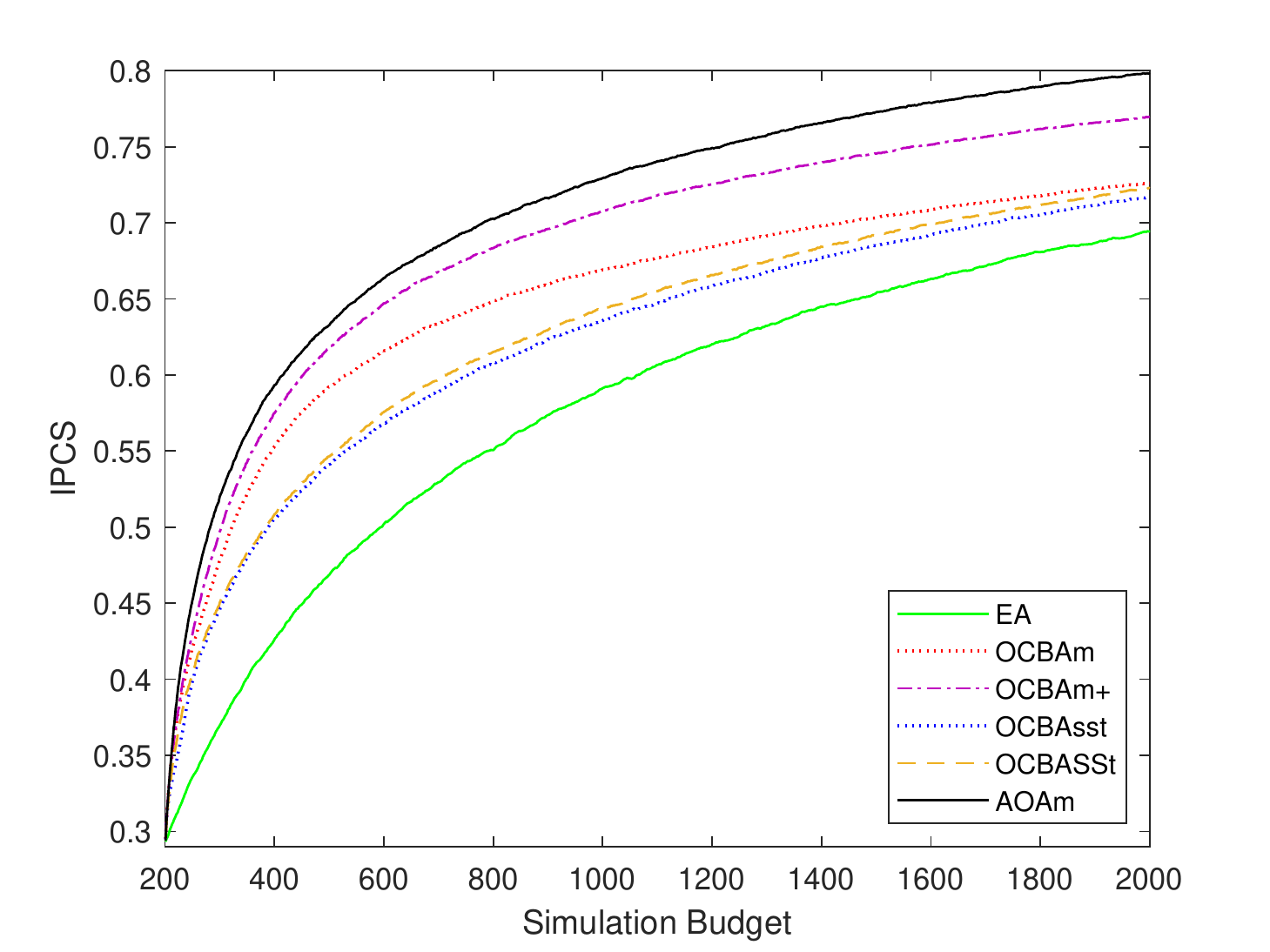}}
\subfigure[EOC]{\includegraphics[width=0.49\textwidth]{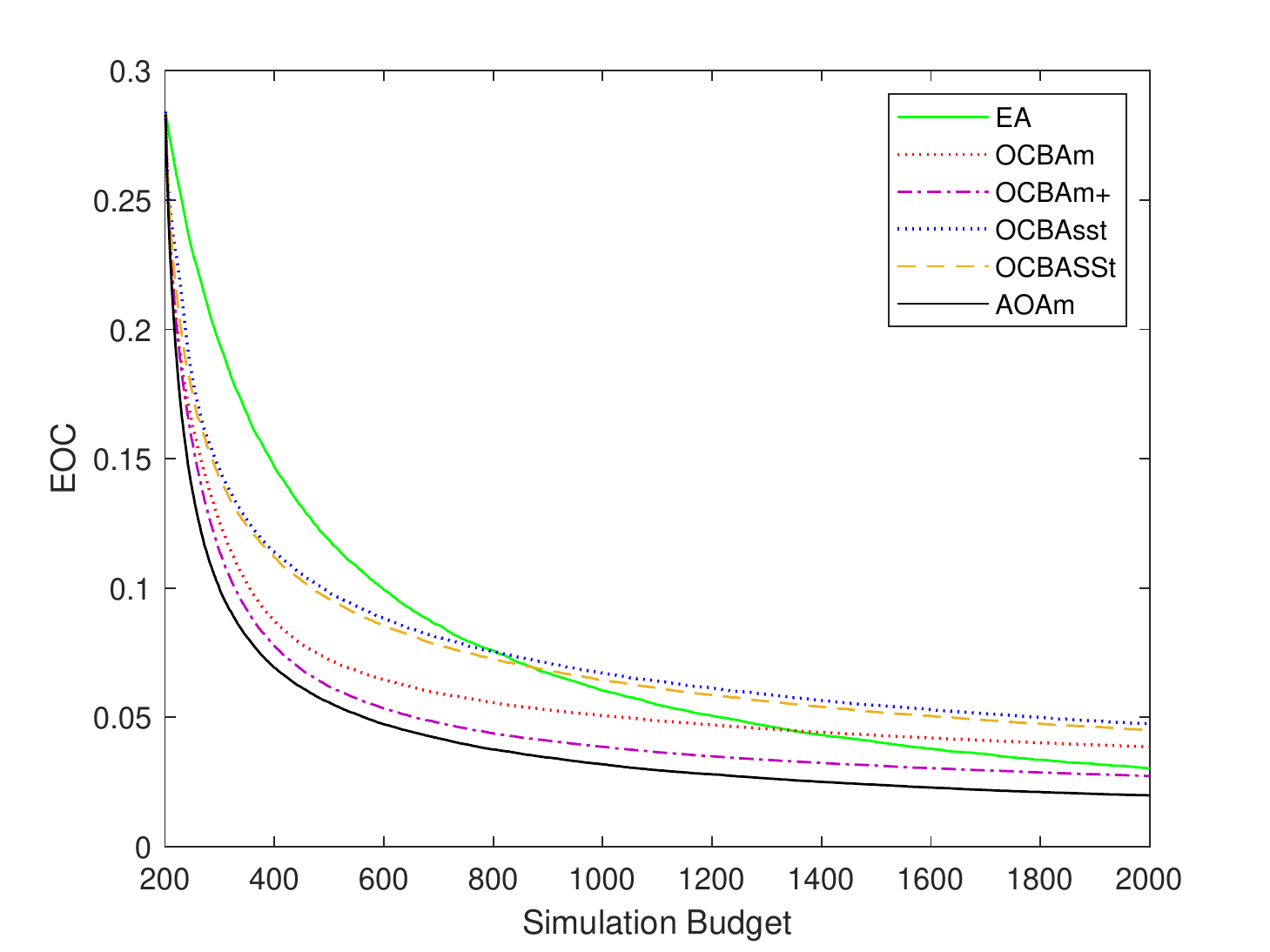}}
\caption{Comparison of IPCSs and EOCs for 6 sampling allocation procedures in Experiment A.3.1.}
\label{figA24}
\end{figure}

From Figure~\ref{figA24} (a), we can see that the IPCS of EA increases at a slow pace and it performs the worst among all allocation procedures. OCBASSt and OCBAsst have a comparable performance at the beginning, and the former has a slight edge over the latter as the simulation budget grows. OCBAm has an edge over OCBAsst and OCBASSt. AOAm, which performs the best among all sampling procedures, has an edge over OCBAm+. The performances of different allocation procedures in terms of EOC shown in Figure~\ref{figA24} (b) exhibit similar behaviors as measured by IPCS.

Compared with Figure~\ref{fig2}, we find that OCBAss and OCBASS allocation procedures lead to much poorer performances when the true parameters are given as inputs. Therefore, the desirable asymptotic property derived from the static optimization problem is inadequate to explain the good performances of OCBAss and OCBASS when simulation budget is finite, which highlights the difference between the asymptotic property and finite-sample performance for allocation procedures. Sequentially implementing the asymptotic optimality conditions properly can result in a good sampling policy in some circumstances, where the objective is to select the alternatives with the optimal mean from alternatives following normal sampling distributions, since the mean-variance trade-off happens to lead to a desirable allocation policy in certain circumstances.

\subsubsection*{A.3.2 Details of emergency evacuation planning problem} 
\begin{figure}[htbp]
\center{\includegraphics[width=0.8\textwidth]{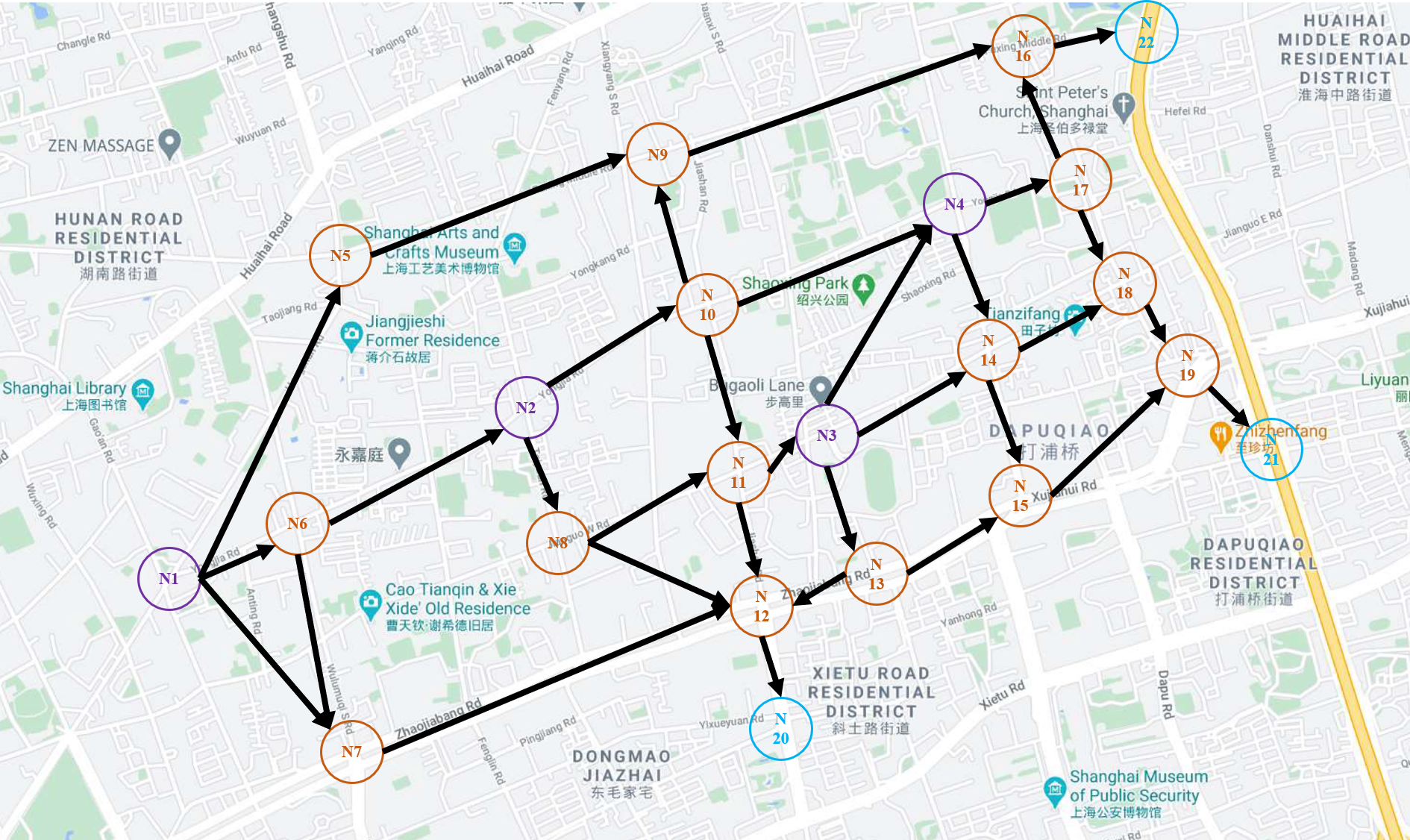}}
\caption{A transportation network of a region of the city of Shanghai, China.}
\label{figshanghai}
\end{figure}

In this section, we show some details of the emergency evacuation planning example in the main text. The evacuation network is based on a transportation network of a region in Shanghai, China, shown in Figure~\ref{figshanghai}.

In order to produce a log-normal distribution with desired mean and variance on the top of each edge in Figure \ref{figrealnetwork}, we can solve the parameters in log-normal distribution by
$$\mu _{ij}^c = \log \left( {\frac{{{{\left( {\mathbb{E}\left( {{c_{ij}}} \right)} \right)}^2}}}{{\sqrt {Var\left( {{c_{ij}}} \right) + {{\left( {\mathbb{E}\left( {{c_{ij}}} \right)} \right)}^2}} }}} \right),\quad \mu _{ij}^t = \log \left( {\frac{{{{\left( {\mathbb{E}\left( {{t_{ij}}} \right)} \right)}^2}}}{{\sqrt {Var\left( {{t_{ij}}} \right) + {{\left( {\mathbb{E}\left( {{t_{ij}}} \right)} \right)}^2}} }}} \right)~,$$
and
$${\left( {\sigma _{ij}^c} \right)^2} = \log \left( {\frac{{Var\left( {{c_{ij}}} \right)}}{{{{\left( {\mathbb{E}\left( {{c_{ij}}} \right)} \right)}^2}}} + 1} \right),\quad {\left( {\sigma _{ij}^t} \right)^2} = \log \left( {\frac{{Var\left( {{t_{ij}}} \right)}}{{{{\left( {\mathbb{E}\left( {{t_{ij}}} \right)} \right)}^2}}} + 1} \right)~,$$
i.e., $c_{ij} \mathop  \sim \limits^{i.i.d.} Lognormal ( \mu _{ij}^c,{\left( {\sigma _{ij}^c} \right)^2} )$ and $t_{ij} \mathop  \sim \limits^{i.i.d.} Lognormal (\mu _{ij}^t,{\left( {\sigma _{ij}^t} \right)^2})$, $i,j\in \left\{1,2,\cdots,\left| {V'} \right|\right\}$.

We prepare 3 alternative evacuation paths for each source node, summarized in Table~\ref{table3}, which are determined by the top-3 minimum value of $\mathbb{E}\left( {T_s^r} \right) = \sum\nolimits_{i = 0}^{h_{sd}^r} {\mathbb{E}( {{t_{n_i^r,n_{i + 1}^r}}} )}$ among all feasible evacuation paths starting from the source node $s \in \left\{1,2,3,4\right\}$ to any destination nodes $d = 20,21,22$. In the example, an evacuation plan contains 2 evacuation paths for each source node, and then there are ${\left( {C_3^2} \right)^4} = 81$ alternatives by combinatorics.

\begin{table}[htbp]
\caption{Alternative Evacuation Paths for Each Node in Experiment 5}
\label{table3}
\centering
\begin{tabular}{ccccc}
\hline
   Nodes & Labels & Paths & \tabincell{c}{Expected Maximum \\ Flow $\mathbb{E}\left(F_s^r\right)$ (persons)} & \tabincell{c}{Expected Travel \\ Time $\mathbb{E}\left(T_s^r\right)$ (mins)} \\
\hline
   $N_{1}$ & \tabincell{c}{$P_{11}$\\$P_{12}$\\$P_{13}$} & \tabincell{l}{1 $\to$ 7 $\to$ 12 $\to$ 20 \\ 1 $\to$ 6 $\to$ 7 $\to$ 12 $\to$ 20 \\ 1 $\to$ 5 $\to$ 9 $\to$ 16 $\to$ 22} & \tabincell{c}{20 \\ 10 \\ 20} & \tabincell{c}{6 \\ 9 \\ 12} \\
   \hline
   $N_{2}$ & \tabincell{c}{$P_{21}$\\$P_{22}$\\$P_{23}$} &  \tabincell{l}{2 $\to$ 8 $\to$ 12 $\to$ 20 \\ 2 $\to$ 10  $\to$ 9 $\to$ 16 $\to$ 22 \\ 2 $\to$ 10 $\to$ 11 $\to$ 12 $\to$ 20} & \tabincell{c}{20 \\ 10 \\ 20} & \tabincell{c}{9 \\ 11 \\ 12} \\
 \hline
   $N_{3}$ & \tabincell{c}{$P_{31}$\\$P_{32}$\\$P_{33}$} & \tabincell{l}{ 3 $\to$ 14 $\to$ 18 $\to$ 19 $\to$ 21 \\ 3 $\to$ 13 $\to$ 12 $\to$ 20 \\ 3 $\to$ 14 $\to$ 15 $\to$ 19 $\to$ 21} & \tabincell{c}{20 \\ 10 \\ 20} & \tabincell{c}{8 \\ 10 \\ 11} \\
   \hline
   $N_{4}$ & \tabincell{c}{$P_{41}$\\$P_{42}$\\$P_{43}$} & \tabincell{l}{4 $\to$ 14 $\to$ 18 $\to$ 19 $\to$ 21 \\ 4 $\to$ 17 $\to$ 16 $\to$ 22 \\ 4 $\to$ 17 $\to$ 18 $\to$ 19 $\to$ 21} & \tabincell{c}{20 \\ 10 \\ 20} & \tabincell{c}{8 \\ 9 \\ 10} \\
\hline
\end{tabular}
\end{table}

From Table~\ref{table3}, we can see that there are some paths having common routes, e.g., $P_{12}$ and $P_{21}$ have common route $\left(12,20\right)$; $P_{31}$ and $P_{41}$ have common routes $\left(14,18\right)$, $\left(18,19\right)$ and $\left(19,21\right)$. Let $t_s^r\left(n_i^r\right)$, $s=1,2,3,4$, $r=1,2$, $i=1,2,\cdots,h_{sd}^r$ be the travel time from source $s$ to intermediate node $n_i^r$ on path $P_{sr}$. Route congestion may exist when evacuation paths having common routes. For example, evacuees from source nodes $1$ and $2$ first arrive at intermediate node 12 at time $t_1^2\left(12\right)$ and $t_2^1\left(12\right)$, respectively, and we assume that $t_1^2\left(12\right) < t_2^1\left(12\right)$. For ${t_1^2\left(12\right)} \le t < {t_2^1\left(12\right)}$, no route congestion exists, whereas for $t \ge {t_2^1\left(12\right)}$, evacuees from both source nodes $1$ and $2$ need to travel on route $\left(12,20\right)$ and route congestion exists.

The performance of each alternative is estimated using sample averages with a large number of sample sizes, i.e., $10^{6}$, and we treat this estimate as the true performance. The top-$5$ evacuation plans with least expected clearance time are shown in Table \ref{table2}.

\begin{table}[htbp]
\caption{Top-$5$ Evacuation Plans in the Emergency Evacuation Planning Problem}
\label{table2}
\centering
\begin{tabular}{ccc}
\hline
   Alternatives & Evacuation Paths & \tabincell{c}{ Expected \\ Clearance Time} \\
\hline
   1 & \tabincell{l}{${N_1}$: 1 $\to$ 6 $\to$ 7 $\to$ 12 $\to$ 20; 1 $\to$ 5 $\to$ 9 $\to$ 16 $\to$ 22 \\ ${N_2}$: 2 $\to$ 8 $\to$ 12 $\to$ 20; 2 $\to$ 10 $\to$ 11 $\to$ 12 $\to$ 20 \\ ${N_3}$: 3 $\to$ 14 $\to$ 18 $\to$ 19 $\to$ 21; 3 $\to$ 14 $\to$ 15 $\to$ 19 $\to$ 21 \\ ${N_4}$: 4 $\to$ 17 $\to$ 16 $\to$ 22; 4 $\to$ 17 $\to$ 18 $\to$ 19 $\to$ 21 } & \tabincell{c}{10.8043 \\ mins}\\
   \hline
   2 & \tabincell{l}{${N_1}$: 1 $\to$ 6 $\to$ 7 $\to$ 12 $\to$ 20; 1 $\to$ 5 $\to$ 9 $\to$ 16 $\to$ 22 \\ ${N_2}$: 2 $\to$ 8 $\to$ 12 $\to$ 20; 2 $\to$ 10 $\to$ 11 $\to$ 12 $\to$ 20 \\ ${N_3}$: 3 $\to$ 14 $\to$ 18 $\to$ 19 $\to$ 21; 3 $\to$ 13 $\to$ 12 $\to$ 20 \\ ${N_4}$: 4 $\to$ 17 $\to$ 16 $\to$ 22; 4 $\to$ 17 $\to$ 18 $\to$ 19 $\to$ 21 } & \tabincell{c}{10.8505 \\mins} \\
   \hline
   3 & \tabincell{l}{${N_1}$: 1 $\to$ 6 $\to$ 7 $\to$ 12 $\to$ 20; 1 $\to$ 5 $\to$ 9 $\to$ 16 $\to$ 22 \\ ${N_2}$: 2 $\to$ 8 $\to$ 12 $\to$ 20; 2 $\to$ 10 $\to$ 11 $\to$ 12 $\to$ 20 \\ ${N_3}$: 3 $\to$ 14 $\to$ 18 $\to$ 19 $\to$ 21; 3 $\to$ 14 $\to$ 15 $\to$ 19 $\to$ 21 \\ ${N_4}$: 4 $\to$ 14 $\to$ 18 $\to$ 19 $\to$ 21; 4 $\to$ 17 $\to$ 16 $\to$ 22 } & \tabincell{c}{10.9206 \\ mins} \\
   \hline
   4 & \tabincell{l}{${N_1}$: 1 $\to$ 6 $\to$ 7 $\to$ 12 $\to$ 20; 1 $\to$ 5 $\to$ 9 $\to$ 16 $\to$ 22 \\ ${N_2}$: 2 $\to$ 8 $\to$ 12 $\to$ 20; 2 $\to$ 10 $\to$ 11 $\to$ 12 $\to$ 20 \\ ${N_3}$: 3 $\to$ 14 $\to$ 18 $\to$ 19 $\to$ 21; 3 $\to$ 13 $\to$ 12 $\to$ 20 \\ ${N_4}$: 4 $\to$ 14 $\to$ 18 $\to$ 19 $\to$ 21; 4 $\to$ 17 $\to$ 16 $\to$ 22 } & \tabincell{c}{10.9633 \\ mins} \\
   \hline
   5 & \tabincell{l}{${N_1}$: 1 $\to$ 6 $\to$ 7 $\to$ 12 $\to$ 20; 1 $\to$ 5 $\to$ 9 $\to$ 16 $\to$ 22 \\ ${N_2}$: 2 $\to$ 8 $\to$ 12 $\to$ 20; 2 $\to$ 10 $\to$ 11 $\to$ 12 $\to$ 20 \\ ${N_3}$: 3 $\to$ 14 $\to$ 18 $\to$ 19 $\to$ 21; 3 $\to$ 14 $\to$ 15 $\to$ 19 $\to$ 21 \\ ${N_4}$: 4 $\to$ 14 $\to$ 18 $\to$ 19 $\to$ 21; 4 $\to$ 17 $\to$ 18 $\to$ 19 $\to$ 21 } & \tabincell{c}{11.2995 \\ mins} \\
\hline
\end{tabular}
\end{table}

\end{document}